\documentclass{article}
\usepackage{amsmath}
\usepackage{amssymb}
\usepackage{mathrsfs}
\usepackage{graphics}
\usepackage{pdfpages}
\usepackage{caption}
\usepackage{float}
\usepackage{amsthm}
\usepackage{graphicx}
\usepackage{subcaption}
\usepackage[english]{babel}
\usepackage{hyperref}
\usepackage{float}
\usepackage[T1]{fontenc}
\usepackage[utf8]{inputenc}
\usepackage{authblk}

\newtheorem{definition}{Definition}
\newtheorem{properties}{Property}
\DeclareMathOperator\erf{erf}
\newtheorem{theorem}{Theorem}

\begin{document}

\title{An Introduction to Partial differential equations}
\author{Per Kristen Jakobsen}
\affil{Department of Mathematics and Statistics, the Arctic University of Norway, 9019 Troms\o, Norway}
\maketitle
\tableofcontents

\section{Introduction}
The field of partial differential
equations (PDEs) is vast in size and
diversity. The basic reason for this is
that essentially all fundamental laws of physics are formulated in terms of PDEs.
In addition, approximations to these
fundamental laws, that form a patchwork
of mathematical models covering the range from the smallest to the largest observable space-time scales, are also formulated in terms of PDEs.
The diverse applications of PDEs in
science and technology testify to the
flexibility and expressiveness of the
language of PDEs, but it also makes
it a hard topic to teach right. Exactly because of the diversity of applications, there
are just so many different points of view when it comes to PDEs. These lecture notes view the subject through the lens of applied mathematics. From this point of view, the physical context for basic equations like the heat equation, the wave equation and the Laplace equation are introduced early on, and the focus of the lecture notes are on methods, rather than precise mathematical definitions and proofs. With respect to methods, both analytical and numerical approaches are discussed. 

These lecture notes has been succesfully used as the text for  a master class in partial differential equations for several years. The students attending this class are assumed to have  previously attended a standard beginners class in ordinary differential equations and a standard beginners class in numerical methods. It is also assumed that they are familiar with programming at the level of a beginners class in informatics at the university level.

While writing these lecture notes we have been influenced by the writings of some of the many authors that previously have written testbooks on  partial differential equations \cite{Folland,John,Courant,Stakgold,Keener,Weinberger,Strauss}. However, for the students that these lecture notes are aimed at, books like \cite{Folland},\cite{John},\cite{Courant},\cite{Stakgold} and \cite{Keener} are much too advanced either mathematically \cite{Folland},\cite{John},\cite{Stakgold},\cite{Stakgold} or technically \cite{Courant},\cite{Keener}. The books \cite{Weinberger},\cite{Strauss} would be accessible to the students we have in mind, but in \cite{Weinberger} there is too much focus on precise mathematical statements and proofs and less on methods, in particular numerical methods. The book \cite{Strauss} is a better match than \cite{Weinberger} for the students these lecure notes has been written for, but it is a little superficial and does not reach far enought. With respect to reach and choise of topics, the book  ``Partial Differential Equations of Applied Mathematics'', written by Erich Zauderer\cite{textbook} would be a perfect match for the students we have mind.  However, the students taking the class covered by these notes does not have any previous exposure to PDEs, and the book \cite{textbook} is for the most part too difficult for them to follow. The book is also very voluminous and contain far to much material for a class covering one semester with five hours of lecture each week. The lecture notes  for the most part follow the structure of \cite{textbook}, but simplify the language and makes a selection of topics that can be covered in the time available. The lecture notes deviate from \cite{textbook} at several points, in particular in the section covering  physical modeling examples, integral transforms and in the treatment of numerical methods. After mastering these lecture notes, the students should be able to use \cite{textbook} as a reference for more advanced methods that are not covered by the notes. In order for students to pass the master class these notes has been written for, the students must complete the three projects included in the last section of the lecture notes and pass an oral exam where they must defend their work and also answer questions based on the content of these lecture notes.

\section{First notions}
A partial differential equation(PDE), is an equation involving one or more functions of two or more variables, and their partial derivatives, up to some finite order. Here are some examples
\begin{eqnarray}
 u_x    +u_y     &=0,&\label{linear1}\\
 u_x    +yu_y    &=0,&\label{linear2} \\
 (u_x)^2+(u_y)^2 &=1,& \label{nonlinear1}
 \end{eqnarray}
 \noindent  where
 \[u_x \equiv \partial_x u \equiv \frac{\partial u}{\partial x}.\]
The highest derivative that occur in the equation is the {\it order} of the equation. All equations (\ref{linear1},\ref{linear2},\ref{nonlinear1}) are of order one.

A PDE is {\it scalar} if it involves only one unknown function. In this course we will mostly work with scalar PDEs of two independent variables. This is hard enough as we will see!

Most important theoretical and numerical issues can be illustrated using such equations. It is also a fact that many physical systems can be modeled by  equations of this type. We will also for the most part restrict our attention on equations of order two or less. This is the type of equations that most frequently occur in modeling situations. This fact is linked to the use of second order Taylor expansions and to the use of Newton's equations in one form or another.

The most general first order scalar PDE, in two independent variables, is of the form.

\begin{equation}
F(x,y,u(x,y),u_x(x,y),u_y(x,y))=0, \label{eq1}
\end{equation}

\[(\textup{or short} \quad F(x,y,u,u_x,u_y)=0)\]
where $F$ is some function of five variables.

Let $D \in \mathbb{R}^2$ be some region in the xy-plane. Then $u(x,y)$, defined on $D$, is a {\it solution} of equation \eqref{eq1} if $u_x$,\,$u_y$ exists and are continuous on $D$ and if $u$ satisfy the equation for all $(x,y)\in D$. PDEs typically have many different solutions. The different solutions can be defined on different, even disjoint domains.

We will later relax this notion of solutions somewhat. In more advanced textbooks what we call a solution here is called a {\it classical} solution.

The most general second order scalar PDE of two independent variables is of the form

\begin{equation}\label{eq2}
F(x,y,u,u_x,u_y,u_{xx},u_{xy},u_{yy})=0.
\end{equation}
A solution in $D\in \mathbb{R}^2$ is now a function where all second order partial derivation are continuous and where the function satisfy equation \eqref{eq2} for all $(x,y)\in D$.

Here are some examples
\begin{eqnarray}
u_{xx}+u_{yy}&=0, \label{linear3}\\
u_{tt}-u_{xx}&=0, \label{linear4}\\
u_t+uu_x-u_{xx}&=0. \label{nonlinear2}
\end{eqnarray}

Equations (\ref{linear1}),(\ref{linear2}),(\ref{linear3}) and (\ref{linear4}) are {\it linear} equations whereas equations (\ref{nonlinear1}) and (\ref{nonlinear2}) are not linear. Equations that are not linear are said to be {\it nonlinear}. The importance of the linearity/nonlinearity distinction in the theory of PDEs can hardly be overstated. We have already seen the importance of this notion in the theory of ODEs.

The precise meaning of linearity is best stated using the language of linear algebra. Let us consider equation (\ref{linear1}).

Define an {\it operator} $\mathscr{L}$ by
\[\mathscr{L} = \partial_x+\partial_y.\]
Equation (\ref{linear1}) can then be written as
\[\mathscr{L}\,u=0.\]
The operator $\mathscr{L}$ has the following important properties:

\begin{description}
\item[1)] For all differentiable functions $u$,$v$
\begin{align*}
\mathscr{L}(u+v)&=\partial_x(u+v)+\partial_y(u+v)\\
             &=\partial_x u+\partial_x v+\partial_y u+\partial_y v\\
             &=\partial_x u+\partial_y u+\partial_x u+\partial_y v\\
             &=\mathscr{L}u+\mathscr{L}v,
\end{align*}
\item[2)] For all differentiable functions $u$ and $c \in \mathbb{R}$ (or $\mathbb{C})$
\begin{align*}
\mathscr{L}(cu)&=\partial_x(cu)+\partial_y(cu)\\
             &=c\partial_xu+c\partial_yu\\
             &=c(\partial_xu+\partial_yu)\\
             &=c\mathscr{L}(u).
\end{align*}
\end{description}
\noindent Properties $1$ and $2$ show us that $\mathscr{L}$ is a {\it linear} differential operator.

\begin{definition}
A scalar PDE is linear if it can be written in the form
\[\mathscr{L}(u)=g,\]
for some linear differential operator $\mathscr{L}$ and given function $g$. The equation is {\it homogeneous} if $g=0$ and {\it inhomogeneous} if $g\neq0$.
\end{definition}
\noindent Let us consider equation (\ref{nonlinear1}) from page 1. Define a differential operator $\mathscr{L}$ by
\[\mathscr{L}(u)=(\partial_xu)^2+(\partial_yu)^2.\]
Then the equation can be written as
\[\mathscr{L}(u)=1.\]
Note that
\begin{align*}
\mathscr{L}(cu)&=(\partial_x(cu))^2+(\partial_y(cu))^2\\
               &=c^2((\partial_xu)^2+(\partial_yu)^2)\\
               &=c^2\mathscr{L}(u).
\end{align*}
$\mathscr{L}$ fails property $2$ and is thus not a linear operator (It also fails property $1$). Thus the equation
\[(\partial_xu)^2+(\partial_yu)^2=1,\]
is nonlinear. It is not hard to see that, basically, a PDE is linear if it can be written as a first order polynomial in $u$ and its derivatives. All other equations are nonlinear.

The most general first order linear scalar PDE of two independent variables is of the form
\[a(x,y) u_x + b(x,y) u_y + c(x,y) u = d(x,y).\]

The most important property of homogeneous linear equations is that they satisfy the {\it  superposition principle}.

\noindent Let $u,v$ be solutions of a linear homogeneous scalar PDE. Thus 
\begin{align*}
\mathscr{L}u&=0,\\
\mathscr{L}v&=0,
\end{align*}
where $\mathscr{L}$ is the linear differential operator defining the equation.
Let $w=au+bv$ where $a,b\in\mathbb{R}$ (or $\mathbb{C}$). Then
\begin{align*}
\mathscr{L}(w)&=\mathscr{L}(au+bv)\\
              &\overset{\text{1}}{=}\mathscr{L}(au)+\mathscr{L}(bu)\\
              &\overset{\text{2}}{=}a\mathscr{L}(u)+b\mathscr{L}(u)\\
              &=0.
\end{align*}
Thus $w$ is also a solution for all choices of $a$ and $b$. This is the superposition principle. We have already seen this principle at work in the theory of ODEs.

Generalizing the argument just given, we obviously have that if $u_1,\dots,u_n$ are solutions then
\[w=\sum_{i=1}^n a_i u_i, \quad a_i \in \mathbb{R}\,(\text{or}\,\mathbb{C}),\]
is also a solution for all choices of the constants $\{a_i\}$.

In the theory of ODEs we found that any solution of a homogeneous linear ODE of order $n$ is of the form 
\[y=c_1 y_1 + \dots + c_n y_n,\]
where $y_1,\dots,y_n$ are $n$ linearly independent solutions of the ODE and $c_1,\dots,c_n$ are some constants. Thus the general solution has $n$ free {\it constants}. For linear PDEs, the situation is not that simple.

Let us consider the PDE
\begin{equation}\label{eq3}
u_{xx}=0.
\end{equation}
Integrating once with respect to $x$ we get 
\[u_x=f(y),\]
where $f(y)$ is an arbitrary function of $y$. Integrating one more time we get the general solution
\[u=xf(y)+g(y),\]
where $g(y)$ is another arbitrary function of $y$. Thus the general solution of equation (\ref{eq3}) depends on two arbitrary {\it functions}. The corresponding ODE
\[y''=0,\]
has a general solution $y=ax+b$ that depends on two arbitrary {\it constants}. The equation (\ref{eq3}) is very special but the conclusion we can derive from it is very general. (But not universal!\,)
\begin{quote}
``\,The general solution of PDEs depends on one or more arbitrary functions.''
\end{quote}
This fact makes it much harder to work with PDEs than ODEs. 

Let us consider another example that will be important later in the course.
The PDE that we want to solve is
\begin{equation}\label{eq4}
u_{xy}=0.
\end{equation}
Integrating with respect to $x$ we get 
\begin{equation}\label{eq5}
u_y=f(y),
\end{equation}
where $f(y)$ is arbitrary. Integrating equation (\ref{eq5}) with respect to $y$ gives us the general solution to equation (\ref{eq4}) in the form
\[u=F(x)+G(y),\]
where $G(y)=\int f(y) \mathrm{d}y$ and $F(x)$ are arbitrary functions of $y$ and $x$.

\section{PDEs as mathematical models of physical systems}
Methods for solving PDEs is the focus of this course. Deriving approximate description of physical systems in terms of PDEs is less of a focus and is in fact better done in specialized courses in physics, chemistry, finance etc.

It is however important to have some insight into the physical modeling context for some of the main types of equations we discuss in the course. This insight will make the methods we introduce to solve PDEs more natural and easy to understand.

\subsection{A traffic model}
We will introduce a model for cars moving on a road. In order for the complexity not to get out of hand, we start by making some assumptions about the constructions of the road.
\renewcommand\theenumi{\roman{enumi})}
\renewcommand\labelenumi{\theenumi}
\begin{enumerate}
\item There is only one lane, so that all cars move in the same direction (no passing is allowed).
\item There are no intersections or off/on ramps. Thus cars can not enter or leave the road.
\end{enumerate}
\begin{figure}[htbp]
\centering
\includegraphics{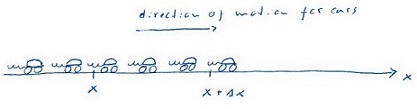}
\label{fig1}
\captionsetup{labelsep=none}
\caption{: A simple traffic model for a one-lane road}
\end{figure}
Because of \romannumeral1) and \theenumi, we can represent the road by a line that we coordinatize by a variable $x$. 
Let
\begin{align*}
\Delta N (x,t) &= \text{\# of  cars between}\,x \,\text{and}\, x+\Delta x \,\text{at time t}.\\
f(x,t) &= \text{\# of cars passing a point}\,x\,\text{from left to right},\text{per unit time, at time}\,t.
\end{align*}
From the definition of $\Delta N$ and $f$ and assumption \romannumeral2) we get the identity
\begin{equation}\label{eq6}
\Delta N(x,t+\Delta t)=\Delta N(x,t)-f(x+\Delta x, t)\Delta t +f(x,t) \Delta t.
\end{equation}
Equation \ref{eq6} express the conservation of the number of cars on the road. 

Both $\Delta N$ and $f$ are integer valued and thus change by whole numbers. After all, they are counting variables! Define the {\it density} of cars, $n(x,t)$, by
\[\Delta N(x,t)=n(x,t)\Delta x.\]
We now make the fundamental continuum assumption: $n(x,t)$ and f(x,t) are smooth functions. This is a reasonable approximation if the number of cars in each interval $\left[x,x+\Delta x\right]$ is very large and change slowly when $x$ and $t$ vary. $f(x,t)$ is also assumed to be slowly varying. (Imagine observing a crowded road from a great height, say from a helicopter).
Rewrite equation (\ref{eq6}) as
\begin{equation}\label{eq7}
\frac{n(x,t+\Delta t)- n(x,t)}{\Delta t} = -\frac{f(x+\Delta x,t)-f(x,t)}{\Delta x}.
\end{equation}
Taylor expanding to first order we get
\begin{align*}
  \frac{1}{\Delta t}&\{n(x,t)+\partial_t n(x,t)\Delta t-n(x,t)\}\\
                 &\approx-\frac{1}{\Delta x}\{f(x,t)+\partial_xf(x,t)\Delta x-f(x,t),\},\\
                 \\
                 &\Downarrow\nonumber\\
                 \\
                & \partial_t n \approx - \partial_x f.
\end{align*}
Letting $\Delta t$ and  $\Delta x $ approach zero, we get a PDE
\begin{equation}\label{eq8}
\partial_t n+\partial_x f=0.
\end{equation}
This type of PDE is called a {\it conservation law}. It express the conservation of cars in differential form. The notion of conservation laws is fundamental in the mathematical description of nature. Many of the most fundamental mathematical models of nature express conservation laws.

Note that equation (\ref{eq8}) is a single first order PDE for two unknown functions, $n$ and $f$.
In order to have a determinate mathematical model we must have the same number of equations and unknowns.

Completing the model involves making further physical assumptions. We will assume that the {\it flux}, $f$, of cars at $(x,t)$ depends on $x$ and $t$ only through the density of cars, $n(x,t)$. Thus
\begin{equation}\label{eq9}
f(x,t)=\alpha(n(x,t)),
\end{equation}
where $\alpha$ is some function whose precise specification requires further modeling assumptions. 

As noted above, the quantity $f$ in this model is called a flux. Such quantities are very common in application and typically measures the flow of some quantity through a point, surface or volume.

In equation (\ref{eq9}), we assume that the flux of cars through a point only depends on the density of cars. This is certainly not a universal law, but it is a reasonable assumption (is it?). The chain rule now gives 
\[\partial_x f=\partial_x \alpha(n) = \alpha'(n) \partial_x n.\]
Thus our model of the traffic flow is given by the PDE
\begin{equation}\label{eq10}
\partial_t n+c(n)\partial_x n =0.  \qquad(c(n)\equiv\alpha'(n))
\end{equation}
This equation describes transport of cars along a road and is an example of a {\it transport equation}.
Transport equations occur often as mathematical models because conserved quantities are common in our description of the world.

The model is still not fully specified because the function $c(n)$ is arbitrary at this point.  If we assume that $c(n)$ has a convergent power series expansion, we have
\begin{equation}\label{eq11}
c(n)=c_0+c_1n+c_2n^2+\cdots\;\;.
\end{equation}
If the density of cars is sufficiently low we can use the approximation
\begin{equation}\label{eq12}
c(n)\approx c_0,
\end{equation}
and our transport equation is a {\it linear} equation
\begin{equation}\label{eq13}
\partial_t n+c_0\partial_x n=0.
\end{equation}
For somewhat higher densities we might have to use the more precise approximation
\begin{equation}\label{eq14}
c(n)\approx c_0+c_1n,
\end{equation}
which gives us a nonlinear transport equation
\[\partial_t n+(c_0+c_1 n)\partial_x n=0.\]
For the linear transport equation (\ref{eq13}), let us look for solutions of the form
\[n(x,t)=\varphi(x-c_0t).\]
The chain rule gives
\begin{align*}
   \partial_t n&=-c_0 \varphi',\\
   \partial_x n&=\varphi',\\
               &\Downarrow\\
   \partial_t n+c_0 \partial_x n&=-c_0\varphi'+c_0\varphi'=0.
\end{align*}
Thus
\begin{equation}\label{eq15}
n(x,t)=\varphi(x-c_0t),
\end{equation}
is a solution of equation (\ref{eq13}) for any function $\varphi$. Let us assume that the density of cars is known at one time, say at $t=0$. Thus
\begin{equation}\label{eq16}
n(x,0)=f(x),
\end{equation}
where $f$ is a given function.
Using equation (\ref{eq16}) and (\ref{eq15}) we have
\[n(x,0)=f(x),\]
\[\Updownarrow\]
\[\varphi(x-c_0\cdot0)=f(x),\]
\[\Updownarrow\]
\[\varphi(x)=f(x).\]
Thus the function
\begin{equation}\label{eq17}
n(x,t)=f(x-c_0 t),
\end{equation}
is a solution of
\begin{eqnarray}
&\partial_t n+c_0 \partial_x n =0,\nonumber\\
&n(x,0)=f(x).\label{eq18}
\end{eqnarray}
(\ref{eq18}) is an example of an initial value problem for a PDE. Initial value problems for PDEs are often, for historical reasons, called {\it Cauchy problems}.

Note that at this point we have not shown that (\ref{eq17}) is the only solution to the Cauchy problem (\ref{eq18}). Proving that there are no other solutions is obviously important because a deterministic model (and this is!) should give us exactly one solution, not several.

The solution (\ref{eq17}) describe an interesting behavior that is common to {\it very} many physical systems. It describe a {\it wave}. 
From (\ref{eq17}), it is evident that
\[ n(x,t)=n(x-c_0t,0).\]
Thus we get the values of $n$ at time $t$ by translating the values of $n$ at $t=0$ along the $x$-axis.
\begin{figure}[h]
\centering
\includegraphics{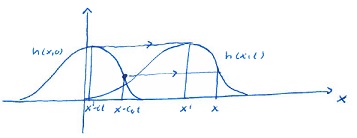}
\label{fig2}
\captionsetup{labelsep=none}
\caption{: A traffic disturbance propagating like a wave.}
\end{figure}

\noindent We are used to observing waves when looking at the surface of the sea, but might not have expected to find it in the description of traffic on a one-lane road.

Traffic waves are in fact common and play no small part in the occurrence of traffic jams.

\subsection{Diffusion}
Consider a narrow pipe oriented along the $x$-axis containing water that is not moving. Assume there is a material substance, a dye for example, of mass density $u(x,t)$ dissolved in the water. Here we assume that the mass density is uniform across the pipe so that the density only depends on $x$ which by choice is the coordinate along the length of the pipe.
\begin{figure}[h]
\centering
\includegraphics{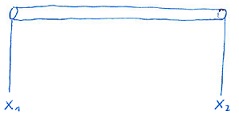}
\label{fig8}
\captionsetup{labelsep=none}
\caption{}
\end{figure}

\noindent The total mass of the dye between $x_1$ and $x_2$ is 
\[M(t)=\int_{x_1}^{x_2}u(x,t)\mathrm{d}x.\]
If the pipe is not leaking, the mass inside the section of the pipe between $x_1$ and $x_2$ can only change by dye entering or leaving through the points $x_1$ and $x_2$. This is the law of {\it conservation of mass} which is the first great pillar of classical physics.

Let $f(x,t)$ be the amount of mass flowing through a point $x$ at time $t$ pr.unit time. This is the {\it mass flux}. By convention $f$ positive means that mass is flowing to the right.

Mass conservation then gives the balance equation
\begin{equation}\label{eq22}
\frac{\mathrm{d}M}{\mathrm{d}t}=f_1-f_2,\quad f_i=f(x_i,t),\,i=1,2,\cdots\;\; .
\end{equation}
{\it Fick's law} states that we have the following relation between mass flux and mass density.
\begin{equation}\label{eq23}
f(x,t)=-ku_x(x,t),\quad k>0
\end{equation}
Thus mass flow from high density areas towards low density areas.

  Fick's law is only a phenomenological approximate relation based on the empirical observations that mass flux is for the most part a function of $(x,t)$, thorugh the density gradient, thus $f(x,t)=f(u_x(x,t))$. Assuming that the flux is a smooth function of the density gradient, we can represent the flux as a power series in the density gradient
\begin{equation}\label{eq24}
f(u_x(x,t))=f_0+f_1u_x(x,t)+f_2(u_x(x,t))^2+\cdots
\end{equation}
Since there is no mass flux if the density gradient is zero we must have $f_0=0$. Furthermore, since mass flows from high density domains to low density domains we must have $f_1<0$. We can therefore write $f_1=-k$, where $k>0$. Clearly, if the density gradient is large enough, the third term in the power series \ref{eq24} will be significant and must be included. Like for the traffic model, this will lead to a nonlinear equation instead of the linear equation that will be the result of truncating the power series at the second term.
In general we will assume that the parameter $k$ in Fick's law depends on $x$. This corresponds to the assumption that the relation between flux and density gradient depends on position. This is not uncommon. Using Fick's law we can write the mass balance equation as
\[\left.\frac{\mathrm{d}M}{\mathrm{d}t}=ku_x\right|_{x_1}^{x_2},\]
\[\Updownarrow\]
\[\int_{x_1}^{x_2}u_t\mathrm{d}x=\int_{x_1}^{x_2}(ku_x)_x\mathrm{d}x,\]
\[\Updownarrow\]
\[\int_{x_1}^{x_2}(u_t-(ku_x)_x)\mathrm{d}x=0,\quad \forall x_1,x_2\]
\[\Updownarrow\]
\begin{equation}\label{eq25}
u_t-(ku_x)_x=0.
\end{equation}
This is the 1D {\it diffusion equation}. By using similar modeling in 2D and 3D, we get the 2D and 3D diffusion equations
\begin{equation*}
u_t=\nabla\cdot (k\nabla u) \quad
\begin{cases}
\nabla=(\partial_{x},\partial_{y}) &2D\\
\nabla^2=(\partial_{x},\partial_{y},\partial_{z}) &3D
\end{cases}.
\end{equation*}
If the mass flow coefficient, $k$, is constant in space, the diffusion equation takes the form 
\begin{equation}\label{eq26}
u_t=k\nabla^2u.
\end{equation}

These equations and their close cousins describe a huge range of physical, chemical, biological and economic phenomena.

If the mass flow is stationary(time independent) and the mass flow coefficient independent of space, we get the simplified equation
\begin{equation}\label{eq27}
\nabla^2 u=0, \quad u=u(\bold x).
\end{equation}
This is the {\it Laplace equation} and occur in the description of many stationary phenomena.

\subsection{Heat conduction}
Let us consider a narrow bar of some material, say a metal. We assume that the temperature of the bar, $T$, only depends on the coordinate along the bar which we call $x$.
\begin{figure}[h]
\centering
\includegraphics{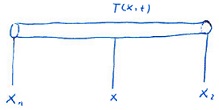}
\label{fig9}
\captionsetup{labelsep=none}
\caption{}
\end{figure}

\noindent Let the temperature of the bar at a point $x$ at time $t$ be $T(x,t)$. Let $e(x,t)$ be the energy density at $x$ at time $t$. The total energy in the bar between $x_1$ and $x_2$ is
\[E(t)=\int_{x_1}^{x_2}e(x,t)\mathrm{d}x.\]
The second great pillar of classical physics is the law of {\it conservation of energy}.
Assuming no energy is gained or lost along the length of the bar, the conservation of energy gives us the balance equation
\begin{equation}\label{eq28}
\frac{\mathrm{d}E}{\mathrm{d}t}=f_1-f_2,
\end{equation}
where $f_i=f_i(x_i,t)$ and $f(x,t)$ is the {\it energy flux} at $x$ at time $t$. Here we use here the convention that positive energy flux means that energy flow to the right. In order to get a model for the bar, we must relate $e(x,t)$ and $f(x,t)$ to $T(x,t)$.

For many materials we have the following approximate identity
\begin{equation}\label{eq29}
e(x,t)=C\rho T(x,t),
\end{equation}
where $\rho$ is the mass density of the bar and $C$ is the heat capacity at constant volume. For (\ref{eq29}) to hold the temperature must not vary too quickly.
Energy is always observed to flow from high to low temperature regions. If the temperature gradient is not too large, {\it Fourier's law} hold for the energy flux
\begin{equation}\label{eq30}
f(x,t)=-q T_x, \quad q>0
\end{equation}
The energy balance equation then become
\[\frac{\mathrm{d}E}{\mathrm{d}t}=f_1-f_2,\]
\[\Updownarrow\]
\[\left.\int_{x_1}^{x_2}C\rho T_t \mathrm{d}x=q T_x\right|_{x_1}^{x_2}=\int_{x_1}^{x_2}(q T_x)_x\mathrm{d}x,\]
\[\Updownarrow\]
\[\int_{x_1}^{x_2}(C\rho T_t-(q T_x)_x)\mathrm{d}x=0, \quad \forall x_1,x_2\]
\[\Updownarrow\]
\[C\rho T_t-(q T_x)_x=0.\]
If the heat conduction coefficient, $q$, does not depend on $x$, we get
\begin{equation}\label{eq31}
T_t=D T_{xx},\quad D=\frac{q}{C\rho},
\end{equation}
Again we get the diffusion equation. In this context the equation is called the {\it equation of heat conduction}. Energy flow in 2D and 3D similarly leads to 2D and 3D version of the heat conduction equations
\begin{equation}\label{eq32}
T_t=D \nabla^2 T,
\end{equation}
and for stationary heat flow, $T_t=0$, we get the Laplace equation
\begin{equation}\label{eq33}
\nabla^2 T=0.
\end{equation}
The wave equation, diffusion equation and Laplace equations and their close cousins will be the main focus in this course.

In the next example we will show that the diffusion equation also arise from a seemingly different context. This example is the beginning of a huge application that have ramifications in all areas of science, technology, economy, etc. I am talking about {\it stochastic processes}.

\subsection{A random walk}
Let $x$ be the position of a particle restricted to moving along a line. We assume that the particle, if left alone, will not move. At random intervals we assume that the particle experience a collision that makes it shift its position by a small amount $\delta$. The collision is equally likely to occur from the left or the right. Thus the shifts in position caused by a collision is equally likely to be $-\delta$ or $+\delta$. The collisions occur at random times, but let us assume that during some small time period, $\tau$, a collision will occur for sure.

Let $P(x,t)$ be the probability of finding the particle at a point $x$ at time $t$. We must then have
\begin{equation}\label{eq34}
P(x,t+\tau)=\frac12 P(x-\delta,t)+ \frac12 P(x+\delta,t).
\end{equation}
This is just saying that {\it probability is conserved}. From this point of view (\ref{eq34}) is a balance equation derived from a conservation law just like equations (\ref{eq5}),\,(\ref{eq19}),\,(\ref{eq22}) and (\ref{eq28}).
Using Taylor expansions
\begin{align*}
P(x,t+\tau) &=P(x,t)+\partial_t P(x,t) \tau + O(\tau^2),\\
P(x\pm \delta,t) &=P(x,t)\pm \partial_x P(x,t) \delta\\
                 &+\frac12 \partial_{xx}P(x,t)\delta^2 +O(\delta^3),
\end{align*}
we get from equation (\ref{eq34})       
\begin{align*}
&\quad \, \,P(x,t)+\partial_t P \tau+O(\tau^2)\\
&=\frac12 P(x,t)-\frac12 \partial_x P \delta+\frac14 \partial_{xx}P\delta^2+\frac12 P(x,t)\\
&+\frac12\partial_x P \delta+\frac14 \partial_{xx}P \delta^2+ O(\delta^4).
\end{align*}
Thus by truncating the expansions we get
\[\partial_t P=D \partial_{xx}P, \quad D=\frac{\delta^2}{2 \tau}.\]
This is again the diffusion equation. A more general and careful derivation of the random walk equation can be found in section 1.1 in \cite{textbook}.

\subsection{The vibrating string}
Let us consider a piece of string of length $l$ that is fixed at $x=0$ and $x=l$. We make the following modeling assumption
\begin{enumerate}
\item The string can only move in a single plane, and in this plane it can only move transversally.
\end{enumerate}
Using \romannumeral1) the state of the string can be described by a single function $u(x,t)$ measuring the height of the deformed string over the $x$-axis at time $t$. Negative values of $u$ means that the string is below the $x$-axis.
\begin{figure}[h]
\centering
\includegraphics{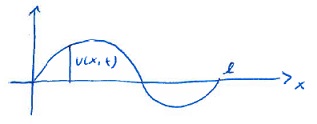}
\label{fig3}
\captionsetup{labelsep=none}
\caption{}
\end{figure}
\begin{enumerate}
\addtocounter{enumi}{1}
\item Any real string has a finite thickness, and bending of the string result in a tension (negative pressure) on the outer edge of the string and a compression(pressure) on the inner edge of the bend.
\end{enumerate}
\begin{figure}[h]
\centering
\includegraphics{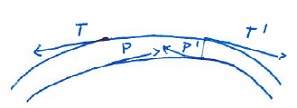}
\label{fig4}
\captionsetup{labelsep=none}
\caption{}
\end{figure}
We model the string as infinitely thin and thus assume no compression force. The tension, $\mathbf T$, in a deformed string is always tangential to the string and assumed to be of constant length and time independent
\[T=|\mathbf T|=\text{const}.\]
\begin{enumerate}
\addtocounter{enumi}{2}
\item The string has a mass density given by $\rho$ which is independent of position and  also time independent.
\end{enumerate}
\begin{description}
\item[Note:] \romannumeral2) and \romannumeral3) are not totally independent because a typical cause of nonuniform tension, $T$, is that the string has a variable material composition and this usually imply varying mass density.
\end{description}
\begin{figure}[h]
\centering
\includegraphics{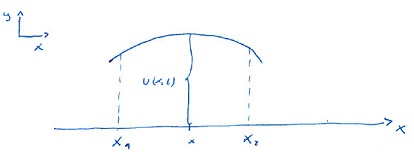}
\label{fig5}
\captionsetup{labelsep=none}
\caption{: A piece of the string above the interval $[x_1,x_2]$.}
\end{figure}
Let us consider a piece of string above the interval $[x_1,x_2]$. The total vertical momentum (horizontal momentum is by assumption \romannumeral1) equal to zero), $P(t)$, in this piece of string at time $t$ is 
\begin{equation}\label{eq19}
P(t)=\int_{x_1}^{x_2}\mathrm{d}x \rho u_t.
\end{equation}
Recall that momentum is conserved in classical (nonrelativistic) physics. This is the third great pillar of classical physics.
Note that with the $y$-axis oriented like in figure 5, a piece of the string moving vertically up gives a positive contribution to the total momentum $P(t)$.
Recall also that force is by definition change in momentum per unit time. Thus, if we assume there are no forces acting on the string inside the interval $[x_1,x_2]$, then momentum inside the interval can only change through the action of the tension force at the bounding points $x_1$ and $x_2$. If we orient the coordinated axes as in figure 5, a positive vertical component, $f_1$,  of the tension force at $x_1$, will lead to an increase the total momentum $p(t)$. A similar convention apply for $f_2$ at $x_2$. 

Conservation of momentum for the string is then expressed by the following balance equation
\[\frac{\mathrm{d}P}{\mathrm{d}t}=f_1+f_2.\]
\begin{figure}[h]
\begin{minipage}[b]{.5\textwidth}
\centering
\includegraphics{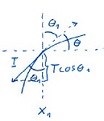}
\caption{}
\end{minipage}%
\begin{minipage}[b][4cm][c]{.5\textwidth}
\centering
\begin{align*}
f_1&=-T \cos \theta_1=-T\cos(\frac{\pi}{2}-\theta)\\
   &=-T\sin \theta \\
   &=-T \left.\frac{u_x}{\sqrt{1+u_x^2}}\right|_{x=x_1}.
\end{align*}
\end{minipage}
\end{figure}

\begin{figure}[h]
\begin{minipage}[b]{.5\textwidth}
\centering
\includegraphics{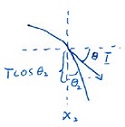}
\caption{}
\end{minipage}%
\begin{minipage}[b][4cm][c]{.5\textwidth}
\centering
\begin{align*}
f_2&=-T \cos \theta_2=-T\cos(\frac{\pi}{2}+\theta)\\
   &=T\sin \theta \\
   &=T \left.\frac{u_x}{\sqrt{1+u_x^2}}\right|_{x=x_2}.
\end{align*}
\end{minipage}
\end{figure}

\noindent Our last modelling assumption is
\begin{enumerate}
\addtocounter{enumi}{3}
\item The string perform small oscillations
\[|u_x|\ll 1.\]
\end{enumerate}
Applying \romannumeral4) gives us
\begin{equation}\label{eq20}
\begin{split}
f_1 &\approx -Tu_x\Big|_{x=x_1},\\
f_2 &\approx Tu_x\Big|_{x=x_2}.
\end{split}
\end{equation}
Thus our fundamental momentum balance equation is well approximated by
\[\frac{\mathrm{d}P}{\mathrm{d}t}=-Tu_x\Big|_{x_1}+Tu_x\Big|_{x_2}=Tu_x\Big|_{x_1}^{x_2},\]
\[\Updownarrow\]
\[\int_{x_1}^{x_2}(\rho u_t)_t\mathrm{d}x = \int_{x_1}^{x_2}(Tu_x)_x\mathrm{d}x,\]
\[\Updownarrow\]
\[\int_{x_1}^{x_2}(\rho u_{tt}-Tu_{xx})\mathrm{d}x=0, \qquad \forall x_1,x_2\]
\[\Updownarrow\]
\[\rho u_{tt}-Tu_{xx}=0,\]
\[\Updownarrow\]
\begin{equation}\label{eq21}
  u_{tt}-c^2u_{xx}=0, \qquad c=\sqrt{\frac T\rho} >0.
\end{equation}
This is the {\it wave equation}. Similar modeling of small vibrations of membranes (2D) and solids\,(3D) give wave equations
\begin{equation*}
u_{tt}-c^2\nabla^2u=0\quad
\begin{cases}
\nabla^2 =\partial_{xx}+\partial_{yy}&2D\\
\nabla^2 =\partial_{xx}+\partial_{yy}+\partial_{zz}&3D
\end{cases}.
\end{equation*}
Wave equations describe small vibrations and wave phenomena in all sorts of physical systems, from stars to molecules.

\section{Initial and boundary conditions}
PDEs typically have a very large number of solutions (free functions, remember?). We usually select a unique solution by specifying extra conditions that the solutions need to satisfy. Picking the right conditions is part of the modeling and will be guided by the physical system that is under consideration.

There are typically two types of such conditions:
\begin{enumerate}
\item[\romannumeral1)] initial conditions
\item[\romannumeral2)] boundary conditions
\end{enumerate}
An initial condition consists of specifying the state of the system at some particular time $t_0$. For the diffusion equation, we would specify
\[u(\mathbf{x}, t_0)=\varphi(\mathbf{x}),\]
for some given function $\varphi$. This would for example describe the distribution of dye at time  $t_0$. For the wave equation, we need {\it two} initial conditions since the equation is of second order in time (same as for ODEs)
\begin{align*}
u(\mathbf{x},t_0)&=\varphi(\mathbf{x}),\\
\partial_tu(\mathbf{x},t_0)&=\varphi(\bold{x}).
\end{align*}
The modelling typically leads us to consider the PDE on some domain $D$ (bounded or unbounded).

For the vibrating string it is natural to assume that the string has a finite length, $l$, so $D=[0,l]$. It is then intuitively obvious that in order to get a unique behavior we must specify what the conditions on the string are at the endpoints $x=0$ and $x=l$. 
After all, a clamped string certainly behave different from a string that is loose at one or both ends! Here are a couple  of natural conditions for the string:
\begin{enumerate}
\item[\romannumeral1)] $u(0,t)=0, u(l,t)=0$ ,
\item[\romannumeral2)] $u(0,t)=g(t), u_x(l,t)=b(t)$.
\end{enumerate}
For case i) the string is clamped at both ends and for case ii) the string is moved in a specified way at $x=0$, and is acted upon by a force $Tb(t)$ at the right end, $x=l$.

The three types of boundary conditions that most often occur are
\begin{enumerate}
\item[D)] $u$ specified at the boundary.
\item[N)] $\partial_{\mathbf{n}}u$ specified at the boundary.
\item[R)] $\partial_{\mathbf{n}}u+au$ specified at the boundary.
\end{enumerate}
The first is a {\it Dirichlet condition}, the second a {\it Neumann condition} and the third a {\it Robin condition}.
Here $\mathbf{n}$ is a normal to the curve (2D) or surface (3D) that defines the boundary of our domain and $a$ is a function defined on the boundary.
\captionsetup{labelsep=none}
\begin{figure}[h]
\begin{minipage}[b]{.5\textwidth}
\centering
\includegraphics{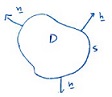}
\caption{}
\end{minipage}
\begin{minipage}[b][3cm][c]{.5\textwidth}
$\mathbf{n}$ is by convention the {\it outward pointing} normal.
\end{minipage}
\end{figure}
Let us consider the temperature profile of a piece of material covering a domain $D$.
\noindent We want to model heat conduction in this material.

\begin{enumerate}
\item[\romannumeral1)] Assume that the material is submerged in a large body of water kept at a temperature $T=0$. The problem we must solve is then
\begin{align*}
\partial_t u=D\nabla^2u,\qquad &\text{in}\,D,\\
u(\mathbf{x},0)=\varphi(\mathbf{x}),\qquad & \mathbf{x}\in D,\\
u(\mathbf{x},t)=0, \qquad & \forall\,\mathbf{x}\in \partial D,
\end{align*}
where $\varphi(\mathbf{x})$ is the initial temperature profile.
\item[\romannumeral2)] If the material is perfectly isolated (no heat loss), we must solve the problem
\begin{align*}
\partial_t u=D\nabla^2u,\qquad&\text{in } D,\\
u(\mathbf{x},0)=\varphi(\mathbf{x}), \qquad&\mathbf{x}\in D,\\
\partial_{\mathbf{n}}u(\mathbf{x},t)=0,\qquad& \forall\,\mathbf{x}\in \partial D.
\end{align*}
\end{enumerate}
The general morale here is that the choice of initial and boundary conditions are part of the modelling process just like the equation themselves.

\subsection{Well posed Problems}
Let us assume that the modelling of some physical system has produced a PDE, initial conditions and/or boundary conditions.

\noindent How can we know whether or not the model is any good? After all, the model is always a caricature of the real system. There are always many simplifying assumptions and approximations involved.

We could of course try to solve it and compare solutions to actual observation of the system. However, if this should even make sense to try, there are some basic requirements that {\it must} be satisfied.
\renewcommand\theenumi{\roman{enumi})}
\renewcommand\labelenumi{\theenumi}
\begin{enumerate}
\item Existence:  There is at least one solution that satisfy all the requirements of the model.
\item Uniqueness:  There exists no more than one solution that satisfy all the requirements of the model.
\item Stability. The unique solution depends in a stable manner on the data of the problem: Small changes in data (boundary conditions, initial conditions, parameters in equations, etc, etc) lead to small changes in the solution.
\end{enumerate}
Producing a formal mathematical proof that a given model is well posed can be a hard problem and belongs to the area of pure mathematics. In applied mathematics we derive the models using our best physical intuition and sound methods of approximation and assume as a ``working hypothesis'' that they are well posed. If they are not, the models will eventually break down and in this way let us know. The way they break down is often a hint we can use to modify them such that the modified models {\it are} well posed.

We thus use an ``empirical'' approach to the question of well-posedness in applied mathematics. We will return to this topic later in the class.

\section{Numerical methods}
In most realistic modelling, the equations and/or the geometry of the domains are too complicated for solution ``by hand'' to be possible. In the ``old days'' this would mean that the problems were unsolvable. Today, the ever increasing memory and power of computers offer a way to solve previously unsolvable problems. Both exact methods and approximate numerical methods are available, often combined in easy to use systems, like Matlab, Mathematica, Maple etc.

In this section we will discuss numerical methods. The use of numerical methods to solve PDEs does not make more theoretical methods obsolete. On the contrary, heavy use of mathematics is often required to reformulate the PDEs and write them in a form that is well suited for the architecture of the machine we will be running it on, which could be serial, shared memory, cluster etc.

Solving PDEs numerically has two sources of error.
\begin{enumerate}
\item[1)] Truncation error: This is an error that arises when the partial differential equation is approximated by some (multidimensional) difference equation. Continuous space and time can not be represented in a computer and as a consequence there will {\it always} be truncation errors. The challenge is to make them as small as possible and at the same time stay within the available computational resources.
\item[2)] Rounding errors: Any computer memory is finite and therefore only a finite set of real numbers can be represented. After any mathematical operation on the available numbers we must select the ``closest possible'' 
available number to represent the result. In this way we accumulate errors after every operation. In order for the accumulated errors not to swamp the result we are trying to compute, we must be careful about {\it how} the mathematical operations are performed and {\it how many} operations we perform.
\end{enumerate}
1)\,and 2)\,apply to {\it all} possible numerical methods. Heavy use of mathematical methods are often necessary to build confidence that some numerical solution is an accurate representation of the true mathematical solution. The accuracy of numerical solutions are often tested by running them against known exact solutions to special initial and/or boundary conditions and to ``stripped down'' versions of the equations.

The simples truncation method to describe is {\it the method of finite differences}. This is the only truncation method we will discuss in this class.  We will introduce the method by applying it to the 1D heat equation, the 1D wave equation and the 2D Laplace equation

\subsection{The finite difference method for the heat equation}
Consider the following initial/boundary value problem for the heat equation
\begin{align*}
u_t&=u_{xx}, \qquad \qquad 0<x<1,t>0,\\
u(0,t)&=u(1,t)=0,\\
u(x,0)&=\varphi(x).
\end{align*}
We introduce a uniform spacetime grid $(x_i,t_n)$, where 
\begin{align*}
x_i&=hi,\quad i=0,\dots,N,\\
t_n&=k n, \qquad n\in \mathbb{N},
\end{align*}
where $x_0=0,x_N=1$ and thus $h=\frac{1}{N}$.
The solution will only be computed at the grid points.
\begin{figure}[h]
\begin{centering}
\includegraphics{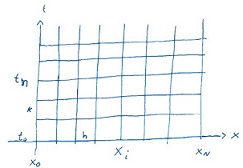}
\end{centering}
\caption{}
\end{figure}

 Evaluating the equation and the initial/boundary conditions at the grid points we find
\begin{equation}\label{eq35}
\begin{split}
(u_t)_i^n=(u_{xx})_i^n,\\
u_0^n=0,\, u_N^n=0,\\ 
u_i^0=\varphi_i=\varphi(x_i),
\end{split}
\end{equation}
where $u_i^n=u(x_i,t_n)$.

We now must find expressions for $u_t$ and $u_{xx}$, using {\it only} function values at the grid points.
Our main tool for this is Taylors theorem.
\begin{align*}
u(x_i,t_n+k)&=u_i^n+\left(\frac{\partial u}{\partial t}\right)_i^nk+\frac12\left(\frac{\partial^2u}{\partial t^2}\right)_i^n k^2+O(k^3),\\
u(x_i,t_n-k)&=u_i^n-\left(\frac{\partial u}{\partial t}\right)_i^nk+\frac12\left(\frac{\partial^2u}{\partial t^2}\right)_i^n k^2+O(k^3).
\end{align*}
Using the Taylor expansions we easily get the following expressions
\begin{align*}
\frac{u(x_i,t_n+k)-u_i^n}{k}&=\left(\frac{\partial u}{\partial t}\right)_i^n+O(k),\\\label{eq35-2}
\frac{u_i^n-u(x_i,t_n-k)}{k}&=\left(\frac{\partial u}{\partial t}\right)_i^n+O(k)),\\
\frac{u(x_i,t_n+k)-u(x_i,t_n-k)}{2k}&=\left(\frac{\partial u}{\partial t}\right)_i^n+O(k^2)).\\
\end{align*}
Since $t_n\pm k=t_{n\pm1}$, we get the following three approximations to the first derivative with respect time 

\begin{equation}\label{eq36}
\begin{split}
 &\left.\begin{aligned}
       \left(\frac{\partial u}{\partial t}\right)_i^n&\approx \frac{u_i^{n+1}-u_i^n}{k} \quad \text{forward difference}\\
        \left(\frac{\partial u}{\partial t}\right)_i^n&\approx \frac{u_i^n-u_i^{n-1}}{k} \quad \text{backward difference}
        \end{aligned}
 \right\}
 \text{approx. has an error of order } k,\\
  & \!\!\!\left(\frac{\partial u}{\partial t}\right)_i^n\approx\frac{u_i^{n+1}-u_i^{n-1}}{2k} \quad\text{center difference}\qquad \,\,\,\text{order }k^2.
\end{split}
\end{equation}
The second derivative in space is also approximated using Taylor theorem
\begin{align*}
u(x_i+h,t_n)&=u_i^n+\left(\frac{\partial u}{\partial x}\right)_i^nh+\frac12\left(\frac{\partial^2u}{\partial x^2}\right)_i^n h^2+\frac16\left(\frac{\partial^3u}{\partial x^3}\right)_i^n h^3+O(h^4),\\
u(x_i-h,t_n)&=u_i^n-\left(\frac{\partial u}{\partial x}\right)_i^nh+\frac12\left(\frac{\partial^2u}{\partial x^2}\right)_i^n h^2-\frac16\left(\frac{\partial^3u}{\partial x^3}\right)_i^n h^3+O(h^4).
\end{align*}
From these two expressions we get
\begin{align*}
\frac{u(x_i+h,t_n)-2u_i^n+u(x_i-h,t_n)}{h^2}&=\left(\frac{\partial^2u}{\partial x^2}\right)_i^n +O(h^2).\\
\end{align*}
Thus we get the following approximation to the second derivative in space
\begin{equation}\label{eq37}
\left(\frac{\partial^2u}{\partial x^2}\right)_i^n\approx\frac{u_{i+1}^n-2u_i^n+u_{i-1}^n}{h^2}\quad(\text{center difference}) \qquad \text{order } h^2.
\end{equation}
Using the forward difference in time and center difference in space in (\ref{eq35}) gives
\begin{equation*}
\frac{u_i^{n+1}-u_i^n}{k}=\frac{u_{i+1}^n-2u_i^n+u_{i-1}^n}{h^2}.
\end{equation*}

\noindent Define $s=\frac{k}{h^2}$. Then the problem we must solve is
\begin{equation}\label{eq38}
\begin{split}
u_i^{n+1}&=s(u_{i+1}^n+u_{i-1}^n)+(1-2s)u_i^n,\quad i=1,\dots,N-1,\\
u_0^n&=u_N^n=0,\\
u_i^0&=\varphi_i.  
\end{split}
\end{equation}
This is now a (partial) difference equation.
\begin{figure}[h]
\centering
\includegraphics{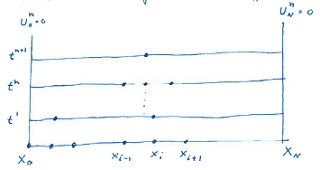}
\caption{}
\end{figure}

\noindent Let us use the initial condition
\begin{equation*}
\varphi(x)=\begin{cases}
2x, &0<x\leqslant\frac12\\
2(1-x), &\frac12\leqslant x<1
\end{cases}.
\end{equation*}
The exact solution is found by separation of variables which we will discuss later in this class.
\[u(x,t)=\frac8{\pi^2}\sum_{m=1}^{\infty}\frac1{m^2}\sin \frac12 m \pi\hspace{1 mm}  e^{-m^2\pi^2t}\sin m \pi x. \]
In order to compute the numerical solution, we must choose values for the spacetime discretization parameters $h$ and $k$.
In the figures 13 and 14, we plot the initial condition, the exact solution and the numerical solution for $h=0.05$, thus $N=50$, and $k=h^2s=0.0004s$ for $s=0.49$ and $0.51$. We use 100 iterations in $n$.

It is clear from the figure 14 that the numerical solution for $s=0.51$ is very bad. This is an example of a {\it numerical instability}. Numerical instability is a constant danger when we compute numerical solutions to PDEs.

\begin{figure}
\centering
\includegraphics{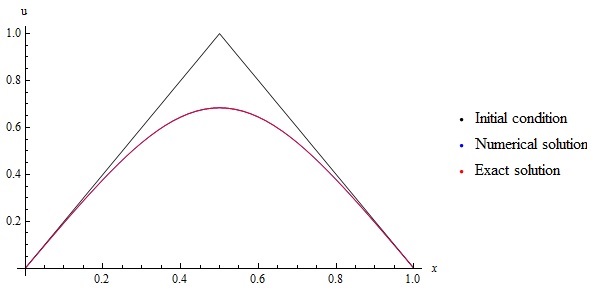}
\caption{: In this figure we have s=0.49}
\label{fig3-1}
\end{figure}
\begin{figure}
\centering
\includegraphics{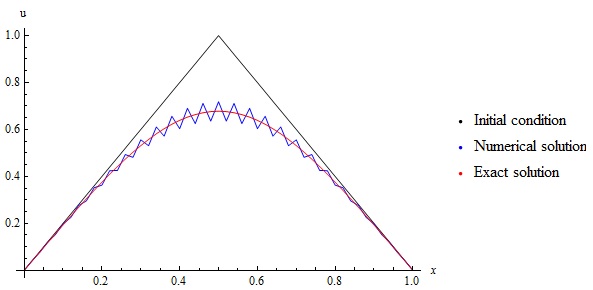}
\caption{: In this figure s=0.51}
\label{fig3-2}
\end{figure}

For a simple equation like the heat equation we can understand this instability, why it appears and how to protect against it.

Let us try to solve the difference equation in (\ref{eq38}) by separating variables. We thus seek a solution of the form
\[u_j^n=X_j T_n,\]
\[\Downarrow\]
\[X_jT_{n+1}=s(X_{j+1}T_n+X_{j-1}T_n)+(1-2s)X_jT_n,\]
\[\Downarrow\]
\[\frac{T_{n+1}}{T_n}=(1-2s)+s\,\frac{X_{j+1}+X_{j-1}}{X_j}.\]
Both sides must be equal to the same constant $\xi$. This gives us two difference equations coupled through the separation constant $\xi$.
\begin{align*}
&1)\quad T_{n+1}=\xi T_n.\\
&2)\quad s\,\frac{X_{j+1}+X_{j-1}}{X_j}+(1-2s)=\xi.
\end{align*}
The solution of equation 1) is simply
\[T_n=\xi^n T_0.\]
The second equation is a boundary value problem
\begin{equation}\label{eq3-1}
s(X_{j+1}+X_{j-1})+(1-2s)X_j=\xi X_j,
\end{equation}
\[X_0=X_N=0,\]
We look for solutions of the form $X_j=z^j$. This implies that 
\begin{align*}
X_{j-1}&=z^{-1}X_j,\\
X_{j+1}&=zX_j,
\end{align*}
Thus (\ref{eq3-1}) implies that
\[\xi=s(z+z^{-1})+(1-2s),\]
\[\Updownarrow\]
\begin{equation}\label{eq3-2}
z^2+(\frac1s-2-\frac{\xi}{s})z+1=0. 
\end{equation}
This equation has two solutions $z_1$ and $z_2$ which we in the following assume to be different.  (Show that the case when $z_1=z_2$ does not lead to a solution of the boundary value problem (\ref{eq3-1}).)

  The general solution of (\ref{eq3-1}) is thus
\begin{equation*}
X_j=Az_1^j+Bz_2^j.
\end{equation*}
The general solution must satisfy the boundary conditions. This implies that
\begin{align*}
X_0=0  &\Leftrightarrow  A+B=0,\\
X_N=0&\Leftrightarrow Az_1^N+Bz_2^N=0.
\end{align*}
Nontrivial solutions exists only if
\begin{equation}\label{eq3-3}
\left(\frac{z_1}{z_2}\right)^N=1.  
\end{equation}
This condition can not be satisfied if $z_1\neq z_2$ are real. Since equation (\ref{eq3-2}) has real coefficient, we must have
\begin{align*}
z_1&=z=\rho e^{i\theta},\\
z_2&=\bar{z}=\rho e^{-i\theta}.
\end{align*}
But equation (\ref{eq3-2}) $\Rightarrow$ $z_1z_2=1$ $\Leftrightarrow$ $\rho^2=1$ $\Leftrightarrow$ $\rho=1$. 
Thus $z_1=e^{i\theta}$, $z_2=e^{-i\theta}$. The condition (\ref{eq3-3}) then implies that
\[e^{2i\theta N}=1,\quad \Leftrightarrow \quad \theta_r=\frac{\pi r}{N}, \quad r=1,2,\dots,N-1.\]
The corresponding solutions of the boundary value problem (\ref{eq3-1}) are
\[X_j^{(r)}=A(e^{i\theta_rj}-e^{-i\theta_rj})=2iA\sin \theta_r j.\]
A basis of solutions  can be choosen to be
\[X_j^{(r)}=\sin \frac{r \pi j}{N},\quad r=1,2,\dots, N-1.\]
It is easy to verify that 
\[c_1X_j^{(1)}+\dots+c_{N-1}X_j^{(N-1)}=0,\quad \forall j\]
\[\Downarrow\]
\[c_1=c_2=\dots=c_{N-1}=0.\]
The vectors $\{X_j^{(r)}\}_{r=1}^{N-1}$ are thus  linearly independent and indeed form a basis because $\mathbb{R}^{N-1}$ has dimension $N-1$.

Each vector $X_j^{(r)}$  gives a separated solution.
\[{u_r}_j^n=\xi_r^n\sin \frac{r\pi j}{N},\]
to the boundary value problem
\begin{equation*}
\begin{split}
u_i^{n+1}&=s(u_{i+1}^n+u_{i-1}^n)+(1-2s)u_i^n,\quad i=1,\dots,N-1,\\
u_0^n&=u_N^n=0.
\end{split}
\end{equation*}
where we have defined
\[\xi_r=(1-2s)+s(e^{i\theta_r}+e^{-i\theta_r})=1-2s(1-\cos \frac{r\pi}{N}).\]
These solutions give exponential growth if $|\xi_r|>1$ for some $r$. Thus the requirement for {\it stability} is
\[|\xi_r|\leqslant1, \qquad \forall r.\]
If this holds,  the numerical scheme is stable.

  Now
\[-1\leqslant\cos \frac{r\pi}{N}\leqslant1,\quad \forall\,r,\]
\[\Updownarrow\]
\[0\leqslant1-\cos \frac{r\pi}{N}\leqslant2,\quad \forall\,r.\]
Therefore $-1\leqslant \xi^r\leqslant1$ if
\[1-4s\eqslantgtr-1,\]
\[\Updownarrow\]
\[s\leqslant\frac12.\]
This is the stability condition indicated by the numerical test illustrated in figures (\ref{fig3-1}) and (\ref{fig3-2}).

\subsection{The Crank-Nicholson scheme}
It is possible to avoid the stability condition by modifying the simple scheme discussed earlier.
Define
\[(\delta^2u)_j^n=\frac{u_{j+1}^{n}-2u_j^n+u_{j-1}^n}{h^2}.\]
and let $0\leqslant Q \leqslant1$. For each such number, we define a ``$Q$-scheme''
\begin{equation}\label{qscheme}
\frac{u_j^{n+1}-u_j^n}{k}=(1-Q)(\delta^2u)_j^n+Q(\delta^2u)_j^{n+1}.
\end{equation}
When $Q=0$, we get the explicit scheme we have discussed. For $Q>0$ the scheme is implicit; At each step we must solve a linear system.

We analyze the stability of the scheme by inserting a separated solution of the form
\begin{equation}\label{eq3-4}
u_j^n=(\xi_r)^ne^{i\frac{2r\pi j}N}.
\end{equation}
\begin{enumerate}
\item[Note:] This separated solution will not satisfy Dirichlet boundary conditions but {\it periodic} boundary conditions, $u_{j+N}^n=u_j^n$ In this type of stability analysis we disregard the actual boundary conditions of the problem, using instead periodic boundary conditions. A numerical scheme is {\it Von Neumann Stable} if
\[|\xi_r|\leqslant1,\quad \forall r\]
There are other notions of stability where the actual boundary conditions are taken into account, but we will not discuss them in this course. Usually Von Neumann Stability gives a good indicator of stability, also when the actual boundary conditions are used.
\end{enumerate}
If we insert (\ref{eq3-4}) into the $Q$-scheme we get
\[\xi_r=\frac{1-2(1-Q)s(1-\cos \frac{2r\pi}{N})}{1+2Qs(1-\cos \frac{2r\pi}{N})}, \quad s=\frac k{h^2}.\]
Evidently $\xi_r\leqslant1$. The condition $\xi_r\eqslantgtr -1$ requires that
\[s(1-2Q)(1-\cos\frac{2r\pi}{N})\leqslant1.\]
This is {\it always} true for $1-2Q\leqslant0$. Thus if $\frac12\leqslant Q\leqslant1$ the numerical scheme is stable for {\it all} choices of grid parameters $k$ and $h$. The price we pay is that the scheme is implicit.

For $\frac12\leqslant Q\leqslant1$ we say that the $Q$-scheme is {\it unconditionally stable}. For $Q=\frac12$, the $Q$-scheme is called the Crank-Nicolson scheme.

Observe that when $N$ is large, the angles 
\begin{equation*}
\theta_r=\frac{2\pi r}{N},r=1,2,\cdots N,
\end{equation*}
form a very dense set of points in the angular interval $[0,2\pi]$. Thus, in the argument above we might for simplicity take $\theta$ to be a continuous variable in the interval $[0,2\pi]$. We will therefore in the following  investigate Von Neuman stability by inserting 
\begin{equation*}
u_j^n=\xi^n\eta^j,
\end{equation*}
where $\eta=e^{i\theta}$, into the numerical scheme of interest. Verify that using this simplified approach leads to the same stability conditions for the two schemes (\ref{eq38}) and (\ref{qscheme}).

\subsubsection{Finite difference method for the heat equation: Neumann boundary conditions}
Let us consider the heat equation as before but now with Neumann boundary conditions.
\begin{align*}
u_t&=u_{xx}, \quad\,\,\,\,0<x<1,\\
u_x(0,t)&=f(t), \quad \,u_x(1,t)=g(t).
\end{align*}
One possible discretization for the boundary conditions is to use forward difference at $x=0$ and backward difference at $x=1$.
\begin{align*}
\frac{u_1^n-u_0^n}{h}&=f^n,\\
\frac{u_N^n-u_{N-1}^n}{h}&=g^n.
\end{align*}
These approximations are only of order $h$. The space derivative in the equation was approximated to second order $h^2$. We get second order also at the boundary by using center difference. However we then need to extend our grid by ``ghost points'' $x_{-1}$ and $x_{N+1}$
\begin{figure}[h]
\centering
\includegraphics{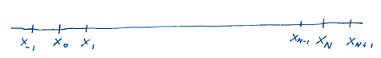}
\caption{}
\end{figure}
\[(u_x)_0=\frac{u_1-u_{-1}}{2h},\qquad (u_x)_N=\frac{u_{N+1}-u_{N-1}}{2h}.\]
Thus our numerical scheme is
\begin{align*}
u_i^{n+1}&=s(u_{i+1}^{n}+u_{i-1}^{n})+(1-2s)u_i^n,\quad i=0,\dots,N,\\
u_{-1}^n&=u_1^n-2hf^n,\\
u_{N+1}^n&=u_{N-1}^n+2hg^n.
\end{align*}
Introduction of ghost points is the standard way to treat derivative boundary conditions for 1D,2D and 3D problems.

\subsection{Finite difference method for the wave equation}
We consider the following initial and boundary value problem for the wave equation
\begin{align*}
u_{tt}&-u_{xx}=0, \qquad \qquad 0<x<L,t>0,\\
u(0,t)&=u(L,t)=0,\\
u(x,0)&=\varphi(x),\\
u_t(x,0)&=\psi(x),
\end{align*}
on the real axis. We introduce a uniform spacetime grid $(x_j,t_n)$ where 
\begin{align*}
x_j&=hj,\quad j=0,\dots,N,\\
t_n&=k n, \qquad n\in \mathbb{N},
\end{align*}
where $x_0=0,x_N=L$ and thus $h=\frac{L}{N}$.
\begin{figure}[h]
\centering
\includegraphics{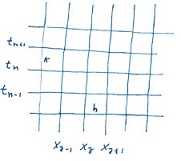}
\caption{}
\end{figure}

\noindent The wave equation is required to hold at the grid points only.
\[(u_{tt})_j^n=(u_{xx})_j^n.\]
The partial derivatives are approximated by centered differences (order $h^2$, $k^2$).
\[\frac{u_j^{n+1}-2u_j^n+u_j^{n-1}}{k^2}=\frac{u_{j+1}^n-2u_j^n+u_{j-1}^n}{h^2},\]
\[\Updownarrow\]
\begin{equation}\label{eq3-5}
u_j^{n+1}=s(u_{j+1}^n+u_{j-1}^n)+2(1-s)u_j^n-u_j^{n-1},
\end{equation}
where $s=\frac{k^2}{h^2}$. The boundary conditions evaluated on the grid simply become $u_0^n=0,u_N^n=0$.

\begin{figure}[!h]
\centering
\includegraphics{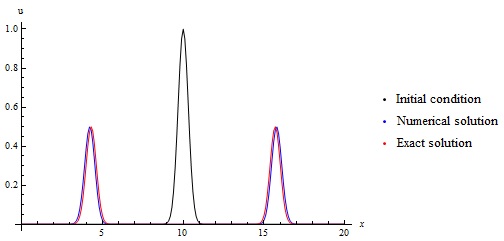}
\caption{}
\label{fig3-3}
\end{figure}

\begin{figure}[!h]
\centering
\includegraphics{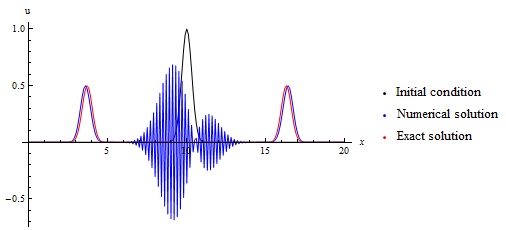}
\caption{}
\label{fig3-4}
\end{figure}

Note that the difference equation (\ref{eq3-5}) is of order 2 in $n$ and $j$. In order to solve the difference equation and compute the solution at time level $n+1$, we need to know the values on time level $n$ and $n-1$. Thus to compute $u_j^2$, we need $u_j^1$ and $u_j^0$. These values are provided by the initial conditions
\begin{align*}
u(x,0)&=\varphi(x),\\
u_t(x,0)&=\psi(x).
\end{align*}
Evaluated on the grid these are
\begin{align*}
u_j^0&=\varphi_j,\\
(u_t)_j^0&=\psi_j.
\end{align*}
We now need a difference approximation for $u_t$ at timelevel $n=0$.
Since the scheme (\ref{eq3-5}) is of order 2 in $t$ and $x$, we should use a 2-order approximation for $u_t$. If we use a 1-order approximation for $u_t$ in the initial condition the whole numerical scheme is of order 1.
Using a centered difference in time, we get
\begin{equation}\label{wbd1}
\frac{u_j^1-u_j^{-1}}{2k}=\psi_j.
\end{equation}
We have introduced a ``ghost point'' $t_{-1}$ in time with associated $u$ value $u_j^{-1}$. If we insert $n=0$ in (\ref{eq3-5}), we get a second identity involving $u_j^{-1}$
\begin{equation}\label{wbd2}
u_j^1=s(u_{j+1}^0+u_{j-1}^0)+2(1-s)u_j^0-u_j^{-1}.
\end{equation}
We can use (\ref{wbd1}) and (\ref{wbd2}) to eliminate the ghost point and find in conjunction with $u_j^0=\varphi_j$ that 
\begin{align*}
u_j^0&=\varphi_j,\\
u_j^1&=\frac12s(\varphi_{j+1}+\varphi_{j-1})+(1-s)\varphi_j+k\psi_j.
\end{align*}
These are the initial values for the scheme (\ref{eq3-5}).
We now test the numerical scheme using initial conditions 
\begin{align*}
\varphi(x)&=e^{-a(x-\frac{1}{2}L)^2},\\
\psi(x)&=0.
\end{align*}
The exact solution 
\begin{equation*}
u(x,t)=\frac{1}{2}(\varphi(x+t)+\varphi(x-t)),
\end{equation*}
to this problem is found using a method that we will discuss later in this class.
 For the test we use $L=20$ and  chose the parameter $a$ so large that $\varphi(x)$ to a good approximation satisfy the boundary conditions for the wave equation.
  In figures (\ref{fig3-3}) and (\ref{fig3-4})  we plot the initial condition, the exact solution and the numerical solution on the interval $(0,L)$. The space-time grid parameters are $h=0.1$, thus $N=200$, and $k=h\sqrt{s}=0.0004\sqrt{s}$ for $s=0.9$ in figure \ref{fig3-3} and $s=1.1$ in figure \ref{fig3-4}. We use 60 iterations in $n$.

\vfill \newpage
\noindent It appears as if the scheme is stable for $s\lesssim1$ and unstable for $s>1$. Let us prove this using the Von Neumann Approach.
\[u_j^n=\eta^j\xi^n,\qquad \eta=e^{i\theta}.\]
Inserting this into the scheme we get
\[\eta^j\xi^{n+1}=s(\eta^{j+1}\xi^n+\eta^{j-1}\xi^n)+2(1-s)\eta^j\xi^n-\eta^j\xi^{n-1},\]
\[\quad\qquad\qquad\qquad\qquad\Downarrow\qquad (\text{divided by }\xi^n\eta^j)\]
\[\xi+\frac1{\xi}-2=s(\eta+\frac1{\eta}-2)=2s(\cos \theta-1).\]
Define $p=s(\cos\theta-1)$. Then the equation can be written as
\[\xi^2-2(1+p)\xi+1=0.\]
The roots are
\[\xi=1+p\pm\sqrt{p^2+2p}.\]
Clearly $p\leqslant0$. If $p<-2$ then $p^2+2p=p(p+2)>0$ and there are two real roots. One of these must be $<-1 \Rightarrow |\xi|>1 \Rightarrow$ scheme unstable. If $p>-2$ then $p^2+2p<0$ and we have two complex conjugate roots
\[\xi_{\pm}=1+p\pm i\sqrt{-p^2-2p}.\]
These two roots satisfy
\begin{align*}
|\xi_{\pm}|^2&=(1+p)^2+(-p^2-2p)\\
            &=1+2p+p^2-p^2-2p=1.
\end{align*}
we therefore have stability if $p>-2$. It is easy to verify that the scheme is also stable for $p=-2$. Taking into account the definition of $p$ we thus have stability if  
\[s(\cos\theta-1)\geq-2,\]
\[s\leqslant\frac2{1-\cos\theta}.\]
This holds for all $\theta$ if $s\leqslant1$ because clearly $\frac2{1-\cos\theta}\eqslantgtr1$. Thus we have stability if 
\[s\leqslant1.\]
which is what our numerical test indicated.

\subsection{Finite difference method for the Laplace equation}
Let us consider the problem
\begin{align*}
\Delta u&=0,\qquad(x,y)\in(0,1)\times(0,1)\equiv D,\\
u&=f\qquad(x,y)\in\partial D.
\end{align*}
We solve the problem by introducing a grid
\[(x_i,y_j)\qquad x_i=i h, y_j=jp, \quad i=0,\dots,N, j=0,\dots,M\]
where $h=\frac 1N, p=\frac 1M$.
\begin{figure}[h]
\centering
\includegraphics{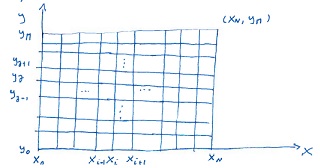}
\caption{}
\label{fig3-7}
\end{figure}

\begin{figure}[!h]
\centering
\includegraphics{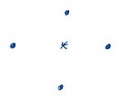}
\caption{}
\label{fig3-8}
\end{figure}

\noindent Using center difference for $u_{xx}$ and $u_{yy}$, we get the scheme
\[\frac{u_{i+1,j}-2u_{i,j}+u_{i-1,j}}{h^2}+\frac{u_{i,j+1}-2u_{i,j}+u_{i,j-1}}{p^2}=0.\]
This scheme is 2-order in $h$ and $p$. Because of the symmetry there is no reason to choose $h$ and $p$ differently. For the case $h=p$, which corresponds to $M=N$, we get
\begin{equation}\label{eq3-7}
u_{i,j}=\frac14(u_{i+1,j}+u_{i-1,j}+u_{i,j+1}+u_{i,j-1}).
\end{equation}
 Observe that the value in the center is simply the mean of the value of all closest neighbors. The {\it stencil} of the scheme i shown in figure (\ref{fig3-8}).

There is no iteration process in the scheme (\ref{eq3-7}) (There is no time!). It is a linear system of variables for the $(N-1)(N-1)$ interior gridpoints. We are thus looking at a linear system of equations
\begin{equation}\label{eq3-8}
A\mathbf{x}=\mathbf{b},
\end{equation}
where $\mathbf{b}$ comes from the boundary values.
Observe however, that the linear system (\ref{eq3-7}) is not actually written in the form (\ref{eq3-8}). In order to apply linear algebra to solve (\ref{eq3-7}), we must write (\ref{eq3-7}) in the form (\ref{eq3-8}). The problem that arise at this point is that the  interior gridpoints in (\ref{eq3-7})  are not {\it ordered}, whereas the vector components in (\ref{eq3-8}) are ordered. We must thus {\it choose} an ordering of the gridpoints. Each choice will lead to a different {\it structure} for the matrix A (for instance which elements are zero and which are nonzero). Since the efficiency of iterative methods depends on the structure of the matrix, this is a point where we can be clever (or stupid!).

For the problem at hand it is natural to order the interiour gridpoints {\it lexiographical}.
\begin{figure}[h]
\centering
\includegraphics{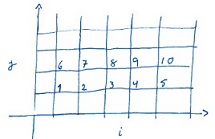}
\caption{: Ordering the grid points.}
\end{figure}

\noindent Thus we define the vector $\mathbf{x}$ as
\[\mathbf{x}=(u_{1,1},u_{2,1},\dots,u_{N-1,1},u_{1,2},u_{2,2},\dots,u_{N,1},\dots, u_{N-1,N-1}).\]
The $\mathbf{b}$ vector is obtained from the values of $u$ at the boundary.
For example when $j=1$, we have 
\[u_{i,1}=\frac14(u_{i+1,1}+u_{i-1,1}+u_{i,2}+u_{i,0}),\quad 1\leqslant i\leqslant N-1.\]
The boundary condition gives
\begin{align*}
u_{i,0}&=f_{i,0} \equiv f(x_i,y_0),\\
u_{i,N}&=f_{i,N},\\
u_{0,j}&=f_{0,j},\\
u_{N,j}&=f_{N,j},
\end{align*}
$i=1,$
\[4u_{1,1}=u_{2,1}+u_{0,1}+u_{1,2}+u_{1,0},\]
\[\Updownarrow\]
\[-4u_{1,1}+u_{1,2}+u_{2,1}=-f_{0,1}-f_{1,0},\]
$1<i<N-1,$
\[4u_{i,1}=u_{i+1,1}+u_{i-1,1}+u_{i,2}+u_{i,0},\]
\[\Updownarrow\]
\[u_{i-1,1}-4u_{i,1}+u_{i+1,1}+u_{i,2}=-f_{i,0},\]
$i=N-1,$
\[4u_{N-1,1}=u_{N,1}+u_{N-2,1}+u_{N-1,2}+u_{N-1,0},\]
\[\Updownarrow\]
\[u_{N-2,1}-4u_{N-1,1}+u_{N-1,2}=-f_{N,1}-f_{N-1,0}.\]
Thus the matrix for the system has a structure that starts off as
\begin{footnotesize}
\[\setcounter{MaxMatrixCols}{30}
\begin{matrix}
   &   &        &       &      &      &\textrm{\tiny{$N-1$}}&  &     &      &      &      &  & &\textrm{\tiny{$(N-1)(N-1)$}}\\
-4 & 1 & 0      &\cdots &\cdots&\cdots &0 &1 &0    &\cdots&\cdots&\cdots&\cdots&\cdots&0\\
1 & -4 & 1      &0      &\cdots &\cdots&0 &0 &1    &0     &\cdots&\cdots&\cdots&\cdots&0\\
\vdots  &    &  &       &       &      &  &  &     &      &      &      &      &      &\vdots\\
0 & \cdots&\cdots&\cdots&\cdots &1     &-4&0 &\cdots&\cdots&\cdots&1     &0     &\cdots&0
\end{matrix}.\]
\end{footnotesize}
 By writing down a larger part of the linear system we eventually realize that the matric of the system has the form of a $(N-1)\times (N-1)$ block matrix
 \vspace{2mm}

\[ A=\begin{pmatrix}
B & I & 0 &\cdots &\cdots &0\\
I & B & I &0      &\cdots &0\\
  &   &\vdots&\vdots &       & \\
0 & \cdots&\cdots &\cdots&I&B
\end{pmatrix},\]
where $I$ is the $(N-1)(N-1)$ identity and where $B$ is the following $(N-1)\times (N-1)$ tri-diagonal matrix

\[ B=\begin{pmatrix}
-4 & 1 & 0 &\cdots &\cdots&\cdots &0\\
1 & -4 & 1 &0      &\cdots &\cdots&0\\
0 & 1 &-4  &1     &0       &\cdots&\cdots \\
  &   & \vdots&   &\vdots  &     &       \\
0 & \cdots&\cdots &\cdots&\cdots &1     &-4
\end{pmatrix}.\]

\noindent We observe that most elements in the matrices $A$ and $B$  are zero, we say the matrices are {\it sparse}. The structural placing of the nonzero elements represents the
{\it sparsity structure} of the matrix. Finite difference methods for PDEs usually leads to sparse matrises with a simple sparsity structure.
      
 When $N$ is large, as it often is, the linear system (\ref{eq3-8}) represents an enormous computational challenge.

In 3D with 1000 points in each direction, we get a linear system with $10^9$ independent variables and the matrix consists of $10^{18}$ entries!

There are two broad classes of methods available for solving large linear systems of $P$ equations in $P$ unknowns.
\begin{enumerate}
\item[1)] Direct methods: These are variants of elementary Gaussian elimination. If the matrix is sparse (most elements are zero) direct methods can solve large linear systems fast. For dense matrices, the complexity of Gaussian elimination $(O(P^3))$ makes it impossible to solve more than few thousand equations in reasonable time.
\item[2)] Iterative methods: These methods try to find approximate solutions to (\ref{eq3-7}) by defining an iteration
\[\mathbf{x}^{n+1}=M\mathbf{x}^n,\]
that converges to a solution of the linear system of interest.
Since each multiplication cost at most $O(P^2)$ operations, fast convergence($n\ll P$) will give an approximate solution to (\ref{eq3-8}) in $O(P^2)$ operations, even for a dense matrix.
\end{enumerate}

\noindent The area of iterative solution methods for linear systems is vast and we can not possibly give an overview of the available methods and issues in a introductory class on PDEs. However, it is appropriate to discuss some of the central ideas.

\subsubsection{Iterations in $\mathbb{R}^n$}
Let $\mathbf{u}^0\in\mathbb{R}^n$ be a given vector in $\mathbb{R}^n$ and consider the following iteration
\[\mathbf{u}^{l+1}=M\mathbf{u}^l+\mathbf{c},\]
where $M$ is a $n\times n$ matrix and $\mathbf{c}$ is a fixed vector. Let us assume that $\mathbf{u}^l\to \mathbf{u}^{*}$ when $l\to\infty$. Then we have
\[\mathbf{u}^{l+1}=M\mathbf{u}^l+\mathbf{c},\]
\[\Downarrow\]
\[\mathbf{u}^*=M\mathbf{u}^*+\mathbf{c},\]
\[\Downarrow\]
\[(I-M)\mathbf{u}^*=\mathbf{c}.\]
Thus $\mathbf{u}^*$ solves a linear system. The idea behind  iterative solution of linear system can be described as follows:

 Find a matrix $M$ and a vector $\mathbf{c}$ such that
\begin{enumerate}
\item \qquad \qquad \qquad\qquad \qquad \qquad\qquad$A\mathbf{u}=f,$
\[\Updownarrow\]
\[(I-M)\mathbf{u}=\mathbf{c},\]
\item The iteration
\[\mathbf{u}^{l+1}=M\mathbf{u}^l+\mathbf{c},\]
converges.
\end{enumerate}
The condition \romannumeral1) can be realized in many ways and one seeks choises $M$ and $\mathbf{c}$ such that the convergence is {\it fast}.

 In order get anywhere, it is at this point evident that we  need to be more precise about convergence in $\mathbb{R}^n$. This we do by introducing a vector norm on $\mathbb{R}^n$ that measure distance between vectors. We will also need a corresponding norm on $n\times n$ matrices.

A $n\times n$ matrix is a linear operator on $\mathbb{R}^n$. Let us assume that $\mathbb{R}^n$ comes with a norm $\lVert\cdot\rVert$ defined. There are many such norms in used for $\mathbb{R}^n$. The following three norms are frequently used
\begin{align*}
&\lVert\mathbf{x}\rVert_2=\left(\sum_{i=1}^n\lvert x_i\rvert^2\right)^{\frac12},\\
&\lVert\mathbf{x}\rVert_1=\sum_{i=1}^n\lvert x_i\rvert,\\
&\lVert\mathbf{x}\rVert_\infty=\mathop{sup}\limits_i\left\{\lvert x_i\rvert\right\}.
\end{align*}
Given a norm, $\lVert\cdot\rVert$, on $\mathbb{R}^n$ there is a corresponding norm on linear operators (matrices) defined by
\[\lVert A \rVert=\mathop{sup}\limits_{\lVert\mathbf{x}\rVert=1}\lVert  A\mathbf{x}\rVert.\]
For matrix norms corresponding to vector norms we have the following two important inequalities.
\[\lVert A \mathbf{x}\rVert \leqslant\lVert A \rVert\lVert \mathbf{x}\rVert,\]
\[\lVert A B\rVert\leqslant\lVert A \rVert\lVert B \rVert.\]
For $\lVert\cdot\rVert_1$ and $\lVert\cdot\rVert_\infty$, one can show that
\begin{align*}
\lVert A\rVert_\infty&=\mathop{max}\limits_{j} \sum_i\lvert a_{i,j}\rvert,\\
\lVert A\rVert_1&=\mathop{max}\limits_{i}\sum_j\lvert a_{i,j} \rvert.
\end{align*}
For any matrix, define the {\it spectral radius} $\rho(A)$ by 
\[\rho(A)=\mathop{max}\limits_s|\lambda_s|,\]
where $\lambda_s$ are the eigenvalues of $A$. For the $\lVert \cdot\rVert_2$ norm we have the following identity
\[\lVert A \rVert_2=\sqrt{\rho(AA^*)}, \]
where $A^*$ is the adjoint of $A$. For the special case when $A$ is selfadjoint, $A=A^*$, we have
\begin{properties}
Let $A=A^*$. Then
\[\lVert A \rVert_2=\rho(A).\]
\end{properties}

\begin{proof}
Let $\lambda$ be an eigenvalue of $A$ with associated eigenvector $u$. Then we have
\[A^2(u)=A(Au)=A(\lambda u)=\lambda A(u)=\lambda^2u,\]
\[\Downarrow\]
\[\text{spectrum}(A^2)=\left\{ \lambda^2|\lambda\in\text{spectrum}(A)\right\},\]
\[\Downarrow\]
\begin{align*}
\lVert A\rVert_2&=(\rho(AA^*))^{\frac12}=(\rho(A^2))^{\frac12}=(\mathop{max}\limits_s|\lambda_s|^2)^{\frac12}\\
&=\mathop{max}\limits_s\sqrt{|\lambda_s|^2}=\mathop{max}\limits_s|\lambda_s|=\rho(A). \qedhere
\end{align*}
\end{proof}
\noindent Let us now consider the iteration
\[\mathbf{u}^{l+1}=M\mathbf{u}^l+\mathbf{c},\]
and let us use the $\lVert\cdot\rVert_2$ norm on $\mathbb{R}^n$ with the associated norm $\lVert\cdot\rVert_2$ on $n\times n$ matrices. The iteration will then converge to $\mathbf{u}^*$ in $\mathbb{R}^n$ iff
\[\lim_{l\to\infty}\lVert\mathbf{u}^l-\mathbf{u}^*\rVert_2=0.\]
But
\[\lVert\mathbf{u}^{l}-\mathbf{u}^*\rVert_2=\lVert M\mathbf{u}^{l-1}+\mathbf{c}-(M\mathbf{u}^*+c)\rVert_2\]
\[=\lVert M(\mathbf{u}^{l-1}-\mathbf{u}^*)\rVert_2\leqslant\lVert M\rVert_2\lVert\mathbf{u}^{l-1}-\mathbf{u}^*\rVert_2.\]
By induction we get 
\[\lVert\mathbf{u}^l-\mathbf{u}^*\rVert_2\leqslant(\lVert M\rVert_2)^l\lVert\mathbf{u}^0-\mathbf{u}^*\rVert_2.\]
Let us assume that $M$ is symmetric, then
\[\lVert\mathbf{u}^l-\mathbf{u}^*\rVert_2\leqslant(\rho(M))^l\lVert\mathbf{u}^0-\mathbf{u}^*\rVert_2.\]
This gives us the following sufficient condition for convergence.
\begin{theorem}
Let $M$ be symmetric and assume that for all eigenvalues $\lambda$ for $M$ we have 
\[|\lambda|<1.\]
Then the iteration
\[\mathbf{u}^{l+1}=M\mathbf{u}^l+\mathbf{c},\]
converges to the unique solution of
\[(I-M)\mathbf{u}=\mathbf{c},\]
for all initial values $\mathbf{u}^0$.
\end{theorem}
Note that even though all initial values converges to the solution, the number of iterations we need to get a predetermined accuracy {\it will} depend on the choice of $\mathbf{u}^0$.

We are now ready to describe the two simplest iteration methods for solving $A\mathbf{u}=\mathbf{b}$.

\subsubsection{Gauss-Jacobi}
Any matrix can be written in the form 
\[A=L+D+U,\]
where $L$ is lower triangle, $U$ is upper triangular and $D$ is diagonal. Using this decomposition we have
\[A\mathbf{u}=\mathbf{b},\]
\[\Updownarrow\]
\[L\mathbf{u}+D\mathbf{u}+U\mathbf{u}=\mathbf{b},\]
\[\Updownarrow\]
\[D\mathbf{u}=-(L+U)\mathbf{u}+\mathbf{b}.\]
Assume that $D$ is invertible (no zeroes on the diagonal). Then we get
\[\mathbf{u}=M_J\mathbf{u}+\mathbf{c},\]
where $M_J=-D^{-1}(L+U)$, $\mathbf{c}=D^{-1}\mathbf{b}$.
This gives the Gauss-Jacobi iteration.
\[\mathbf{u}^{l+1}=M_J\mathbf{u}^l+\mathbf{c}.\]
In component form the Gauss-Jacobi iteration can be written in the form
\[u_i^{l+1}=\frac1{a_{ii}}(b_i-\sum_{j\neq i}a_{ij}u_j^l),\quad i=1,\dots,n\;\;.\]
Observe that
\[M_J=-D^{-1}(L+U),\]
\[\Downarrow\]
\[DM_J=-L-U=-A+D,\]
\[\Downarrow\]
\[M_J=-D^{-1}A+I,\quad \text{and}\quad D^{-1}=-\frac14I,\]
\[\Downarrow\]
\[M_J=\frac14A+I.\]
Let $\lambda$ be an eigenvalue of $A$ with corresponding eigenvector $v$. Then we have
\[M_J\mathbf{v}=\frac14Av+Iv=(\frac14\lambda+1)v.\]
This means that
\begin{equation}\label{Gspectrum}
\text{spectrum}(M_J)=\left\{\left.\left(\frac14\lambda+1\right)\right|\lambda\in\text{spectrum}(A)\right\}.
\end{equation}
In order to find the spectrum of $A$, we will start from a property of the spectrum for block-matrices.

Let $C$ be a matrix on the form
\[C=\begin{bmatrix}
C_{11} &\cdots &\cdots &C_{1m}\\
\vdots &       &       &\vdots\\
C_{m1} &\cdots &\cdots &C_{mm}
\end{bmatrix}.\]
where $C_{ij}$ are $n\times n$ matrices. Thus $C$ is a $(nm)\times(nm)$ matrix. Let us assume that there is a common set of eigenvectors for the matrices $C_{ij}$. Denote this set by $\{v^k\}^n_{k=1}$. Thus we have
\[C_{ij}v^k=\lambda_{ij}^kv^k,\]
where $\lambda_{ij}^k$ is the eigenvalue of the matrix $C_{ij}$ corresponding to the eigenvector $v^k$.

Let $\alpha_1,\dots,\alpha_m$ be real numbers. We will seek eigenvectors of $C$ on the form
\[w^k=\begin{bmatrix}
\alpha_1v^k\\
\alpha_2v^k\\
\vdots\\
\alpha_mv^k
\end{bmatrix}.\]
Then $w^k$ is an eigenvector of $C$ with associated eigenvalue $\mu^k$ only if
\[Aw^k=\mu^kw^k,\]
\[\Updownarrow\]
\[\begin{bmatrix}
(\lambda_{11}^k\alpha_1+\dots+\lambda_{1m}^k\alpha_m)v^k\\
\vdots\\
(\lambda_{m1}^k\alpha_1+\dots+\lambda_{mm}^k\alpha_m)v^k
\end{bmatrix}
\begin{bmatrix}
\mu^k\alpha_1v^k\\
\vdots\\
\mu^k\alpha_mv^k
\end{bmatrix}.\]
The eigenvalues of $C$ are thus given by the eigenvalues of the matrix
\[\Lambda=\begin{bmatrix}
\lambda_{11}&\dots&\lambda_{1m}\\
\vdots&&\vdots\\
\lambda_{m1}&\dots&\lambda_{mm}
\end{bmatrix}.\]
Let us now return to the convergence problem for Gauss-Jacobi. 

 $A$ is a $(N-1)\times (N-1)$ block matrix given by 

\[A=\begin{bmatrix}
B&I&0    &\dots&\dots&0\\
I&B&I    &0    &\dots&0\\
\vdots&&&&&\vdots\\
0&0&\dots&\dots&I    &B
\end{bmatrix},\]
where $B$ is a $(N-1)\times (N-1)$ tri-diagonal matrix given by

\[B=\begin{bmatrix}
-4   &1&0&\dots&\dots&0\\
1    &-4&1&0&\dots&0\\
\vdots&&&&&\vdots\\
0&\dots&\dots&\dots&1&-4
\end{bmatrix}.\]

\noindent  $B$ obviously commute with $I$ and the zero matrix and they therefore has common eigenvectors. Since all vectors in $\mathbb{R}^{N-1}$ are eigenvectors of the identity matrix and the zero matrix, with eigenvalues 1 and zero, respectively, in order to set up the $\Lambda$ matrix we only need to find the eigenvalues of $B$.

Let $D$ be a $n\times n$ tri-diagonal matrix of the form

\[D=\begin{bmatrix}
a     &b    &0    &\dots&\dots&\dots&0\\
b     &a    &b    &0    &\dots&\dots&0\\
\vdots&     &     &\vdots&     &     &\vdots\\
0     &\dots&\dots&\dots&b    &a    &b\\
0     &\dots&\dots&\dots&0    &b    &a
\end{bmatrix}.\]

\noindent Let $\lambda$ be an eigenvalue with corresponding eigenvector $\mathbf{v}=(v_1,\dots,v_n)$
Thus
\[Dv=\lambda v,\]
\[\Updownarrow\]
\[bv_{i-1}+av_i+bv_{i+1}=\lambda v_i,\quad i=1,\dots,n\;,\]
and $v_0=v_{n+1}=0$. This is a boundary-value problem for a difference equation of order 2. We write the equation on the form
\[bv_{i-1}+(a-\lambda)v_i+bv_{i+1}=0,\]
and look for solutions of the form
\[v_i=r^i\Rightarrow v_{i-1}=r^{-1}v_i,\quad v_{i+1}=rv_i.\]
Thus we get the algebraic condition
\begin{equation}\label{charpoly}
br^2+(a-\lambda)r+b=0,
\end{equation}
\[\Downarrow\]
\[r_{\pm}=\frac1{2b}(-(a-\lambda)\pm((a-\lambda)^2-4b^2)^{\frac12}).\]
If $r_{+}\neq r_-$, the general solution to the difference equation is 
\[v_i=Ar_+^i+Br_-^i,\]
the boundary conditions give
\begin{align*}
v_0=0,\quad &\Leftrightarrow\quad A+B=0,\\
v_{n+1}=0, \quad&\Leftrightarrow\quad Ar_+^{n+1}+Br_-^{n+1}=0,
\end{align*}
nontrival solutions exits only if
\begin{eqnarray}\label{trieqn}
r_-^{n+1}-r_-^{n+1}=0,\nonumber\\
\quad\Updownarrow\quad\nonumber\\
\left(\frac{r_+}{r_-}\right)^{n+1}=1.
\end{eqnarray}
But \quad $r_+r_-=1 \Rightarrow r_-=\frac1{r_+}$ and therefore from equation (\ref{trieqn}) 
\begin{eqnarray*}
(r_+^2)^{(n+1)}&=1,\\
&\Updownarrow\\
 r_+&=e^{\pi i\frac{k}{n+1}},\\
 &\Downarrow\\
  r_-&=e^{-\pi i\frac{k}{n+1}}.
 \end{eqnarray*}
 Furthermore, from equation (\ref{charpoly}) we have 
 \begin{eqnarray}\label{tracecondition}
\frac{a-\lambda}{b}&=-(r_++r_-),\nonumber\\
&\Updownarrow\nonumber\\
\lambda&=a+b(r_++r_-),\nonumber\\
&\Updownarrow\nonumber\\
\lambda&=a+2b\cos(\pi \frac {k}{n+1}).
\end{eqnarray}
The eigenvalues of the tri-diagonal matrix $D$ are thus
\[\lambda_k=a+2b\cos\left(\frac{\pi k}{n+1}\right),\quad k=1,\dots,n\;.\]
Returning to our analysis of the matrix $A$, we can now assemble our matrix $\Lambda^k$
\[\Lambda^k=\begin{bmatrix}
\lambda_k     &1    &0    &\dots&\dots&0\\
1     &\lambda_k    &1    &0    &\dots&0\\
\vdots&     &     &     &          &\vdots\\
0     &\dots&\dots&\dots    &1    &\lambda_k
\end{bmatrix},\]
where $\lambda_k=-4+2\cos\left(\frac{\pi k}{N}\right),\quad k=1,\dots,N-1$. $\Lambda^k$ is also tri-diagonal and the previous theory for tri-diagonal matrices give us the eigenvalues
\[\lambda_{kl}=\lambda_k+2\cos\left(\frac{\pi l}{N}\right),\quad l=1,\dots,N-1\;,\]
\[\Downarrow\]
\[\lambda_{kl}=-4+2\cos\left(\frac{\pi k}{N}\right)+2\cos\left(\frac{\pi l}{N}\right),\]
or using standard trigonometric identities
\[\lambda_{kl}=-4\left(\sin^2\left(\frac{\pi k}{2N}\right)+\sin^2\left(\frac{\pi l}{2N}\right)\right),\quad k=1,\dots,N-1,\quad l=1,\dots,N-1\;,\]
and finally the spectrum of $M_J$ in the Gauss-Jacobi iteration method is from (\ref{Gspectrum})
\[\mu_{kl}=1-\sin^2\left(\frac{\pi k}{2N}\right)-\sin^2\left(\frac{\pi l}{2N}\right).\]
Using Theorem 1, the stability condition for Gauss-Jacobi is
\[| \mu_{kl}|<1,\]
\[\Updownarrow\]
\[\left|1-\sin^2\left(\frac12\pi x_k\right)-\sin^2\left(\frac12\pi y_l\right)\right|<1,\quad x_k=\frac kN,\quad y_l=\frac lN\;.\]
Since both $x_k,y_l\in(0,1)$, it is evident that the condition holds and that we have convergence. We get the largest eigenvalue by choosing $k=l=1$. The spectral radius can then for large $N$ be approximated as 
\[\rho(M)=1-2\sin^2\left(\frac12 \pi h\right)\approx1-\frac12(\pi h)^2,\]
where we have recall that $h=\frac{1}{N}$ grid parameter for the grid defined in figure (\ref{fig3-7}). The relative error is defined as
\[\frac{\lVert\mathbf{u}^l-\mathbf{u}^*\rVert}{\lVert \mathbf{u}^0-\mathbf{u}^*\rVert}\leqslant(\rho(M))^l\approx\left(1-\frac12 \pi^2h^2\right)^l.\]
The discretization error of our scheme is of order $h^2$. It is natural to ask how many iterations we need in order for the relative error to be of this size. There is no point in iterating more than this.
\[\left(1-\frac12 \pi^2h^2\right)^l\approx h^2,\]
\[\Updownarrow\]
\[e^{-\frac12l\pi^2h^2}\approx h^2,\]
\[\Updownarrow\]
\[l\approx\frac4{\pi^2}N^2\ln N.\]
Gauss elimination requires in general $O(N^3)$ and since for large $N$ we have $N^2\ln N\ll N^3$ we observe that Gauss-Jacobi iteration is much faster that Gauss elimination.

\subsubsection{Gauss-Seidel}
We use the same decomposition of $A$ as for the Gauss-Jacobi method, $A=L+D+U$
\[(L+D+U)\mathbf{u}=\mathbf{b},\]
\[\Downarrow\]
\[-(L+D)\mathbf{u}=U\mathbf{u}-\mathbf{b},\]
\[\Downarrow\]
\[\mathbf{u}=M_S\mathbf{u}+\mathbf{c},\]
where $M_S=-(L+D)^{-1}U$, $\mathbf{c}=(L+D)^{-1}\mathbf{b}$. In component form we have the iteration
\[u_i^{l+1}=\frac1{a_{ii}}(b_i-\sum_{j<i}a_{ij}u_j^{l+1}-\sum_{j>i}a_{ij}u_j^l),\quad i=1,\dots,n\;.\]

\noindent We will now see that Gauss-Seidel converges faster than Gauss-Jacobi. It also use the memory in a more effective way.
For the stability of this method we must find the spectrum of the matrix
\[M_S=-(L+D^{-1})U.\]
Note that the iteration matrix for the Gauss-Seidel method can be rewritten into the form
\begin{equation*}
M_S=-(I+D^{-1}L)^{-1} D^{-1} U.
\end{equation*}
 Using straight forward matrix manipulations it is now easy to verify that $\lambda_s$ is an eigenvalue to the matrix $M_S$ iff 
 \begin{equation}
 \text{det}\left(\lambda_s D+\lambda_s L+U\right)=0.\label{SeideEigenValueEquation}
 \end{equation}
  easy to show that 
 Let now $\alpha$ be any number, and 
define a matrix $E$ by
\[E=\begin{bmatrix}
 I   &0        &\dots&\dots&0\\
0     &\alpha  I  &0        &\dots&0\\
\vdots&     &     &              &\vdots\\
0     &\dots&\dots&\dots        &\alpha^{N-2}I
\end{bmatrix}.\]
Doing the matrix multiplication, it is easy to verify that we have the identity
\[E(D+L+U)E^{-1}=D+\alpha L+\frac1{\alpha} U,\]
from which we immediately can conclude
\begin{equation}
\text{det}\left(D+\alpha L+\frac{1}{\alpha} U\right)=\text{det}\left(D+L+U\right).\label{AlphaDeterminantRule}
\end{equation}
This derivation only require $D$ to be a diagonal matrix, $L$ to be a lower triangular matrix and $U$ to be a upper triangular matrix.
Using (\ref{SeideEigenValueEquation}) and (\ref{AlphaDeterminantRule}), we have
\begin{align*}
0&=\text{det}\left(\lambda_s I-M_S\right)\\
&=\text{det}\left(\lambda_s D+\lambda_s L+U\right)\\
&=\text{det}\left(\lambda_s D+\alpha\lambda_s L+\frac{1}{\alpha}U\right)\\
&=\text{det}\left(\lambda_s D+\sqrt{\lambda_s}(L+U)\right),\;\;\;\alpha=\frac{1}{\sqrt{\lambda_s}}\\
&=\lambda_s^{\frac{(N-1)}{2}}\text{det}\left(\mu D+L+U\right), \;\;\;\mu=\sqrt{\lambda_s}\\
&=\lambda_s^{\frac{(N-1)}{2}}\text{det}(D)\text{det}\left(\mu +D^{-1}(L+U)\right)\\
&=\lambda_s^{\frac{(N-1)}{2}}\text{det}(D)\text{det}\left(\mu I-M_J\right).
\end{align*}
Thus, we have proved that the nonzero eigenvalues $\lambda_s$ of the Gauss-Seidel iteration matrix are determined by the eigenvalues $\mu$ of the Jacobi iteration matrix  through the formula
\begin{equation}
\lambda_s=\mu^2\label{MuAplhaRelation}
\end{equation}
Since all  eigenvalues of the Jacobi iteration matrix $M_J$ are inside the unit circle, this formula implies that the same is true for the eigenvalues of the Gauss-Seidel iteration matrix $M_S$. Thus Gauss-Seidel also converges for all initial conditions. As for the {\it rate} of convergence, using the results from the analysis of the convergence of the Gauss-Jacobi method, we find that
 \[\rho(M_S)\approx\left(1-\frac12\pi^2h^2\right)^2\approx1-\pi^2h^2,\]
 and therefore
 \[\frac{\lVert\mathbf{u}^l-\mathbf{u}^*\rVert}{\lVert \mathbf{u}^0-\mathbf{u}^*\rVert}=h^2\Leftrightarrow (1-\pi^2h^2)^l=h^2 \Leftrightarrow l\approx\frac2{\pi^2}N^2\ln N.\]
 For a fixed accuracy the Gauss-Seidel only needs half as many iterations as Gauss-Jacobi.
 
 These two examples have given a taste of what one needs to do in order to define iteration methods for linear systems and study their convergence.
 
 As previously mentioned there is an extensive set of iteration methods available for solving large linear systems of equation, and many has been implemented in numerical libraries available for most common programming languages. The choice of method depends  to a large extent on the structure of your matrix. Here are two common choices:
\[
 \begin{tabular}{ll}
\multicolumn{1}{c}{Matrix property}&\multicolumn{1}{c}{Method name}\\
\vspace{1mm}\\
{Sparse, symmetric, positive:}&{ \hspace{2cm} Conjugate gradient}\\
\vspace{1mm}\\
{Dense, nonsymmetric:}&{Generalized minimal 
residual method CGMRE}
 \end{tabular}
\]

\noindent The conjugate gradient method has also been extended to more general matrices. A standard reference for iterative methods is ``Iterative methods for sparse linear systems'' by Yoosef Saad.

\section{First order PDEs}
So far we have only discussed numerical methods for solving PDEs. It is now time to introduce the so-called {\it analytical methods}. The aim here is to produce some kind of formula(s) for the solution of PDEs. What counts as a formula is not very precisely defined, but in general it involves infinite series and integrals that can not be solved by quadrature. Thus, if the ultimate aim is to get a numerical answer, and it usually is, numerical methods must be used to approximate the said infinite series and integrals. There is thus no clear-cut line separating analytical and numerical methods.

We started our exposition of the finite difference method by applying it to second order equations. This is because the simplest finite difference methods work best for such equations.

For our exposition of analytical methods we start with first order PDEs because the simplest such methods work best for these types of equations.

\subsection{Higher order equations and first order systems}
It turns out that higher order PDEs and first order systems of PDEs are not really separate things. Usually one can be converted into the other. Whether this is useful or not depends on the context.

Let us recall the analog situation from the theory of ODEs. The ODE for the Harmonic Oscillator is
\[x''+w^2x=0.\]
Define two new functions by
\begin{align*}
x_1(t)&=x(t),\\
x_2(t)&=x'(t).
\end{align*}
Then
\begin{align*}
x'_1&=x'=x_2,\\
x'_2&=x''=-w^2x=-w^2x_1.
\end{align*}
Thus we get a first order system of ODEs.
\begin{align*}
x'_1&=x_2,\\
x'_2&=-w^2x_1.
\end{align*}
In the theory of ODEs this is a very useful device that can be used to transform any ODE of order $n$ into a first order system of ODEs in $\mathbb{R}^n$.

Note that the transition from one high order ODE, to a system of first order ODEs is not by any means unique. If we for the Harmonic Oscillator rather define
\begin{align*}
x_1&=x+x',\\
x_2&=x-x',
\end{align*}
we get the first order system
\begin{align*}
x'_1&=\alpha x_1-\beta x_2,\\
x'_2&=\beta x_1-\alpha x_2,\qquad \alpha=\frac12(1-w^2),\,\beta=\frac12(1+w^2).
\end{align*}
Sometimes it is also possible to convert a first order system into a higher order scalar equation. This is however not possible in general but can always be done for linear first order systems of both ODEs and PDEs.

Let us consider the wave equation
\[u_{tt}-\gamma^2 u_{xx}=0.\]
Observe that
\begin{align*}
&\quad\,\,(\partial_t+\gamma\partial_x)(\partial_t-\gamma\partial_x)u=(\partial_t+\gamma\partial_x)(u_t-\gamma u_x)\\
&=(\partial_t+\gamma\partial_x)u_t-(\partial_t+\gamma\partial_x)\gamma u_x\\
&=u_{tt}+\gamma u_{tx}-\gamma u_{xt}-\gamma^2u_{xx}\\
&=u_{tt}-\gamma^2u_{xx}=0.
\end{align*}
Inspired by this let us define a function
\[v=(\partial_t-\gamma\partial_x)u.\]
Then the wave equation implies
\[(\partial_t+\gamma\partial_x)v=0.\]
Thus the wave equation can be converted into the system
\begin{align*}
u_t-\gamma u_x&=v,\\
v_t+\gamma v_x&=0.
\end{align*}
Observe that the second equation is decoupled from the first. We will later in this class see that this can be used to solve the wave equation.

Writing a higher order PDE as a first order system is not always useful. Observe that
\[(\partial_x+i\partial_y)(\partial_x-i\partial_y)u=\partial_{xx}u+\partial_{yy}u.\]
Thus we might write the Laplace equation as the system
\begin{align*}
u_x-iu_y&=v,\\
v_x+iv_y&=0.
\end{align*}
However this is not a very useful system since we most of the time are looking for {\it real} solutions to the Laplace equation.

If we rather put
\begin{align*}
v&=u_x,\\
w&=u_y.
\end{align*}
Then using the Laplace equation and assuming equality of mixed partial derivatives we have
\begin{align*}
v_y&=u_{xy}=u_{yx}=w_x,\\
v_x&=u_{xx}=-u_{yy}=-w_y.
\end{align*}
Thus we get the first order system
\begin{equation}\label{eq4-1}
\begin{split}
v_y-w_x&=0,\\
v_x+w_y&=0.
\end{split}
\end{equation}
These are the famous {\it Cauchy-Riemann equations}. If you have studied complex function theory you know that (\ref{eq4-1}) is the key to the whole subject; they determine the complex differentiability of the function
\[f(z)=v(x,y)+iw(x,y),\quad z=x+iy.\]
(\ref{eq4-1}) is a coupled system and does not on its own make it any simpler to find solutions to the Laplace equation. However, by using the link  between complex function theory and the Laplace equation, great insight into it's set of solutions can be gained.

Going from first order systems to higher order equations is not always possible but sometimes {\it very} useful.

A famous case is the Maxwell equations for electrodynamics, which in a region containing no sources of currents and charge are
\begin{equation}\label{eq4-2}
\begin{split}
&\nabla\times\mathbf{E}+\partial_t\mathbf{B}=0,\\
&\nabla\times\mathbf{B}-\mu_0\varepsilon_0\partial_t\mathbf{E}=0,\\
&\nabla\cdot\mathbf{E}=0,\\
&\nabla\cdot\mathbf{B}=0,
\end{split}
\end{equation}
where $\mathbf{E}=(E_1,E_2,E_3)$ and $\mathbf{B}=(B_1,B_2,B_3)$ are the electric and magnetic fields. We can write this first order system for $\mathbf{E}$ and $\mathbf{B}$ as a second order system for $\mathbf{E}$ alone. We have, using (\ref{eq4-2}), that
\[\mu_0\varepsilon_0\partial_t\mathbf{E}=\nabla\times\mathbf{B},\]
\[\Downarrow\]
\[\mu_0\varepsilon_0\partial_{tt}\mathbf{E}=\nabla\times\partial_t\mathbf{B}=\nabla\times(-\nabla\times \mathbf{E}),\]
\[\Downarrow\]
\[\mu_0\varepsilon_0\partial_{tt}\mathbf{E}+\nabla\times\nabla\times\mathbf{E}=0,\]
and
\[\nabla\times\nabla\times \mathbf{E}=\nabla(\nabla\cdot\mathbf{E})-\nabla^2\mathbf{E}=-\nabla^2\mathbf{E},\]
\[\Downarrow\]
\[\varepsilon_0\mu_0\partial_{tt}\mathbf{E}-\nabla^2\mathbf{E}=0.\]
This is the 3D wave equation and it predicts, as we will see later in this course, that there are solutions describing electromagnetic disturbances traveling at a speed $v=\frac1{\sqrt{\varepsilon_0\mu_0}}$ in source free regions (for example in a vacuum). The parameters $\varepsilon_0, \mu_0$ are constants of nature that was known to Maxwell, and when he calculated $v$, he found that the speed of these disturbances was very close to the known speed of light.

He thus conjectured that light is an electromagnetic disturbance. This was a monumental discovery whose ramifications we are still exploring today.

\subsection{First order linear scalar PDEs}
A general first order linear scalar PDE with two independent variables, $(x,t)$, is of the form
\begin{equation}\label{eq4-3}
a\partial_xv+b\partial_tv=cv+d,
\end{equation}
where in general, $a,b,c$ and $d$ depends on $x$ and $t$. The solution method for (\ref{eq4-3}) is based on a geometrical interpretation of the equation.

Observe that the two functions $a(x,t)$ and $b(x,t)$ can be thought of as the two components of a vector field $\mathbf{f}(x,t)$ in the plane.
\[\mathbf{f}(x,t)=(a(x,t),b(x,t)).\]
Let $\boldsymbol{\gamma}(s)=(x(s),t(s))$ be a parametrized curve in the plane whose velocity at each point $(x,t)$ is equal to $\mathbf{f}(x,t)$.
\begin{figure}[h]
\centering
\includegraphics{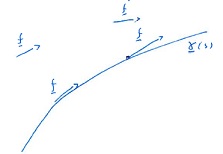}
\caption{}
\end{figure}

\noindent Thus we have
\begin{equation}\label{eq4-4}
\begin{split}
\frac{\mathrm{d}\boldsymbol{\gamma}}{\mathrm{d}s}&=\mathbf{f}(\boldsymbol{\gamma}(s)),\\
&\Updownarrow\\
\frac{\mathrm{d}x}{\mathrm{d}s}&=a,\\
\frac{\mathrm{d}t}{\mathrm{d}s}&=b.
\end{split}
\end{equation}
Let $v(x,t)$ be a solution of (\ref{eq4-3}) and define 
\[v(s)=v(\boldsymbol{\gamma}(s))=v(x(s),t(s)).\]
Then 
\begin{align*}
\frac{\mathrm{d}v}{\mathrm{d}s}&=\partial_x v\frac{\mathrm{d}x}{\mathrm{d}s}+\partial_tv\frac{\mathrm{d}t}{\mathrm{d}s}\\
{}&=a\partial_xv+b\partial_tv=cv+d.
\end{align*}
Thus $v(s)$ satisfy an ODE!
\begin{equation}\label{eq4-5}
\frac{\mathrm{d}v(s)}{\mathrm{d}s}=c(x(s),t(s))v(s)+d(x(s),t(s)).
\end{equation}
From (\ref{eq4-4}) and (\ref{eq4-5}) we get a coupled system of ODEs
\begin{equation}\label{eq4-6}
\begin{split}
\frac{\mathrm{d}x}{\mathrm{d}s}&=a,\\
\frac{\mathrm{d}y}{\mathrm{d}s}&=b,\\
\frac{\mathrm{d}v}{\mathrm{d}s}&=cv+d.
\end{split}
\end{equation}
The general theory of ODEs tells us that through each point $(x_0,t_0,v_0)$, there goes exactly one solution curve of (\ref{eq4-6}) if certain regularity conditions on $a,b,c$ and $d$ are satisfied.
\begin{figure}[h]
\centering
\includegraphics{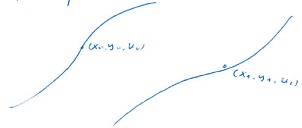}
\caption{: Characteristic curves for the PDE (\ref{eq4-3}).}
\end{figure}
The family of curves that solves (\ref{eq4-6}) are called {\it the characteristic curves} for the PDE (\ref{eq4-3}).

Observe that the system (\ref{eq4-6}) is partly decoupled, the two first equations can be solved independently of the third. The plane curves solving
\begin{align*}
\frac{\mathrm{d}x}{\mathrm{d}s}&=a,\\
\frac{\mathrm{d}t}{\mathrm{d}s}&=b,
\end{align*}
are called {\it characterictic base curves}. They are also, somewhat imprecisely, called characteristic curves.

Most of the time we are not interested in the general solution of (\ref{eq4-3}), but rather in special solutions that in addition to (\ref{eq4-3}) has given values on some curve $C$ in the $(x,t)$ plane.
\begin{figure}[h]
\centering
\includegraphics{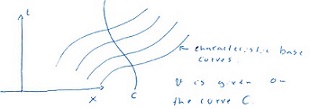}
\caption{: Characteristic base curves for the PDE (\ref{eq4-3}).}
\end{figure}
Any such problem is called an initial value problem for our equation (\ref{eq4-3}).

This is a generalization of the notion of initial value problem for a PDE from what we have discussed earlier. This earlier notion consisted of giving the solution at $t=0$ (or $t=t_0$)
\[v(x,0)=\varphi(x),\]
where $\varphi(x)$ is some given function of $x$.

Let $(x(\tau),y(\tau))$ be a parametrization of $C$. Then our generalized notion of initial value problem for (\ref{eq4-3}) consists of imposing the condition
\[v(x(\tau),t(\tau))=\varphi(\tau),\]
where $\varphi(\tau)$ is a given function. Our previous notion of initial value problem corresponds to the choice.
\[(x(\tau),t(\tau))=(\tau,0).\]
Since this is a parametrization of the x-axis in the $(x,t)$ plane. Thus in summary, the (generalized) initial value problem for (\ref{eq4-3}) consists of finding functions
\begin{equation}\label{eq4-7}
\begin{split}
x&=x(s,\tau),\\
t&=t(s,\tau),\\
v&=v(s,\tau),
\end{split}
\end{equation}
such that
\begin{equation}\label{eq4-8}
\begin{split}
\frac{\mathrm{d}x}{\mathrm{d}s}&=a,\,\quad\quad \quad x(0,\tau)=x(\tau),\\
\frac{\mathrm{d}t}{\mathrm{d}s}&=b,\,\quad\quad \quad t(0,\tau)=t(\tau),\\
\frac{\mathrm{d}v}{\mathrm{d}s}&=cv+d,\quad v(0,\tau)=\varphi(\tau).
\end{split}
\end{equation}
Where $(x(\tau),t(\tau))$ parametrize the initial curve and where $\varphi(\tau)$ are the given values of $v$ on the initial curve. The reader might at this point object and ask how we can possibly consider (\ref{eq4-7}), solving (\ref{eq4-8}), to be a solution to our original problem (\ref{eq4-3}). Where is my $v(x,t)$ you might ask? This is after all what we started out to find.

In fact, the function $v(x,t)$, solving our equation (\ref{eq4-3}), is implicitly contained in (\ref{eq4-7}).

We can in principle use the two first equations in (\ref{eq4-7}) to express $s$ and $\tau$ as functions of $x$ and $t$.
\begin{equation}\label{eq4-9}
\begin{split}
x&=x(s,\tau),\\
t&=t(s,\tau),\\
&\quad\Updownarrow\\
s&=s(x,t),\\
\tau&=\tau(x,t).\\
\end{split}
\end{equation}
It is known that this can be done locally close to a point $(s_0,\tau_0)$ if the determinant of the Jacobian matrix at $(s_0,\tau_0)$, this is what is known as the {\it Jacobian} in calculus, is nonzero.
\[\left.(\partial_sx\,\partial_{\tau} t-\partial_st\,\partial_{\tau}x)\right|_{(s_0,\tau_0)}\neq0.\]
Thus if the Jacobian is nonzero for all points on the initial curve the initial value problem for (\ref{eq4-3}) has a unique solution for any specification of $\varphi$ on the initial curve. This solution is
\begin{equation*}
v(x,t)=v(s(x,t),\tau(x,t)).
\end{equation*}

\noindent Actually finding a formula for the solution $v(x,t)$ using this method, which is called {\it the method of characteristics}, can be hard. This is because in general it is hard to find explicit solutions to the coupled ODEs (\ref{eq4-8}) and even if you succeed with this, inverting the system (\ref{eq4-9}) analytically might easily be impossible.

However, even if we can not find an explicit solution, the method of characteristics can still give great insight into what the PDE is telling us. Additionally, the method of characteristic is the starting point for a numerical approach to solving PDEs of the type (\ref{eq4-3}) and of related types as well.

In our traffic model, we encountered the equation
\begin{equation}\label{traffic}
\partial_tn+c_0\partial_xn=0.
\end{equation}
This equation is of the type (\ref{eq4-3}) and can be solved using the method of characteristic.
\noindent Let the initial data be given on the x-axis.
\[h(x,0)=\varphi(x).\]
Our system of ODEs determining the base characteristic curves for the traffic equation (\ref{traffic}) is thus
\begin{align*}
\frac{\mathrm{d}t}{\mathrm{d}s}&=1, \quad t(0,\tau)=0,\\
\frac{\mathrm{d}x}{\mathrm{d}s}&=c_0, \quad x(0,\tau)=\tau,\\
\frac{\mathrm{d}n}{\mathrm{d}s}&=0,\quad n(0,\tau)=\varphi(\tau).
\end{align*}
Observe that the last ODE tell us that the function $n$ is constant along the base characteristic curves.
\begin{align*}
\frac{\mathrm{d}t}{\mathrm{d}s}&=1\,\Rightarrow\, t=s+t_0,\,\,t(0,\tau)=0 \,\Rightarrow\, t_0=0\,\Rightarrow\,t(x,\tau)=s,\\
\frac{\mathrm{d}x}{\mathrm{d}s}&=c_0\,\Rightarrow\,x=c_0 s+x_0,\,x(0,\tau)=\tau\,\Rightarrow\,x_0=\tau\\
{ }&\Rightarrow\,x(s,\tau)=c_0s+\tau,\\
\frac{\mathrm{d}n}{\mathrm{d}s}&=0\,\Rightarrow\,n=n_0,n(0,\tau)=\varphi(\,\tau)\Rightarrow\,n_0=\varphi(\tau)\\
&\Rightarrow n(s,\tau)=\varphi(\tau).
\end{align*}
Now we must express $s$ and $\tau$ in terms of $x$ and $t$
\begin{align*}
t&=s\,\quad\quad\quad\Rightarrow s=t,\\
x&=c_0s+\tau\,\,\,\,\Rightarrow \tau=x-c_0s=x-c_0t,
\end{align*}
and our solution to the initial value problem is 
\[n(x,t)=\varphi(\tau)=\varphi(x-c_0t).\]
This is the unique solution to our initial value problem.
The base characteristic curves are given in parametrized form as
\begin{align*}
t&=s,\\
x&=c_0s+\tau,\quad\tau\in\mathbb{R},
\end{align*}
or in unparametrized form as
\[x-c_0t=\tau.\]
In either case they are a family of straight lines that is parametrized by $\tau$.
\begin{figure}[!h]
\centering
\includegraphics{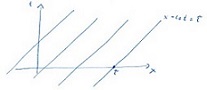}
\caption{: Base characteristic curves for the traffic equation (\ref{traffic}).}
\end{figure}

\noindent The solution in general depends on both the choice of initial curve and which values we specify for the solution on this curse. In fact, if we are not careful here, there might be many solutions to the initial value problem or none at all.

Since $n$ is constant along a base characteristic curve, its value is determined by the initial value on the point of the initial curve where the base characteristic curve and the initial curve coincide. If this occurs at more than one point there is in general no solution to the initial value problem.

\noindent Let for example $C$ be given through a parametrization
\begin{align}\label{parabola}
x(\tau)&=\tau,\\
t(\tau)&=\tau(1-\tau).
\end{align}
The initial values is some function $\varphi(\tau)$ defined on $C$. 
The initial curve is thus a parabola and the base characteristic curves form a family of parallel lines
\[x-c_0t=\text{const}\]
\begin{figure}[!h]
\centering
\includegraphics{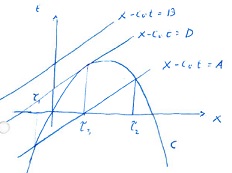}
\caption{: Base characteristic curves and parabolic initial curve for the traffic equation (\ref{traffic}).}
\end{figure}

\noindent Since $x-c_0t=A$ intersect the initial curve at two points corresponding to $x=\tau_1$ and $x=\tau_2$ and since $n=\varphi(\tau)$ is constant on the base characteristic curves, we must have 
\[\varphi(\tau_1)=\varphi(\tau_2).\]
This only holds for very special initial data. At a certain point $x=\tau_3, \,t=\tau_3(1-\tau_3)$, the characteristic base curve is tangent to the initial curve.

At this point the equations
\begin{align*}
x&=x(s,\tau),\\
t&=t(x,\tau),
\end{align*}
can not be inverted, even locally.

 In order to see this, we observe that the base characteristic curves for the traffic equation (\ref{traffic}), with initial data prescribed on the parabola (\ref{parabola}), are given by
\begin{align*}
x&=c_0s+\tau,\\
t&=s+\tau(1-\tau).
\end{align*}
The Jacobian is 
\[\Delta=\partial_s x\,\partial_{\tau}t-\partial_st\,\partial_{\tau}x=c_0 (1-2\tau)-1,\]
and at the point of tangency we have
\[\frac1c_0=\frac{\mathrm{d}}{\mathrm{d}\tau}(\tau(1-\tau))=(1-2\tau),\]
\[\Downarrow\]
\[c_0(1-2\tau)=1,\]
\[\Downarrow\]
\[\Delta=0.\]
A point where the base characteristic curves and the initial curve are tangent is called a {\it characteristic point}. The initial curve is said to be {\it characteristic} at this point. Thus the initial curve $C$ is characteristic at the point $(x,t)$ where
\begin{eqnarray*}
x=\frac{1}{2}(1-\frac{1}{c_0}),\\
t=\frac{1}{4}(1-\frac{1}{c_0^2}).
\end{eqnarray*}
 If the initial curve is characteristic at all of  its points, we say that the initial curve is characteristic. It in fact {\it is} a base characteristic curve. The initial value problem then in general has no solution, and if it does have a solution, the solution is not unique. A characteristic initial value problem for our traffic equation (\ref{traffic}) would be
\[\partial_tn+c_0\partial_xn=0,\]
\[n(c_0t,t)=\varphi(t).\]
General solution of the equation is
\begin{equation}\label{nonunique} 
n(x,t)=f(x-c_0t).
\end{equation}
The initial values require that
\[n(c_0t,t)=\varphi(t),\]
\[\Updownarrow\]
\[f(0)=\varphi(t),\quad \forall t.\]
In order for a solution to exist, $\varphi(t)=A$ must be a constant. And  if $\varphi(t)=A$ then (\ref{nonunique}), with {\it any} function $f(y)$ where $f(0)=A$, would solve the initial value problem. Thus we have nonuniqueness.

Let us next consider the initial value problem
\begin{eqnarray}
t\partial_tv+\partial_xv&=&0,\nonumber\\
v(x,1)&=&\varphi(x)\label{prob2}.
\end{eqnarray}
The initial curce is the horizontal line $t=1$ which can be parametrized using
\begin{eqnarray*}
x(\tau)&=\tau,\\
t(\tau)&=1.
\end{eqnarray*}
Applying the method of characteristics, we get the following system of ODEs
\begin{align*}
\frac{\mathrm{d}t}{\mathrm{d}s}&=t,\quad t(0,\tau)=1,\\
\frac{\mathrm{d}x}{\mathrm{d}s}&=1,\quad x(0,\tau)=\tau,\\
\frac{\mathrm{d}v}{\mathrm{d}s}&=0,\quad v(0,\tau)=\varphi(\tau),
\end{align*}
\noindent and
\begin{align*}
\frac{\mathrm{d}t}{\mathrm{d}s}&=t\,\Rightarrow\,t=t_0e^s,\,t(0,\tau)=1\,\Rightarrow\,t_0=1\,\Rightarrow\,t=e^s,\\
\frac{\mathrm{d}x}{\mathrm{d}s}&=1\,\Rightarrow\,x=s+x_0,\,x(0,\tau)=\tau\,\Rightarrow\,x_0=\tau\,\Rightarrow\,x=s+\tau,\\
\frac{\mathrm{d}v}{\mathrm{d}s}&=0\,\Rightarrow\,v=v(\tau),\,v(0,\tau)=\varphi(\tau)\,\Rightarrow\,v=\varphi(\tau).
\end{align*}
We now solve for $s$ and $\tau$ in terms of $x$ and $t$
\begin{align*}
t&=e^s,\\
x&=s+\tau.
\end{align*}
The first equation only have solutions for $t>0$
\begin{align*}
s&=\ln t,\\
\tau&=x-s=x-\ln t,
\end{align*}
\[\Downarrow\]
\[v(x,t)=\varphi(\tau)=\varphi(x-\ln t).\]
\pagebreak
This solution is constant along the base characteristic curves
\[\tau=x-\ln t,\]
\[\Updownarrow\]
\[t=e^{{x-\tau}}=e^{-\tau}e^x,\quad \tau\in\mathbb{R}.\]
\noindent The solution we have found exists only for $t>0$. This is its domain of definition
\begin{figure}[!h]
\centering
\includegraphics{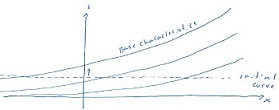}
\caption{: Base characteristic curves and initial curve for the initial value problem (\ref{prob2}). }
\label{4-6}
\end{figure}

As a final example, let us consider the initial value problem for the wave equation
\begin{eqnarray}
v_{tt}-\gamma^2 v_{xx}&=&0,\nonumber\\
v(x,0)&=&f(x),\nonumber\\
v_t(x,0)&=&g(x).\label{prob3}
\end{eqnarray}
We have previously seen that this can be written as a first order system
\begin{equation}\label{4-1}
v_t-\gamma v_x=u,
\end{equation}
\begin{equation}\label{4-2}
u_t+\gamma u_x=0.
\end{equation}
The initial conditions are 
\begin{align*}
v(x,0)&=f(x),\\
u(x,0)&=(v_t-\gamma v_x)(x,0)=g(x)-\gamma f'(x).
\end{align*}
We start by solving the homogeneous equation (\ref{4-2}). We will use the method of characteristics on the problem
\[u_t+\gamma u_x=0,\]
\[u(x,0)=g-\gamma f'.\]
The corresponding system of ODEs is
\begin{align*}
\frac{\mathrm{d}x}{\mathrm{d}s}&=\gamma,\quad x(0,\tau)=\tau,\\
\frac{\mathrm{d}t}{\mathrm{d}s}&=1,\quad t(0,\tau)=0,\\
\frac{\mathrm{d}u}{\mathrm{d}s}&=0,\quad u(0,\tau)=g(\tau)-\gamma f'(\tau).
\end{align*}
Solving we find
\begin{align*}
\frac{\mathrm{d}x}{\mathrm{d}s}&=\gamma\,\Rightarrow\,x=\gamma s+x_0,\,x(0,\tau)=\tau\,\Rightarrow\,x_0=\tau\,\Rightarrow\,x=\gamma s+\tau,\\
\frac{\mathrm{d}t}{\mathrm{d}s}&=1 \,\Rightarrow\,t=s+t_0,\,t(0,\tau)=0\,\Rightarrow\,t_0=0\,\Rightarrow\,t=s,\\
\frac{\mathrm{d}u}{\mathrm{d}s}&=0\,\Rightarrow\,u=u_0,\,u(0,\tau)=g(\tau)-\gamma f'(\tau)\,\Rightarrow\,u_0=g(\tau)-\gamma f'(\tau)\\
&\Rightarrow\,u(s,\tau)=g(\tau)-\gamma f'(\tau).
\end{align*}
We now express $\tau$ and $s$  as function of $x$ and $t$  
\begin{equation*}
\begin{split}
x&=\gamma s+\tau,\\
t&=s,
\end{split}
\,\,\Leftrightarrow\,\,
\begin{split}
\tau&=x-\gamma t,\\
s&=t,
\end{split}
\end{equation*}
from which it follows that the function $u(x,t)=u(s(x,t),\tau(x,t))$ is given by
\[u(x,t)=g(x-\gamma t)-\gamma f'(x-\gamma t).\]
Let us next apply the method of characteristic to the initial value problem for equation (\ref{4-1}). This problem is 
\[v_t-\gamma v_x=g(x-\gamma t)-\gamma f'(x-\gamma t),\]
\[v(x,0)=f(x).\]
The corresponding system of ODEs is
\begin{align*}
\frac{\mathrm{d}x}{\mathrm{d}s}&=-\gamma, &x(0,\tau)=\tau,\\
\frac{\mathrm{d}t}{\mathrm{d}s}&=1,&t(0,\tau)=0,\\
\frac{\mathrm{d}v}{\mathrm{d}s}&=g(x(s)-\gamma t(s))-\gamma f'(x(s)-\gamma t(s)),&v(\tau,0)=f(\tau).
\end{align*}
The first two equations give
\begin{equation*}
\begin{split}
x&=-\gamma s+\tau,\\
t&=s,
\end{split}
\,\,\Leftrightarrow\,\,
\begin{split}
\tau&=x+\gamma t,\\
s&=t.
\end{split}
\end{equation*}
The last equation is solved by
\[v(s,\tau)=\int_0^s \mathrm{d}\sigma(g(x(\sigma)-\gamma t(\sigma))-\gamma f'(x(\sigma)-\gamma t(\sigma)))+v_0(\tau),\]
and using the initial data we have
\[v(0,\tau)=f(\tau)\,\Rightarrow\,v_0(\tau)=f(\tau),\]
\[\Downarrow\]
\[v(s,\tau)=f(\tau)+\int_0^s\mathrm{d}\sigma(g(x(\sigma)-\gamma t(\sigma))-\gamma f'(x(\sigma)-\gamma t(\sigma))).\]
Observing that
\[x-\gamma t=-\gamma s+\tau-\gamma s=-2\gamma s+\tau,\]
the formula for $v$ become
\[v(s,\tau)=f(\tau)+\int_0^s\mathrm{d}\sigma(g(-2\gamma\sigma+\tau)-\gamma f'(-2\gamma\sigma+\tau)).\]
Changing variable in the integral using  $y=-2\gamma \sigma+\tau$, we get
\[v(s,\tau)=f(\tau)-\frac1{2\gamma}\int_\tau^{-2\gamma s+\tau}\mathrm{d}y(g(y)-\gamma f'(y)),\]
and $s=t,\,\tau=x+\gamma t$ from which it follows that
\[-2\gamma s+\tau=-2\gamma t+x+\gamma t=x-\gamma t.\]
Therefore, the function $v(x,t)=v(s(x,t),\tau(x,t))$ is given by the formula
\begin{align*}
v(x,t)&=f(x+\gamma t)-\frac1{2\gamma}\int_{x+\gamma t}^{x-\gamma t}\mathrm{d}y(g(y)-\gamma f'(y))\\
&=f(x+\gamma t)+\frac1{2}\int_{x+\gamma t}^{x-\gamma t}\mathrm{d}yf'(y)+\frac1{2\gamma}\int_{x-\gamma t}^{x+\gamma t}\mathrm{d}yg(y),
\end{align*}
\[\Downarrow\]
\[v(x,t)=\frac12(f(x+\gamma t)+f(x-\gamma t))+\frac1{2\gamma}\int_{x-\gamma t}^{x+\gamma t}\mathrm{d}yg(y).\]
This is the famous D'Alembert formula for the solution of the initial value problem for the wave equation on an infinite domain. It was first derived by the French mathematician D'Alembert in 1747.

Let $G(\xi)=\int g(y)\mathrm{d}y$. Then the formula can be rewritten as 

\[v(x,t)=\left[\frac12f(x+\gamma t)+\frac12G(x+\gamma t)\right]+\left[\frac12f(x-\gamma t)-\frac1{2\gamma}G(x-\gamma t)\right].\]

\noindent The solution is clearly a linear superposition of a left moving and a right moving wave.

\subsection{First order Quasilinear scalar PDEs}
A first order, scalar, quasilinear PDE is an equation of the general form
\begin{equation}\label{eq6-1}
a\partial_xu+b\partial_tu=c,
\end{equation}
where
\begin{align*}
a&=a(x,t,u),\\
b&=b(x,t,u),\\
c&=c(x,t,u).
\end{align*}
This equation can also be solved through a geometrical interpretation of the equation.

Let $u(x,t)$ be a solution of (\ref{eq6-1}).The {\it graph} of $u$ is a surface in $(x,t,u)$ space defined by the equation.
\[F(x,t,u)\equiv u(x,t)-u=0.\]
This surface is called an {\it integral surface} for the equation.

From calculus we know that the gradient of $F$
\[\nabla F=(\frac{\partial u}{\partial x},\frac{\partial u}{\partial t},-1),\]
is normal to the integral surface $F$. Observe that the equation (\ref{eq6-1}) can be rewritten as 
\[a\partial_x u+b\partial_tu-c=0,\]
\[\Updownarrow\]
\[(a,b,c)\cdot\nabla F(x,t,u)=0.\]
Thus the vector $(a,b,c)$ must always be in the tangent plane to the integral surface corresponding to a solution $u(x,t)$.

At each point $(x,t,u)$ the vector $(a,b,c)$ determines a direction that is the {\it characteristic direction} at this point. Thus the vector $(a,b,c)$ determines a {\it characteristic direction field} in $(x,t,u)$ space.

The characteristic curves are the curves that has $(a,b,c)$ as velocity at each of its points. Thus, if the curve is parametrized by $(x(s),t(s),u(s))$, we must have
\begin{align*}
\frac{\mathrm{d}x}{\mathrm{d}s}&=a(x,t,u),\\
\frac{\mathrm{d}t}{\mathrm{d}s}&=b(x,t,u),\\
\frac{\mathrm{d}u}{\mathrm{d}s}&=c(x,t,u).
\end{align*}
Observe that for this quasilinear case the two first equations does not decouple and the 3 ODEs has to be solved simultaneously. This can obviously be challenging.

For sufficiently smooth $a,b$ and $c$, there is a solution curve, or {\it characteristic curve}, passing through any point $(x_0,t_0,u_0)$ in $(x,t,u)$ space.

The initial value problem is specified by giving the value of $u$ on a curve in $(x,t)$ space. This curve together with the initial data determines a curve, $C$, in $(x,t,u)$ space that we call {\it the initial curve}. To solve the initial value problem we pass a characteristic curve through every point of the initial curve. If these curves generate a smooth surface, the function, whose graph is equal to the surface, is the solution of the initial value problem.
\begin{figure}[h]
\centering
\includegraphics{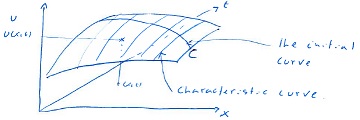}
\caption{}
\end{figure}

\noindent Let the initial curve $C$ be parametrized by 
\[(x(\tau),t(\tau),\varphi(\tau)).\]
The characteristic system of ODEs with initial conditions imposed is
\begin{align*}
\frac{\mathrm{d}x}{\mathrm{d}s}&=a,\,x(0,\tau)=x(\tau),\\
\frac{\mathrm{d}t}{\mathrm{d}s}&=b,\,t(0,\tau)=t(\tau),\\
\frac{\mathrm{d}u}{\mathrm{d}s}&=c,\,u(0,\tau)=\varphi(\tau).
\end{align*}
The theory of ODEs gives a unique solution
\begin{equation}\label{eq6-2}
\begin{split}
x&=x(s,\tau),\\
t&=t(s,\tau),\\
u&=\varphi(s,\tau).
\end{split}
\end{equation}
Observe that 
\[\partial_s x=a,\quad \partial_s t=b,\]
so the Jacobian is 
\[\Delta(s,\tau)=
\begin{array}{|cc|}
\partial_s x&\partial_\tau x\\
\partial_s t&\partial_\tau t
\end{array}
=(a\partial_\tau t-b\partial_\tau x)(s,\tau).\]
Let us assume that the Jacobian is nonzero along a smooth initial curve. We also assume that $a,b$ and $c$ are smooth. Then $\Delta(s,\tau)$ is continuous in $s$ and $\tau$, and therefore nonzero in an open set, surrounding the initial curve.

By the inverse function theorem, the equations
\begin{align*}
x&=x(s,\tau),\\
t&=t(s,\tau),
\end{align*}
can be inverted in an open set surrounding the initial curve. Thus we have 
\begin{align*}
s&=s(x,t),\\
\tau&=\tau(x,t).
\end{align*}
Then
\[u(x,t)=\varphi(s(x,t),\tau(x,t))\]
is the unique solution to the initial value problem
\begin{eqnarray*}
a\partial_xu+b\partial_tu&=&c,\\
u(x(\tau),t(\tau))&=&\varphi(\tau).
\end{eqnarray*}

\noindent In general, going through this program and finding the solution $u(x,t)$ can be very hard. However deploying the method can give great insight into a given PDE.

Let us consider the following quasilinear equation
\begin{eqnarray}
u_t+uu_x&=&0,\nonumber\\
u(x,0)&=&\varphi(x).\label{eq6-3}
\end{eqnarray}
\noindent This equation occurs as a first approximation to many physical phenomena. We have seen how it appears in the modelling of traffic.

  It also describe waves on the surface of the ocean. If we assume that the waves moves in one direction only, then (\ref{eq6-3}) appears as a description of the deviation of the ocean surface from a calm, flat state.
\begin{figure}[h]
\centering
\includegraphics{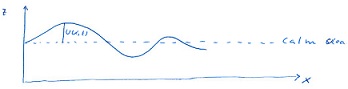}
\caption{: Ocean surface waves.}
\end{figure}

\noindent It is called the (inviscid) Burger's equation. Let us solve Burger's equation using the method of characteristics.
Our initial value problem is 
\begin{align*}
u_t+uu_x&=0,\\
u(x,0)&=\varphi(x).
\end{align*}
The initial curve is parametrized by
\[x=\tau,\quad t=0,\quad u=\varphi(\tau).\]
Here 
\begin{align*}
a(x,t,u)&=u,\\
b(x,t,u)&=1,\\
c(x,t,u)&=0.
\end{align*}
Let us compute the Jacobian along the initial curve. If the characteristic curves are
\begin{align*}
x&=x(s,\tau),\\
t&=t(s,\tau),\\
u&=u(s,\tau),
\end{align*}
where
\begin{equation}\label{eq6-4}
\begin{split}
\frac{\mathrm{d}x}{\mathrm{d}s}&=u,\,\,x(0,\tau)=\tau,\\
\frac{\mathrm{d}t}{\mathrm{d}s}&=1,\,\,t(0,\tau)=0,\\
\frac{\mathrm{d}u}{\mathrm{d}s}&=0,\,\,u(0,\tau)=\varphi(\tau).
\end{split}
\end{equation}
Then 
\[\partial_\tau t=0,\quad \partial_\tau x=1\,\,\,\text{at}\,\,\,(0,\tau),\]
\[\Downarrow\]
\[\Delta(0,\tau)=\left.(a\partial_\tau t-b\partial_\tau x)\right|_{(0,\tau)}=-1\neq0,\]
so $\Delta\neq0$ along the entire initial curve. Let us next solve (\ref{eq6-4}). From the last equation in (\ref{eq6-4}) we get
\[u=u_0.\]
Thus $u$ is a constant independent of $s$. Since $u$ is a constant, the second equation in (\ref{eq6-4}) give us 
\[x=u_0s+x_0,\]
and 
\begin{align*}
u(0,\tau)&=\varphi(\tau)\quad \Rightarrow u_0=\varphi(\tau),\\
x(0,\tau)&=\tau\,\, \quad\quad \Rightarrow x_0=\tau,
\end{align*}
\[\Downarrow\]
\begin{align*}
u&=\varphi(\tau),\\
x&=\varphi(\tau)s+\tau.
\end{align*}
Finally for $t$ we get in the usual way
\[t=s.\]
If the conditions for the inverse function theorem holds, the equations
\begin{align*}
x&=\varphi(\tau)s+\tau,\\
t&=s,
\end{align*}
can be inverted to yield
\[\tau=\tau(x,t),\quad s=t,\]
and the solution to our problem is 
\[u(x,t)=\varphi(\tau(x,t)).\]
Without specifying $\varphi$, we can not invert the system and find a solution to our initial value problem. We can however write down an {\it implicit} solution.
We have
\begin{align*}
u&=\varphi(\tau),\\
x&=\varphi(\tau)s+\tau,\\
t&=s,\\
&\Downarrow\\
\tau&=x-us,\\
&\Downarrow\\
u&=\varphi(x-tu).
\end{align*}
Even without precisely specifying $\varphi$ we can still gain great insight by considering the way the characteristics depends on the initial data. Projection of the characteristics in $(x,t,u)$ space down into $(x,t)$ space give us base characteristics of the form
\[x=\varphi(\tau)t+\tau.\]
This is a family of straight lines in $(x,t)$ space parametrized by the parameter $\tau$. Observe that the slope of the lines depends on $\tau$, thus the lines in the family are {\it not} parallel in general. This is highly significant. 
The slope of the base characteristics is 
\[\frac{\mathrm{d}t}{\mathrm{d}x}=\frac1{\varphi(\tau)},\]
and is evidently varying as a function of $\tau$. This implies the possibility that the base characteristics will intersect. If they do, we are in trouble because the solution, $u$, is constant along the base characteristics. Thus if the base characteristics corresponding to parameters $\tau_1$ and $\tau_2$ intersect at a point $(x^*,t^*)$ in $(x,t)$ space then at this point $u$ will have two values
\begin{align*}
u_1&=\varphi(\tau_1),\\
u_2&=\varphi(\tau_2),
\end{align*}
that are different unless it happened to be the case that 
\[\varphi(\tau_1)=\varphi(\tau_2).\]
But in this case the two characteristics have the same sloop and would not have intersected.
The conclusion is that base characteristics can intersect and when that happens $u$ is no longer a valid solution.

\noindent There is another way of looking at this phenomenon. Recall that the implicit solution is
\begin{equation}\label{eq6-5}
u=\varphi(x-ut).
\end{equation}
We know that a formula
\[v=\varphi(x-cv),\]
describe a wave moving at speed $c$ to the right.
The formula (\ref{eq6-5}) can then be interpreted to say that each point on the wave has a ``local'' speed $u(x,t)$ that vary along the wave. Thus the wave does not move as one unit but can deform as it moves.

If $u(x,t)>0$, the point $(x,u(x,t))$ on the wave move to the right and if $u(x,t)<0$, it moves to the left. If $u(x,t)$ is increasing, points on the wave move at greater and greater speed as we move to the right. The wave will stretch out. If $u(x,t)$ is decreasing the wave will compress. In a wave that has the shape of a localized ``lump'',
\begin{figure}[h]
\centering
\includegraphics{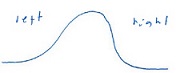}
\caption{: Initial conditions in the form of a symmetric lump.}
\end{figure}

\noindent there are points to the left that move faster than points further to the right. Thus points to the left can catch up with points to the right and lead to a steepening of the wave
\begin{figure}[h]
\centering
\includegraphics{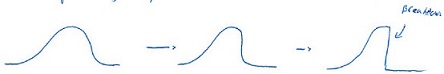}
\caption{: Steepning of initial lump over time.}
\end{figure}

\noindent Eventually the steepening can create a wave with a vertical edge. At this point $|u_x|=\infty$ and the model no longer describe what it is supposed to model. We thus have a {\it breakdown}. 

It is evidently important to find out if and when such a breakdown occurs. There are two ways to work this out. The first approach use the implicit solution.
\begin{equation}\label{eq6-6}
u=\varphi(x-ut).
\end{equation}
We ask when the slope $u_x$ become infinite. From (\ref{eq6-6}) we get
\[u_x=\varphi'(x-ut)(1-u_xt),\]
\[\Downarrow\]
\[u_x=\frac{\varphi'(x-ut)}{1+t\varphi'(x-ut)}.\]
Clearly $|u_x|$ become infinite if the denominator become zero
\[1+t\varphi'(x-ut)=0,\]
\[\Updownarrow\]
\[t=-\frac1{\varphi'(\tau)},\]
where we recall that $\tau=x-ut$. What is important is to find the smallest time for which infinite sloop occurs. The solution is not valid for times larger than this. Thus we must determine
\[t^*=\mathop{\mathrm{min}}_\tau\left(-\frac1{\varphi'(\tau)}\right).\]
When initial data is given, $\varphi(\tau)$ is known and $t^*$ can be found. 

The second method involve looking at the geometry of the characteristics. It leads to the same breakdown time $t^*$. 

We can immediately observe that if $\varphi$ is increasing for all $x$ so that
\[\varphi'(\tau)>0,\]
then $t^*<0$ and no breakdown occur. For breakdown to occur we must have $\varphi'(\tau)<0$ for some $\tau$.
Let us look at some specific choices of initial conditions.
\begin{enumerate}
\item  $u(x,0)=\varphi(x)=-x$
\[\!\!\!\!\!\!\!\!\!\!\Rightarrow \quad t(\tau)=-\frac1{\varphi'(\tau)}=1,\]
       \[\Downarrow\]
       \[t^*=\mathop{\mathrm{min}}_\tau \varphi(\tau)=1.\]
\end{enumerate}
Thus we have breakdown at $t=1$. The solution is 
\begin{eqnarray*}
&u=\varphi(x-ut)=-(x-ut)=-x+ut,\\
&\Downarrow\\
&u(1-t)=-x\\
&\Downarrow\\
&u(x,t)=\frac{x}{t-1}
\end{eqnarray*}
\begin{figure}[h]
\centering
\includegraphics{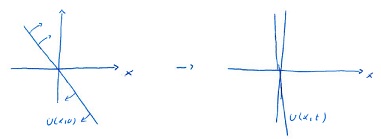}
\caption{: The solution $u(x,t)=\frac{x}{t-1}$ at $t=0$ and $t=t_1>0$.}
\label{fig6-5}
\end{figure}
\noindent The whole graph becomes vertical at the breakdown time $t^*=1$. This is illustrated in figure \ref{fig6-5}.
Let us describe this same case in terms of base characteristics. 
They are 
\[x+\tau t=\tau.\]
They all intersect the t-axis at $t=1$ and the base characteristics corresponding to $\tau$ intersect the x-axis at $x=\tau$. This case is illustrated in figure \ref{fig6-6}.
\begin{figure}[h]
\centering
\includegraphics{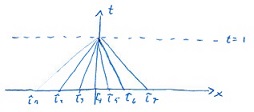}
\caption{: Intersecting base characteristics for the initial condition $u(x,0)=-x$}
\label{fig6-6}
\end{figure}
\begin{enumerate}
\addtocounter{enumi}{1}
\item $u(x,0)=\varphi(x)=1-x^2$.
\[\!\!\!\!\!\!\!\!\!\!\Rightarrow \quad t(\tau)=-\frac1{\varphi'(\tau)}=\frac{1}{2\tau-1},\]
       \[\Downarrow\]
       \[t^*=\mathop{\mathrm{min}}_\tau \varphi(\tau)=0.\]
Thus we have breakdown already at $t=0$!   
The implicit formula for the solution now is 
\[u=1-(x-ut)^2,\]
\[\Downarrow\]
\[u=1-x^2+2xut-u^2t^2,\]
\[\Downarrow\]
\[u^2t^2+(1-2xt)u+x^2-1,\]
\[\Downarrow\]
\[u=\frac1{2t^2}\left(-(1-2xt)\pm\left((1-2xt)^2-4t^2(x^2-1)^{\frac12}\right)\right),\]
\[\Downarrow\]
\[u=\frac{x}{t}-\frac{1\pm(1+4t(t-x))^{\frac12}}{2t^2}.\]
\end{enumerate}
The initial condition requires that 
\[\lim_{t\to 0^{+}} u(x,t)=1-x^2,\]
and for small $t$ we have
\[(1+4t(t-x))^{\frac14}\approx1+2t(t-x)-2t^2(t-x)^2.\]
Using the `\,-\,' solution, we have for small $t$
\begin{align*}
u&\approx \frac{x}{t}-\frac{1-1-2t(t-x)+2t^2(t-x)^2}{2t^2},\\
&=\frac{x}{t}-\frac1{2t^2}(-2t^2+2tx+2t^4-4t^3x+2t^2x^2),\\
&=\frac{x}{t}+1-\frac{x}{t}-x^2-2t^2x-t^2.
\end{align*}
Thus
\[u(x,t)\to1-x^2,\qquad \text{when}\,t\to 0.\]
Our unique solution of the initial value problem is therefore
\[u(x,t)=\frac{x}{t}-\frac{1-(1+4t(t-x))^{\frac12}}{2t^2}.\]
We have
\[u_x=\frac1t(1-(1+4t(t-x))^{-\frac12}),\]
and thus infinite sloop occurs at all space-time points $(x,t)$ determined by 
\[1+4t(t-x)=0,\]
\[\Updownarrow\]
\[x=t+\frac1{4t}.\]
This curve is illustarted in figure (\ref{fig6-7}).
There clearly are space-time points with arbitrarily small times on this curve and we therefore conclude that the breakdown time is $t^*=0$.
\begin{figure}[h]
\centering
\includegraphics{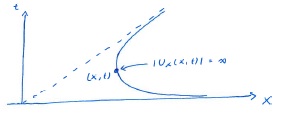}
\caption{: Set of space-time points where $u(x,t)$ has infinite sloop.}
\label{fig6-7}
\end{figure}

\begin{figure}[h]
\centering
\includegraphics{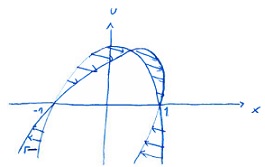}
\caption{: The solution $u(x,t)$ with $u(x,0)=1-x^2$ at $t=0$ and $t=t_1>0$ }
\label{fig6-8}
\end{figure}

\noindent The actual wave propagation is illustrated in figure \ref{fig6-8}.
The steepening is evident.
\begin{enumerate}
\addtocounter{enumi}{2}
\item $u(x,0)=\varphi(x)=\sin x.$
\end{enumerate}
The implicit solution is now
\[u=\sin(x-ut)\]
and thus
\begin{eqnarray*}
&\varphi(\tau)=\sin \tau,\\
&\Downarrow\\
& \varphi'(\tau)=\cos \tau,\\
&\Downarrow\\
&t^*=\mathop{\mathrm{min}}_\tau\left(-\frac1{\cos \tau}\right).
\end{eqnarray*}
From this formula it is evident that
\[t^*=1,\]
and this minimum breakdown time is achieved for each $\tau_n=(2n+1)\pi$.

\begin{figure}[h]
\centering
\includegraphics{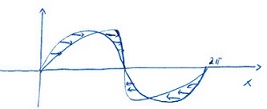}
\caption{: The solution $u(x,t)$ with $u(x,0)=\sin x$ at $t=0$ and $t=t_1>0$}
\label{fig6-9}
\end{figure}

  Recall that the characteristic curves are parametrized by $\tau$ through
\begin{equation*}
x=sin(\tau)t+\tau.
\end{equation*}
Thus we can conclude that the breakdown occurs at the space-time points $(t^*,x_n^*)$ where
\begin{equation*}
x^*=sin(\tau_n)t^*+\tau_n=(2n+1)\pi.
\end{equation*}
The breakdown is illustrated in figure \ref{fig6-9}.

\section{Classification of equations: Characteristics}
We will start our investigations by considering second order equations. More specifically we will consider scalar linear second order PDEs with two independent variables. The most general such equation is of the form
\begin{equation}\label{eq7-1}
Au_{xx}+2Bu_{xy}+Cu_{yy}+Du_x+Eu_y+Fu=G,
\end{equation}
where $u=u(x,y), A=A(x,y)$ etc.
The wave equation, diffusion equation and Laplace equation are clearly included in this class.

What separates one equation in this class from another are determined entirely by the functions $A\rightarrow G$. Thus the exact form of these functions determines all properties of the equation and the exact form of the set of solutions (or {\it space} of solutions as we usually say). In fact if $G=0$ the space of solutions is a {\it vectorspace}. This is because the differential operator
\[\mathscr{L}=A\partial_{xx}+2B\partial_{xy}+C\partial_{yy}+D\partial_x+E\partial_y+F,\]
is a {\it linear operator}. (Prove it!)
Thus even if we at this point knows {\it nothing} about the functions forming the space of solutions of the given equation, we do know that it has the {\it abstract property} of being a vectorspace.

Studying abstract properties of the space of solutions of PDEs is a central theme in the modern theory of PDEs. Very rarely will we be able to get an analytical description of the solution space for a PDE.

Knowing that a given PDE has a solution space that has a set of abstract properties is often of paramount importance. This knowledge can both inform us about expected qualitative properties of solutions and act as a guide for selecting which numerical or analytical approximation methods are appropriate for the given equation.

The aim of this section is to describe three abstract properties that the equation (\ref{eq7-1}) might have or not. These are the properties of being {\it hyperbolic}, {\it parabolic} or {\it elliptic}.

Let us then begin: We first assume that $A(x,y)$ is nonzero in the region of interest. This is not a very restrictive assumption and can be achieved for most reasonable $A$. We can then divide (\ref{eq7-1}) by $A$
\begin{equation}\label{eq7-2}
u_{xx}+\frac{2B}{A}u_{xy}+\frac CAu_{yy}+\frac DAu_x+\frac EAu_y+\frac FAu=\frac GA.
\end{equation}
We will now focus on the terms in (\ref{eq7-2}) that are of second order in the derivatives. These terms form the {\it principal part} of the equation. The differential operator defining the principal part is
\begin{equation}\label{eq7-3}
\mathscr{L}=\partial_{xx}+\frac{2B}{A}\partial_{xy}+\frac{C}{A}\partial_{yy}.
\end{equation}
We now try to factorize the differential operator by writing it in the form
\begin{equation}\label{eq7-4}
\mathscr{L}=(\partial_x-a\partial_y)(\partial_x-b\partial_y)+c\partial_y,
\end{equation}
where $a,b,c$ are some unknown functions. We find these functions by expanding out the derivatives and comparing with the expression (\ref{eq7-3})
\begin{align*}
\mathscr{L}&=\partial_{xx}-b_x\partial_y-b\partial_{xy}-a\partial_{xy}+ab_y\partial_y+ab\partial_{yy}+c\partial_y\\
&=\partial_{xx}-(a+b)\partial_{xy}+ab\partial_{yy}+(c-b_x+ab_y)\partial_y.
\end{align*}
Comparing with (\ref{eq7-3}) we get
\begin{align*}
a+b&=-\frac{2B}{A},\\
ab&=\frac CA,\quad\quad c=b_x-ab_y.
\end{align*}
From the first two equations we find that both $a$ and $b$ solves the equation
\begin{equation}\label{eq7-5}
Az^2+2Bz+C=0.
\end{equation}
Let us call these solutions $\omega^{\pm}$ so that 
\[a=\omega^+, b=\omega^-.\]
We thus have
\begin{equation}\label{eq7-6}
\omega^{\pm}=\frac1A\left\{-B\pm(B^2-AC)^{\frac12}\right\},
\end{equation}
and then
\begin{align*}
&\partial_{xx}+\frac{2B}A\partial_{xy}+\frac CA\partial_{yy}\\
=&(\partial_x-\omega^+\partial_y)(\partial_x-\omega^-\partial_y)+(\omega_x^--\omega^+\omega_y^-)\partial_y.
\end{align*}
We are in general interested in real solutions to our equations and aim to use this factorization to rewrite and simplify our equation. Thus it is of importance to know how many real solutions equation (\ref{eq7-5}) has. From formula (\ref{eq7-6}) we see that this is answered by studying the sign of $B^2(x,y)-A(x,y)C(x,y)$ in the domain of interest. We will here assume that this quantity is of one sign in the region of interest. This can always be achieved by restricting the domain, if necessary.

The equation (\ref{eq7-1}) is called
\begin{enumerate}
\item {\it Hyperbolic}\quad if\quad $B^2-AC>0,$
\item {\it Parabolic}\quad \,\,\,\,if\quad $B^2-AC=0,$
\item {\it Elliptic}\quad\quad \,\,\,\,if\quad $B^2-AC<0.$
\end{enumerate}

\noindent For the wave equation
\[u_{tt}-c^2u_{xx}=0,\]
we have $A=1,B=0,C=-c^2$ ($x\to t,y\to x$ in the general form of the equation (\ref{eq7-1}))
\[\Rightarrow B^2-AC=-1\cdot(-c^2)=c^2>0.\]
Thus the wave equation is {\it hyperbolic}. 

 For the diffusion equation we have
\[u_t-Du_{xx}=0,\]
we have $A=0,B=0,C=-D$
\[\Rightarrow B^2-AC=0.\]
Thus the diffusion equation is {\it parabolic}

 For the Laplace equation we have
\[u_{xx}+u_{yy}=0,\]
we have $A=1,B=0,C=1$
\[\Rightarrow B^2-AC=-1<0.\]
Thus the Laplace equation is {\it elliptic}.

We will now start deriving consequences of the classification into hyperbolic, parabolic and elliptic.

\subsection{Canonical forms for equations of hyperbolic type}
For this type of equations $\omega^+\neq \omega^-$ and both $\omega^+$ and $\omega^-$ are real. Take $\omega^+$ and write down the following ODE 
\begin{equation}\label{eq7-7}
\frac{\mathrm{d}y}{\mathrm{d}x}=-\omega^+.
\end{equation}
Let the solution curves of this equation be parametrized by the parameter $\xi$. Thus the curves can be written as
\[\xi=\xi(x,y).\]
In a similar way we can write the solution curves for the ODE
\begin{equation}\label{eq7-8}
\frac{\mathrm{d}y}{\mathrm{d}x}=-\omega^-,
\end{equation}
in the form
\[\eta=\eta(x,y),\]
where $\eta$ is the parameter parameterizing the solution curves for (\ref{eq7-8}). We now introduce $(\xi,\eta)$ as new coordinates in the plane. The change of coordinates is thus 
\begin{equation}\label{eq7-9}
\begin{split}
\xi&=\xi(x,y),\\
\eta&=\eta(x,y).
\end{split}
\end{equation}
For this to really be a coordinate system the Jacobian has to be nonzero.
\[\Delta=\xi_x\eta_y-\xi_y\eta_x.\]
The solution curves, or {\it characteristic curve}, of the equation (\ref{eq7-1}) are by definition determined by
\begin{align*}
\xi(x,y)&=\text{const},\\
\eta(x,y)&=\text{const}.
\end{align*}
Differentiate implicit with respect to $x$. This gives, upon using (\ref{eq7-7}) and (\ref{eq7-8}),
\begin{equation}\label{eq7-10}
\begin{split}
0&=\xi_x+y'\xi_y\,=\xi_x-\omega^+\xi_y \Rightarrow \xi_x=\omega^+\xi_y,\\
0&=\eta_x+y'\eta_y=\eta_x-\omega^-\eta_y\Rightarrow \eta_x=\omega^-\eta_y.
\end{split}
\end{equation}
Thus
\[\Delta=(\omega^+\xi_y)\eta_y-\xi_y(\omega^-\eta_y)=(\omega^+-\omega^-)\eta_y\xi_y\neq0,\]
because, by assumption, $\omega^+\neq \omega^-$ in our domain. Thus the proposed change of coordinates (\ref{eq7-9}) does in fact define a change of coordinates.
\begin{figure}[h]
\centering
\includegraphics{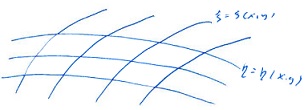}
\caption{}
\end{figure}

\noindent Recall that $z=\omega^+$ is a solutions to the quadratic equation (\ref{eq7-5}), thus 
\begin{equation}\label{eq7-11}
A(\omega^{+})^2+2B\omega^{+}+C=0.
\end{equation}
Inserting $\omega^+=\frac{\xi_x}{\xi_y}$ from equation (\ref{eq7-10}) we get, after multiplication by $\xi_y^2$
\[A(x,y)\xi_x^2(x,y)+2B(x,y)\xi_x(x,y)\xi_y(x,y)+C(x,y)\xi_y^2(x,y)=0.\label{eq7-11.1}\]
In a similar way the root $z=\omega^-$ show that the function $\eta(x,y)$ satisfy the equation
\[A(x,y)\eta_x^2(x,y)+2B(x,y)\eta_x(x,y)\eta_y(x,y)+C(x,y)\eta_y^2(x,y)=0.\label{eq7-11.2}\]
Thus both $\xi(x,y)$ and $\eta(x,y)$, that determines the characteristic coordinate system are solutions to the {\it characteristic PDE} corresponding to our original PDE (\ref{eq7-1})
\[A(x,y)\varphi_x^2(x,y)+2B(x,y)\varphi_x(x,y)\varphi_y(x,y)+x(x,y)\varphi_y^2(x,y)=0.\]
\noindent Using the chain rule we have
\begin{align*}
u_x&=\xi_x u_\xi+\eta_x u_\eta,\\
u_y&=\xi_y u_\xi +\eta_y u_\eta,
\end{align*}
\begin{equation}\label{eq7-12}
\begin{split}
\Downarrow\\
u_x-\omega^+u_y&=\xi_x u_\xi+\eta_x u_\eta-\omega^+\xi_yu_\xi-\omega^+\eta_yu_\eta\\
&=(\xi_x-\omega^+\xi_y)u_\xi+(\eta_x-\omega^+\eta_y)u_\eta\\
&=(\eta_x-\omega^+\eta_y)u_\eta=\eta_y(\omega^{-}-\omega^{+})u_\eta,
\end{split}
\end{equation}
because 
\[\xi_x-\omega^+\xi_y=\xi_x+y'(x)\xi_y=\frac{\mathrm{d}}{\mathrm{d}x}\xi(x,y(x))=0.\]
Since 
\[\xi(x,y)=\text{const},\]
describe the solution curves to the equation (\ref{eq7-7}). 
In a similar way we find
\begin{equation}\label{eq7-13}
 u_x-\omega^-u_y=\xi_y(\omega^{+}-\omega^{-})u_\xi.
\end{equation}
Therefore, using equations (\ref{eq7-10}), we have
\begin{equation}\label{eq8-14}
\begin{split}
\partial_x-\omega^+\partial_y&=-\eta_y(\omega^+-\omega^-)\partial_\eta,\\
\partial_x-\omega^-\partial_y&=\xi_y(\omega^+-\omega^-)\partial_\xi,
\end{split}
\end{equation}
and both operators are nonzero because $\omega^+\neq \omega^-$. But then we have
\begin{equation}\label{eq7-15}
\begin{split}
&(\partial_x-\omega^+\partial_y)(\partial_x-\omega^-\partial_y)=\left[-\eta_y(\omega^+-\omega^-)\partial_\eta\right]\left[\xi_y(\omega^+-\omega^-)\partial_\xi\right] \\
&=-\xi_y\eta_y(\omega^+-\omega^-)^2\partial_{\xi\eta}^2-\eta_y(\omega^+-\omega^-)\partial_\eta(\xi_y(\omega^+-\omega^-))\partial_\xi.
\end{split}
\end{equation}
Substituting this into (\ref{eq7-1}) and expressing all first order derivatives with respect to $x$ and $y$ in terms of $\partial_\xi$ and $\partial_\eta$, we get an equation of the form
\begin{equation}\label{eq7-16}
u_{\xi\eta}+au_\xi+bu_\eta+cu=d,
\end{equation}
where we have divided by the nonzero quantity $\xi_y\eta_y(w^+-w^-)$ and where now $u=u(\xi,\eta)$, $a=a(\xi,\eta)$ etc. For any particular equation $a,b,c$ and $d$ will be determined explicitely in terms of $A,B,C$ and $D$.

\noindent Since (\ref{eq7-16}) results from (\ref{eq7-1}) by a change of coordinates, (\ref{eq7-16}) and are {\it equivalent} PDEs. Their spaces of solutions will consist of functions that are related through the change of coordinates (\ref{eq7-9}).

(\ref{eq7-16}) is called {\it the canonical form} for equations of hyperbolic type. 

There is a second common form for hyperbolic equation. We find this by doing the additional transformation

\begin{equation*}
\begin{split}
\xi&=\alpha+\beta,\\
\eta&=\alpha-\beta,
\end{split}
\end{equation*}
\[\Downarrow\]
\begin{equation}\label{eq7-18}
u_{\alpha\alpha}-u_{\beta\beta}+\widetilde{a}u_\alpha+\widetilde{b}u_\beta+\widetilde{c}u=\widetilde{d},
\end{equation}
for some new functions $\widetilde{a}=\widetilde{a}(\alpha,\beta)$ etc. This is also called the canonical form for hyperbolic equations. We observe that the wave equation
\[u_{tt}-u_{xx}=1,\]
is exactly of this form.

\subsection{Canonical forms for equations of parabolic type}
For parabolic equations we have
\[B^2(x,y)-A(x,y)C(x,y)=0,\]
in the domain of interest. 
For this case we have
\[\omega^+=\omega^-=\omega,\]
 and there exists only one family of characteristic curves
\[\xi=\xi(x,y),\]
where the curves
\begin{equation*}
 \xi(x,y)=\text{const},
 \end{equation*}
  by definition are solutions to the ODE
\[\frac{\mathrm{d}y}{\mathrm{d}x}=-\omega.\]
Let 
\[\eta=\eta(x,y),\]
be another family of curves chosen such that 
\begin{equation}\label{eq7-19}
\begin{split}
\xi&=\xi(x,y),\\
\eta&=\eta(x,y),
\end{split}
\end{equation}
defines a new coordinate system $(\xi,\eta)$ in the domain of interest.
For example if $A,B$ and $C$ are constants we have
\[\xi=y+\omega x,\]
and we can let 
\[\eta=\omega y+cx,\]
where $c\neq \omega^2$. Then (\ref{eq7-19}) is a coordinate system because
\[\Delta=\xi_x\eta_y-\xi_y\eta_x=\omega^2-c\neq 0.\]

We now express all derivatives $u_x,u_{xx},u_{xy},u_y,u_{yy}$ in terms of $u_\xi,u_\eta,u_{\xi\eta},u_{\xi\xi},u_{\eta\eta}$. using the chain rule. For example we have
\begin{align*}
u_x&=\xi_xu_\xi+\eta_xu_\eta,\\
u_y&=\xi_yu_\xi+\eta_yu_\eta,\\
u_{xx}&=\xi_{xx}u_\xi+\eta_{xx}u_\eta+\xi_x^2u_{\xi\xi}+2\xi_x\eta_xu_{\xi\eta}+\eta_x^2u_{\eta\eta},
\end{align*}
etc. 

Inserting these expressions into the equation we get a principal part of the form
\begin{align}\label{principalpart}
&\left[A\xi_x^2+2B\xi_x\xi_y+C\xi_y^2\right]u_{\xi\xi}
+2\left[A\xi_x\eta_x+B\xi_x\eta_y+B\xi_y\eta_x+C\xi_y\eta_y\right]u_{\xi\eta}\nonumber\\
&+\left[A\eta_x^2+2B\eta_x\eta_y+C\eta_y^2\right]u_{\eta\eta},
\end{align}
Since \[\xi_x=\omega\xi_y,\] we have 
\begin{eqnarray*}
A\xi_x^2+2B\xi_x\xi_y+C\xi_y^2&=A(\omega\xi_y)^2+2B(\omega\xi_y)\xi_y+C\xi_y^2\\
&=\xi_y^2(A\omega^2+2B\omega+C)=0.
\end{eqnarray*}
Thus, the term $u_{\xi\xi}$  drops out of the equation. Also recall that for the parabolic case we have $\omega=-\frac{B}{A}$. Therefore we have
\begin{eqnarray*}
A\xi_x\eta_x+B\xi_x\eta_y+B\xi_y\eta_x+C\xi_y\eta_y\\
&=A(-\frac{B}{A}\xi_y)\eta_x+B(-\frac{B}{A}\xi_y)\eta_y+B\xi_y\eta_x+C\xi_y\eta_y\\
&=-B\xi_y\eta_x-\frac{B^2}{A}\xi_y\eta_y+B\xi_y\eta_x+C\xi_y\eta_y\\
&=-\frac{1}{A}\xi_y\eta_y(B^2-AC)\\
&=0.
\end{eqnarray*}
Therefore the term $u_{\xi\eta}$ also drops out of the equation. Also, observe that the term 
\begin{equation}\label{uksiksi}
A\eta_x^2+2B\eta_x\eta_y+C\eta_y^2,
\end{equation}
must be nonzero. This is  because if it was zero, $\eta=\eta(x,y)$ would define a second family of characteristic curves different from $\xi=\xi(x,y)$ and we know that no such family exists for the parabolic case.

We can therefore divide the equation by the term (\ref{uksiksi}) and we get
\begin{equation}\label{eq7-20}
u_{\eta\eta}+au_{\xi}+bu_{\eta}+cu=d.
\end{equation}
This is the canonical form for equations of parabolic type. We observe that the diffusion equation
\[u_t-Du_{xx}=0,\]
is of this form.

\subsection{Canonical forms for equations of elliptic type}
For elliptic equations we have
\[B^2-AC<0,\]
so both $w^+$ and $w^-$ are complex and $\xi$ and $\eta$ become complex coordinates. In \cite{textbook} all functions are extended  into the complex domain, and an argument like the one for the hyperbolic case lead to the canonical form
\begin{equation}\label{eq7-21}
u_{\alpha\alpha}+u_{\beta\beta}+au_\alpha+bu_\beta+cu=d,
\end{equation}
where $a=a(\alpha,\beta)$ etc. The same canonical form can be found using only real variables. The key step in the reduction to canonical form is to simplify the principal part of the equation using a change of variables. Let us therefore do a general change of variables
\begin{align*}
\xi&=\xi(x,y),\\
\eta&=\eta(x,y).
\end{align*}
From equation (\ref{principalpart}) we have the following expression for the principal part of the transformed equation
\begin{align*}
&\left[A\xi_x^2+2B\xi_x\xi_y+C\xi_y^2\right]u_{\xi\xi}
+2\left[A\xi_x\eta_x+B\xi_x\eta_y+B\xi_y\eta_x+C\xi_y\eta_y\right]u_{\xi\eta}\\
&+\left[A\eta_x^2+B\eta_x\eta_y+C\eta_y^2\right]u_{\eta\eta}.
\end{align*}
The functions $\xi(x,y)$ and $\eta(x,y)$ are at this point only constrained by the Jacobi condition
\[\Delta(x,y)=(\xi_x\eta_y-\xi_y\eta_x)(x,y)\neq0.\]

\noindent In the elliptic case, both solutions $\xi,\eta$ to the characteristic equation are complex. Thus, for the elliptic case the coefficients in front of $u_{\xi\xi}$ and $u_{\eta\eta}$ can {\it not} be made to vanish. However if we let $\eta(x,y)$ be arbitrary and let $\xi=\xi(x,y)$ be the solution family of the system
\begin{equation}\label{eta}
\frac{\mathrm{d}y}{\mathrm{d}x}=\frac{B\eta_x+C\eta_y}{A\eta_x+B\eta_y}.
\end{equation}
Then along the solution curves
$\xi(x,y)=\text{const},$
we have
\[\xi_x+\frac{\mathrm{d}y}{\mathrm{d}x}\xi_y=0,\]
\[\Downarrow\]
\[\xi_x=-\xi_y\frac{B\eta_x+C\eta_y}{A\eta_x+B\eta_y}.\]
Thus
\begin{align*}
\Delta&=\xi_x\eta_y-\xi_y\eta_x=-\xi_y\eta_y\frac{B\eta_x+C\eta_y}{A\eta_x+B\eta_y}-\xi_y\eta_x\\
&=\frac{-\xi_y}{A\eta_x+B\eta_y}\left\{\eta_y(B\eta_x+C\eta_y)+\eta_x(A\eta_x+B\eta_y)\right\}\\
&=\frac{-\xi_y}{A\eta_x+B\eta_y}\left\{B\eta_x\eta_y+C\eta_y^2+A\eta_x^2+B\eta_x\eta_y\right\}\\
&=\frac{-\xi_y}{A\eta_x+B\eta_y}\left\{A\eta_x^2+2B\eta_x\eta_y+C\eta_y^2\right\}\neq0,
\end{align*}
because the characteristic equation has no real solutions.

  Therefore
\begin{align*}
\xi&=\xi(x,y),\\
\eta&=\eta(x,y),
\end{align*}
{\it is} a change of coordinates.

  Furthermore
\begin{align*}
A\xi_x\eta_x+B\xi_x\eta_y+B\xi_y\eta_x+C\xi_y\eta_y
&=(A\eta_x+B\eta_y)(\xi_x+\frac{B\eta_x+C\eta_y}{A\eta_x+B\eta_y}\xi_y)=0.
\end{align*}
Thus, the coefficient in front of $u_{\xi\eta}$ is zero. Since the quantity
\[A\varphi_x^2+2B\varphi_x\varphi_y+C\varphi_y^2\neq0,\quad \forall\,\,\text{real}\,\,\varphi,\]
the coefficients in front of $u_{\xi\xi}$ and $u_{\eta\eta}$ are nonzero and of the same sign. Dividing, we find that the principal part takes the form
\begin{equation}\label{elliptic1}
u_{\xi\xi}+h(\xi,\eta)u_{\eta\eta},
\end{equation}
where the function $h(x,y)$ is a positive function.

Let us at this point simplify our considerations by assuming that the coefficients of the equation (\ref{eq7-1}) are in fact constant. If that is the case,  we choose the function $\eta(x,y)$ to be linear in $x$ and $y$. Then the derivatives of $\eta$ are constant and therefore the function $\xi(x,y)$ defined in (\ref{eta}) is also a linear function. For this particular case, the function $h(\xi,\eta)$ is a positive constant which we denote by $h_0$. Thus equation (\ref{elliptic1}) takes the form
\begin{equation}\label{elliptic2}
u_{\xi\xi}+h_0u_{\eta\eta}.
\end{equation}
Introduce a change of variables defined by 
\begin{eqnarray}
\alpha=&\xi,\\
\beta=&\frac{1}{\sqrt{h_0}}\eta.
\end{eqnarray}
Using the chain rule we find that the principal part of an elliptic equation with constant coefficients, takes the form
\begin{equation*}
u_{\alpha\alpha}+u_{\beta\beta}.
\end{equation*}
Thus the canonical form of elliptic equations with constant coefficients  is 
\begin{equation}\label{CanonicalFormElliptic}
u_{\alpha\alpha}+u_{\beta\beta}+au_\alpha+bu_\beta+cu=d.
\end{equation}
One can show that also elliptic equations with variable coefficients can be put into canonical form (\ref{CanonicalFormElliptic}) by a suitable change of coordinates. This argument however either relies on complex analysis, or some rather intricate real variable manipulations.

\subsection{Characteristic curves}
One consequence of the classification of equations into hyperbolic, parabolic and elliptic, is the possibility to reduce each class into a canonical form where the principal part has the form
\begin{align*}
&u_{\xi\xi}-u_{\eta\eta} \quad\quad\text{hyperbolic},\\
&u_{\xi\xi}\quad\quad\quad\quad\,\,\,\,\text{parabolic},\\
&u_{\xi\xi}+u_{\eta\eta}\quad\quad\text{elliptic}.
\end{align*}
In these reductions, certain special curve families played an important role. 
We will now see that these special curve families reappear when we consider the possibility of solutions that contain singularities. These are solutions whose existence,or not, is a consequence of which class our equation belong to. This will give us a new window into what the abstract classification of equations into the three classes hyperbolic, parabolic and elliptic actually means.

Let our domain, $D$, be separated into two parts $D_1$ and $D_2$ by a curve.
\begin{figure}[h]
\centering
\includegraphics{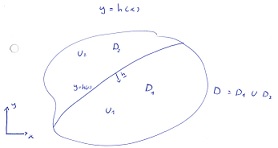}
\caption{}
\end{figure}

\noindent Let $u_1$ and $u_2$ be solutions of the equation
\begin{equation}\label{eq7-22}
Au_{xx}+2Bu_xy+Cu_{yy}+Du_x+Eu_y+Fu=G,
\end{equation}
in $D_1$ and $D_2$. Define a function $U$ on the whole domain $D$ by
\begin{equation*}
u(x,y)=\begin{cases}
u_1(x,y) & (x,y)\in D_1,\\
u_2(x,y) & (x,y)\in D_2.
\end{cases}
\end{equation*}
We will assume that $u,u_x,u_y$ are continuous on $D$ but that $u_{xx}, u_{xy}, u_{yy}$ can be discontinuous across the curve separating $D_1$ and $D_2$. Thus $u$ is not a solution in the usual sense to (\ref{eq7-22}) on the whole domain $D$.

\noindent The curve is defined by
\[\xi(x,y)=y-h(x)=0\]
Introduce a curve family
\[\xi(x,y)=\xi,\quad \xi\in\mathbb{R}\]
Our curve $y=h(x)$ corresponds to $\xi=0$.
Introduce a second curve family determined $\eta=\eta(x,y)$ in such a way that 
\begin{align*}
\xi&=\varphi(x,y),\\
\eta&=\psi(x,y),
\end{align*}
define a change of coordinates. For this to be the case the Jacobian must be nonzero everywhere.\\
We now  transform our equation into the new coordinate system  in the same way as on page 79, and get
\begin{equation}\label{eq7-23}
\begin{split}
&(A\xi_x^2+2B\xi_x\xi_y+C\xi_y^2)u_{\xi\xi}+2(A\xi_x\eta_x+B(\xi_x\eta_y+\xi_y\eta_x)\\
&+C\xi_y\eta_y)u_{\xi\eta}+(A\eta_x^2+2B\eta_x\eta_y+C\eta_y^2)u_{\eta\eta}+\cdots=G,
\end{split}
\end{equation}
where we only has displayed the principal part of the equation.
What we precisely mean when we say that $u_{xx}, u_{xy}$ and $u_{yy}$ can be discontinous across $y=h(x)$ is that
\[u,u_\xi,u_\eta,u_{\xi\eta},u_{\eta\eta},\]
are continuous but that $u_{\xi\xi}$ can be discontinuous. Evaluating (\ref{eq7-23}) above and below the curve $\xi=0 \Leftrightarrow y=h(x)$, taking the difference and going to the limit as the points approach $y=h(x)$ from above and below we get
\[\left.(A\xi_x^2+2B\xi_x\xi_y+C\xi_y^2)\left[u_{\xi\xi}\right]\right|_{\xi=0}=0,\]
where $\left.\left[u_{\xi\xi}\right]\right|_{\xi=0}$ is the jump in $u_{\xi\xi}$ across $y=h(x)$. Since $\left.\left[u_{\xi\xi}\right]\right|_{\xi=0}\neq0$, we {\it must} have
\[A\xi_x^2+2B\xi_x\xi_y+C\xi_y^2=0.\]
Thus $\xi(x,y)$ must satisfy the characteristic PDE and consequently
\[\xi=0\Leftrightarrow y=h(x),\]
{\it must} be a {\it characteristic curve} for the PDE (\ref{eq7-22}).

Thus, discontinuities in the second derivatives can only occur on characteristic curves. 

Let us see how this works out for the wave equation, diffusion equation and Laplace equation. We have seen that these three equations are respectively of hyperbolic type, parabolic type and elliptic type.

The wave equation is
\[u_{tt}-\gamma^2u_{xx}=0,\]
\[A=1, B=0, C=\gamma.\]
The characteristic PDE is
\[\varphi_t^2-\gamma^2\varphi_x^2=(\varphi_t+\gamma\varphi_x)(\varphi_t-\gamma\varphi_x)=0,\]
and the characteristic curves are
\begin{enumerate}
\item $\varphi_t-\gamma\varphi_x=0 \Leftarrow \varphi(x,t)=x+\gamma t,$
\item $\varphi_t+\gamma\varphi_x=0 \Leftarrow \varphi(x,t)=x-\gamma t.$
\end{enumerate}
The two families of characteristic curves are
\begin{align*}
x+\gamma t&=\text{const},\\
x-\gamma t&=\text{const}.
\end{align*}
\begin{figure}[h]
\centering
\includegraphics{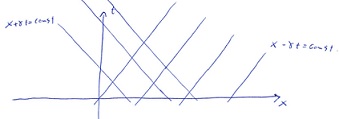}
\caption{}
\end{figure}
 Singularities, that is discontinuities in $u_{xx}, u_{xy},u_{yy}$, can only occur along these curves. Singularities given in initial data at $t=0$ is propagated along the characteristic curves and persists forever. This behavior is typical for all hyperbolic equations. 

The diffusion equation is
\[u_t=Du_{xx},\]
\[\Rightarrow A=D, B=C=0.\]
So the characteristic PDE is 
\[\varphi_x^2=0, \Rightarrow \varphi=\varphi(t),\]
and characteristic curves $\varphi(x,t)=\text{const}$ must be of the form
\[t=\text{const}.\]
\begin{figure}[h]
\centering
\includegraphics{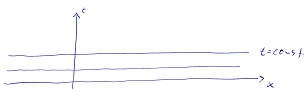}
\caption{}
\end{figure}

\noindent Singularities in the initial data given on the x-axis can not leave the x-axis and will not be seen at any point $(x,t)$ for $t>0$. Thus, all singularities in initial data are washed out instantly for $t>0$. 
This behaviour is typical for all parabolic equations.

\noindent Finally, the Laplace equation is 
\[u_{xx}+u_{yy}=0,\]
and the characteristic equation is 

\[\varphi_x^2+\varphi_y^2=0,\Leftrightarrow \varphi_x=\varphi_y=0,\Leftrightarrow \varphi(x,y)=\text{const}.\]
There are no characteristic curves that can support singularities. The solutions of elliptic equations are very smooth. This is typical for all elliptic equations.

\subsection{Characteristic curves again: A different point of view}
We first discussed characteristic curves in the context of linear and quasilinear first order scalar PDEs in two independent variables. There they played the star role in a method for solving the initial value problem for such equations. This was the {\it method of characteristics}. 
Next, we met them in the context of second order linear scalar PDEs in two independent variables. Here we have seen that they played a dual role. On the one hand they were instrumental in the reduction of hyperbolic and parabolic equations to canonical form, and on the other hand they where curves along which solutions to the PDEs could have discontinuities in their second derivatives. This second role is also played by the characteristic of first order linear scalar PDEs.

\noindent We will now discuss a third way in which characteristic curves play a role.

Let us start  by considering the case of first order linear scalar PDEs in two independent variables. The initial value problem in assumed to be given in terms of an initial curve that is the graph of a function $h(x)$
\[y=h(x).\]
On this curve the solution $u$ is assumed to be given by some function $f(x)$. Thus we have
\begin{equation}
au_x+bu_y=cu+d,\label{eq7-24}
\end{equation}
\begin{equation}
u(x,h(x))=f(x). \label{eq7-25}
\end{equation}
Differentiating (\ref{eq7-25}) we get 
\begin{equation}\label{eq7-26}
u_x+h'u_y=f',
\end{equation}
and evaluating (\ref{eq7-24}) along the initial curve we get
\begin{equation}\label{eq7-27}
au_x+bu_y=cf+d,
\end{equation}
(\ref{eq7-26}) and (\ref{eq7-27}) is a linear system of equations for determining $u_x$ and $u_y$
\begin{equation}\label{eq7-28}
\begin{bmatrix}
1 & h' \\
a & b
\end{bmatrix}
\begin{bmatrix}
u_x\\
u_y
\end{bmatrix}=
\begin{bmatrix}
f'\\
cf+d
\end{bmatrix}.
\end{equation}
For points close to the initial curve we can express the solution of the initial value problem (\ref{eq7-24}),(\ref{eq7-25}) as a Taylor series.
\begin{figure}[h]
\centering
\includegraphics{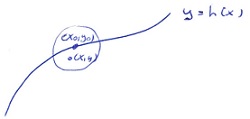}
\caption{}
\end{figure}
\[u(x,y)=u(x_0,y_0)+u_x(x_0,y_0)(x-x_0)+u_y(x_0,y_0)(y-y_0)+\cdots\]
Clearly if the solution $u(x,y)$ is to exist and be unique, then at least $u_x(x_0,y_0)$ and $u_y(x_0,y_0)$ must exist as a unique solution of the system (\ref{eq7-28}).

This is ensured if
\[\text{det}\left.
\begin{bmatrix}
1 & h'\\
a & b
\end{bmatrix}\right|_{(x_0,y_0)}\neq0.\]
Even if this is true it is not certain that the initial value problem is uniquely solvable, there is the question of the coefficients of the rest of the Taylor series and of course there is the matter of convergence for such series \dots.

However, if on the other hand, 
\begin{equation}\label{eq7-29}
\text{det}\left.
\begin{bmatrix}
1 & h'\\
a & b
\end{bmatrix}\right|_{(x_0,y_0)}=0,
\end{equation}
$u_x(x_0,y_0), u_y(x_0,y_0)$ either does not exist and if they do exist they are not unique. This follows from the elementary theory of linear systems of equations. So if (\ref{eq7-29}) holds we can know for sure that we do not have existence and uniqueness for the initial value problem (\ref{eq7-24}),(\ref{eq7-25}).

Condition (\ref{eq7-29}) can be written as
\[ b-ah'=0,\]
\[\Updownarrow\]
\begin{equation}\label{eq7-30}
h'=\frac ba.
\end{equation}
This equation tells us that $y=h(x)$ must be a characteristic curve! This follow, because the family of base characteristic curves $(x(s),y(s))$ are determined by being solutions to 
\begin{equation}\label{eq7-31}
\frac{\mathrm{d}x}{\mathrm{d}s}=a,\frac{\mathrm{d}y}{\mathrm{d}s}=b.
\end{equation}
If the curves in question are graphs of functions we can parametrise them using $x$. So that they are given by a function $y(x)$. For this function (\ref{eq7-31}) together with the chain rule gives
\[\frac{\mathrm{d}y}{\mathrm{d}x}=\frac{\mathrm{d}y}{\mathrm{d}s}\cdot\frac{\mathrm{d}s}{\mathrm{d}x}=\frac ba.\]
and (\ref{eq7-30}) say that $h(x)$ solves this equation so that $y=h(x)$ is a base characteristic curve.

Let us next consider the case of second order linear scalar PDEs in two independent variables. 

The characteristic curves have preciously been found to satisfy
\begin{equation}\label{eq7-32}
A\varphi_x^2+2B\varphi_x\varphi_y+C\varphi_y^2=0.
\end{equation}
The curves on which the PDE could have singularities, $y=h(x)$, were related to $\varphi(x,y)$ through
\[\varphi(x,y)=y-h(x)=0,\]
and (\ref{eq7-32}) then give the equation
\begin{equation}\label{eq7-33}
Ah'^2-2Bh'+C=0
\end{equation}
The solutions of this nonlinear ODE gives the characteristic curves for PDEs of the form
\begin{equation}\label{eq7-34}
Au_{xx}+Bu_{xy}+Cu_{yy}=F(x,y,u,u_x,u_y),
\end{equation}
where
\[F(x,y,u,u_x,u_y)=Du_x+Eu_y+Fu+G.\]
Recall that all parameters are functions of $x$ and $y$.

Let $y=h(x)$ be a initial curve where we have given $u,u_x$ and $u_y$. Thus 
\begin{eqnarray}\label{eq7-37}
u(x,h(x))=p(x),\nonumber\\
u_x(x,h(x))=q(x),\nonumber\\
u_y(x,h(x))=r(x).
\end{eqnarray}
\begin{figure}[h]
\centering
\includegraphics{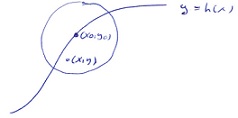}
\caption{}
\end{figure}

\noindent For points close to the initial curve the solution can be expand as a Taylor series.
\begin{align*}
u(x,y)&=u(x_0,y_0)+u_x(x_0,y_0)(x-x_0)+u_y(x_0,y_0)(y-y_0)\\
&+\frac12u_{xx}(x_0,y_0)(x-x_0)^2+u_{xy}(x_0,y_0)(x-x_0)(y-y_0)\\
&+\frac12u_{yy}(x_0,y_0)(y-y_0)^2.
\end{align*}
Differentiating the last two equations in the initial data (\ref{eq7-37}) with respect to $x$,  we get
\begin{equation}\label{eq7-38}
u_{xx}+h'u_{xy}=q',
\end{equation}
\begin{equation}\label{eq7-39}
u_{xy}+h'u_{yy}=r'.
\end{equation}
Evaluating (\ref{eq7-34}) along $y=h(x)$ we get the equation
\begin{equation}\label{eq7-40}
Au_{xx}+2Bu_{xy}+Cu_{yy}=F(x,h(x),p,q,r).
\end{equation}
Equations (\ref{eq7-38}),(\ref{eq7-38}) and (\ref{eq7-40}) define a linear system of equations for $u_{xx},u_{xy}$ and $u_{yy}$
\begin{equation}\label{eq7-41}
\begin{bmatrix}
1 & h' & 0\\
0 & 1  & h'\\
A & 2B & C
\end{bmatrix}
\begin{bmatrix}
u_{xx}\\
u_{xy}\\
u_{yy}
\end{bmatrix}=
\begin{bmatrix}
q'\\
r'\\
F(x,h(x),p,q,r)
\end{bmatrix}.
\end{equation}
The Taylor expansion show that the the problem (\ref{eq7-34}) with initial data  (\ref{eq7-37}) can have a unique solution only if the determinant of the matrix in the linear system (\ref{eq7-41}) is nonzero.

If, on the other hand, the determinant is zero, the problem is sure to either have no solution or to have a nonunique solution. This condition is 
\begin{eqnarray}
&\text{det}\begin{bmatrix}
1 & h' & 0 \\
0 & 1 & h'\\
A & 2B & C
\end{bmatrix}=0,\nonumber\\
&\Updownarrow\\
&C-2Bh'+Ah'^2=0,\label{eq7-42}
\end{eqnarray}
and this is exactly the equation for the characteristic curves (\ref{eq7-33}).

Thus characteristic curves {\it are}
curves for which the (generalized) initial value problem either can not be solved or for which the solution is non-unique.

There is something awkward about the presentation given and equation (\ref{eq7-42}). Since we are considering the general initial value problem as giving $u,u_x$ and $u_y$ along some curve in the $(x,y)$ plane, there is really no point in restricting to curves that is the graph of a function $y=h(s)$. Why not consider curves that are graph of a function $x=k(y)$ or indeed a curve that is neither a graph of a function of $x$ or a function of $y$?
\begin{figure}[h]
\centering
\includegraphics{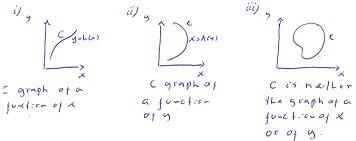}
\caption{}
\end{figure}

\noindent Let us consider case \romannumeral3). Let $\varphi(x,y)$ be a function such that $C$ is a level curve of $\varphi$
\[C:\varphi(x,y)=0.\]
Pick any point  $(x_0,y_0)$ on $C$. Then the piece of the curve close to $(x_0,y_0)$ will either be the graph of a function of $x$ or the graph of a function of $y$
\begin{figure}[h]
\centering
\includegraphics{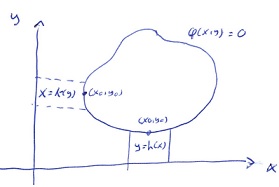}
\caption{}
\end{figure}

\noindent Let us assume, without loss of generality, that the piece of $C$ close to $(x_0,y_0)$ is the graph of a function of x, $y=h(x)$. But then we must have
\[\varphi(x,h(x))=0,\]
\[\Downarrow\]
\[\varphi_x+h'\varphi_y=0,\]
\[\Downarrow\]
\[h'=-\frac{\varphi_x}{\varphi_y}.\]
Inserting this into equation (\ref{eq7-42}) gives
\[C+2B\frac{\varphi_x}{\varphi_y}+A\left(\frac{\varphi_x}{\varphi_y}\right)^2=0,\]
\[\Updownarrow\]
\begin{equation}\label{eq7-43}
A\varphi_x^2+2B\varphi_x\varphi_y+C\varphi_y^2=0.
\end{equation}
The conclusion of this is that curves, $C$, for which the generalized initial value problem $(\ref{eq7-34}), (\ref{eq7-37})$ either has no solution or has a nonunique solution, must be level curves
\[\varphi(x,y)=0,\]
where $\varphi$ is a solution of (\ref{eq7-43}). And all solutions of (\ref{eq7-43}) determines such special initial curves.

We have met equation (\ref{eq7-43}) before on page 9 of these notes where it was called the characteristic equation for $(\ref{eq7-34})$.

Thus the special initial curves discussed in this section are in fact {\it characteristic} {\it curves} as defined previously.

The  characteristic curves are seen to play three different roles for PDEs of type (\ref{eq7-34}).
\begin{enumerate}
\item Level curves for special coordinate system used to reduce the PDE to canonical form.
\item Curves on which solutions of the PDE can have singularities.
\item Curves for which the initial value problem for the PDE either have no solution or have a non-unique solution.
\end{enumerate}
So far the statements \romannumeral1), \romannumeral 2) and \romannumeral3) only applies to linear second order scalar PDEs in two independent variables. Can we extend these results to more general equations?

The answer is yes we can! Many of these ideas can be extended far beyond the class of equations we have discussed so far. In \cite{textbook} classification and associated characteristic curves and surfaces for large classes of equations are discussed on pages 137-151.

In this course we will not discuss all these various cases but will rather restrict to linear scalar second order equations in $n$ independent variables where $n$ can be any positive integer.

Equations of this type takes the form
\begin{equation}\label{eq9-44}
\sum_{i=1}^n\sum_{j=1}^na_{ij}\partial_{x_i}\partial_{x_j}u+\sum_{i=1}^nb_i\partial_{x_i}u+cu+d=0,
\end{equation}
where all parameters $\{a_{ij}\}$ and $\{b_i\}, c, d$ are functions of $(x_1,\dots,x_n)$
 
Remark on notation: Here and in many other areas of applied mathematics and theoretical physics it is common to use the so-called {\it Einstein summation} {\it convention}.

In this convention repeated indices in any expression are implicitly summed over. Thus the expressions
\[x_iy_i,\]
\[a_{ij}x_j,\]
\[x_ia_{ij}y_j,\]
means
\[\sum_{i=1}^nx_iy_j,\]
\[\sum_{j=1}^na_{ij}x_j,\]
\[\sum_{i=1}^n\sum_{j=1}^nx_ia_{ij}y_j.\]
Using this convention, equation (\ref{eq9-44}) is written as
\begin{equation}\label{eq9-45}
a_{ij}\partial_{x_i}\partial_{x_j}u+b_i\partial_{x_i}u+cu+d=0,
\end{equation}
Observe that
\begin{align*}
a_{ij}\partial_{x_i}\partial_{x_j}u&=\frac12a_{ij}\partial_{x_i}\partial_{x_j}u+\frac12a_{ij}\partial_{x_i}\partial_{x_j}u\\
                                   &=\frac12a_{ij}\partial_{x_i}\partial_{x_j}u+\frac12a_{ji}\partial_{x_j}\partial_{x_i}u\\
                                   &=\frac12a_{ij}\partial_{x_i}\partial_{x_j}u+\frac12a_{ji}\partial_{x_i}\partial_{x_j}u\\
                                   &=\left( \frac12a_{ij}+\frac12a_{ji}\right) \partial_{x_i}\partial_{x_j}u\\
                                   &=\widetilde{a}_{ij}\partial_{x_i}\partial_{x_j}u,\\
\end{align*}
where $\widetilde{a}_{ij}=\frac12(a_{ij}+a_{ji})$.

\noindent $\widetilde{a}_{ij}$ are evidently the components of a symmetric matrix 
\[\widetilde{a}_{ij}=\widetilde{a}_{ji}.\]
Thus, we can without loss of generality assume that the matrix $A=(a_{ij})$ in equation (\ref{eq9-45}) is a symmetrix matrix.

Let us first assume that all $a_{ij}$ are constants. Let $\lambda_k$ be the eigenvalues of $A$ and $R=(r_{ki})$ be the orthogonal matrix whose rows are the correponding eigenvectors. In terms of matrices we have 
\[RAR^t=\begin{pmatrix}
\lambda_1  & & 0\\
           &\ddots &\\
0 &         &\lambda_n
\end{pmatrix}.\]
This can be written, using the Einstein convention, in the form
\[r_{pi}a_{ij}r_{qj}=\Lambda_{pq},\]
where $\Lambda_{pq}$ is defined by
\[\Lambda_{pq}=\begin{cases}
\lambda_p,&p=q\\
0,&p\neq q
\end{cases}.\]
Introduce new coordinates $\xi_1,\dots,\xi_n$ through
\[\xi_p=r_{pk}x_k.\]
This is evidently a linear change of coordinates. Observe that
\begin{align*}
\partial_{x_i}&=\partial_{x_i}\xi_p\partial_{\xi_p}=r_{p_i}\partial_{\xi_p},\\
\partial_{x_j}&=\partial_{x_j}\xi_q\partial_{\xi_q}=r_{q_j}\partial_{\xi_q},
\end{align*}
\[\Updownarrow\]
\begin{align*}
a_{ij}\partial_{x_i}\partial_{x_j}u&=a_{ij}r_{p_i}r_{q_j}\partial_{\xi_p}\partial_{\xi_q}u\\
&=\Lambda_{pq}\partial_{\xi p}\partial_{\xi q}u\\
&=\lambda_p\partial_{\xi p}\partial_{\xi_p}u.
\end{align*}
Thus the second order mixed partial derivatives all vanish. This is the canonical form for (\ref{eq9-45}).
\noindent If we return to the $n=2$ case we observe that
\begin{align*}
\lambda_1,\lambda_2\,\,\text{same sign}&-\,\,\text{elliptic case},\\
\lambda_1,\lambda_2\,\,\text{opposite sign}&-\,\,\text{hyperbolic case},\\
\text{one of}\,\,\lambda_1,\lambda_2\,\,\text{is zero}&-\,\,\text{parabolic case}.
\end{align*}
Based on this we now make the following classification of equations of type (\ref{eq9-45}).
\begin{align*}
&\begin{cases}
\lambda_i>0,&\forall i\\
\text{or}& \\
\lambda_i<0,&\forall i
\end{cases}\qquad \qquad\quad\qquad\qquad\,\,\,-\text{elliptic type},\\
&\begin{cases}
\text{one of}\,\,\lambda_i>0\,\,\text{or} \lambda_i<0\\
\text{All other}\,\,\lambda_i \,\,\text{opposite sign}
\end{cases}\qquad\,\, - \text{hyperbolic type},\\
&\text{one or more}\,\,\lambda_i=0 \qquad\quad\qquad\qquad\,\,\,- \text{parabolic type}.
\end{align*}
For the elliptic case we can make the further change of coordinates
\begin{equation}\label{eq9-46}
\alpha_i=\frac1{|\lambda_i|}\xi_i,
\end{equation}
leading to a principal part of the form
\begin{equation}\label{eq9-47}
\pm\partial_{\alpha_i}\partial_{\alpha_i}u.
\end{equation}
Assuming $\lambda_1>0$ and $\lambda_2 \dots \lambda_n<0$, the hyperbolic case is reduced, using (\ref{eq9-46}), to the form
\begin{equation}\label{eq9-48}
\partial_{x_1}^2-\sum_{i=2}^n\partial_{x_i}^2.
\end{equation}
Similar reductions are possible for the parabolic case.
Some $n=3,4$ basic examples of elliptic, hyperbolic and parabolic equations are 
\begin{align*}
\nabla^2u&=u_{xx}+u_{yy}+u_{zz}=0 \qquad\qquad\,\,\, \text{elliptic},\\
u_{tt}-\nabla^2u&=u_{tt}-u_{xx}-u_{yy}-u_{zz}=0\qquad \text{hyperbolic},\\
u_t-\nabla^2u&=u_t-u_{xx}-u_{yy}-u_{zz}=0\qquad\,\,\text{parabolic}.
\end{align*}
The classification of equations of type (\ref{eq9-45}) into elliptic, hyperbolic and parabolic based on eigenvalues of the matrix $A$ is still possible when $a_{ij}$ depends on $(x_1,\dots,x_n)$ and are thus not constant. The eigenvalues and the classification will however in general depend on which point we are at. Also, it is in general not possible to reduce it to a canonical form where all terms with mixed partial derivatives vanish. This is seen to be true by a simple counting argument.

\noindent Let
\begin{equation}\label{9-49}
\begin{split}
\xi_1&=\xi_1(x_1,\dots,x_n),\\
\vdots& \\
\xi_n&=\xi_n(x_n,\dots,x_n),
\end{split}
\end{equation}
be a change of coordinates in $\mathrm{R}^n$. Using the chain rule the principal part of equation (\ref{eq9-45}) can be written as
\begin{equation}
\sum_{ij}F_{ij}(\xi_1,\dots,\xi_n)\partial_{\xi_i}\partial_{\xi_j}u,
\end{equation}
where $F_{ij}$ is some expression containing derivatives of the functions $\xi_j$.
Let us try to eliminate all mixed partials. 
Then for this to happened, we must have
\begin{equation}\label{eq9-51}
F_{ij}(\xi_1,\dots,\xi_n)=0,\quad \forall i\neq j.
\end{equation}
There are exactly $\frac12n(n-1)$ terms with mixed partial derivatives. Thus there are $\frac12n(n-1)$ equations in (\ref{eq9-51}). On the other hand we only have $n$ unknown functions $\xi_1,\dots,\xi_n$. If
\[\frac12 n(n-1)>n,\]
the system (\ref{eq9-51}) in general has no solution. This happened for $n\eqslantgtr4$.

For equations of type (\ref{eq9-45}), surfaces supporting singularities of solutions, and on which given initial data does not give a unique solution for the initial value problem, are given as solutions to
\begin{equation}\label{eq9-52}
a_{ij}\varphi_{x_i}\varphi_{x_j}=0.
\end{equation}
These surfaces are called {\it characteristic} surfaces and equation (\ref{eq9-52}) is called the {\it characteristic equation} for equation (\ref{eq9-45}).
For example the 3D wave equation
\[u_{tt}-u_{xx}-u_{yy}-u_{zz}=0,\]
has a characteristic equation
\[\varphi_t^2-\varphi_x^2-\varphi_y^2-\varphi_z^2=0.\]
One important type of solution is
\[\varphi=\omega t-k_1x-k_2y-k_3z+c,\]
where $\omega^2=k_1^2+k_2^2+k_3^2$.
For each value of $t$ the characteristic surface
\[\varphi=0,\]
\[\Updownarrow\]
\[k_1x+k_2y+k_3z=\omega t+c,\]
is a plane normal to the vector $\boldsymbol{k}=(k_1,k_2,k_3)$. Such surfaces play a pivotal role in the theory of wave equations where they appear through special solutions called {\it plane waves}.

\subsection{Stability theory, energy conservation and dispersion}
So far we have used the highest order terms, the so-called principal part of the equation, to classify our equations. The lower order terms however also play an important role.

In this section of the notes we will investigate the effects of lower order terms by considering the class of equations of the form
\begin{equation}\label{eq9-53}
Au_{xx}+2Bu_{xt}+Cu_{tt}+Du_x+Eu_t+Fu=0,
\end{equation}
where $A,B,C,D,E$ and $F$ are real constants. Certain aspects of our conclusions for this class can be extended to a wider class of equations\cite{textbook}.

We start our investigations by looking for special solutions to (\ref{eq9-53}) of exponential form
\begin{equation}\label{eq9-54}
u(x,t)=a(k)e^{(ik x+\lambda(k)t)}.
\end{equation}
This is a complex solution but both the real and imaginary part of $u$ are real solutions of (\ref{eq9-53}). (Verify this).
We will in this section interpret the variable $t$ as time and study the evolution of the solutions (\ref{eq9-54}) for $t>0$. The special solutions (\ref{eq9-54}) are called {\it normal modes} for (\ref{eq9-53}). By the principle of linear super-position we can add together several normal modes corresponding to different values of $k$ to create more general solutions. Taking an infinite set of discrete values $\{k_j\}$ we get a solution of the form
\begin{equation}\label{eq9-55}
u(x,t)=\sum_{j=-\infty}^{\infty}a(k_j)e^{(ik_jx+\lambda (k_j) t)}.
\end{equation}
This is a {\it Fourier series}, and the sense in which such series solves (\ref{eq9-53}) will be discussed extensively later in the class. There we will also discuss a much more general notion of normal mode that applies to a much wider class of equations. We will see that the normal modes (\ref{eq9-54}) can be combined in a ``continuous sum'' by which we mean an integral
\begin{equation}\label{eq9-56}
u(x,t)=\int_{-\infty}^{\infty}\mathrm{d}k\,a(k)e^{(ik x+\lambda(k)t)},
\end{equation}
to generate even more general solutions to (\ref{eq9-53}). We recognize (\ref{eq9-56}) as a {\it Fourier transform}. More general transforms are associated with other types of normal modes. These are discussed in chapter 5 in \cite{textbook}.

At this point, the only thing to note is that through the general representations (\ref{eq9-55}) and (\ref{eq9-56}) the behavior of the normal modes (\ref{eq9-54}) will determine what any solution to equation (\ref{eq9-53}) will do for $t>0$.

So, let us now start analyzing the normal modes.
From (\ref{eq9-54}) we have
\begin{eqnarray}
&u(x,t)=a(k)e^{(ik x+\lambda(k)t)},\nonumber\\
&\Downarrow\nonumber\\
&\begin{split}
u_t=\lambda(k)u,\,\,u_x=ik u,\,\, u_{tt}=\lambda^2(k)u,\nonumber\\
u_{xt}=ik\lambda(k)u,\,\,u_{xx}=-k^2u.\label{eq9-57}
\end{split}
\end{eqnarray}
Inserting (\ref{eq9-57}) into (\ref{eq9-53}) we get
\[(-Ak^2+2iBk\lambda+C\lambda^2+iDk+E\lambda+F)u=0.\]
A solution of this type thus only exists if 
\begin{equation}\label{eq9-58}
C\lambda^2+(2ikB+E)\lambda+(-Ak^2+iDk+F)=0.
\end{equation}
This equation has two complex solutions. For each of these we can write
\begin{equation}\label{eq9-59}
\lambda(k)=Re[\lambda(k)]+iIm[\lambda(k)],
\end{equation}
and the corresponding normal mode can be written as
\begin{equation}\label{eq9-60}
u(x,t)=a(k)e^{Re[\lambda(k)]t}e^{i(k x+Im[\lambda(k)]t)}.
\end{equation}
The magnitude of $u(x,t)$ is
\begin{equation}\label{eq9-61}
\left|u(x,t)\right|=\left|a(k)\right|e^{Re[\lambda(k)]t}.
\end{equation}
We will assume that the expression $\left|a(k)\right|$ is bounded for all $k$. The growth of $|u(x,t)|$ as a function of $k$ for $t>0$ is then determined by $Re[\lambda(k)]$.

There are now two distinct possibilities. Either $Re[\lambda(k)]$ is bounded above for all $k\in\mathbb{R}$, or it is not.

\noindent Define
\begin{equation}\label{eq9-62}
\Omega=\mathop{\mathrm{sup}}_{k\in\mathbb{R}}Re[\lambda(k)].
\end{equation}
If $Re[\lambda(k)]$ is unbounded we let $\Omega=+\infty$. Let us first consider the case $\Omega=+\infty$. Previously in this class we introduced the notion of being well-posed for initial and or boundary value problems for PDEs. In order to be well-posed, three criteria had to be satisfied
\begin{enumerate}
\item[1)] The problem must have at least one solution.
\item[2)] The problem must have at most one solution.
\item[3)] The solution must depend in a continuous way on the data of the problem: Small data must produce a small solution.
\end{enumerate}
Data here means parameters occurring in the PDEs, description of the boundary curves or surfaces, if any, and functions determining boundary and/or initial conditions

For the case $\Omega=+\infty$, let $a(k)=\frac1{\lambda(k)^2}$ and consider the initial value problem for the normal mode with data at $t=0$. The data are from (\ref{eq9-54})
\begin{align*}
u(x,0)&=\frac1{\lambda^2(k)}e^{ik x},\\
u_t(x,0)&=\frac1{\lambda(k)}e^{ik x}.
\end{align*}
Since $Re[\lambda(k)]$ are unbounded $(\Omega=+\infty)$ there are values of $k$ such that $|u(x,0)|$ and $|u_t(x,0)|$ are as small as we like. Thus the data in this initial value problem can be made arbitrarily small. However
\[|u(x,t)|=\frac1{\lambda^2(k)}e^{Re[\lambda(k)]t},\]
grows exponentially as a function of $k$  and can for any fixed $t>0$ be made arbitrarily large. This is because an algebraic decay $\frac1{\lambda^2(k)}$ is always overwhelmed by an exponential growth $e^{Re[\lambda(k)]t}$ (Basic calculus).
Thus arbitrarily small data can be made arbitrarily large for any $t>0$.
Condition 3) on page 52 is not satisfied and consequently the problem is {\it not} well posed.
If the case $\Omega=+\infty$ is realized, it is back to the drawing board, something is seriously wrong with the model!
If $\Omega<+\infty$ one can show that the problem is well posed.

Let us look at some specific examples.

The equation
\[u_{xx}+u_{tt}+\rho u=0,\]
where $\rho$ is constant, is an elliptic equation. Our equation for $\lambda(k)$ for this case is
\[\lambda^2+(-k^2+\rho)=0,\]
\[\Downarrow\]
\[\lambda(k)=\pm\sqrt{k^2-\rho}.\]
$\lambda(k)$ is clearly unbounded when $k$ ranges over $\mathbb{R}$. Thus $\Omega=+\infty$ and the Cauchy problem for the equation with data at $t=0$ is not well posed.

The equation
\[\rho u_{xx}+u_t=0,\]
where $\rho$ is constant is a parabolic equation. In fact, we can also write it as 
\[u_t=-\rho u_{xx}.\]
For $\rho<0$ it is the heat equation. For $\rho>0$ it is called {\it the reverse heat equation}. We get the reverse heat equation from the heat equation by making the transformation $t'=-t$. Increasing values of $t'$ then corresponds to decreasing values of $t$. Solving the reverse heat equation thus corresponds to solving the heat equation backward in time.

\noindent We now get the following equation for $\lambda(k)$
\begin{equation}
-\rho k^2+\lambda=0 \Rightarrow \lambda=\rho k^2.\label{CharEq}
\end{equation}
Clearly $\lambda(k)$ is unbounded for $\rho>0\Rightarrow \Omega=+\infty$ which implies that the reverse heat equation is not well posed. This say something profound about the process of heat diffusion or any other process described by the diffusion equation. Such processes are {\it irreversible}, they have a built-in direction of time. It is a well known fact that no fundamental physical theory is irreversible, they have no direction of time. So how come the diffusion equation has one? Where and how did the direction of time enter into the logical chain connecting microscopic reversible physics to macroscopic irreversible diffusion processes? This question has been discussed for more than a hundred years and no consensus has been reached. Nobody really knows.

For $\rho<0$, $\lambda(k)$ is bounded above by zero $\Rightarrow \Omega=0$ and the Cauchy problem for the heat equation is well-posed.

The equation
\[u_{tt}-u_{xx}+\rho u=0,\]
where $\rho$ is a constant, is a hyperbolic equation. The equation for $\lambda(k)$ is 
\[\lambda^2+k^2+\rho=0,\]
\[\Downarrow\]
\[\lambda=\pm i\sqrt{k^2+\rho}.\]
For $\rho<0$ the real part of $\lambda$ vanish for $|k|$ large enough and for $\rho>0$ the real part vanish for all $k$. Thus in either case $\Omega<+\infty$ and the Cauchy problem for this equation is well-posed.

The constant $\Omega$ is called the {\it stability index} for the equation. 
If $\Omega>0$ there is some set or range of $k$ values whose corresponding normal modes experience exponential growth.
If this is the case the equation is said to be {\it unstable}. If $\Omega<+\infty$ the Cauchy problem is well posed, but the unbounded growth can still invalidate the PDE as a mathematical model for some $t>0$. The reason for this is that PDE models are usually derived using truncated Taylor expansions or expansions of some other sort. In order for these truncations to remain valid the solutions must remain small. This assumption will break down in the case of unbounded growth.
If $\Omega<0$ all normal modes decay exponentially. If this is the case the PDE is {\it strictly stable}. If $\Omega=0$ there may be some modes that does not decay and will continue to play a role for all $t>0$. If this is the case the PDE is {\it neutrally stable}.

Consider the parabolic equation
\[u_t-c^2u_{xx}+au_x+bu=0.\]
For this equation the equation for $\lambda(k)$ is 
\[\lambda+c^2k^2+iak+b=0,\]
\[\Downarrow\]
\[\lambda(k)=-c^2k^2-b-iak,\]
\[Re[\lambda(k)]=-b-c^2k^2.\]
Clearly $\Omega=-b$. It is stable if $b>0$ and unstable if $b<0$ and neutrally stable for $b=0$.

The equation
\[u_{xx}+2u_{xt}+u_{tt}=0,\]
is parabolic and has an equation for $\lambda$ given by 
\[-k^2+2ik\lambda+\lambda^2=0,\]
\[\Updownarrow\]
\[(\lambda+ik)^2=0,\]
\[\Updownarrow\]
\[\lambda(k)=-ik,\]
and thus  $\Omega=0$. The equation is therefore neutrally stable. It  has a normal mode solution
\[u_1(x,t)=ae^{ik(x-t)}.\]
It is straight forward to verify that it also have a solution
\[u_2(x,t)=bte^{ik(x-t)}.\]
By linear superposition it has a solution
\[u(x,t)=ae^{ik(x-t)}+bte^{ik(x-t)}.\]
The solution to the Cauchy problem with start data
\begin{align*}
u(x,0)&=0\\
u_t(x,0)&=1,
\end{align*}
is then easily shown to be
\[u(x,t)=te^{-ikx}.\]
Thus, even if the data for this case is bounded, the solution grows without bound  as $t$ approach infinity.  However, the growth is algebraic rather than exponential, which was the case for an unstable equation. This kind of behavior is generally true for the neutrally stable case; if there is an instability, then the solution grows algebraically rather than exponentially.

A special situation occur when $Re[\lambda(k)]=0,\,\,\forall k$. Then
\begin{align*}
|u(x,t)|&=\left|a(k)e^{i(k x+Im[\lambda(k)]t)}\right|\\
&=|a(k)|=|u(x,0)|.
\end{align*}
The amplitude of the normal mode is preserved for all time. In many cases of importance, $|u|^2$ is a measure of the {\it energy} of the normal mode and we have that the energy in the system is conserved in time.

In general, an equation where $Re[\lambda(k)]=0$, is called {\it conservative}.
For the conservative case we {\it define} the positive quantity
\[\omega=|Im[\lambda(k)]|, \]
to be the {\it frequency} of the mode. Given this definition we have
\[\lambda(k)=\pm i\omega(k),\]
depending on the sign of the imaginary part of the complex root $\lambda(k)$ of equation (\ref{CharEq}).

The normal mode then takes the form
\begin{equation}\label{eq9-63}
u(x,t)=a(k)e^{i(k x\pm \omega(k)t)}.
\end{equation}
The equation
\begin{equation}\label{eq9-64}
\omega=\omega(k),
\end{equation}
is called the {\it dispersion relation} for the equation.
If \begin{equation}\label{eq9-65}
\omega''(k)\neq0,
\end{equation}
and thus $\omega(k)$ is not a linear function of $k$, the equation is said to be of {\it dispersive type}.
The normal mode can, for the conservative case, be written as
\begin{equation}\label{eq9-66}
u(x,t)=a(k)e^{i\theta},
\end{equation}
where $\theta=k x\pm \omega(k)t$ is called the {\it phase} of the normal mode. It is clear that for $k>0$, (\ref{eq9-66}) represents a wave moving from right  to left(+) or from left to right(-), with speed
\begin{equation}\label{eq9-67}
v_f(k)=\frac{\omega(k)}{k}.
\end{equation}
The quantity $v_f$  is called the {\it phase speed} of the wave. If $\omega''(k)=0$ then
\[\omega(k)=ck\]
where $c$ is some constant. For this case the phase speed is independent of $k$
\[v_f=\frac{\omega(k)}{k}=\frac{ck}{k}=c.\]
Thus all normal modes moves at the same speed.

For the case when $w''(k)\neq0$, each normal mode moves at its own speed, and the linear superpositions (\ref{eq9-55}) and (\ref{eq9-56}) represents spreading or {\it dispersing} waves. For such solutions the phase speed is not the relevant quantity to consider. We will see later that for such waves the {\it group velocity}.
\begin{equation}\label{eq9-68}
v_g=\frac{\mathrm{d}\omega}{\mathrm{d}k},
\end{equation}
is the relevant quantity.

If the stability index $\Omega\leqslant0$ and $Re[\lambda(k)]<0$ for all, except a finite number of $k$-values, the equation is said to be of {\it dissipative type}.

Let us consider the so-called Telegrapher's equation
\[u_{tt}-\gamma^2u_{xx}+2\widehat{\lambda}u_t=0,\quad \widehat{\lambda}>0.\]
This is a hyperbolic  equation and the polynomial determining $\lambda(k)$ is
\[\lambda^2+\gamma^2k^2+2\widehat{\lambda}\lambda=0,\]
\[\Downarrow \]
\[\lambda(k)=-\widehat{\lambda}\pm\sqrt{\widehat{\lambda}^2-\gamma^2k^2}.\]
Since $\widehat{\lambda}>0, Re[\lambda(k)]<0,\,\forall k\neq0$. For $k=0$, either $\lambda(0)=0$ or $\lambda(0)=-2\widehat\lambda$ so that
\[\Omega=0.\]
The telegraphers equation is thus neutrally stable. It is also of dissipative type. When $\widehat{\lambda}=0$ we get the regular wave equation which is hyperbolic and conservation. For this reason the telegrapher's equation is also called the {\it damped wave} {\it equation}.

The {\it Klein-Gordon} equation takes the form
\[u_{tt}-\gamma^2u_{xx}+c^2u=0.\]
The equation for $\lambda(k)$ is
\[\lambda^2+\gamma^2k^2+c^2=0,\]
\[\Downarrow\]
\[\lambda(k)=\pm i\sqrt{\gamma^2k^2+c^2}.\]
The real part of $\lambda(k)$ is zero so the equation is of conservative type, also $\Omega=0$ so it is neutrally stable. The dispersion relation is 
\[\omega=\sqrt{\gamma^2k^2+c^2},\]
and clearly $\omega''(k)=0$ so the Klein-Gorden equation is of dispersive type.

\noindent If $\gamma=0$ we get
\[\omega=ck.\]
and for this case the equation is not of dispersive type.

The Klein-Gordon equation was proposed as a relativistic extension of the Schrodinger equation of quantum mechanics. It was first found by Schrodinger but was rediscovered by Klein and Gordon. In both cases it was rejected as an equation for a relativistic electron on physical grounds. It predicted negative probability! Later it however reappeared in quantum field theory and is now one of the fundamental equations of nature. It for example describes the Higgs Boson.

\section{Adjoint differential operators}
Let us consider the differential operator
\begin{equation}\label{eq9-69}
L=\frac{\mathrm{d}^2}{\mathrm{d}x^2}.
\end{equation}
For any pair of functions $u,w$ we have
\begin{equation}\label{eq9-70}
\begin{split}
uLw-wLu&=uw''-wu''\\
       &=uw''+u'w'-wu''-u'w'\\
       &=(uw'-wu')'.
\end{split}
\end{equation}
Thus $uLw-wLu$ is a total derivative. Identities like (\ref{eq9-69}) are of outmost importance in the theory of PDEs. Here is another example. Let
\begin{equation}\label{eq9-71}
L=\frac{\mathrm{d}}{\mathrm{d}x},
\end{equation}
and define a operator $L^*$ by
\begin{equation}\label{eq9-72}
L^*=-\frac{\mathrm{d}}{\mathrm{d}x}.
\end{equation}

Then
\[uLw-wL^*u=uw'+wu'=(uw)'.\]
Thus we get a total derivative. 
Next consider an operator
\begin{equation}\label{eq9-73}
L=\nabla^2.
\end{equation}
Then from vector calculus we have 
\begin{align*}
uLw-wLu&=u\nabla^2w-w\nabla^2u\\
       &=\nabla\cdot(u\nabla w-w\nabla u).
\end{align*}
Thus we get a total derivative. In genera,l if $L$ is a differential operator and if we have, for some vector field $\mathbf{F}$, the identity
\[uLw-wL^*u=\nabla\cdot\mathbf{F},\]
then $L^*$ is the {\it formal adjoint} of $L$. If $L=L^*$ the operator is {\it formally self-adjoint}. The operators (\ref{eq9-69}), (\ref{eq9-73}) are formally self-adjoint. The operator $L^*=-\frac{\mathrm{d}}{\mathrm{d}x}$ is the formal adjoint of $L=\frac{\mathrm{d}}{\mathrm{d}x}$.

Let $D$ be some domain in $\mathbb{R}^n$ and let $L$ be a differential operator defined for functions whose domain include $D$. Let $L^*$ be the formal adjoint of $L$. Then by the divergence theorem we have
\begin{equation}\label{eq9-74}
\int_D\mathrm{d}v\left\{uLw-wL^*u\right\}=\int_D \mathrm{d}v\nabla\cdot \mathbf{F}=\int_S\mathrm{d}A\, \mathbf{n}\cdot\mathbf{F}.
\end{equation}
\begin{figure}[!h]
\centering
\includegraphics{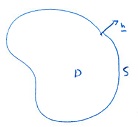}
\caption{}
\end{figure}

\noindent It is the possibility of reducing a integral over a domain to a integral over the boundary of the domain that is the reason why the formal adjoint of differential operators are important. We can use this possibility in many important ways. We will return to this topic many times in this and the next class (Mat-3202) in the spring.

\noindent For now observe the following: 

 Let $L=\nabla^2\Rightarrow L^*=\nabla^2$. Then (\ref{eq9-74}) is
\begin{equation}\label{eq9-75}
\begin{split}
\int_D\mathrm{d}v\left\{uLw-w Lu\right\}&=\int_S\mathrm{d}A\,\mathbf{n}\cdot\left\{u\nabla w-w\nabla u\right\}\\
           &=\int_S\mathrm{d}A\left\{u\partial_{\mathbf{n}}w-w\partial_{\mathbf{n}}u\right\},
\end{split}
\end{equation}
where $\partial_{\mathbf{n}}=\mathbf{n}\cdot\nabla$ is the normal derivative on the surface $S$. Let us assume that the domain of $L$ consists only of functions that satisfy the boundary conditions
\begin{eqnarray}
 u=&0\hspace{1mm}\text{on}\hspace{1mm} S,\label{eq9-76}\\
 \text{or}\nonumber\\
 \partial_\mathbf{n}u=&0\hspace{1mm}\text{on}\hspace{1mm} S.\nonumber
 \end{eqnarray}
Then (\ref{eq9-75}) implies that
\begin{equation}\label{eq9-77}
\int_D\mathrm{d}v\left\{uLw-wLu\right\}=0,
\end{equation}
which can be rewritten into the form
\begin{equation}\label{eq9-78}
\int_{D}\mathrm{d}v\,uLw=\int_{D}wLu.
\end{equation}
Let $V$ be the vector space of smooth functions defined on $D$ and that satisfy the boundary conditions (\ref{eq9-76}) (Prove that this is a vector space). On the vector space $V$, we can define and inner product using the formula
\[(u,w)=\int_{D}\mathrm{d}v\,uw.\]
Then (\ref{eq9-78}) show that
\begin{equation}\label{eq9-79}
(u,Lw)=(Lu,w).
\end{equation}
So, $L$ is a self-adjoint operator on $V$. If we choose different boundary conditions so that the right hand side of (\ref{eq9-75}) does not vanish then (\ref{eq9-79}) will not hold and the operator is not self-adjoint in the sense of linear algebra. It is however still formally self-adjoint according to the definition on the previous page.

\noindent We know from linear algebra that self-adjoint operators have many good properties, we will discuss many of these in the next section of these lecture notes.

An important class of differential operators in the theory of PDEs consists of operators of the form 
\begin{equation}\label{eq9-80}
Lu=a_{ij}\partial_{x_i}\partial_{x_j}u+b_i\partial_{x_i}u+cu,
\end{equation}
where we use the Einstein convention.

\noindent The formal adjoint of this operator is
\begin{equation}\label{eq9-81}
L^*u=\partial_{x_i}\partial_{x_j}(a_{ij}u)-\partial_{x_i}(b_iu)+cu.
\end{equation}
To see why this is so take the single term $wa_{11}\partial_{x_1}\partial_{x_1}u$. We have
\begin{eqnarray*}
wa_{11}\partial_{x_1}\partial_{x_1}u=\partial_{x_1}(wa_{11}\partial_{x_1}u)-\partial_{x_1}(wa_{11})\partial_{x_1}u\\
=\partial_{x_1}(wa_{11}\partial_{x_1}u)-\partial_{x_1}(u\partial_{x_1}(wa_{11}))+u\partial_{x_1}\partial_{x_1}(wa_{11}),\\
\Downarrow\\
wa_{11}\partial_{x_1}\partial_{x_1}u-\partial_{x_1}\partial_{x_1}(wa_{11})u=\partial_{x_1}(wa_{11}\partial_{x_1}u-u\partial_{x_1}(wa_{11})).
\end{eqnarray*}
The other terms are treated in the same way and (\ref{eq9-81}) appears.

\noindent By carrying out all differentiations in (\ref{eq9-81}) and comparing with (\ref{eq9-80}) it is easy to verify that $L$ is formally self-adjoint only if
\[b_i=\partial_{x_j}a_{ij}.\]
If this is the case we have
\begin{equation}\label{eq9-82}
\begin{split}
&\quad\,\, a_{ij}\partial_{x_i}\partial_{x_j}u+b_i\partial_{x_i}u\\
&=a_{ij}\partial_{x_i}\partial_{x_j}u+\partial_{x_j}a_{ij}\partial_{x_i}u\\
&=a_{ij}\partial_{x_i}\partial_{x_j}u+\partial_{x_i}a_{ji}\partial_{x_j}u\\
&=a_{ij}\partial_{x_i}\partial_{x_j}u+\partial_{x_i}a_{ij}\partial_{x_j}u\\
&=\partial_{x_i}(a_{ij}\partial_{x_j}u),
\end{split}
\end{equation}
where we have used the fact that $a_{ij}$ is a symmetric matrix.
From (\ref{eq9-82}) we see that a self-adjoint operator of the form(\ref{eq9-80}) can be rewritten as
\begin{equation}\label{eq9-83}
Lu=\partial_{x_i}(a_{ij}\partial_{x_j}u)+cu.
\end{equation}
If the coefficients are constants (\ref{eq9-83}) imply that a self-adjoint operator can not have any first derivative terms.

  Let us look at some examples.

\noindent The operator
\[Lu=-\nabla\cdot(p\nabla u)+qu,\quad p>0,\]
is an elliptic operator (Show this).

 Observe that
\begin{align*}
wLu-uLw&=-w\nabla\cdot(p\nabla u)+quw+u\nabla\cdot(p\nabla w)-quw\\
       &=\nabla\cdot(pu\nabla w-pw\nabla u),
\end{align*}
so $L$ is formally self-adjoint.

The operator
\[\widetilde{L}u=\rho u_{tt}+Lu,\]
is a hyperbolic operator (Show this).

\noindent Assuming that $\rho$ is independent of time we observe that
\begin{align*}
w\widetilde{L}u-u\widetilde{L}w&=\rho wu_{tt}+wLu-\rho uw_{tt}-uLw\\
&=\partial_t(\rho wu_t-\rho uw_t)+\nabla\cdot(pu\nabla w-pw\nabla u)\\
&=\widetilde{\nabla}\cdot(pu\nabla w-pw\nabla u, \rho w u_t-\rho uw_t),
\end{align*}
where $\widetilde{\nabla}=(\nabla\cdot,\partial_t)$ is a 4D divergence operator. Thus $\widetilde{L}$ is formally self-adjoint.

  Finally define the operator
\[\widehat{L}u=\rho u_t+Lu,\]
where we still assume that $\rho$ is independent of time. This operator is parabolic (Show this).
Define another parabolic operator
\[\widehat{L}^*u=-\rho u_t+Lu.\]
Then
\begin{align*}
w\widehat{L}u-u\widehat{L}^*w&=\rho wu_t+wLu+\rho uw_t-uLw\\
&=\partial_t(\rho uw)+\nabla\cdot(\rho u\nabla w-\rho w \nabla u)\\
&=\widetilde{\nabla}\cdot(\rho u \nabla w-\rho w \nabla u,\rho uw).
\end{align*}
Thus $\widehat{L}^*$ is the formal adjoint of $\widehat{L}$. Since $\widehat{L}^*\neq\widehat{L}$ the operator $\widehat{L}$ is not formally self-adjoint.

The notion of a formal adjoint can be extended to a much larger class of operators\cite{textbook} than (\ref{eq9-80}). 
\section{Initial and boundary value problems in bounded domains}
In this section we will introduce the standard approach for solving PDEs in bounded domains. This approach consists of expanding solutions into infinite series of eigenfunctions. The eigenfunctions are determined by the PDE and the boundary values. Eigenfunctions are generalizations of the normal modes introduced in the previous section. In fact eigenfunctions are often called normal modes.

This approach is very general and work for a large range of linear PDEs. We will see, at the end of this section, that they are also very useful when we investigate nonlinear PDEs.

In these lecture notes, we choose a certain level of generality. 

\noindent Our first restriction is to consider only scalar PDEs. Within the class of scalar PDEs we chose a subclass that is quite general and cover situations that occur in many applications in science and technology. In chapter 8 of  \cite{textbook} a much larger class of PDEs is investigated using similar methods to the ones we discuss here.

Let $L$ be a second order differential operator defined by
\begin{equation}\label{eq11-1}
Lu=-\nabla\cdot(p(\mathbf{x})\nabla u)+q(\mathbf{x})u,
\end{equation}
where $\mathbf{x}\in \mathbb{R}^2$ or $\mathbb{R}^3$. For one space-dimension we define $L$ by
\begin{equation}\label{eq11-2}
Lu=-\partial_x(p(x)\partial_x u)+q(x)u.
\end{equation}
We will in general assume that $p>0$ and $q\eqslantgtr0$. The type of equations we want to consider are

\noindent Elliptic case:
\begin{equation}\label{eq11-3}
Lu=\rho(x)F(\mathbf{x}),\quad \rho>0.
\end{equation}
Parabolic case:
\begin{equation}\label{eq11-4}
\rho(\mathbf{x})\partial_tu+Lu=\rho(\mathbf{x})F(\mathbf{x},t),\quad \rho>0.
\end{equation}
Hyperbolic case:
\begin{equation}\label{eq11-5}
\rho(\mathbf{x})\partial_{tt}u+Lu=\rho(\mathbf{x})F(\mathbf{x},t),\quad \rho>0.
\end{equation}
It is an important point here that the functions $p,\rho$ and $q$ does not depend on time.

These equations, of elliptic, parabolic and hyperbolic type appear as the first approximation to a large set of mathematical modelling problems in science and technology.

Let us consider the case of mass flow. Let $D$ be some domain with boundary $S$.
\begin{figure}[h]
\centering
\includegraphics{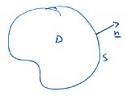}
\caption{}
\end{figure}

\noindent In $D$ we have a concentration of some material substance given by $u(\mathbf{x},t)$. The total mass of the material inside $D$ is
\begin{equation}\label{eq11-6}
M(t)=\int_D\mathrm{d}v\,u(\mathbf{x},t)\rho,
\end{equation}
where $\rho$ have the interpretation of mass density.
\[\text{Mass of substance}\,\,=\,\,\text{concentration of substance}\,\,\times\,\,\text{mass density}\]
This requirement explains the assumption $\rho>0$ in equations (\ref{eq11-3}),(\ref{eq11-4}) and (\ref{eq11-5}). 

We have the balance equation for mass 
\begin{equation}\label{eq11-7}
\frac{\mathrm{d}M}{\mathrm{d}t}=-\int_S\mathrm{d}S \mathbf{j}\cdot\mathrm{n}+\int_D\mathrm{d}V\,H.
\end{equation}
where $\mathbf{j}$ is the mass flux and $H$ is a source of mass. In general, $\mathbf{j}=\mathbf{j}(\mathbf{x},u,\nabla u)$, but for low concentration and small gradients we can Taylor expand and find
\begin{equation}\label{eq11-8}
\mathbf{j}\approx-p(\mathbf{x})\nabla u,\quad p>0.
\end{equation}
Here a constant term and a term only depending on $u$ has been dropped on physical grounds: All mass flux is driven by gradients, and we have introduced a minus sign and assumed $p>0$ in order to accommondate the fact that mass always flow from high to low concentration. In general $H=H(\mathbf{x},t,u)$ but for low concentration we can also here Taylor expand and find
\begin{equation}\label{eq11-9}
H=-q(\mathbf{x})u+\widetilde{F}(\mathbf{x},t),\quad q\eqslantgtr0.
\end{equation}
For notational convenience we introduce $F$ by
\begin{equation}\label{eq11-10}
\widetilde{F}=\rho F.
\end{equation}
The minus sign in the first term is introduced to model a situation where high concentration leads to high loss. Imagine a sink hole where the amount leaving increase in proportion to the mass present.

Using (\ref{eq11-8}), (\ref{eq11-9}) and (\ref{eq11-10}) in (\ref{eq11-7}) we get
\[\int_V\mathrm{d}v\,\rho \partial_tu=\int_S\mathrm{d}S\,p\nabla u\cdot\mathbf{n}+\int_D\mathrm{d}V\left\{-qu+\rho F\right\},\]
\[\Downarrow\quad\text{divergence theorem}\]
\[\int_V\mathrm{d}V\left\{\rho\partial_tu-p\nabla\cdot(p\nabla u)+qu-\rho F\right\}=0.\]
This is assumed to hold for all domains $D$. Therefore we must have
\[\rho\partial_tu+Lu=\rho F,\]
where $L$ is the operator defined in (\ref{eq11-1}).
In a similar way, using conservation of momentum, equation (\ref{eq11-5}) can be shown to model vibrations in homogeneous solids, membranes and bars.

  Our basic model equations must be subject to boundary conditions since the domains in this section are assumed to be bounded. In two and three spatial dimensions we pose the conditions

\begin{equation}\label{eq11-11}
\left.\alpha(\mathbf{x})u+\beta(\mathbf{x})\partial_{\mathbf{n}}u\right|_{\partial G}=B(\mathbf{x},t),
\end{equation}

\begin{figure}
\centering
\includegraphics{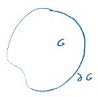}
\caption{}
\end{figure}

\noindent where $\alpha,\beta$ and $B$ are given functions defined on the boundary, $\partial G$, of $G$.

\noindent In the elliptic case, (\ref{eq11-3}), there is no time dependence and $B=B(\mathbf{x})$. For the parabolic and hyperbolic case $G$ is a space-time domain of the form $G=G_{\mathbf{x}}\times[t_0,t_1]$
\begin{figure}[!h]
\centering
\includegraphics{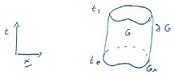}
\caption{}
\end{figure}

\noindent For most common situations occurring in applications, both $\alpha$ and $\beta$ are positive. They can not both be zero at the same point since then the boundary condition (\ref{eq11-11}) would not actually be a constraint!

There are 3 special cases of (\ref{eq11-11}).
\begin{enumerate}
\item $\alpha\neq0,\beta=0$\quad--\,Diricelet condition or condition of the first kind.
\item $\alpha=0,\beta\neq0$\quad--\,Neuman condition or condition of the second kind.
\item $\alpha\neq0,\beta\neq0$\quad--\,Robin condition or condition of the third kind.
\end{enumerate}
It is not uncommon that conditions of the first, second and third kind are specified at different parts of the boundary.
\begin{figure}[h]
\begin{minipage}[2cm]{0.5\textwidth}
\centering
\includegraphics{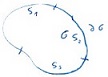}
\caption{}
\end{minipage}%
\begin{minipage}[c][2cm][c]{0.5\textwidth}
Conditions of type $j$ are specified on $S_j, j=1,2,3$.
\end{minipage}
\end{figure}

\noindent For heat conduction problems, part of the boundary could have a given temprature (first kind) whereas other parts could be insulated (second kind).

\noindent For the one dimensional case, $\partial G$ consists of the two end points of the interval $G=(0,l)$. So $\partial G=\{0,l\}$.
For the hyperbolic and parabolic case the general conditions (\ref{eq11-11}) are
\begin{equation}\label{eq11-12}
\begin{split}
\alpha_1u(0,t)-\beta_1u_x(0,t)&=g_1(t),\\
\alpha_2u(l,t)+\beta_2u_x(l,t)&=g_2(t),\quad \alpha_j,\beta_j>0.
\end{split}
\end{equation}
The minus sign at $x=0$ occur because the unit normal at the point $x=0$ is given by the number (-1).

\noindent If $F$ and $B$ in (\ref{eq11-3}), (\ref{eq11-4}), (\ref{eq11-5}) and (\ref{eq11-11}) are zero the corresponding boundary value problem is  {\it homogeneous}.

\subsection{Separation of variables}
We will now apply the method of separation of variables to the equations and boundary conditions introduced in the previous section.

For the {\it hyperbolic case} we consider the homogeneous PDE
\begin{equation}\label{eq11-13}
\rho(\mathbf{x})\partial_{tt}u+Lu=0,\quad\quad\mathbf{x}\in G,t>0,
\end{equation}
and homogeneous boundary conditions
\begin{equation}\label{eq11-14}
\left.\alpha(\mathbf{x})u+\beta(\mathbf{x})\partial_{\mathbf{n}}u\right|_{\partial G}=0,\quad t>0,
\end{equation}
in 2D and 3D.

 In 1D, the homogeneous condition is
\begin{equation}\label{eq11-15}
\begin{split}
\alpha_1u(0,t)-\beta_1u_x(0,t)&=0,\\
\alpha_2u(l,t)+\beta_2u_x(l,t)&=0,\quad t>0.
\end{split}
\end{equation}
In order to get a completely determined problem we also must specify initial conditions 
\begin{equation}\label{eq11-16}
u(\mathbf{x},0)=f(\mathbf{x})\quad u_t(\mathbf{x},0)=g(\mathbf{x})\quad \mathbf{x}\in G.
\end{equation}
For the {\it parabolic case} we consider the homogeneous equation
\begin{equation}\label{eq11-17}
\rho\partial_tu+Lu=0,\mathbf{x}\in G,t>0,
\end{equation}
with the boundary conditions (\ref{eq11-14}) or (\ref{eq11-15}) and the single initial condition
\begin{equation}\label{eq11-18}
u(\mathbf{x},0)=f(\mathbf{x})\quad\mathbf{x}\in G.
\end{equation}
In order to give a uniform treatment of the hyperbolic, parabolic and elliptic case, we introduce for the elliptic case a function $u=u(\mathbf{x},y)$ where now $\mathbf{x}$ is in $\mathbb{R}$ or $\mathbb{R}^2$. The variable $y$ is scalar and restricted to the domain $0<y<\widehat{l}$. The elliptic PDE can in terms of these variables be written as
\begin{equation}\label{eq11-19}
\rho(\mathbf{x})\partial_{yy}u-Lu=0,\quad \mathbf{x}\in G,0<y<\widehat{l},
\end{equation}
where $\rho,p,q$ are functions of $\mathbf{x}$ alone. The equation (\ref{eq11-19}) is supplied with the boundary conditions (\ref{eq11-14}) for the 2D case or (\ref{eq11-15}) for the 1D case. The analog of the initial conditions in the time dependent case are the following conditions at $y=0$ and $y=l$. 
\begin{equation}\label{eq11-20}
u(\mathbf{x},0)=f(\mathbf{x}), u(\mathbf{x},\widehat{l})=g(\mathbf{x}),\mathbf{x}\in G.
\end{equation}
This assumed form for the elliptic equation certainly does not cover all possibilities that occur in applications. We observe for example that the domain always is a generalized cylinder
\begin{figure}[!h]
\centering
\includegraphics{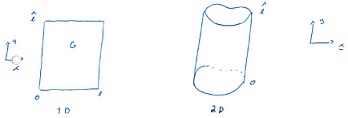}
\caption{}
\end{figure}

\noindent However, the assumed form makes a unified presentation for the hyperbolic, parabolic and elliptic case possible.

After all this preparatory work, let us start separating variables!

\noindent The method of separation of variables starts by looking for special solutions
\begin{equation}\label{eq11-21}
u(\mathbf{x},t)=M(\mathbf{x})N(t),
\end{equation}
for (\ref{eq11-13}) and (\ref{eq11-14}) and
\begin{equation}\label{eq11-22}
u(x,y)=M(\mathbf{x})N(y),
\end{equation}
for (\ref{eq11-19}). Inserting (\ref{eq11-21}) and (\ref{eq11-22}) into (\ref{eq11-13}), (\ref{eq11-17}) and (\ref{eq11-19}) gives, upon division by $NM$
\begin{align}\label{eq11-23} 
\frac{N''(t)}{N(t)}&=-\frac{LM(\mathbf{x})}{\rho(\mathbf{x})M(\mathbf{x})}, \quad \quad\text{hyperbolic case},\\  \label{eq11-24}
\frac{N'(t)}{N(t)}&=-\frac{LM(\mathbf{x})}{\rho(\mathbf{x})M(\mathbf{x})},\quad\quad\text{parabolic case},\\ \label{eq11-25}
-\frac{N''(y)}{N(y)}&=-\frac{LM(\mathbf{x})}{\rho(\mathbf{x})M(\mathbf{x})},\quad\quad\text{elliptic case}. 
\end{align}
Since the two sides of equations (\ref{eq11-23}) $\rightarrow$ (\ref{eq11-25}) depends on different variables, both sides must be equal to the same constant. This constant we write, by convention, in the form $-\lambda$. This give us in all three cases the following equation for $M$.
\begin{equation}\label{eq11-26}
LM(\mathbf{x})=\lambda\rho(\mathbf{x})M(\mathbf{x}),
\end{equation}
and
\begin{align}\label{eq11-27}
N''(t)+\lambda N(t)&=0,\quad \text{hyperbolic case},\\ \label{eq11-28}
N'(t)+\lambda N(t)&=0,\quad \text{parabolic case},\\ \label{eq11-29}
N''(t)-\lambda N(y)&=0,\quad \text{elliptic case}.
\end{align}
The function $M(\mathbf{x})$ in equation (\ref{eq11-26}) is subject to boundary conditions (\ref{eq11-14}) or (\ref{eq11-15}). Since both conditions are homogeneous, as are equation (\ref{eq11-26}), $M(\mathbf{x})=0$ is a solution. This solution, also called the {\it trivial solution}, gives $u=0$ and thus plays no role.

\noindent Defining $A=\frac1\rho L$, we observe that the structure of equation (\ref{eq11-26}) is
\[AM=\lambda M\]
Thus $M$ is an eigenvector of $A$ and $\lambda$ the corresponding eigenvalue. The vector space is the set of smooth functions defined on $G$ and satisfying the homogeneous boundary conditions (\ref{eq11-14}) or (\ref{eq11-15}) on the boundary of $G$. (Show that this is a vector space!) In this vector space, vectors are functions, so solutions of (\ref{eq11-26}) are often called {\it eigenfunctions}. The big difference between the eigenvalue problems that are studied in elementary linear algebra and the current one, is that here, the vector space is infinite dimensional. This makes the problem much harder, much more interesting and also much more useful.

\subsection{Self-adjoint and positive operators}
The operator
\begin{equation}\label{eq11-30}
L=-\nabla\cdot(p(\mathbf{x})\nabla u)+q(\mathbf{x})u,
\end{equation}
occurring in equation (\ref{eq11-26}) and determining the functions $M(\mathbf{x})$ in the separation
\[u=NM,\]
for the hyperbolic case (\ref{eq11-13}), the parabolic case (\ref{eq11-17}) and the elliptic case (\ref{eq11-19}), is very special. We will now determine some of the properties of the {\it spectrum} of $L$. The spectrum of a linear operator is the collection of eigenvalues and corresponding eigenspaces of the operator. Recall that the eigenspace corresponding to a given eigenvalue consists of all eigenvectors of the operator that has the same eigenvalue.

\noindent Using vector calculus we observe that
\begin{align*}
w\nabla \cdot(p\nabla u)&=\nabla\cdot(pw\nabla u)-p\nabla w\cdot \nabla u,\\
u\nabla\cdot(p\nabla w)&=\nabla\cdot(pu\nabla w)-p\nabla u\cdot \nabla w,
\end{align*}
\[\Rightarrow
w \nabla\cdot(p\nabla u)-u\nabla\cdot(p\nabla w)=\nabla\cdot(pw\nabla u-pu\nabla w),\]
\[\Downarrow\]
\begin{equation}\label{eq11-31}
wLu-uLw=\nabla\cdot(pu\nabla w-pw\nabla u).
\end{equation}
Thus $L$ is formally self-adjoint.

  Let $G$ be a domain in $\mathbb{R}^2$ or $\mathbb{R}^3$. Integrating (\ref{eq11-31})
over $G$ gives
\begin{equation}\label{eq11-32}
\int_G \mathrm{d}V\,\{wLu-uLw\}=\int_{\partial G}\mathrm{d}S\,p(u\partial_{\mathbf{n}}w-w\partial_\mathrm{n}u).
\end{equation}
Let us assume that both $u$ and $w$ satisfy the boundary conditions (\ref{eq11-11}) with $B=0$. Thus
\begin{equation}\label{eq11-33}
\begin{split}
\alpha(\mathbf{x})u+\beta(\mathbf{x})\partial_{\mathbf{n}}u\big|_{\partial G}&=0,\\
\alpha(\mathbf{x})w+\beta(\mathbf{x})\partial_{\mathbf{n}}w\big|_{\partial G}&=0.
\end{split}
\end{equation}
For a given $u,w$ that satisfy (\ref{eq11-33}) we can think of (\ref{eq11-33}) as a linear system for $\alpha$ and $\beta$


\begin{equation} \label{eq11-34}
\begin{bmatrix}
u&\partial_{\mathbf{n}} u \\ 
w & \partial_{\mathbf{n}} w
\end{bmatrix} \bigg|_{\partial G} 
\begin{bmatrix} 
\alpha \\ \beta
\end{bmatrix} = 0.
\end{equation}
But our assumptions on the functions $\alpha$ and $\beta$ imply that they can not both be zero at any point. Thus, the system (\ref{eq11-34}) must have a nonzero solution. This imply that the determinant of the matrix defining the linear system (\ref{eq11-34}), must be zero or 
\begin{equation} \label{eq11-35}
(u\,\partial_{\mathbf{n}}w-w\,\partial_{\mathbf{n}}w)|_{\partial G}=0. 
\end{equation}
Thus, the right hand side of (\ref{eq11-32}) is zero 
\begin{equation} \label{eq11-36}
\Rightarrow
\int\limits_{G}  dV \lbrace wLu-uLw \rbrace =0.
\end{equation}
On the vector space, $V$, consisting of real valued functions defined on $G$ and satisfying the homogeneous boundary conditions (\ref{eq11-33}), we define an inner product 
\begin{equation}
(f,g)=\int\limits_{G} dV \rho(\textbf{x})f(\textbf{x})g(\textbf{x}). 
\end{equation}
Given this inner product we define the \emph{norm} of functions in $V$ as
\begin{equation}
\|f\|=\left( (f,f) \right)^{1/2} =\sqrt{\int_{G}dV\,\rho(\mathbf{x})f^2(\mathbf{x})}.
\end{equation} 
Assuming $f$ is continuous we have 
\begin{equation}
\|f\|=0 \Leftrightarrow f=0. 
\end{equation}
Thus, what we have defined is a norm in the usual sense of linear algebra. Using this inner product (\ref{eq11-36}) can be written as $$(w, \frac{1}{\rho}Lu)=(\frac{1}{\rho}Lw,u). $$ Thus the operator $\tilde{L}=\frac{1}{\rho}L$ is self-adjoint in the sense of linear algebra. In the most common cases when $\rho=1$, $L$ itself is self-adjoint.

Let now $M_m$, $M_l$ be eigenfunctions of $L$ corresponding to \emph{different} eigenvalues, $\lambda_k\neq\lambda_l$. Since both $M_k$ and $M_l$ satisfy the homogeneous boundary conditions, (\ref{eq11-33}), we have $$\int\limits_{G}dV\,\lbrace M_k\,L\,M_l-M_l\,L\,M_k\rbrace =0, $$
but 
\begin{align*}
LM_k & =\rho\lambda_k M_k,\\
LM_l & =\rho\lambda_l M_l.
\end{align*} 
Thus
\begin{align*}0=\int\limits_{G}dV\,\lbrace M_k\,L\,M_l-M_l\,L\,M_k\rbrace = \int\limits_{G} dV\,\rho(\lambda_l-\lambda_k)M_k M_l \\
=(\lambda_l-\lambda_k)(M_k,M_l). 
\end{align*}
But $\lambda_k\neq\lambda_l \Rightarrow \lambda_k-\lambda_l\neq 0$, and therefore we must have $(M_k, M_l)=0$.\\

\emph{"Eigenfunctions corresponding to different eigenvalues are orthogonal."}\\ 

\noindent{Actually we have in our choice of inner product assumed that the eigenfunctions are real valued and the spectrum real. This is in fact true:}

Let $\varphi$ be a complex valued function. Then, redoing the argument leading up to (\ref{eq11-32}) using $w=u^*$ (the complex conjugate of $u$), we get 
\begin{equation}\label{eq11-38}
\int\limits_{G} dV\,\lbrace u(Lu)^*-(Lu)u^*\rbrace=\int\limits_{\partial G} dS\, \rho\lbrace u^*\partial_{\mathbf{n}}u-u \partial_{\mathbf{n}}u^*\rbrace,
\end{equation}
and (\ref{eq11-35}) becomes

\[ \lbrace u^*\partial_{\mathbf{n}} u- u\partial_{\mathbf{n}}u^*\rbrace\big|_{\partial G}=0,\label{eq11-39}\]
\[\Downarrow\]
\[\int\limits_{G}dV\,\lbrace uLu^* -Luu^*\rbrace=0.\]

\noindent Let now $\lambda$ be a complex eigenvalue and $M$ the corresponding complex valued eigenfunction. Then (\ref{eq11-39}) imply that
\begin{align*}
0 &=  \int\limits_{G}dV\,\lbrace M(\rho\lambda M)^*- \rho\lambda M M^*\rbrace\\
&=\int\limits_{G}dV\,\lbrace \rho\lambda^* M M^*-\rho\lambda M M^*\rbrace\\
&= (\lambda^*-\lambda)\int\limits_{G}dV\,\rho M M^* .
\end{align*}
But $$\int\limits_{G} dV\,\rho |M|^2=0\Rightarrow M=0,$$ so we must have $$\lambda^*=\lambda,$$ which means that the eigenvalues are \emph{real}. The eigenfunctions can then always be chosen to be real. 

Let now $u$ be a (real valued) function that satisfies the boundary condition (\ref{eq11-11}) with $B=0$. Thus $u\in V$, as defined on page 108. 
Then we have 
\begin{align*}
\int\limits_{G}dV\, uLu =& -\int\limits_{G}dV\,\lbrace\nabla\cdot(pu\nabla u)-p(\nabla u)^2\rbrace +\int\limits_{G}dV\, q u^2 \\
=& \int\limits_{G} dV\,\lbrace p(\nabla u)^2 +qu^2 \rbrace -\int\limits_{\partial G} dS\, pu\partial_\mathbf{n} u. 
\end{align*}
Split the boundary $\partial G$ into three pieces $S_1, S_2$ and $S_3$, where on $S_i$ the boundary conditions are of  kind $i$, where $i=1,2,3$ (see page 7). Then on $S_1\cup S_2, pu\partial_{\mathbf{n}}u=0$ and on $S_3$, $pu\partial_{\mathbf{n}}u=-p\frac{\alpha}{\beta}u^2$
\begin{equation}\label{eq11-40}
\Rightarrow \int\limits_{G}dV\, uLu = \int\limits_{G}dV\,\lbrace p(\nabla u)^2 +qu^2\rbrace + \int\limits_{S_3}dS\,p\frac{\alpha}{\beta}u^2. 
\end{equation}
The assumptions we have put on $p,\alpha,\beta$ and $q$ earlier implies that the right hand side of (\ref{eq11-40}) is positive. Thus, using the inner product defined on page 16, we have found that 
\begin{equation}\label{eq11-41}
(\frac{1}{\rho}Lu,u)\geq 0.
\end{equation}
This means, by definition, that $\frac{1}{\rho}L$ is a \emph{positive} operator.

Let now $M_k$ be an eigenfunction corresponding to the eigenvalue $\lambda_k$. From (\ref{eq11-41}) we have 
\begin{equation*}
(M_k,\frac{1}{\rho}LM_k)=\int\limits_{G}dV\, M_kLM_k=\int\limits_{G}dV\, M_k\rho\lambda_kM_k=\lambda_k\int\limits_{G}dV\, \rho M_k^2\geq 0,
\end{equation*}
and 
\begin{equation*}
\int\limits_{G}dV\, \rho M^2_k> 0. 
\end{equation*}
Thus we must have $$\lambda_k\geq 0.$$
Therefore, the spectrum of $L$ is not only real but non-negative. 

\begin{figure}[h!]
	\centering
	\includegraphics[width=10cm]{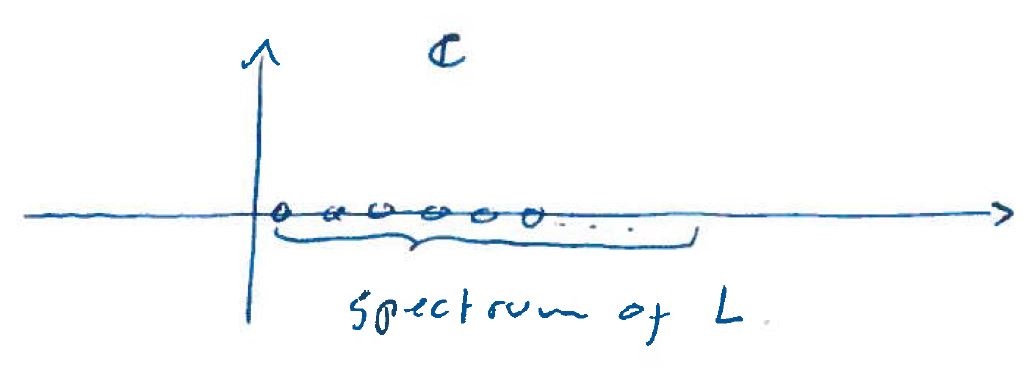}
	\caption{}
\end{figure}

One can prove\cite{textbook} that there are countably infinite many eigenvalues $\lbrace\lambda_k\rbrace_{k=1}^{\infty}$ and that $\infty$ is the only limit point of the sequence of eigenvalues. Recall that a limit point for a sequence $\lbrace a_n\rbrace$ is a point $a^*\in\mathbb{R}$ such that $\forall \epsilon> 0$ there exists infitely many $n$ such that $$|a_n-a^*|<\epsilon.$$ Thus the eigenvalues only "accumulate" around $\lambda=\infty$. 

Let now $\lbrace\lambda_k\rbrace$ be the countably infinite set of eigenvalues. One can also prove that there are at most a finite set of linearly independent eigenvectors corresponding to each $\lambda_k$. We can assume that these has been orthogonalized using the Gram-Schmidt orthogonalization process. 

Let us in the following for simplicity only consider the case when each eigenspace is of dimension one and where all $k$,  $\lambda_k >0$. Thus $\lambda=0$ is not an eigenvalue. Then the equations (\ref{eq11-27}-\ref{eq11-29}) for the function $N$ gives upon solution the following:
\begin{eqnarray}
N_k(t)=& a_k \cos{\sqrt{\lambda_k}t}+b_k \sin{\sqrt{\lambda_k}t},&\quad\text{hyperbolic case,}\label{eq11-42}\\
N_k(t)=& a_k e^{\lambda_k t},&\quad\text{parabolic case,}\label{eq11-43}\\ 
N_k(y)=& a_k e^{\sqrt{\lambda_k}y}+b_k e^{-\sqrt{\lambda_k}y},&\quad\text{elliptic case.}\label{eq11-44}
\end{eqnarray}
Here, $a_k$, $b_k$ are arbitrary constants, one for each eigenvalue $\lambda_1, \lambda_2,\dots$. We now have infinitely many solutions $$u_k=N_kM_k,$$ that solves equations (\ref{eq11-13}), (\ref{eq11-17}) or (\ref{eq11-19}) and the boundary conditions (\ref{eq11-14}) or (\ref{eq11-15}). These functions will however not in general satisfy the initial conditions. In order to also satisfy the initial conditions (\ref{eq11-10}) for the hyperbolic case, (\ref{eq11-18}) for the parabolic case and (\ref{eq11-20}) for the elliptic case we consider formal superposition 
\begin{equation}\label{eq11-45}
u=\sum\limits_{k=1}^{\infty}u_k=\sum\limits_{k=1}^{\infty}M_k\,N_k.
\end{equation}
For the hyperbolic case we get the conditions 
\begin{equation} \label{eq11-46}
u(\mathbf{x},0)=\sum\limits_{k=1}^{\infty}M_k(\mathbf{x})N_k(0)=\sum\limits_{k=1}^{\infty}a_kM_k(\mathbf{x})=f(\mathbf{x}),
\end{equation}

\begin{equation} \label{eq11-47}
u_t(\mathbf{x},0)=\sum\limits_{k=1}^{\infty}M_k(\mathbf{x})N'_k(0)=\sum\limits_{k=1}^{\infty}\sqrt{\lambda_k}b_kM_k(\mathbf{x})=g(\mathbf x).
\end{equation}
(\ref{eq11-46}) and (\ref{eq11-47}) require us to find \emph{eigenfunction expansions} for the the functions $f$ and $g$ in terms of the eigenfunctions $\lbrace M_k \rbrace^{\infty}_{k=1}$. In what  sense $f$ and $g$ can be represented by their eigenfunction expansions is an important, and in general very hard, question that will be discussed more fully only for the 1D case in this class.

You might recall from elementary calculus that for infinite series of functions, several different notions of convergence are available. The most important ones are pointwise convergence, uniform convergence and mean-square convergence. We will see that the last type of convergence will hold for very mild conditions on $f$ and $g$. However, mean-square convergence is not enough to ensure that (\ref{eq11-45}) is an actual, honest to God, solution to our initial-boundary value problem (\ref{eq11-13}),(\ref{eq11-14}),(\ref{eq11-16}). For this to be the case we need a stronger convergence than mean-square. This stronger form of convergence will allow us to differentiate the formal superposition series (\ref{eq11-45}) term by term. If this is possible, we have a \emph{classical} solution. If it is not possible, (\ref{eq11-45}) is a \emph{generalized solution} that might still carry important information about the problem. 

In order to express the coefficients $\lbrace a_k \rbrace $ and $\lbrace b_k \rbrace $ in terms of $f$ and $g$, we use the fact that the eigenfunctions $\lbrace M_k \rbrace $ are orthogonal for $k\neq k'$ 
\begin{align*}
(f,M_j)&=\sum\limits_{k=1}^{\infty}a_k(M_k, M_j) = a_j (M_j, M_j),
\end{align*}
\begin{equation}\label{eq11-48}
\Rightarrow\medspace a_j=\frac{(f, M_j)}{(M_j, M_j)}.
\end{equation}
In the same way we find
\begin{equation}\label{eq11-49}
b_j=\frac{(g, M_j)}{\sqrt{\lambda_j}(M_j, M_j)}.
\end{equation}
The \emph{formal} solution to the hyperbolic initial-boundary value problem (\ref{eq11-13}-\ref{eq11-15}) is then 

\begin{equation}\label{eq11-50}
u(\mathbf x,t)=\sum\limits_{k=1}^{\infty}(a_k\cos \sqrt{\lambda_k t}+b_k\sin{\sqrt{\lambda_k}}t)\,M_k(\mathbf x), 
\end{equation}
where $\lbrace a_k\rbrace$ and $\lbrace b_k\rbrace$ are given by (\ref{eq11-48}) and (\ref{eq11-49}). 

In a similar way, the solution of the initial-boundary problem for the parabolic case is 
\begin{equation}\label{eq11-51}
u(\mathbf x,t)=\sum\limits_{k=1}^{\infty}a_k e^{-\lambda_k t}\,M_k(\mathbf x),
\end{equation}
where  $\lbrace a_k\rbrace$ is given by (\ref{eq11-48}). 

For the elliptic case we get the conditions
\begin{align*}
u(\mathbf x, 0)&= \sum\limits_{k=1}^{\infty} (a_k+b_k)M_k(\mathbf x)=f(\mathbf x), \\
u(\mathbf x, \tilde l)&= \sum\limits_{k=1}^{\infty} (a_k e^{\lambda_k \tilde l}+b_k e^{-\lambda_k \tilde l})M_k(\mathbf x)=g(\mathbf x).
\end{align*}
Proceeding as for the hyperbolic case we get
\begin{eqnarray}\label{eq11-52}
a_k+b_k=(f, M_k),\nonumber \\
e^{\lambda_k \tilde l}a_k+ e^{-\lambda_k \tilde l}b_k=(g, M_k).
\end{eqnarray}
These equations can be solved to give unique values for $\lbrace a_k \rbrace$ and $\lbrace b_k \rbrace$. Our formal solution for the elliptic case is then

\begin{equation}\label{eq11-53}
u(\mathbf x, y)= \sum\limits_{k=1}^{\infty} (a_k e^{\lambda_k y}+b_k e^{-\lambda_k y})\,M_k(\mathbf x).
\end{equation}
We will now take a closer look at the formal solutions (\ref{eq11-50}), (\ref{eq11-51}) and (\ref{eq11-53}) for the 1D case. 

\section{The Sturm-Liouville problem and Fourier series}

The 1D version of the eigenvalue problem (\ref{eq11-26}) + (\ref{eq11-14}) is 
\begin{equation}\label{eq11-52b} 
L[u(x)]=-\frac{\mathrm d}{\mathrm{d}x}\left[p(x)\,\frac{\mathrm{d}u}{\mathrm{d}x}\right] + q(x)u(x)=\lambda \rho(x)u(x),
\end{equation}
for $0<x<l$ and boundary conditions 
\begin{eqnarray}\label{eq11-53b} 
&\alpha_1 u(0)-\beta_1u'(0)=0,\nonumber \\
&\alpha_2 u(l)+\beta_2u'(l)=0.
\end{eqnarray}
The problem defined by (\ref{eq11-52b}) and (\ref{eq11-53b}) is known as the \emph{Sturm-Liouville problem}. If we require the following additional conditions

\begin{align}\label{eq11-54}
& i) \quad  p(x)>0,\, \rho(x)>0,\,q(x)\geq 0,\quad 0\leq x\leq l, \nonumber\\
& ii)\quad  p(x),\, p'(x),\,\rho(x),\,q(x)\text{ are continous for }0\leq x\leq l,\\
& iii) \quad \alpha_i\geq 0,\, \beta_i\geq 0, \,\alpha_i+\beta_i>0,\nonumber
\end{align}

\noindent to hold, (\ref{eq11-52b}), (\ref{eq11-54}) define a \emph{regular Sturm-Liouville} problem. If we relax one or more condition on the coefficient functions $p, q$ and $\rho$ we get more generalized versions of the Sturm-Liouville problem. For example, one can allow the functions $p$ and/or $\rho$ to vanish at the endpoints $x=0, x=l$. We then get a \emph{singular Sturm-Liouville} problem. This generalization do occur in important applications of this theory.

Before we start developing the theory, let us introduce some important constructions. These constructions have already been introduced in the more general setting discussed earlier, but will be repeated here in the 1D context. 

For real valued functions on $0\leq x\leq l$ which are bounded and  integrable, we introduce the \emph{inner product} 

\begin{equation}\label{eq11-55}
(\varphi, \psi) = \int_{0}^{l} \rho(x) \varphi(x)\psi(x)\,\mathrm{d}x.
\end{equation}
Observe that (\ref{eq11-55}) is evidently bilinear and symmetric. For continuous functions we have 
\begin{equation*}
\int_{0}^{l} \rho(x) \varphi(x)^2\,\mathrm{d}x=0, \Rightarrow \varphi(x)=0\quad \forall x.
\end{equation*}
 Thus, (\ref{eq11-55}) really defines an inner product in the sense of linear algebra on the linear vector space consisting of continuous functions on $0 \leq x \leq l$. 

Introduce the norm of a function by 
\begin{equation}\label{eq11-56}
\|\varphi\|=\left(\int_0^l\rho(x)\varphi^2(x)\right)^{1/2}.
\end{equation}
A function such that $\|\varphi\|<\infty$ is said to be \emph{square-integrable}. Furthermore, a square-integrable function is \emph{normalized to unity} if 
\begin{equation}\label{eq11-57}
\|\varphi\|=1.
\end{equation}
If $\varphi$ is any square-integrable function we can \emph{normalize} this function by defining a new function

\begin{gather*}
\hat{\varphi} =\frac{1}{\|\varphi\|}\varphi,\\
\Downarrow\\
\|\hat{\varphi}\|=\|\frac{1}{\|\varphi\|}\varphi\|=\frac{1}{\|\varphi\|}\|\varphi\|=1.
\end{gather*}

\noindent Two functions $\varphi$ and $\psi$ are \emph{orthogonal} on the interval $0<x<l$ if

\begin{equation}\label{eq11-58}
(\varphi,\psi)=0.
\end{equation}
A set of functions $\lbrace \varphi_k(x)\rbrace,\, k=1,2,\dots$ is said to be an \emph{orthogonal set} if

\begin{equation}\label{eq11-59}
(\varphi_k,\varphi_j)=0,\quad k\neq j.
\end{equation}
If also $\|\varphi_k\|=1$ for all $k$, the set is \emph{orthonormal}. An orthogonal set can always be turned into an orthonormal set by using the approach described above.

If we allow the functions to assume complex values, we replace the inner product (\ref{eq11-55}) by

\begin{equation}\label{eq11-60}
(\varphi,\psi)=\int_0^l\rho(x)\varphi(x)\overline{\psi(x)}\,\mathrm{d}x,
\end{equation}
where for any complex number $z$, $\overline{z}$ is its complex conjugate. Observe that we have bilinearity and 

\begin{equation}\label{eq11-61}
	(\varphi,\psi)=(\overline{\psi,\varphi})
\end{equation}
Furthermore for any $\varphi$ 

\begin{equation}\label{eq11-62}
(\varphi,\varphi)=\int_0^l\rho(x)|\varphi(x)|^2\,\mathrm{d}x\geq 0,
\end{equation}
and if $\varphi$ is, say, continuous then 

\begin{equation}\label{eq11-63}
(\varphi,\varphi)=0\Rightarrow\varphi(x)=0,\quad\forall x.
\end{equation}
Thus (\ref{eq11-60}) defines a \emph{Hermitian} inner product in the sense of linear algebra.

  Let us return to the real valued case. Let an orthonormal set of functions $\lbrace \varphi_k \rbrace $ be given and let $\varphi$ be a square-integrable function on $0<x<l$. Then the infinite set of numbers $(\varphi,\varphi_k)$ are called the \emph{Fourier coefficients} of $\varphi$ with respect to $\lbrace\varphi_k\rbrace$. The formal series 

\begin{equation}\label{eq11-64}
\sum\limits_{k=1}^{\infty}(\varphi,\varphi_k)\,\varphi_k,
\end{equation}
is called the \emph{Fourier series} of $\varphi$ with respect to the set $\lbrace \varphi_k \rbrace$. Note that we use the summation sign in (\ref{eq11-64}) in a formal sense - we have at this point not assumed anything about the convergence of the infinite series of functions (\ref{eq11-64}). 

The name, Fourier series, is taken from the special case when $\lbrace \varphi_k \rbrace$ is a set of trigonometric functions. Often in the literature the term Fourier series refer \emph{only} to this trigonometric case, but we will use the term Fourier series for all formal series (\ref{eq11-64}).

So, in what sense does (\ref{eq11-64}) converge? And if it converge, does it converge to the function $\varphi$ that we used to construct the formal series?

 \noindent  Let us start by considering a finite sum,
 \begin{equation*}
  S_N=\sum\limits_{k=1}^{N}(\varphi,\varphi_k)\,\varphi_k.
  \end{equation*}
  For this sum we have
\begin{align*}
&\|\varphi-\sum_{k=1}^N(\varphi, \varphi_k)\,\varphi_k\|^2= (\varphi-\sum_{k=1}^N(\varphi, \varphi_k)\,\varphi_k,\, \varphi-\sum_{j=1}^N(\varphi, \varphi_j)\,\varphi_j)\\
& = (\varphi,\varphi)-\sum_{k=1}^{N}(\varphi,\varphi_k)^2-\sum_{j=1}^N(\varphi_j,\varphi_j)^2 + \sum\limits_{k,j=1}^N (\varphi,\varphi_k)(\varphi,\varphi_j)(\varphi_k,\varphi_j)\\
& = \|\varphi\|^2-2\sum_{k=1}^N(\varphi, \varphi_k)^2+\sum_{k=1}^N(\varphi, \varphi_k)^2\\
& = \|\varphi\|^2-\sum_{k=1}^N(\varphi, \varphi_k)^2.
\end{align*}
But the norm of any function is non-negative. Thus we must have

\begin{gather*}
\|\varphi\|^2-\sum_{k=1}^N(\varphi, \varphi_k)\geq 0,\\
\Downarrow\\
\sum_{k=1}^N(\varphi, \varphi_k) \leq \|\varphi\|^2.
\end{gather*}
This inequality says that \emph{all} partial sums of the \emph{numerical series} $\sum\limits_{k=1}^N(\varphi, \varphi_k)^2$, consisting entirely of positive terms, are bounded by $\|\varphi\|^2$. It then follows, from the fact that all increasing bounded sequences are convergent that the series $\sum_{k=1}^N(\varphi, \varphi_k)^2$ is \emph{convergent} and that

\begin{equation}\label{eq11-65}
\sum_{k=1}^\infty(\varphi, \varphi_k)^2\leq \|\varphi\|^2,
\end{equation}
holds. (\ref{eq11-65}) is known as \emph{Bessel's inequality}. 

A sequence of square-integrable functions $\lbrace S_N\rbrace_{N=1}^{\infty}$ is said to converge to a function $\varphi(x)$ \emph{in the mean} if 

\begin{equation*}
\lim\limits_{N\rightarrow\infty} \|\varphi-S_N\|=0.
\end{equation*}
This kind of convergence does not in general imply pointwise convergence. It is called \emph{mean square} convergence. Observe that if we have equality in Bessel's inequality (\ref{eq11-65}),

\begin{equation}\label{eq11-66}
\sum_{k=1}^\infty(\varphi, \varphi_k)^2= \|\varphi\|^2,
\end{equation}
then we have 

\begin{align*}
\lim_{N\rightarrow\infty}\|\varphi-S_N\|^2&=\lim\limits_{N\rightarrow\infty}\|\varphi-\sum\limits_{k=1}^N(\varphi,\varphi_k)\,\varphi_k\|^2\\
&= \lim\limits_{N\rightarrow\infty}\lbrace\|\varphi\|^2-\sum\limits_{k=1}^N(\varphi,\varphi_k)^2\rbrace=0.
\end{align*}
(\ref{eq11-66}) is called \emph{Parsevals's equality}. The argument we just did shows that if Parseval's equality holds, for a given $\varphi$ and given orthonormal set $\lbrace\varphi_k\rbrace$, then the Fourier series (\ref{eq11-64}) converges to $\varphi$ in the mean. 

An orthonormal set $\lbrace \varphi_k\rbrace _{k=1}^{\infty}$ is said to be \emph{complete}  if the Parseval equality (\ref{eq11-66}) holds for all square integrable functions. We have proved that if an orthonormal set of functions is complete then the formal Fourier series (\ref{eq11-64}) for any square integrable functions $\varphi$, converge in the mean square square sense to the function $\varphi$.

\subsection{Properties of eigenvalues and eigenfunctions}
The properties of the eigenvalues and eigenspaces are the same as for the 2D and 3D case discussed earlier. After all, they follow from structural properties of $L$ and the boundary conditions that are completely analogous in 1D, 2D and 3D. Let us just record these properties here:

\begin{enumerate}
	\item Eigenfunctions corresponding to different eigenvalues are orthogonal.
	\item The eigenvalues are real and non-negative and the eigenfunctions can be chosen to be real.
	\item Each eigenvalue is simple. This means that each eigenspace is of dimension one.
	\item There is a countable infinity of eigenvalues having a limit point at infinity. The set of eigenvalues can be arranged as $$0\leq\lambda_1<\lambda_2<\lambda_3<\dots,$$ with $\lim\limits_{k\rightarrow\infty}\lambda_k=\infty$.
	\item The set of eigenfunctions $\lbrace u_k(x)\rbrace_{k=1}^{\infty}$ form a complete orthonormal set of square integrable functions. Thus, for any square integrable function on $0<x<l$,  the Fourier series $$u=\sum\limits_{k=1}^{\infty}(u,u_k)\,u_k,$$ converge to $u$ in the mean.
	\item If $u(x)$ is continuous on $0\leq x\leq l$ \emph{and} has a piecewise continuous first derivative on $0\leq x\leq l$ \emph{and} $u(x)$ satisfies the boundary conditions (\ref{eq11-53b}) then the Fourier series (\ref{eq11-64}) converge absolutely and uniformly to $u(x)$ i  $0\leq x\leq l$.
	\item If $u(x)$ has some jump discontinuity at a point $x_0$ in the interior of the interval  $0\leq x\leq l$ then (\ref{eq11-64}) converge to $\frac{1}{2}(u(x_0^-)+u(x_0^+))$ at this point.
\end{enumerate}
Here
\begin{align*}
	& u(x_0^-)=\lim\limits_{x\rightarrow x_0^-} u(x),\quad\text{(approaching from the left)},\\
	& u(x_0^+)=\lim\limits_{x\rightarrow x_0^+} u(x),\quad\text{(approaching from the right)},
\end{align*}
so the only thing that remains to do is to actually find the eigenvalues and eigenfunctions for particular cases of interest. We will do so for a few important cases. 

All regular Sturm-Liouville problems can be simplified if we use the following uniform approach:

Let $v(x;\lambda),\,w(x;\lambda) $ be solutions of the initial value problem 

\begin{equation}\label{eq11-67}
-\frac{\mathrm{d}}{\mathrm d x}\left[p(x)\frac{\mathrm d u}{\mathrm d x}\right]+q(x)u(x)=\lambda+\rho(x)u(x),
\end{equation}
with initial conditions
\begin{align}\label{eq11-68}
& v(0;\lambda)=1,\quad w(0;\lambda)=0,\nonumber\\
& v'(0;\lambda)=0,\quad w'(0;\lambda)=1.
\end{align}
Then it is easy to verify that the function 

\begin{equation}\label{eq11-69}
u(x;\lambda)=\beta_1v(x;\lambda)+\alpha_1w(x;\lambda),
\end{equation}
is a solution to (\ref{eq11-67}) and satisfies the boundary condition

\begin{equation*}
\alpha_1 u(0;\lambda)-\beta_1u'(0;\lambda)=0.
\end{equation*}
In order to also satisfy the boundary condition at $x=l$ we must have

\begin{equation}\label{eq11-70}
\alpha_2\beta_1v(l;\lambda)+\alpha_2\alpha_1w(l;\lambda)+\beta_2\beta_1v'(l;\lambda)+\beta_2\alpha_1w'(l;\lambda)=0.
\end{equation}
The eigenvalues are solutions of equation (\ref{eq11-70}) and the corresponding eigenfunctions are determined by (\ref{eq11-69}).

We will now apply the uniform approach to several cases leading to various trigonometric eigenfunctions and then turn to a couple of important examples of singular Sturm-Liouville problems.

\subsubsection{Trigonometric eigenfunctions}
For this case we choose $p(x)=\rho(x)=1$ and $q(x)=0$. Equation (\ref{eq11-52}) in the Sturm-Liouville problem for this case reduces to 

\begin{equation}\label{eq11-71}
-u''(x)=\lambda u(x).
\end{equation}
It is easy to verify that the functions

\begin{align}\label{eq11-72}
& v(x;\lambda)=\cos(\sqrt{\lambda} x), \\
& w(x;\lambda)=\frac{\sin(\sqrt{\lambda }x)}{\sqrt{\lambda}},
\end{align}
satisfy (\ref{eq11-71}) and the initial values (\ref{eq11-68}). Here we assume $\lambda>0$. We leave the case $\lambda=0$ as an exercise for the reader. The eigenfunction for $\lambda>0$ is found from (\ref{eq11-69}) and is 

\begin{equation}\label{eq11-73}
u(x)=\beta_1\cos(\sqrt{\lambda} x) +\alpha_1\frac{\sin(\sqrt{\lambda} x)}{\sqrt{\lambda}}.
\end{equation}
The eigenvalues are determined from the equation (\ref{eq11-69}), which for this case simplifies into

\begin{equation}\label{eq11-74}
\sqrt{\lambda}(\alpha_1\beta_2+\beta_1\alpha_2)\cos(\sqrt{\lambda} l)+(\alpha_2+\alpha_1-\lambda\beta_2\beta_1)\sin(\sqrt{\lambda} l)=0.
\end{equation}
We will get different sets of eigenvalues and eigenfunctions depending on which boundary conditions (\ref{eq11-53}) we assume for the Sturm-Liouville problem (\ref{eq11-71}).

\paragraph{Fourier Sine series}
Assume $\alpha_1=\alpha_2=1, \beta_1=\beta_2=0$. Thus the boundary conditions are 

\begin{equation}\label{eq11-75}
u(0)=u(l)=0.
\end{equation}
Our equation (\ref{eq11-74}) then reduces to $$\sin(\sqrt{\lambda}l)=0,$$ which gives us the following nonzero eigenvalues

\begin{equation}\label{eq11-76}
\lambda_k=\left(\frac{\pi k}{l}\right)^2,\quad k=1,2,3,\dots\;\;,
\end{equation}
and the corresponding normalized eigenfunctions are from (\ref{eq11-73})

\begin{equation}\label{eq11-77}
u_k(x)=\sqrt{\frac{2}{l}}\sin\left(\frac{\pi k}{l}x\right),\quad k=1,2,\dots\;\;.
\end{equation}
For this case $\lambda=0$ is not an eigenvalue because for $\lambda=0$ equation (\ref{eq11-71}) reduces to

\begin{gather*}
-u''(x)=0,\\
\Downarrow\\
u(x)=Ax+B,\\
\\
\left.
\begin{array}{lcr}
u(0)=0 & \Rightarrow & B=0\\
u(l)=0 & \Rightarrow & A=0
\end{array}
\right\} \Rightarrow u(x)=0.
\end{gather*}

\noindent The expansion of a function in a series

\begin{equation*}
u=\sum\limits_{k=1}^{\infty}(u,u_k)\,u_k,
\end{equation*}
is called the \emph{Fourier Sine series}. Using explicit expressions for the eigenfunction and the inner product we get the following formulas for the Fourier Sine series
\begin{equation*}
u(x)=\sqrt{\frac{2}{l}}\sum_{k=1}^\infty\,a_k\,\sin\frac{\pi\,k}{l}x,
\end{equation*}
where
\begin{equation*}
a_k=\sqrt{\frac{2}{l}}\int_0^ldx\,u(x)\,\sin\frac{\pi\,k}{l}x.
\end{equation*}
\paragraph{Fourier Cosine series}
Assume $\alpha_1=\alpha_2=0,$ $\beta_1=\beta_2=1$. Thus the boundary conditions are 

\begin{equation}\label{eq11-77b} 
u'(0)=u'(l)=0.
\end{equation}
Our equation (\ref{eq11-74}) now reduces to 

\begin{equation*}
\lambda\sin{\sqrt{\lambda}l}=0.
\end{equation*}
$\lambda_0=0$ is clearly a solution. For $\lambda_0=0$ our equation reduces to 

\begin{gather*}
-u''(x)=0,\\
\Downarrow\\
u(x)=Ax+B,\\
\\
\left.
\begin{array}{lcr}
u'(0)=0 & \Rightarrow & A=0\\
u'(l)=0 & \Rightarrow & A=0
\end{array}
\right\} \Rightarrow u(x)=B.
\end{gather*}
Thus, $\lambda=0$ is an eigenvalue for this case. The corresponding normalized eigenfunction is 

\begin{equation}\label{eq11-78}
u_0(x)=\frac{1}{\sqrt{l}}.
\end{equation}
For $\lambda>0$ we can divide the eigenvalue equation by $\lambda$ and get 

\begin{gather}\label{eq11-79}
\sin(\sqrt{\lambda}l)=0,\\
\Downarrow\\
\lambda_k=\left(\frac{\pi k}{l}\right)^2,\quad k=1,2,\dots\;\;.
\end{gather}
The corresponding normalized eigenfunctions are 

\begin{equation}\label{eq11-80}
u_k(x)=\sqrt{\frac{2}{l}}\cos\left(\frac{\pi k}{l}x\right), \quad k=1,2,\dots\;\;.
\end{equation}
The series expansion

\begin{equation*}
u=(u,u_0)u_0+\sum\limits_{k=1}^{\infty}(u,u_k)\,u_k,
\end{equation*}
is the \emph{Fourier Cosine} series.
Using explicit expression for the eigenfunctions and the inner product we get
\begin{equation*}
u(x)=\sqrt{\frac{1}{l}}a_0+\sqrt{\frac{2}{l}}\sum_{k=1}^\infty\,a_k\,\cos\frac{\pi\,k}{l}x,
\end{equation*}
where
\begin{align*}
&a_0=\sqrt{\frac{1}{l}}\int_0^ldx\,u(x),\\
&a_k=\sqrt{\frac{2}{l}}\int_0^ldx\,u(x)\,\cos\frac{\pi\,k}{l}x.
\end{align*}
Note that all the general properties of eigenvalues and eigenfunctions for Sturm-Liouville problems can be verified to hold for the eigenvalues and eigenfunctions we have found for the Fourier Sine and Fourier Cosine cases.

\paragraph{Fourier series}
In this case we consider a slightly generalized Sturm-Liouville problem. We assume that the domain is $(-l,l)$ instead of $(0,l)$. We still make the assumptions $p(x)=\rho(x)=1$, $q(x)=0$, so the eigenvalue equation is 

\begin{equation}\label{eq11-81}
-u''(x)=\lambda u(x),\quad -l<x<l.
\end{equation}
Instead of using the boundary conditions (\ref{eq11-53}) we use \emph{periodic} boundary conditions

\begin{align} \label{eq11-82}
& u(-l)=u(l),\\
& u'(-l)=u'(l).
\end{align}
The general solution of (\ref{eq11-81}) is for $\lambda>0$ 

\begin{align}\label{eq11-83}
& u(x) = a\cos\sqrt{\lambda}x+b\sin\sqrt{\lambda}x,\\
& \Downarrow\,\\
  & u'(x)=-a\sqrt{\lambda}\sin\sqrt{\lambda}x+b\sqrt{\lambda}\cos\sqrt{\lambda}x,\nonumber
\end{align}	 
and

\begin{gather}\label{eq11-84}
u(-l)=u(l)\Rightarrow a\cos\sqrt{\lambda}l-b\sin\sqrt{\lambda}l= a\cos\sqrt{\lambda}l+b\sin\sqrt{\lambda}l\nonumber,\\
\Downarrow\nonumber\\
2b\sin\sqrt{\lambda}l=0.
\end{gather}

\begin{gather}\label{eq11-85}
u'(-l)=u'(l)\Rightarrow a\sqrt{\lambda}\sin\sqrt{\lambda}l+b\sqrt{\lambda}\cos\sqrt{\lambda}l=-a\sqrt{\lambda}\sin\sqrt{\lambda}l+b\sqrt{\lambda}\cos\sqrt{\lambda}l,\nonumber\\
\Downarrow\nonumber\\
2a\sin\sqrt{\lambda}l=0.
\end{gather}
The eigenvalues $\lambda>0$ are then clearly 

\begin{equation*}
\lambda_k=\left(\frac{\pi k}{l}\right)^2,\quad k=1,2,\dots\;\;.
\end{equation*}
For each such eigenvalue $a$ and $b$ are arbitrary. Thus (\ref{eq11-83}) gives us an eigenspace of dimension 2 spanned by $\cos\sqrt{\lambda_k}x$ and $\sin\sqrt{\lambda_k}x.$ 
Introducing an inner product and corresponding norm

\begin{gather*}
(u,v)=\int_{-l}^l\mathrm d x\,u(x)v(x),\\
\|u\|^2=\int_{-l}^l\mathrm d x\,u^2(x),
\end{gather*}
we can normalize the eigenfunctions and find 

\begin{align}\label{eq11-86}
& u_k(x)=\frac{1}{\sqrt{l}}\sin \frac{\pi k}{l} x,\quad k=1,2,\dots\\
& \hat{u_k}(x)=\frac{1}{\sqrt{l}}\cos \frac{\pi k}{l} x.\nonumber
\end{align}
One can show that there are no negative eigenvalues and that $\lambda_0=0$ is an eigenvalue with corresponding normalized eigenfunction

\begin{equation}\label{eq11-87}
\hat{u_0}(x)=\frac{1}{\sqrt{2l}}.
\end{equation}
Using these functions we get a \emph{Fourier} series

\begin{equation}\label{eq11-88}
u=(u,\hat{u_0})\hat{u_0}+\sum\limits_{k=1}^{\infty}\left[(u,\hat{u_k})\hat{u_k}+(u,u_k)u_k.\right]
\end{equation}
Using experssions for the eigenfunctions and inner product we get
\begin{equation*}
u(x)=\frac{1}{\sqrt{2l}}a_0+\frac{1}{\sqrt{l}}\sum_{k=1}^\infty\,a_k\cos\frac{\pi\,k}{l}x+b_k\sin\frac{\pi\,k}{l}x,
\end{equation*}
where
\begin{align*}
&a_0=\frac{1}{\sqrt{2l}}\int_0^ldx\,u(x),\\
&a_k=\frac{1}{\sqrt{l}}\int_0^ldx\,u(x)\,\cos\frac{\pi\,k}{l}x,\\
&b_k=\frac{1}{\sqrt{l}}\int_0^ldx\,u(x)\,\sin\frac{\pi\,k}{l}x.
\end{align*}
If $u$ is continuous, has a piecewise continous derivative in the interval $[-l,l]$ and satisfies the boundary conditions (\ref{eq11-82}), then the series (\ref{eq11-88}) converges absolutely and uniformly to $u(x)$.

\subsubsection{Bessel eigenfunctions}
We are now going to consider a singular Sturm-Liouville problem that arise from the eigenvalue problem for Bessel's equation of order $n$
\begin{equation}\label{eq11-89}
L(u(x))=-\frac{\mathrm d}{\mathrm d x}\left(x\frac{\mathrm du}{\mathrm d x}\right)+\frac{n^2}{x}=\lambda xu(x),\quad  0<x<l,
\end{equation}
where $n\geq 0$ is an integer. This is a Sturm-Liouville problem with $p(x)=\rho(x)=x$, $q(x)=\frac{n^2}{x}$. The boundary conditions of interest for (\ref{eq11-89}) are

\begin{align}
& u(x)\text{ bounded at x=0},\nonumber\\
& u(l)=0.\label{eq11-90}
\end{align}
Requiring that $u(x)$ is bounded at $x=0$ is a significant condition here because (\ref{eq11-89}) is a singular Sturm-Liouville problem. We observe that $p,\,\rho$ and $q$ are positive on $0<x<l$ but that $p$ and $\rho$ vanish at $x=0$ and that $q$ is singular at $x=0$. We might therefore in general expect there to be solutions that are singular at $x=0$ and (\ref{eq11-90}) eliminate such solutions as acceptable eigenfunctions. One might worry that \emph{all} solutions are singular at $x=0$. Then we will get no eigenvalues and eigenfunctions at all! By using the \emph{method of Frobenious} from the theory of ordinary differential equations one can show that (\ref{eq11-89}) has solutions that are bounded at $x=0$. These solutions can be expressed in terms of \emph{Bessel functions} of order $n$, $J_n$, in the following simple way

\begin{equation}\label{eq11-91}
u(x)=J_n(\sqrt{\lambda}x).
\end{equation}
We find the eigenvalues by solving the equation

\begin{equation}\label{eq11-92}
J_n(\sqrt{\lambda}l)=0.
\end{equation}
One can prove that this equation has infinitely many solutions $\lambda_{kn},\,k=1,2,3,\dots$ They are all positive and $\lambda_{kn}\rightarrow\infty$ as $k\rightarrow\infty$. A large number of zeros for the function $J_n$ has been found using numerical methods. Let these zeros be $\alpha_{kn},\,k=1,2,3,\dots$ In terms of these zeros the eigenvalues are 

\begin{equation}\label{eq11-93}
\lambda_{kn}=\left(\frac{\alpha_{kn}}{l}\right)^2,\quad k=1,2,\dots\;\;.
\end{equation}
The corresponding eigenfunctions are 

\begin{equation}\label{eq11-94}
\hat{u}_{kn}(x)=J_n(\sqrt{\lambda_{kn}}x),\quad k=1,2,\dots\;\;.
\end{equation}
Using the inner product

\begin{equation*}
(u,v)=\int_0^l\mathrm d x\,xu(x)v(x),
\end{equation*}
the eigenfunctions can be normalized 

\begin{equation}\label{eq11-95}
u_{kn}=\frac{\sqrt{2}}{l}\frac{J_n(\sqrt{\lambda_{kn}}x)}{|J_{n+1}(\sqrt{\lambda_{kn}}l)|},\quad k=1,2,\dots,\;\;,
\end{equation}
where we have used the formulas

\begin{equation}\label{eq11-96}
\|J_n(\sqrt{\lambda_{kn}}x)\|^2=\int_0^l x J_n^2(\sqrt{\lambda_{kn}}x)\,\mathrm d x=\frac{l^2}{2}J_{n+1}^2(\sqrt{\lambda_{kn}}l).
\end{equation}
As you have seen in calculus there are a large number of formulas linking elementary functions like $\sin(x)$, $\cos(x)$, $\ln(x)$, $e^x$ etc. and their derivatives and anti-derivatives. The Bessel functions $J_n$ is an extension of this machinery. There are a large set of formulas linking the $J_n$ and their derivatives and anti-derivatives to each other and other functions like $\sin(x)$, $\cos(x)$, $\ln(x)$ etc. The functions in calculus ($\sin(x)$, $\cos(x),\dots$) are called \emph{elementary functions}. The function $J_n$ is our first example of a much larger machinery of functions called \emph{special functions}. There are many handbooks that list some of the important relations and formulas involving special functions. Symbolic systems like mathematica, maple, etc. embody most of these formulas and often produce answers to problems formulated in terms of such functions.

\noindent Note that the boundary conditions (\ref{eq11-90}) are not of the type (\ref{eq11-53}) discussed in the general theory of Sturm-Liouville theory. Thus, there is no guarantee at this point that eigenfunctions corresponding to different eigenvalues are orthogonal. This is typical for singular Sturm-Liouville problems; the properties of the eigenvalues and eigenfunctions must be investigated in each separate case. 

Observe that, in general, with $L$ defined as in (\ref{eq11-52}) we have

\begin{align*}
\int_0^l\mathrm d x\,u Lv & = \int_0^l\mathrm d x\,u[-(pv')'+qv]\\
&=-\int_0^l\mathrm d x\,u(pv')'+\int_0^l\mathrm d x\,u qv\\
&= -u(pv')\big|^l_0+\int_0^l\mathrm d x\,u'pv'+\int_0^l\mathrm d x\,quv\\
&= -u(pv')\big|^l_0+(u'p)v\big|^l_0-\int_0^l\mathrm d x\,(pu')'v+\int_0^l\mathrm d x\,u qv\\
&=\lbrace pu'v-puv'\rbrace\big|^l_0+\int_0^l\mathrm d x\,Luv
\end{align*}

\begin{equation}\label{eq11-97}
\Rightarrow\int_0^l\mathrm d x\,\lbrace u Lv-v Lu\rbrace=\lbrace pu'v-pv'u\rbrace\big|_0^l
\end{equation}

\noindent Recall that for our current singular Sturm-Liouville problem, the space of functions we need to expand in terms of the eigenfunctions (\ref{eq11-95}), consists of functions such that

\begin{align}\label{eq11-98}
& u(l)=0,\\
& \lim\limits_{x\to 0^+} u(x)\text{ bounded.}\nonumber
\end{align}
Thus on the right hand side of (\ref{eq11-97}) the evaluation at $x=0$ must be interpreted as a limit as $x$ approach $0$ from inside the interval $(0,l)$. 

Let us apply (\ref{eq11-97}) to a pair of eigenfunctions $u_i, u_j$ corresponding to different eigenvalues $\lambda_i\neq\lambda_j$. Thus 

\begin{align*}
& Lu_i=\rho\lambda_iu_i,\\
& Lu_j=\rho\lambda_ju_j.
\end{align*}
But then  (\ref{eq11-97}) gives us  

\begin{gather}\label{eq11-99}
(\lambda_i-\lambda_j)(u_i,u_j)=(\lambda_iu_i,u_j)-(u_i,\lambda_j u_j)\\
=(\frac{1}{\rho}Lu_i,u_j)-(u_i,\frac{1}{\rho}Lu_j)\nonumber\\
=\int^l_0\mathrm d x \lbrace Lu_iu_j-u_iLu_j\rbrace=[pu_ju_i'-pu_iu_j']\big|^l_0.\nonumber
\end{gather}
In order for eigenvectors corresponding to different eigenvalues to be orthogonal, the right hand side (\ref{eq11-99}) must be zero. The evaluation at $x=l$ is zero since the eigenfunctions are zero at this point. In order for the right hand side of (\ref{eq11-99}) to vanish, it is in general \emph{not} enough to require that the eigenfunctions are bounded at $x=0$. The derivative of the eigenfunctions might still be unbounded at $x=0$ in such a way that the right hand side of (\ref{eq11-99}) is nonzero and thus ruining the orthogonality of the eigenfunctions. However, if we also require that

\begin{equation}\label{eq11-100}
p(x)u_i'(x)\to 0\text{ when }x\to 0^+,
\end{equation}
everything is okay and the eigenfunctions form an orthogonal set. Since $p(x)=x$ we thus have the requirement

\begin{equation*}
\lim\limits_{x\to 0^+}xu_i'(x)=0.
\end{equation*}
From the theory of Bessel functions we have 

\begin{equation*}
\left.
\begin{array}{lccr}
 J_0(\sqrt{\lambda}x)&=& 1+\mathcal{O}(x^2),&\\
 J_n(\sqrt{\lambda}x)&=&\mathcal{O}(x^n),&\,n>0
\end{array}
\right\}\text{ as }x\to 0^+.
\end{equation*}
Thus the functions are bounded at $x=0$. Also

\begin{equation*}
\left.
\begin{array}{lccr}
	J_0'(\sqrt{\lambda}x)&=&\mathcal{O}(x),&\\
	J_n'(\sqrt{\lambda}x)&=&\mathcal{O}(x^{n-1}),&\,n>0
\end{array}
\right\}\text{ as }x\to 0^+,
\end{equation*}
but then 

\begin{align*}
&\lim\limits_{x\to 0^+} xJ_0'(\sqrt{\lambda}x)=0,\\
&\lim\limits_{x\to 0^+} xJ_n'(\sqrt{\lambda}x)=0,\quad n>0.
\end{align*}
Thus, by the argument given, the eigenfunctions 

\begin{equation*}
u_{kn}(x)=\alpha_{kn}J_n(\sqrt{\lambda_{kn}}x),\,\,\,\text{with}\,\,\alpha_{kn}=\frac{\sqrt{2}}{l|J_{n+1}(\sqrt{\lambda_{kn}}l)|},
\end{equation*}
form an orthonormal set. 

One can show that these eigenfunctions form a complete set. Any smooth function defined on the interval $0\leq x\leq l$ that satisfies the boundary condition $u(l)=0$ can be expanded in a convergent series 

\begin{equation}\label{eq11-101}
u=\sum\limits_{k=1}^{\infty}(u,u_{kn})\,u_{kn}.
\end{equation}

\subsubsection{Legendre polynomial eigenfunctions}
Another singular Sturm-Liouville problem that occur frequently is

\begin{equation}\label{eq11-102}
L(u(x))=-\frac{\mathrm d}{\mathrm d x}\left[(1-x^2)\frac{\mathrm du}{\mathrm d x}\right]=\lambda u(x),\quad  -1<x<l,
\end{equation}
with boundary conditions 

\begin{equation}\label{eq11-103}
u(x)\text{ bounded at }x=\pm 1.
\end{equation}
The function $p(x)=1-x^2$ vanishes at both endpoints $x=\pm 1$ so we indeed have a singular Sturm-Liouville problem. 

One can use the method of Frobenius, alluded to earlier, to show that the only bounded solutions to (\ref{eq11-102}) are in fact polynomials and that these exists only for 

\begin{equation}\label{eq11-104}
\lambda_k=k(k+1),\quad k=0,1,2,\dots\;\;.
\end{equation}
The normalized eigenfunctions, $P_k$, are called the \emph{Legendre polynomials} of degree $k$

\begin{equation}\label{eq11-105}
P_k(x)=\left(\frac{1}{2^k k!}\right)\frac{\mathrm d^k}{\mathrm d x^k}(x^2-1)^k.
\end{equation}
It is very easy to verify that 

\begin{equation}\label{eq11-106}
\lim\limits_{x\to\pm 1}(1-x^2)P'_k(x)=0,
\end{equation}
so this singular Sturm-Liouville problem has an orthogonal set of eigenfunctions according to previous arguments. One can show that 

\begin{equation}\label{eq11-107}
\|P_k(x)\|^2=\int_{-1}^1 P_k^2(x)=\frac{2}{2k+1},\quad k=0,1,2,\dots\;\;.
\end{equation}
Thus we get an orthonormal set of eigenfunctions 

\begin{equation}\label{eq11-108}
u_k(x)=\sqrt{\frac{2k+1}{2}}P_k(x),\quad k=0,1,2,\dots\;\;.
\end{equation}
One can show that this set is complete and that any smooth function defined over the interval $-1\leq x\leq 1$ can be expanded in a convergent series

\begin{equation}\label{eq11-109}
u=\sum\limits_{k=0}^{\infty}(u,u_k)\,u_k.
\end{equation}

\subsection{Series solutions to initial boundary value problems}
It is now time for returning to the solution of initial boundary value problems. We will consider a few specific examples and do detailed calculations for these cases. However, before we do that let us recapitulate the separation of variable strategy for the hyperbolic case. As we have seen, the strategy for parabolic and elliptic equations is very similar.

\begin{enumerate}
	\item Look for functions of the form
	
	\begin{equation}\label{eq11-110}
	u(\mathbf x, t)=M(\mathbf x)N(t)
	\end{equation}
	that solves the PDE (\ref{eq11-13}) \emph{and} the boundary conditions. These could be of the form (\ref{eq11-14}), but as we have seen in the previous section on singular Sturm-Liouville problems, other possibilities exist.
	
	\item If the boundary condition has been choosen in an appropriate way we find that $M$ must be one of a complete set of eigenfunctions $\lbrace M_k\rbrace_{k=1}^{\infty}$ of a self-adjoint eigenvalue problem with eigenvalues $\lbrace \lambda_k\rbrace_{k=1} ^{\infty}$ corresponding to the eigenfunctions $\lbrace M_k\rbrace_{k=1} ^{\infty}$.
	
	\item The function $N(t)$ is found to be one of 
	
	\begin{equation}\label{eq11-111}
	N_k(t)=a_k\cos(\sqrt{\lambda}t)+b_k\sin(\sqrt{\lambda}t),
	\end{equation}
	where $\lbrace a_k\rbrace,$ $\lbrace b_k\rbrace$ are constants that can be freely specified.
	
	\item A formal solution to the PDE and the boundary conditions is
	
	\begin{equation}\label{eq11-112}
	u(\mathbf x, t)= \sum\limits_{k=1}^{\infty}\lbrace a_k\cos(\sqrt{\lambda}t)+b_k\sin(\sqrt{\lambda}t)\rbrace\,M_k(\mathbf x).
	\end{equation}
	
	\item Let $a_k$ and $b_k$ be the Fourier coefficients 
	
	\begin{equation}\label{eq11-113}
	a_k=\frac{(f,M_K)}{(M_k,M_k)},\quad b_k=\frac{(g,M_K)}{\sqrt{\lambda_k}(M_k,M_k)},
	\end{equation}
	where $f$ and $g$ are the initial data. Then
	
	\begin{align}\label{eq11-114}
	& u(\mathbf x,0)=\sum\limits_{k=1}^{\infty}a_kM_k(\mathbf x)=f(\mathbf x),\\
	& u_t(\mathbf x,0)=\sum\limits_{k=1}^{\infty}\sqrt{\lambda_k}b_kM_k(\mathbf x)=g(\mathbf x)\nonumber
	\end{align}
	and thus the formal solution (\ref{eq11-112}) satisfies the PDE, the boundary conditions \emph{and} the initial conditions. 
	
	\item Investigate the convergence of the series (\ref{eq11-112}). If the initial data $f$ and $g$ are smooth enough the series will converges so fast that it can be differentiated twice termwise. In this case the formal series solution is a classical (honest to God!) solution to our initial/boundary value problem. For less smooth $f$ and $g$ the convergence is weaker and will lead to generalized (non-classical) solutions.
\end{enumerate} 

\subsubsection{Vibrations of a string with fixed endpoints}
We have seen earlier in this class that a vibrating string with fixed endpoints can in certain cases be modelled using the equation

\begin{equation}\label{eq11-115}
u_{tt}(x,t)-c^2u_{xx}(x,t)=0,\quad 0<x<l,\,t>0,
\end{equation}

\begin{equation}\label{eq11-116}
u(0,t)=0,\, u(l,t)=0,
\end{equation}
where $c=\sqrt{\frac{T}{\rho}}$, and where $T$ and $\rho$ are the, assumed constant, tension and mass density of the string. In order to get a unique motion of the string, we must specify the displacement and movement of the string at $t=0$
	
\begin{align}\label{eq11-117}
& u(x,0)=f(x),\quad 0<x<l,\\
& u_t(x,0)=g(x). \nonumber
\end{align}
These are the initial values. We now separate variables, $u(x,t)=M(x)N(t)$, and get the following regular Sturm-Liouville problem for $M$
	
\begin{align}\label{eq11-118}
& -M''(x)=\lambda M(x),\\
& M(0)=M(l)=0. \nonumber
\end{align}
This is the problem we investigated on page 122. 
We found eigenvalues and normalized eigenfunctions of the form

\begin{equation}\label{eq11-119}
\left.
\begin{array}{lccl}
&\lambda_k &=&\left(\frac{\pi k}{l}\right)^2\\
& M_k(x) &= &\sqrt{\frac{2}{l}}\sin\left(\frac{\pi k}{l}x\right)
\end{array}
\right\}
\quad k=1,2,\dots
\end{equation}
For the function $N(t)$ we find 

\begin{equation}\label{eq11-120}
N_k(t)=a_k\cos\left(\frac{\pi k c}{l}t\right)+b_k\sin\left(\frac{\pi k c}{l}t\right),\quad k=1,2,\dots
\end{equation}
and the corresponding separated solutions $u_k(x,t)$ are 

\begin{equation}\label{eq11-121}
u_k(x,t)=\left[a_k\cos\left(\frac{\pi k c}{l}t\right)+b_k\sin\left(\frac{\pi k c}{l}t\right)\right]\sqrt{\frac{2}{l}}\sin\left(\frac{\pi k}{l}x\right),\quad\text{for } k=1,2,\dots
\end{equation} 
In order to satisfy the initial conditions we write down the formal solution 

\begin{equation}\label{eq11-122}
u(x,t)=\sum\limits_{k=1}^{\infty}u_k(x,t).
\end{equation}
Upon applying the initial conditions (\ref{eq11-117}) we get 

\begin{align*}
& u(x,0)=\sum\limits_{k=1}^{\infty}u_k(x,0)=\sum\limits_{k=1}^{\infty}a_k\sqrt{\frac{2}{l}}\sin\left(\frac{\pi k}{l}x\right)=f(x),\\
& u_t(x,0)=\sum\limits_{k=1}^{\infty}\partial_t u_k(x,0)=\sum\limits_{k=1}^{\infty}\left(\frac{\pi k c}{l}\right)b_k\sqrt{\frac{2}{l}}\sin\left(\frac{\pi k}{l}x\right)=g(x).
\end{align*}
These are formally satisfied if 

\begin{align}\label{eq11-123}
& a_k=(f,M_k)=\sqrt{\frac{2}{l}}\int_0^l\mathrm d x \,f(x)\sin\left(\frac{\pi k}{l}x\right),\\
& b_k = \frac{l}{\pi k c}(g,M_k)=\frac{\sqrt{2l}}{\pi k c}\int_0^l\mathrm d x\,g(x)\sin\left(\frac{\pi k}{l}x\right). \nonumber
\end{align}
With these values for $a_k$ and $b_k$, (\ref{eq11-122}), is the formal solution of the initial value problem for the vibrating string fixed at both ends.

The separated solutions (\ref{eq11-121}) are called \emph{normal modes} for the vibrating string. The formal solution (\ref{eq11-122}) is a sum, or \emph{linear superposition}, of normal modes. The notion of normal modes is a concept that goes way beyond their use for the vibrating string. Separated solutions satisfying a hyperbolic PDE and boundary conditions are in general called normal modes. 

Normal modes is a physical concept that has very wide applicability and deep roots. Physicists in general are not very familiar with eigenvalue problems for linear differential operators and their associated eigenfunctions, but they all know what normal modes are! 

In quantum theory applied to fields, normal modes is a key idea. In this setting, they are in fact the mathematical embodiment of elementary particles. The newly discovered Higgs Boson is a normal mode for a certain field equation closely related to the 3D wave equation. In the quantum theory of materials, so-called many-body theory, the interaction between material normal modes (phonons) and electronic normal modes is responsible for the phenomenon of superconductivity.

Now then, when can we expect the formal solution (\ref{eq11-122}) to be an actual, honest to God, classical solution? In general, finding sharp conditions on the initial data $f(\mathbf x)$ and $g(\mathbf x)$ that will imply the existence of a classical solution is a hard mathematical problem that is best left to specialists in mathematical analysis. In fact, sometimes the problem is so hard that entirely new mathematical ideas have to be developed in order to solve this problem. 

You might know that all of modern mathematics is formulated in terms of, and gets its consistency from, \emph{set theory}. You might not know that set theory was invented by George Cantor in the 1880's while studying the convergence of trigonometric Fourier series! 

So what does mathematical analysis tell us about the convergence of the formal solution (\ref{eq11-122})? The following conditions are \emph{sufficient} to ensure the existence of a classical solution for the vibrating string:

\begin{enumerate}
	\item $f$ must have two continuous derivatives and a third piecewise continuous derivative,
	\item $f(0)=f(l)=0,\quad f''(0)=f''(l)=0$,
	\item $g$ must have one continuous derivative and a second piecewise continuous derivative,
	\item $g(0)=g(l)=0.$
\end{enumerate}


By using standard trigonometric formulas, the normal modes (\ref{eq11-121})  can be rewritten as

\begin{equation}\label{eq11-126}
u_k(x,t)=\alpha_k\cos\left(\frac{\pi k c}{l}t+\delta_ k\right)\sin\left(\frac{\pi k}{l}x\right),\quad k=1,2,\dots,
\end{equation} 
where 

\begin{equation}\label{eq11-127}
\left.
\begin{array}{lr}
&\alpha_k=\sqrt{\frac{2}{l}}\sqrt{a_k^2+b_k^2}\\
&\delta_k=-tg^{-1}\left(\frac{b_k}{a_k}\right)
\end{array}
\right\}
\quad k=1,2,\dots
\end{equation}
The formula (\ref{eq11-126}) shows that for a normal mode each point $x_0$ perform a harmonic vibration with frequency 

\begin{equation}\label{eq11-128}
w_k=\frac{\pi k c}{l}
\end{equation}
and amplitude

\begin{equation}\label{eq11-129}
A_k=\alpha_k\sin\left(\frac{\pi k}{l}x_0\right)
\end{equation}

\noindent Solutions to the wave equation are in general called {\it waves}. The normal mode (\ref{eq11-121}) then is a \emph{standing wave}; each point on the string oscillate in place as $t$ varies. 

Some points on the string remain fixed during the oscillation. These are the points

\begin{equation*}
x_m=\frac{ml}{k},\quad m=1,2,\dots,k-1
\end{equation*}
which are called \emph{nodes} for the standing wave. The standing wave has maximum amplitude at the points 

\begin{equation*}
x_n=\frac{(2n+1)l}{2k},\quad n=0,1,2,\dots,k-1,
\end{equation*}
which are called \emph{antinodes}.

Since each normal mode performs a harmonic oscillation, $u_k(x,t)$ is often called the $k^{th}$ harmonic. Since a vibrating string makes sound, the normal mode $u_k$ will produce a pure tone at frequency $\omega_k$. Thus the decomposition of the string into normal modes corresponds to the decomposition of sound into pure tones. If the tone-spectrum is pleasing we call it music! 

Recall that the energy of a small piece of a vibrating string located between $x$ and $x+\mathrm d x$ is 

\begin{equation*}
E(x,t)=\lbrace \underbrace{\rho(\partial_t u(x,t))^2}_{\text{kinetic energy}}+\underbrace{T(\partial_xu(x,t)^2}_{\text{potential energy}}\rbrace\,\mathrm d x.
\end{equation*}
Thus, the total energy of the vibrating string is 

\begin{equation*}
E(t)=\frac{1}{2}\int_0^l\mathrm d x\,\lbrace \rho(\partial_t u(x,t))^2+T(\partial_xu(x,t))^2\rbrace.
\end{equation*}
Writing 
\begin{equation*}
u_k(x,t)=N_k(t)M_k(x),
\end{equation*}
where
\begin{eqnarray*}
&N_k(t)=\alpha_k\cos\left(\omega_kt+\delta_ k\right),\\
&M_k(x)=\sin\left(\frac{\pi k}{l}x\right),
\end{eqnarray*}
we observe that
\begin{align*}
 \int_0^l\mathrm d x\,\lbrace (\partial_t u(x,t))^2 & = (\partial_t u,\partial_t u) \\
&=\left(\sum\limits_{k} N'_k(t)M_k,\,\sum\limits_{k'} N'_{k'}(t)M_{k'}\right)\\
&=\sum\limits_{kk'} N_k'(t)N'_{k'}(t)(M_k,M_{k'})\\
&=\sum\limits_{k}(N_k'(t))^2(M_k,M_{k})=\sum\limits_{k}(N_k'(t)M_k,N_k'(t)M_{k})\\
&=\sum\limits_{k}\int_0^l\mathrm d x\,(\partial_t u_k(x,t))^2 .
\end{align*}
In a similar way we find that 

\begin{equation*}
\int_0^l\mathrm d x\,(\partial_xu(x,t))^2=\sum\limits_{k}\int_0^l\mathrm d x\,(\partial_x u_k(x,t))^2.
\end{equation*}
Therefore 

\begin{equation*}
E(t)=\sum\limits_kE_k(t),
\end{equation*}
where 

\begin{align*}
E_k(t)& =\frac{1}{2}\int_0^l\mathrm d x\,\lbrace \rho(\partial_t u_k(x,t))^2+T(\partial_xu_k(x,t)^2\rbrace\\
& =\frac{\omega_k^2m(a_k^2+b_k^2)}{2l},
\end{align*}
where $ m=\rho l =\text{ mass of string}$.(Do this calculation!)

Thus each $E_k$ is independent of $t$ and so is the total energy 

\begin{equation*}
E(t)=\frac{m}{2l}\sum\limits_{k=1}^{\infty}\omega_k^2(a_k^2+b_k^2),
\end{equation*}
which is the sum of the energies, $E_k$,  of each normal mode. 

The decomposition of the total energy of a system into the energy of each normal mode plays a key role in the quantum theory of fields and in many-body theory. 

Using trigonometric addition formulas one can show that 
\begin{align}\label{eq11-130}
 2\cos \frac{A-B}{2}\sin\frac{A-B}{2}&=\sin A+\sin B\\
 -2\sin \frac{A-B}{2}\sin\frac{A+B}{2}&=\cos A-\cos. B\nonumber
\end{align}
Using these formulas on the expression (\ref{eq11-121}) for the normal modes we find

\begin{align}\label{eq11-131}
u_k(x,t)&=\frac{1}{\sqrt{2l}}(a_k\sin\left(\frac{\pi k}{l}(x+ct)\right)+a_k\sin\left(\frac{\pi k}{l}(x-ct)\right)\\
& + \frac{1}{\sqrt{2l}}(b_k\cos\left(\frac{\pi k}{l}(x-ct)\right)-b_k\cos\left(\frac{\pi k}{l}(x-ct)\right).
\end{align}
Thus the normal mode can be written as a sum of left and right traveling waves. This is not surprising since we know that \emph{all} solutions of the wave equation can be written as a sum of left and right traveling waves. What is special for the normal modes is that their traveling wave components are arranged in such a way that we get a standing wave. 

\subsubsection{Heat conduction in a finite rod}
As we have seen previously in this class, heat conduction in a finite rod, where there are no heat sources and where the endpoints of the rod are kept at zero temperature, is in certain cases modeled by the boundary value problem 

\begin{align}\label{eq11-132}
& u_t(x,t)-c^2u_{xx}(x,t)=0,\quad 0<x<l,\,t>0,\\
& u(0,t)=u(l,t)=0,
\end{align}
where $u(x,t)$ is the temperature of the rod as a function of time, and where $c^2$ is the coefficient of heat conductivity. In order to get a unique solution we must specify the initial temperature distribution

\begin{equation}\label{eq11-133}
u(x,0)=f(x),\quad 0<x<l.
\end{equation}
The separation of variable approach leads to the same eigenvalue problem as for the string. The corresponding functions $N_k(t)$ now become 

\begin{equation}\label{eq11-134}
N_k(t)=a_ke^{-(\frac{\pi k c}{l})^2t},\quad k=1,2,\dots
\end{equation}
and we get a formal solution 

\begin{equation}\label{eq11-135}
u(x,t)=\sum\limits_{k=1}^{\infty}N_k(t)M_k(x)=\sqrt{\frac{2}{l}}\sum\limits_{k=1}^{\infty} a_ke^{-(\frac{\pi k c}{l})^2t}\sin\left(\frac{\pi k}{l}x\right).
\end{equation}
The initial value (\ref{eq11-132}) implies that the constants $a_k$ must be 

\begin{equation}\label{eq11-136}
a_k=\sqrt{\frac{2}{l}}\int_0^l\mathrm d x\,f(x)\sin\left(\frac{\pi k}{l}x\right).
\end{equation}
(\ref{eq11-135}) with the constants given by (\ref{eq11-136}) is a formal solution to our initial boundary value problem. 

Because of the exponential term in the solution for $t>0$, the series (\ref{eq11-135}) converges so strongly that the series can be differentiated termwise as many times as we like. Thus (\ref{eq11-135}) satisfies the PDE and boundary conditions if we merely assume that $f(x)$ is a bounded function on $0\leq x\leq l$. However, as $t\to 0$ the exponential term vanishes, so in order to ensure that (\ref{eq11-135}) also satisfies the initial condition, we need stronger conditions on $f$. Our friendly mathematical analyst tells us that the following conditions are sufficient: 
\begin{enumerate}
\item  f must be continuous and f' piecewise continuous in $0 \leq x\leq l$
\item  f(0)=f(l)=0
\end{enumerate}

For the wave equation we know that discontinuities in the initial data are preserved by the equation and propagated along the characteristic curves  $x\pm ct=\text{const}$. For the heat conduction problem, the exponential factor will ensure that the solution is infinitely many times differentiable for any $t>0$ even if the initial data is discontinuous!

The heat conduction equation can not support any kind of discontinuities for $t>0$. We have seen this same result appear previously from the fact that the only characteristic curves for the heat conduction equation are of the form $t=\text{const}$. 

\subsubsection{The Laplace equation in a rectangle}
Consider the boundary value problem 

\begin{align}\label{eq11-138}
& u_{xx}(x,y)+u_{yy}(x,y)=0,\quad 0<x<l,\,0<y<\hat{l},\\
& u(x,0)=f(x),\,u(0,y)=0,\nonumber\\
& u(x,\hat{l})=g(x),\,u(l,y)=0.\nonumber
\end{align}
The equation, which is the Laplace equation, and boundary conditions determine the steady state displacement of a stretched membrane on a rectangular domain. The functions $f$ and $g$ specify the stretching of the membrane on the boundary curves $y=0$ and $y=\hat{l}$. 

The separation of variable method leads to the same eigenvalue problem as before but now the corresponding functions $N_k$ are 

\begin{equation}\label{eq11-139}
N_k(y)=\hat{a_k}e^{\frac{\pi k}{l}y}+\hat{b_k}e^{-\frac{\pi k}{l}y}.
\end{equation}
The separated solutions then are 

\begin{equation}\label{eq11-140}
u_k(x,y)=[a_k\sinh\left(\frac{\pi k}{l}y\right)+b_k\sinh\left(\frac{\pi k}{l}(y-\hat{l})\right)]\sqrt{\frac{2}{l}}\sin\left(\frac{\pi k}{l}x\right),
\end{equation}
where we have introduced hyperbolic functions instead of exponential ones. 

The formal solution to our equation and homogeneous boundary conditions are then 

\begin{equation}\label{eq11-141}
u(x,y)=\sum\limits_{k=1}^{\infty}u_k(x,y).
\end{equation}
The boundary conditions at $y=0$ and $y=\hat{l}$ imply that

\begin{align}\label{eq11-142}
& u(x,0)=\sum\limits_{k=1}^{\infty}b_k\sinh\left(-\frac{\pi k}{l}\hat{l}\right)\sqrt{\frac{2}{l}}\sin\left(\frac{\pi k}{l}x\right)=f(x),\\
& u(x,\hat{l})=\sum\limits_{k=1}^{\infty}a_k\sinh\left(\frac{\pi k}{l}\hat{l}\right)\sqrt{\frac{2}{l}}\sin\left(\frac{\pi k}{l}x\right)=g(x).	
\end{align}
These are formally satisfied if

\begin{equation}\label{eq11-143}
\left.
\begin{aligned}
a_k&=\frac{(g,M_k)}{\sinh\left(\frac{\pi k\hat{l}}{l}\right)}\\
b_k&=\frac{(f,M_k)}{\sinh\left(-\frac{\pi k\hat{l}}{l}\right)}
\end{aligned}
\right\}\quad k=1,2,\dots
\end{equation}

(\ref{eq11-141}) with $a_k$ and $b_k$ given by (\ref{eq11-143}) determines the formal solution to the problem (\ref{eq11-138}). 

Assume that

\begin{equation}\label{eq11-144}
\left.
\begin{aligned}
& \int_0^l\mathrm d x\, |f(x)|\leq m\\
& \int_0^l\mathrm d x\, |g(x)|\leq m
\end{aligned}
\right\}\quad\text{for some }m>0.
\end{equation} 
Then 

\begin{align}\label{eq11-145}
|b_k|&=\frac{1}{|\sinh\left(-\frac{\pi k\hat{l}}{l}\right)|}\Big|\int_0^l\mathrm d x\, f(x)\sqrt{\frac{2}{l}}\sin\left(\frac{\pi k}{l}x\right)\Big|\\
& \leq \frac{m\sqrt{\frac{2}{l}}}{|\sinh\left(-\frac{\pi k\hat{l}}{l}\right)|}=
\frac{m\sqrt{\frac{2}{l}}}{\frac{1}{2}e^{\frac{\pi k\hat{l}}{l}}(1-e^{-\frac{2\pi k\hat{l}}{l}})} \nonumber
\end{align}
and using the fact that for $z\geq 0$

\begin{equation}\label{eq11-146}
\sinh z =\frac{1}{2}(e^z-e^{-z})\leq\frac{1}{2}e^z
\end{equation}
we get the estimate

\begin{eqnarray*}
|b_k\sinh\left(\frac{\pi k}{l}(y-\hat{l})\right)|&\leq\frac{m\sqrt{\frac{2}{l}}e^{\frac{\pi k}{l}(\hat{l}-y)}}{e^{\frac{\pi k \hat{l}}{l}(1-e^{-\frac{2\pi k\hat{l}}{l}})}}\\
&=\frac{m\sqrt{\frac{2}{l}}e^{-\frac{\pi k}{l}y}}{(1-e^{-\frac{2\pi k\hat{l}}{l}})}		
\end{eqnarray*}

Thus, this term is bounded by an exponential as $k\to\infty,\medspace\forall\, 0<y<\hat{l}.$ In a similar way $|a_k\sinh\left(\frac{\pi k}{l} y\right)|$ is bounded by a decaying exponential as $k\to\infty$. Therefore the terms in the formal solution (\ref{eq11-141}) decay exponentially as $k\to\infty$ for all $0<x<l, 0<y<\hat{l}$. The formal series (\ref{eq11-141}) therefore converges and defines a function in the inside of the rectangle $[0,l]\times[0,\hat{l}]$ whose series representation can be differentiated termwise as many times as we like. Thus the formal series (\ref{eq11-141}) defines a classical solution to the Laplace equation inside a closed rectangle  $[0,l]\times[0,\hat{l}]$. It, by construction of the separated solutions, satisfies the boundary conditions on the edges $x=0,\,x=l,\,0\leq y\leq\hat{l}$. As we approach $y=0$ or $y=\hat{l}$, the exponential decay is weakened and in order to ensure that the formal solution satisfies the boundary conditions on $y=0,\,\hat{l},\,0\leq x\leq l$ our friendly mathematical analysis tells us that it is sufficient that:
\begin{enumerate}
\item $f$, $g$ are continuous and $f'$, $g'$
        piecewise continuous in the interval $0 \leq x \leq l$
\item $f(0)=f(l)=g(0)=g(l)=0$
\end{enumerate}

The fact that solutions to the Laplace equations are smooth in the interior of the rectangle even if the boundary data has discontinuities is, as we have seen, also a consequence of having no real characteristics at all.


\subsection{Inhomogenous equations - Duhamel's principle}
So far, we have only considered the separation of variable method for homogeneous equations. This means physically that we are modeling source-free problems. It is now time to consider the presence of sources. We will here focus on the hyperbolic and parabolic case:

\begin{align}\label{eq11-149}
& \rho(\mathbf x)u_{tt}(\mathbf x,t)+L[u(\mathbf x,t)]=g(\mathbf x, t),\quad\text{Hyperbolic case},\\
& \rho(\mathbf x)u_{t}(\mathbf x,t)+L[u(\mathbf x,t)]=g(\mathbf x, t),\quad\text{Parabolic case}.\nonumber
\end{align}
The equations (\ref{eq11-149}) are supplied with homogeneous boundary conditions and also homogeneous initial conditions

\begin{align}\label{eq11-150}
u(\mathbf x,0)&=0,\quad u_t(\mathbf x,0)=0,&\quad\text{Hyperbolic case},\\
u(\mathbf x,0)&=0,&\quad \text{Parabolic case.}\nonumber
\end{align}
The initial/boundary value problem (\ref{eq11-149}),(\ref{eq11-150}) is solved using \emph{Duhamel's principle}. We proceed as follows:

Let $v(x,t;\tau)$ be a function that satisfies
\begin{align*}\
& \rho(\mathbf x)v_{tt}(\mathbf x,t)+L[v(\mathbf x,t)]=0,\quad\text{Hyperbolic case},\\
& \rho(\mathbf x)v_{t}(\mathbf x,t)+L[v(\mathbf x,t)]=0,\quad\text{Parabolic case}.
\end{align*}
with initial conditions 

\begin{align}\label{eq11-151}
v(\mathbf x,\tau;\tau)=0,\,v_t(\mathbf x,\tau;\tau)&=\frac{g(\mathbf x,\tau)}{\rho(\mathbf x)},\quad\text{Hyperbolic case},\\
v(\mathbf x,\tau;\tau)=&\frac{g(\mathbf x,\tau)}{\rho(\mathbf x)},\quad\text{Parabolic case}.\nonumber
\end{align}
We assume this solution has been found using separation of variables or some other method. 

Define a function $u(\mathbf x,t)$ by 

\begin{equation}\label{eq11-152}
u(\mathbf x,t)=\int^t_0\mathrm d\tau\,v(\mathbf x,t;\tau).
\end{equation}
Then we have by the rules of calculus

\begin{align}\label{eq11-153}
& u_t(\mathbf x,t)=v(\mathbf x,t;t)+\int_0^t\mathrm d t\,v_t(\mathbf x,t;\tau),\nonumber\\
& v_{tt}=\partial_t(v(\mathbf x,t;t))+v_t(\mathbf x,t;t)+\int^t_0\mathrm d\tau v_{tt}(x,t;\tau),\\
& L[u(\mathrm x,t)]=\int^t_0\mathrm d\tau\,L[v(\mathbf x,t;\tau)] \nonumber
\end{align}
and from the initial conditions

\begin{align}\label{eq11-154}
& v(\mathbf x,t;t)=0\Rightarrow \partial_t(v(\mathbf x,t;t))=0,\quad\text{hyperbolic},\\
& v(\mathbf x,t;t)=\frac{g(\mathbf x,t)}{\rho(\mathbf x)},\quad\text{parabolic}\nonumber
\end{align}
Thus we have for the hyperbolic case
\begin{align}\label{eq11-155}
\rho u_{tt}(\mathbf x,t)+L[u(\mathbf x,t)]=\rho(\mathbf x)\left(\frac{g(\mathbf x,t)}{\rho(\mathbf x)}\right)\nonumber\\
+\int^t_0\mathrm d \tau\,\lbrace\rho v_{tt}(\mathbf x,t;\tau)+L[v(\mathbf x,t;\tau)]=g(\mathbf x,t).
\end{align}
(\ref{eq11-155}) show that the $u(\mathbf x,t)$ constructed in (\ref{eq11-152}) solves the inhomogeneous problem for the hyperbolic case. In a similar
way we can show that the function $v(\mathbf x,t;\tau)$ for the parabolic case also solves the ihomogenous equation for that case. (\ref{eq11-152}) is Duhamel's principle.

\subsubsection{The inhomogeneous wave equation}
Let the following initial/boundary value problem be given 

\begin{align}\label{eq11-156}
& u_{tt}(x,t)-c^2u_{xx}(x,t)=g(x,t),\quad -\infty<x<\infty,\,t>0,\\
& u(x,0)=u_t(x,0)=0,\quad -\infty<x<\infty.\nonumber 
\end{align}
According to Duhamel's principle we should construct a function $v(x,t;\tau)$ which solves 

\begin{align}\label{eq11-157}
& v_{tt}(x,t;\tau)-c^2v_{xx}(x,t;\tau)=0,\quad -\infty<x<\infty,\,t>\tau,\nonumber\\
& v(x,\tau;\tau)=0,\\
& v_t(x,\tau;\tau)=g(x,\tau).\nonumber
\end{align}
We solve (\ref{eq11-157}) using the d'Alembert formula

\begin{equation}\label{eq11-158}
u(x,t;\tau)=\frac{1}{2c}\int^{x+c(t-\tau)}_{x-c(t-\tau)}\mathrm d \theta\,g(\theta,\tau).
\end{equation}
The solution to (\ref{eq11-156}), which we denote by $u_p$, is then according to Duhamel's principle 

\begin{equation}\label{eq11-159}
u_p(x,t)=\frac{1}{2c}\int^t_0 \mathrm d\tau\,\int^{x+c(t-\tau)}_{x-c(t-\tau)}\mathrm d \theta\,g(\theta,\tau).
\end{equation}
Let $u_h(x,t)$ be the solution to the homogeneous wave equation with arbitrary initial data 

\begin{align}\label{eq11-160}
& u(x,0)=F(x),\quad-\infty<x<\infty,\\
& u_t(x,0)=G(x).\nonumber
\end{align}
According to the d'Alembert formula the solution with initial data (\ref{eq11-160}) is 

\begin{equation}\label{eq11-161}
u_h(x,t)=\frac{1}{2}\lbrace F(x+ct)+F(x-ct)\rbrace+\frac{1}{2c}\int^{x+ct}_{x-ct}\mathrm d \theta\,G(\theta).
\end{equation}
Let now

\begin{equation}\label{eq11-162}
u(x,t)=u_h(x,t)+u_p(x,t).
\end{equation}
Verify that $u(x,t)$ satisfies the inhomogeneous wave equation with initial data (\ref{eq11-160}). Inserting the formulas for $u_p$ and $u_h$ we get

\begin{align}\label{eq11-163}
u(x,t)=\frac{1}{2}\lbrace F(x+ct)+F(x-ct)\rbrace+\frac{1}{2c}\int^{x+ct}_{x-ct}\mathrm d \theta\,G(\theta)\\
+\frac{1}{2c}\int^t_0 \mathrm d\tau\,\int^{x+c(t-\tau)}_{x-c(t-\tau)}\mathrm d \theta\,g(\theta,\tau).\nonumber
\end{align}
From this formula it is clear that the value of $u$ at a point $(x_0,t_0$), with $t_0>0$ depends only on the data at points in a \emph{characteristic triangle}.

\begin{figure}[h!]
	\centering
	\includegraphics[scale=0.5]{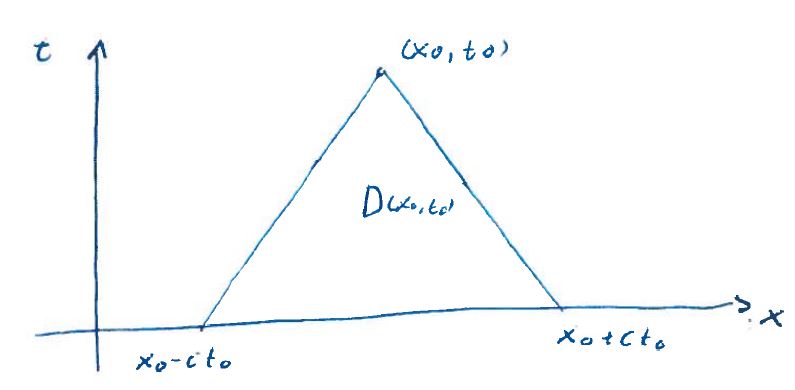}
	\caption{}
\end{figure}

The characteristic triangle is called \emph{the domain of dependence} of the solution at the point $(x_0,t_0)$. 

Similarly, we define the \emph{domain of influence} of a point  $(x_0,t_0)$ to be the set of space-time points $(x,t)$ such that  the domain of dependence of the solution at $(x,t)$ includes the point $(x_0,t_0)$.

\begin{figure}[h!]
	\centering
	\includegraphics[scale=0.5]{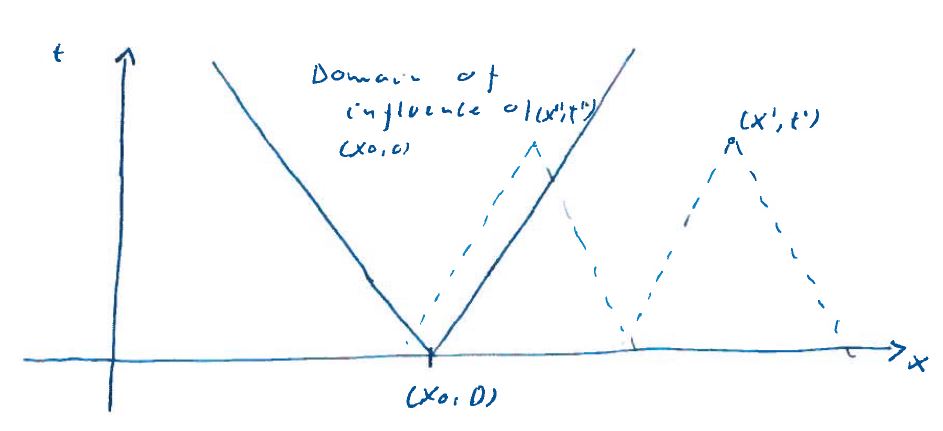}
	\caption{}
\end{figure}
The point $(x_0,t_0)$ will have no causal connection to points beyond its domain of influence. These notions play a fundamental role in the theory of relativity, where the domain of influence of a space-time point $(x,t)$ is called a {\it  light cone}. All material particles and force mediating particles must move within the light cone. This is strictly true only for classical physics. In quantum theory, both material particles and force mediators can, in the form of virtual particles, move beyond the light cone.

\subsubsection{The inhomogeneous heat equation}
Let us consider the 1D heat equation on a finite interval, $0<x<l$ , driven by a heat source $g(x,t)$ located in the interval. Assuming the endpoints are kept at zero temperature, our problem is

\begin{align} \label{eq11-164}
& u_t-c^2u_{xx}=g(x),\quad 0<x<l,\,t>0,\\
& u(0,t)=u(l,t)=0.\nonumber
\end{align} 
The function $u(x,t;\tau)$ appearing in Duhamel's principle satisfies (\ref{eq11-164}) with $g$ replaced by zero and the initial condition 

\begin{equation}\label{eq11-165}
u(x,\tau;\tau)=g(x).
\end{equation}
Separation of variables for this problem gives (page 72) 

\begin{equation}\label{eq11-166}
u(x,t;\tau)=\sqrt{\frac{2}{l}}\sum\limits_{k=1}^{\infty}a_k(\tau)e^{-\left(\frac{\pi k c}{l}\right)^2(t-\tau)}\sin\left(\frac{\pi k}{l}x\right),
\end{equation}
where 

\begin{equation}\label{eq11-167}
a_k(\tau)=\sqrt{\frac{2}{l}}\int_0^l\mathrm d x\,g(x,\tau)\sin\left(\frac{\pi k}{l}x\right),\quad k=1,2,\dots
\end{equation}
Duhamel's principle now gives the solution to (\ref{eq11-164}) in the form 

\begin{equation}\label{eq11-168}
u(x,t;\tau)=\sqrt{\frac{2}{l}}\sum\limits_{k=1}^{\infty}\left\{\int^t_0\mathrm d\tau\, a_k(\tau)e^{-\left(\frac{\pi k c}{l}\right)^2(t-\tau)}\right\}\sin\left(\frac{\pi k}{l}x\right),
\end{equation}
Let us consider the special case 

\begin{equation}\label{eq11-169}
g(x,t)=\sin\left(\frac{\pi}{l}
x\right).
\end{equation}
It is easy to verify that for this choice

\begin{align*}
& a_1=\sqrt{\frac{l}{2}},\\
& a_k=0,\quad k>1,
\end{align*}

\begin{equation}\label{eq11-170}
\Rightarrow\quad u(x,t)=\left(\frac{l}{\pi c}\right)^2\left\{1-e^{-\left(\frac{\pi c}{l}\right)^2 t}\right\}\sin\left(\frac{\pi}{l}x\right).
\end{equation}
Observe that 

\begin{equation}\label{eq11-171}
u(x,t)\rightarrow\left(\frac{l}{\pi c}\right)^2\sin\left(\frac{\pi}{l}x \right)\equiv u(x),
\end{equation}
when $t\to +\infty$. By direct differentiation we find that the limiting function, $u(x)$, is a solution to 

\begin{align}\label{eq11-172}
-c^2u_{xx}(x)&=\sin\left(\frac{\pi}{l}x\right),\\
u(x)&=u(l)=0,\nonumber
\end{align}
which shows that the function $u(x)$ is a time independent, or \emph{stationary} solution to the heat equation.

  In fact, if $g(x,t)=G(x)$ then from (\ref{eq11-168}) we get 

\begin{equation}\label{eq11-174}
u(x,t)=\sqrt{\frac{2}{l}}\sum\limits_{k=1}^{\infty}\left\{\left(\frac{l}{\pi k c}\right)^2 a_k\left(1-e^{-\left(\frac{\pi k c}{l}\right)^2 t}\right)\sin\left(\frac{\pi k}{l}x\right) \right\}.
\end{equation}
For this case the limiting function $u(x)$ is 

\begin{align}\label{eq11-175}
u(x)&=\lim\limits_{t\to +\infty}u(x,t)\\
&=\sqrt{\frac{2}{l}}\sum\limits_{k=1}^{\infty}\left(\frac{l}{\pi k c}\right)^2a_k \sin\left(\frac{\pi k}{l}x\right).
\end{align}
It is easy to verify that $u(x)$ is in fact a stationary solution to the driven heat equation (\ref{eq11-164})

\begin{align}\label{eq11-176}
-c^2u_{xx}(x)&=G(x),\\
u(x)&=u(l)=0.\nonumber
\end{align}

\subsection{The finite Fourier transform}
Duhamel's principle can, as we have seen, be used to construct solutions to inhomogeneous equations. However, the solutions are expressed as an integral over the parameter $\tau$ and are rarely explicit solutions. In order to get actual numbers out we typically have to do numerical integration. This can quickly become awkward and expensive to do in realistic situations. 

There is however another approach that is more general and which is better behaved numerically. This is the \emph{finite Fourier transform}. In this method we expand the solutions in terms of eigenfunctions of the spatial part of the differential operator defining the equation. The expansion coefficients are time dependent and the PDE is reduced to an infinite set of coupled ODEs. We have thus transformed the PDE into a system of ODEs. This process is called the finite Fourier transform. The infinite series of eigenfunctions expressing the solution of the PDE in terms of the time dependent coefficients is called the \emph{inverse finite Fourier transform}.

We have thus the following solution strategy:


\begin{figure}[h!]
	\centering
	\includegraphics[scale=0.5]{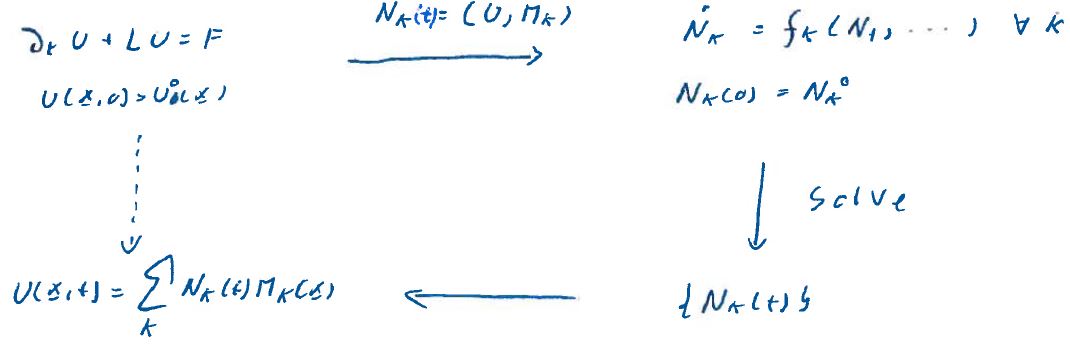}
	\caption{}
\end{figure}

\noindent This \emph{transform strategy} is an important approach to solving complex problems in applied mathematics and theoretical physics. 

We will stick to the restricted class discussed previously in these notes. For that class we have an operator, $L$, given by

\begin{equation}\label{eq11-177}
Lu=-\nabla\cdot(p\nabla u)+qu.
\end{equation}
We can give a unified presentation of the finite Fourier transform if we introduce an operator $K$ which is 

\begin{align}\label{eq11-178}
&K=\partial_t\quad&\text{-- parabolic case,}\nonumber\\
&K=\partial_{tt}\quad&\text{-- hyperbolic case,}\\
&K=-\partial_{yy}\quad&\text{-- elliptic case.}\\
\end{align}
 In fact, if we also allow $K=0$, the general 3D elliptic problem can be included. With this definition of $K$ our problem of interest is 
 
 \begin{equation}\label{eq11-179}
 \rho(\mathbf x)K u+Lu=\rho(\mathbf x)F.
 \end{equation}

 Let $\lbrace M_k\rbrace$ be the complete set of eigenfunctions for the eigenvalue problem
 
 \begin{equation}\label{eq11-180}
 LM_k=\lambda_k\rho M_k,
 \end{equation} 
 where as usual we apply the spatial boundary conditions. We will assume that the eigenfunctions have been normalized so that
 
 \begin{equation}\label{eq11-181}
 (M_k,M_{k'})=\delta_{kk'}.
 \end{equation}
We now expand the solutions of (\ref{eq11-179}) in terms of the eigenfunctions $M_k$

\begin{equation}\label{eq11-182}
u=\sum\limits_k N_k M_k,
\end{equation}
where 

\begin{equation} \label{eq11-183}
N_k=(u,M_k).
\end{equation}
The coefficients $N_k$ are functions of $t$ for the hyperbolic and parabolic case and a function of $y$ for the elliptic case with (generalized) cylinder symmetry. For the general elliptic case the $N_k$'s are constants.

The formula (\ref{eq11-183}) defines the finite Fourier transform and (\ref{eq11-182})  the inverse finite Fourier transform. 

In order to derive an equation for the coefficients $N_k$ we multiply (\ref{eq11-179}) by $M_k(\mathbf x)$ and integrate over the domain for the variable $\mathbf x$. We let this domain be denoted by $G$.
\newpage
 \noindent This integration gives us the equation

\begin{equation}\label{eq11-184}
\int\limits_G\mathrm d V\,\rho(\mathbf x)M_k(\mathbf x)K u(\mathbf x,t)+\int\limits_G\mathrm d V\, M_k(\mathbf x)Lu(\mathbf x,t)=\int\limits_G\mathrm d V\, \rho(\mathbf x)M_k(\mathbf x)F.
\end{equation}
Using the definition of the finite Fourier transform and the inner product, each term in this equation can be reformulated into a more convenient form.

For the first term we have
\begin{eqnarray}\label{eq11-185}
\int\limits_G\mathrm d V \,\rho(\mathbf x)M_k(\mathbf x)K u(\mathbf x,t)&=K\int\limits_G\mathrm d V \,\rho(\mathbf x)u(\mathbf x,t) M_k(\mathbf x)\nonumber\\
&=K(u,M_k(\mathbf x))=K N_k.
\end{eqnarray}
For the second term in(\ref{eq11-184}) we have

\begin{eqnarray}\label{eq11-186}
\int\limits_G\mathrm d V\, M_k(\mathbf x)Lu(\mathbf x,t)=\int\limits_G\mathrm d V\, u(\mathbf x,t)LM_k(\mathbf x)\nonumber\\
 -\int\limits_{\partial G}\mathrm d S\, p(\mathbf x)\lbrace M_k(\mathbf x)\partial_{\mathbf n}u(\mathbf x,t)-u(\mathbf x,t)\partial_{\mathbf n}M_k(\mathbf x)\rbrace\nonumber\\
=\lambda_kN_k-\int\limits_{\partial G}\mathrm d S\, p(\mathbf x)\lbrace M_k(\mathbf x)\partial_{\mathbf n}u(\mathbf x,t)-u(\mathbf x,t)\partial_{\mathbf n}M_k(\mathbf x)\rbrace.
\end{eqnarray}
The second term in (\ref{eq11-186}) can be further simplified using the boundary conditions.

 Recall that the boundary conditions on $u$ and $M_k$ are 

\begin{align}\label{eq11-187}
& \alpha M_k+\beta\partial_{\mathbf n} M_k=0,\\
& \alpha u+\beta\partial_{\mathbf n} u=B,\nonumber
\end{align}
where $u$ and $B$ are functions of $\mathbf x$ and $t$ for the hyperbolic and parabolic case, functions of $\mathbf x$ and $y$ for the restricted elliptic case and a function of $\mathbf x$ for the general elliptic case. As usual  $\alpha=\alpha(\mathbf x)$ and $\beta=\beta(\mathbf x)$.

Using the classification of the boundary conditions into type $1,2,3$ introduced on page 104 
we have for the second term in (\ref{eq11-186})

\begin{align}\label{eq11-188}
& \int\limits_{\partial G}\mathrm d S\, p\lbrace M_k\partial_{\mathbf n}u-u\partial_{\mathbf n}M_k\rbrace.\nonumber\\
& =\int\limits_{\partial G}\mathrm d S\, pM_k\partial_{\mathbf n}u-\int\limits_{\partial G}\mathrm d S\, p\partial_{\mathbf n}M_ku\nonumber\\
& =\int\limits_{S_2\bigcup S_3}\mathrm d S\, pM_k\partial_{\mathbf n}u-\int\limits_{S_1}\mathrm d S\,\frac{p}{\alpha}\partial_{\mathbf n}M_kB\nonumber\\
& -\int\limits_{S_2\bigcup S_3}\mathrm d S\, p\partial_{\mathbf n} M_k u\nonumber\\
& =\int\limits_{S_2\bigcup S_3}\mathrm d S\, \frac{p}{\beta}M_k\beta\partial_{\mathbf n}u+\int\limits_{S_2\bigcup S_3}\mathrm d S\, \frac{p}{\beta}M_k\alpha u\nonumber\\
& -\int\limits_{S_1}\mathrm d S\, \frac{p}{\alpha}\partial_{\mathbf n}M_kB\nonumber\\
& =\int\limits_{S_2\bigcup S_3}\mathrm d S\, \frac{p}{\beta}M_kB-\int\limits_{S_1}\mathrm d S\, \frac{p}{\alpha}\partial_{\mathbf n}M_kB\nonumber\\
&\equiv B_k.
\end{align}
Substituting (\ref{eq11-185}), (\ref{eq11-186}) and (\ref{eq11-188}) into (\ref{eq11-184}) gives us the system 

\begin{equation}  \label{eq11-189}
K N_k+\lambda_kN_k=F_k+B_k,\quad k=1,2,\dots,
\end{equation}
where 

\begin{equation*}
F_k=(F,M_k).
\end{equation*}

\noindent For the general elliptic case $K=0$ so (\ref{eq11-189}) are algebraic equations for the $N_k$. For the restricted elliptic case we have a system of ordinary differential equations

\begin{equation}\label{eq11-190}
-\partial_{yy}N_k(y)+\lambda_kN_k(y)=F_k(y)+B_k(y).
\end{equation}
This system of ODEs must be solved as a boundary value problem. The boundary conditions are 

\begin{align}\label{eq11-191}
N_k(0)&=(u(0),M_k)=(f,M_k),\\
N_k(\hat{l})&=(u(\hat{l}),M_k)=(g,M_k),\nonumber
\end{align}
where we have used the boundary conditions and where $u(0)\equiv u(\mathbf x,t)$ etc.

For the parabolic case $K=\partial_t$ we have a system of ODEs

\begin{equation}\label{eq11-193}
\partial_tN_k(t)+\lambda_kN_k(t)=F_k(t)+B_k(t),\quad k=1,2,\dots,
\end{equation}
subject to the initial condition

\begin{equation}\label{eq11-194}
N_k(0)=(u(0),M_k)=(f,M_k).
\end{equation}
Observe that (\ref{eq11-193}) is an infinite system of ODEs, but they are \emph{uncoupled}. This is a result of the linearity of the underlying PDEs and the fact that the $M_k$'s are eigenfunctions for the operator $L$ defining the equations. The finite Fourier transform can also be used on nonlinear PDEs, as we will see, and for this case the ODE's for the $N_k(t)$ will be \emph{coupled}. The ODEs (\ref{eq11-193}) are of a particular simple type. They are linear first order equations. Such equations can be solved using an integrating factor. The solution is found to be	

\begin{equation}\label{eq11-195}
N_k(t)=N_k(0)e^{-\lambda_k t}+\int\limits_0^t\mathrm d \tau\,[F_k(\tau)+B_k(\tau)]e^{-\lambda_k(t-\tau)}.
\end{equation}
This formula is only useful if we can get some exact or approximate analytical solution of the integral. If a numerical solution is required it might be more efficient to solve the ODEs (\ref{eq11-193})
numerically, using for example Runge-Kutta.

For the hyperbolic case, $K=\partial_{tt}$ and we get a system of second order ODEs

\begin{equation} \label{eq11-196}
\partial_{tt}N_k(t)+\lambda_kN_k(t)=F_k(t)+B_k(t),
\end{equation}
subject to the initial conditions

\begin{align}\label{eq11-197}
N_k(0)&=(u(0),M_k)=(f,M_k),\\
\partial_t N_k(0)&=(u_t(0),M_k)=(g,M_k),\nonumber
\end{align}
This second order equation is also of a simple type which is solvable by the method of variation of parameters. We get the formula

\begin{align}\label{eq11-198}
N_k(t)&=N_k(0)\cos\sqrt{\lambda_k}t+\frac{1}{\sqrt{\lambda_k}}\partial_tN_k(0)\sin\sqrt{\lambda_k}t\\
&+\frac{1}{\sqrt{\lambda_k}}\int^t_0\mathrm d \tau\,[F_k(\tau)+B_k(\tau)]\sin\sqrt{\lambda_k}(t-\tau).\nonumber
\end{align}
Just as for the parabolic case the formula (\ref{eq11-198}) is only useful if some exact or approximate analytical solution can be found. If a numerical solution is required it is more efficient to solve the ODEs (\ref{eq11-196}) using a numerical ODE solver. 

Recall that when the $N_k$s have been computed, the solution to the PDE is constructed using the inverse finite Fourier transform

\begin{equation}\label{eq11-199}
u=\sum\limits_kN_kM_K.
\end{equation}

The representation (\ref{eq11-199})  is problematic when we approach the boundary since all the $M_k$s satisfy homogeneous boundary conditions, whereas $u$ satisfies inhomogeneous boundary conditions. If these conditions are of, say, type 1 \emph{and} if (\ref{eq11-199}) is a classical solution (series converges pointwise) then (\ref{eq11-199}) implies that $u$ is zero on the boundary. But we know that $u$ in general satisfies nonzero boundary conditions. Thus the series (\ref{eq11-199}) can not converge pointwise as we approach the boundary. Therefore, (\ref{eq11-199}) is in general not a classical solution, but rather some form of generalized solution. 

The series might still converge pointwise away from the boundary, but the convergence will typically be slow since the expansion coefficients, $N_k$, "know about" the function over the whole domain and this leads to slow convergence everywhere. 

There is however a way, which we will discuss shortly, to avoid this boundary problem and get series that converge fast enough to get a classical solution.

Let us now look at a couple of examples. 

\subsubsection{Hyperbolic equations}
We consider the following hyperbolic PDE 

\begin{equation}\label{eq11-200}
\rho(\mathbf x)u_{tt}(\mathbf x,t)+Lu(\mathbf x,t)=\rho(\mathbf x)M_i(\mathbf x)\sin \omega t,
\end{equation}
where $M_i(\mathbf x)$ are one of the eigenfunctions of $L$. We assume that initial and boundary conditions are homogeneous

\begin{gather}\label{bc1}
f(\mathbf x)=0,\nonumber\\
g(\mathbf x)=0,\\
B(\mathbf x,t)=0.\nonumber
\end{gather}
For this case we have

\begin{equation}\label{eq11-202}
F_k(t)=(F(t),M_k)=\sin \omega t(M_i,M_k)=\delta_{ik}\sin \omega t,
\end{equation}
and because the initial conditions for the PDE (\ref{eq11-200}) are homogeneous, we have

\begin{align}
N_k(0)&=\partial_tN_k(0)=0,\\
B_k(t)&=0.\nonumber
\end{align}
Thus our system of ODEs simplifies into 

\begin{equation}\label{eq11-204}
\partial_{tt}N_k(t)+\lambda_kN_k(t)=\delta_{ik}\sin \omega t,
\end{equation}
with initial conditions
\begin{align*}
N_k(0)&=0,\\
\partial_tN_k(0)&=0.\nonumber
\end{align*}
We get $N_k=0$, $k\neq i$ and from (\ref{eq11-198})

\begin{align}\label{eq11-205}
N_i(t)&=\frac{1}{\sqrt{\lambda_i}}\int^t_0\mathrm d \tau\sin \omega \tau\sin\sqrt{\lambda_i}(t-\tau)\\
&=\frac{\omega\sin(\sqrt{\lambda_i}t)-\sqrt{\lambda_i}\sin\omega t}{\sqrt{\lambda_i}(\omega^2-\lambda_i)}.\nonumber
\end{align}
The inverse finite Fourier transform gives us

\begin{equation}\label{eq11-206}
u(\mathbf x,t)=\frac{1}{\omega^2-\lambda_i}\left[\frac{\omega}{\sqrt{\lambda_i}}\sin\sqrt{\lambda_i}t-\sin\omega t\right]M_i(\mathbf x).
\end{equation}
This is now the exact solution to this problem. From (\ref{eq11-206}) it is clear that as $\omega\to\sqrt{\lambda_i}$ the amplitude for the solution (\ref{eq11-206}) becomes larger. We can not evaluate (\ref{eq11-206}) at $\omega=\sqrt{\lambda_i}$, since the denominator is then zero, but if we take the limit of (\ref{eq11-206}) as $\omega\to\sqrt{\lambda_i}$, using L'Hopitals rule on the way, we find

\begin{equation}\label{eq11-207}
u(\mathbf x,t)=\frac{1}{2\sqrt{\lambda_i}}\left[\frac{\sin \sqrt{\lambda_i}t}{\sqrt{\lambda_i}}-t\cos \sqrt{\lambda_i} t\right]M_i(\mathbf x).
\end{equation}
The formulas (\ref{eq11-206}) and (\ref{eq11-207}) tell the following story: The general solution of the homogeneous equation

\begin{equation}\label{eq11-208}
\rho(\mathbf x)u_{tt}(\mathbf x,t)+Lu(\mathbf x,t)=0,
\end{equation}
is a sum of normal modes $u_k(\mathbf x,t)$, whose oscillation frequencies are $\omega_k=\sqrt{\lambda_k}$. If we drive (\ref{eq11-205}) with a source that oscillates at a frequency $\omega$ and whose spatial profile is an eigenfunction, $M_i$, the solution is a linear superposition of terms oscillating at $\omega$ and $\omega_i$, with a common spatial profile $M_i(\mathbf x)$. As $\omega$ approaches $\omega_i$, the amplitude of the oscillation grows and in the limit $\omega\to\omega_i$ the system displays linear unbounded growth in time. 

What we have here is the well known phenomena of \emph{resonance} from particle mechanics. Driven continuum system modeled by PDEs can display resonance. Understanding, predicting and preventing resonance is of major importance in applied mathematics and physics. After all we don't want to end up in a plane whose jet engines fall off because resonance has weakened the bolts tying the engine to the plane fuselage! 

\subsubsection{The Poisson equation in a disk}
We consider the boundary value problem for Poisson's equation

\begin{gather}\label{eq11-209}
\nabla^2u(x,y)=-F(x,y),\quad x^2+y^2<R^2,\\
u|_{\tiny{x^2+y^2=R^2}}=f.\nonumber
\end{gather}
For this geometry it makes sense to use polar coordinates

\begin{align}\label{eq11-210}
x=r\cos\theta,\\
y=r\sin\theta.\nonumber
\end{align}
Then 

\begin{align}
\partial_r&=\partial_rx\partial_x+\partial_ry\partial y=\cos\theta\,\partial x+\sin\theta\,\partial_y,\nonumber\\
\partial_\theta&=\partial_{\theta} x\partial_x+\partial_{\theta} y\partial y=-r\sin\theta\,\partial_x+r\cos\theta\,\partial_y,\nonumber\\
&\qquad\Downarrow\nonumber\\
\partial_x&=\cos\theta\,\partial_r-\frac{1}{r}\sin\theta\,\partial_{\theta},\nonumber\\
\partial_y&=\sin\theta\,\partial_r+\frac{1}{r}\cos\theta\,\partial_{\theta},\nonumber\\
&\qquad\Downarrow\nonumber\\
 \partial_{xx}&=(\cos\theta\,\partial_r-\frac{1}{r}\sin\theta\,\partial_{\theta})(\cos\theta\,\partial_r-\frac{1}{r}\sin\theta\,\partial_{\theta})\nonumber\\
&=\cos^2\theta\,\partial_{rr}+\frac{1}{r^2}\sin\theta\,\cos\theta\partial_{\theta}-\frac{1}{r}\sin\theta\cos\theta\,\partial_{r\theta}\nonumber\\
&+\frac{1}{r}\sin^2\theta\,\partial_r-\frac{1}{r}\sin\theta\cos\theta\,\partial_{r\theta}+\frac{1}{r^2}\sin\theta\cos\theta\,\partial_{\theta}\nonumber\\
&+\frac{1}{r^2}\sin^2\theta\partial_{\theta\theta},\nonumber\\
\partial_{yy}&=(\sin\theta\,\partial_r+\frac{1}{r}\cos\theta\,\partial_{\theta})(\sin\theta\,\partial_r+\frac{1}{r}\cos\theta\,\partial_{\theta})\nonumber\\
&=\sin^2\theta\,\partial_{rr}-\frac{1}{r^2}\sin\theta\cos\theta\,\partial_{\theta}+\frac{1}{r}\sin\theta\cos\theta\,\partial_{r\theta}\nonumber\\
&+\frac{1}{r}\cos^2\theta\,\partial_r+\frac{1}{r}\sin\theta\cos\theta\,\partial_{\theta r}-\frac{1}{r^2}\sin\theta\cos\theta\,\partial_{\theta}\nonumber\\
&+\frac{1}{r^2}\cos^2\theta\,\partial_{\theta\theta}\nonumber,\\
&\qquad\Downarrow\nonumber\\
\nabla^2&=\cos^2\theta\,\partial_{rr}+\frac{1}{r^2}\sin\theta\cos\theta\,\partial_{\theta}-\frac{1}{r}\sin\theta\cos\theta\,\partial_{r\theta}\nonumber\\
&+\frac{1}{r}\sin^2\theta\,\partial_{r}-\frac{1}{r}\sin\theta\cos\theta\,\partial_{r\theta}+\frac{1}{r^2}\sin\theta\cos\theta\,\partial_{\theta}\nonumber\\
&+\frac{1}{r^2}\sin^2\theta\,\partial_{\theta\theta}+\sin^2\theta\,\partial_{rr}-\frac{1}{r^2}\sin\theta\cos\theta\,\partial_{\theta}\nonumber\\
&+\frac{1}{r}\sin\theta\cos\theta\,\partial_{r\theta}+\frac{1}{r}\cos^2\theta\,\partial_{r}+\frac{1}{r}\sin\theta\cos\theta\,\partial_{\theta r}\nonumber\\
&-\frac{1}{r^2}\sin\theta\cos\theta\,\partial_{\theta}+\frac{1}{r^2}\cos^2\theta\,\partial_{\theta\theta}\nonumber\\
&=\partial_{rr}+\frac{1}{r}\partial_r+\frac{1}{r^2}\partial_{\theta\theta}.\nonumber
\end{align}
Thus in polar coordinates (\ref{eq11-209}) can be written as

\begin{gather}
\partial_{rr}u(r,\theta)+\frac{1}{r}\partial_ru(r,\theta)+\frac{1}{r^2}\partial_{\theta\theta}u(r,\theta)=F(r,\theta)\nonumber\\
u(R,\theta)=f(\theta). \label{eq11-215}
\end{gather}
By finding the eigenfunctions, $M_k(r,\theta)$, of the 2D elliptic operator

\begin{equation}\label{eq11-216}
L=\partial_{rr}+\frac{1}{r}\partial_r+\frac{1}{r^2}\partial_{\theta\theta},
\end{equation}
we can express the solution $u$ as a Fourier series

\begin{equation}\label{eq11-217}
u(r,\theta)=\sum\limits_ka_kM_k(r,\theta).
\end{equation}
This is done in chapter 8 of our text-book. Here we will solve this problem using the finite Fourier transform. Introduce the operators

\begin{align}\label{eq11-218}
K&=-r^2\partial_{rr}-r\partial_r,\\
\hat{L}&=-\partial_{\theta\theta}.
\end{align}
Then (\ref{eq11-215}) can be written as 

\begin{equation}\label{eq11-219}
K u(r,\theta)+\hat{L}u(r,\theta)=r^2F(r,\theta).
\end{equation}
The idea is now to introduce a finite Fourier transform based on the eigenfunctions of the operator $\hat{L}$. Observe that since $\theta$ is the polar angle, running from zero to $2\pi$, all functions are periodic in $\theta$ with period $2\pi$. Thus we are lead to the eigenvalue problem
 
\begin{gather}
-\partial_{\theta\theta}M(\theta)=\lambda M(\theta),\quad 0<\theta<2\pi,\label{eq11-220}\\
M(\theta+2\pi)=M(\theta). \label{eq11-221}
\end{gather}
This problem is easiest to solve using complex exponentials. The general solution of the equation (\ref{eq11-220})  is for $\lambda>0$

\begin{equation}\label{eq11-222}
M(\theta)=Ae^{i\sqrt{\lambda}\theta}+Be^{-i\sqrt{\lambda}\theta},
\end{equation}
and the periodicity condition (\ref{eq11-221}) implies that

\begin{gather}
Ae^{i\sqrt{\lambda}\theta}+Be^{-i\sqrt{\lambda}\theta}=Ae^{i\sqrt{\lambda}(\theta+2\pi)}+Be^{-i\sqrt{\lambda}(\theta+2\pi)},\nonumber\\
\Updownarrow\nonumber\\
e^{i\sqrt{\lambda}\theta}\left(1-e^{i\sqrt{\lambda}2\pi}\right)A+e^{-i\sqrt{\lambda}\theta}\left(1-e^{-i\sqrt{\lambda}2\pi}\right)B=0. \label{eq11-223}
\end{gather}
The functions $e^{\pm i\sqrt{\lambda}\theta}$ are linearly independent (verify this), so we must conclude from (\ref{eq11-223}) that 

\begin{align}\label{eq11-224}
A\left(1-e^{i\sqrt{\lambda}2\pi}\right)&=0,\\
B\left(1-e^{-i\sqrt{\lambda}2\pi}\right)&=0,\Leftrightarrow B\left(1-e^{i\sqrt{\lambda}2\pi}\right)=0.\nonumber
\end{align}
If $\left(1-e^{-i\sqrt{\lambda}2\pi}\right)\neq 0$ both $A$ and $B$ must be zero. This gives $M(\theta)=0$, so such $\lambda$s are not eigenvalues. The other possibility is that 

\begin{equation}\label{eq11-225}
e^{-i\sqrt{\lambda}2\pi}=1,
\end{equation}
and then $A$, $B$ are arbitrary. The equation (\ref{eq11-225}) is only satisfied if 

\begin{gather}\label{eq11-226}
2\pi\sqrt{\lambda_k}=2\pi k,\quad k=1,2\dots,\nonumber\\
\Updownarrow\nonumber\\
\lambda_k=k^2.
\end{gather}
These are the positive eigenvalues. The corresponding eigenspaces are two-dimensional and spanned by the normalized eigenfunctions

\begin{equation}\label{eq11-227}
\frac{1}{\sqrt{\pi}}\cos k\theta,\quad\frac{1}{\sqrt{\pi}}\sin k\theta.
\end{equation} 
In order to find out if $\lambda=0$ is an eigenvalue we put $\lambda=0$ in (\ref{eq11-220}) and get the equation 

\begin{equation}\label{eq11-228}
\partial_{\theta\theta}M(\theta)=0, 
\end{equation} 
whose general solution is

\begin{equation}\label{eq11-229a} 
M(\theta)=A\theta+B.
\end{equation}
The boundary condition now gives 

\begin{gather*}
M(\theta+2\pi)=M(\theta),\\
\Updownarrow\\
A(\theta+2\pi)+B=A\theta+B,\\
\Updownarrow\\
2\pi A=0\,\Rightarrow\,A=0,
\end{gather*}
and $B$ is arbitrary, so the eigenspace corresponding to $\lambda=0$ is one-dimensional and spanned by the normalized eigenfunction

\begin{equation} \label{eq11-229b} 
\frac{1}{\sqrt{2\pi}}.
\end{equation}
Verify that there are no negative eigenvalues.

  Using (\ref{eq11-227}) and (\ref{eq11-229b}) we define a finite Fourier transform for functions $u(r,\theta)$ by 

\begin{align}\label{eq11-230}
& a_0(r)=\frac{1}{\sqrt{2\pi}}\int_0^{2\pi}\mathrm d\theta\,u(r,\theta),\\
&\left.
\begin{aligned}
 a_k(r)=\frac{1}{\sqrt{\pi}}\int_0^{2\pi}\mathrm d\theta\,u(r,\theta)\cos k\theta\\
 b_k(r)=\frac{1}{\sqrt{\pi}}\int_0^{2\pi}\mathrm d\theta\,u(r,\theta)\sin k\theta
\end{aligned}
\right\}\quad k=1,2,\dots\;\;.
\end{align} 
The corresponding inverse Finite Fourier transform is

\begin{align}\label{eq11-231}
u(r,\theta)&=\frac{1}{\sqrt{2\pi}}a_0(r)\nonumber\\
&+\sum\limits_{k=1}^{\infty}\frac{1}{\sqrt{\pi}}\lbrace a_k(r)\cos k\theta+b_k(r)\sin k\theta\rbrace.
\end{align}
The ODEs for the coefficients $a_k$ and $b_k$ are now found by multiplying (\ref{eq11-219}) by the eigenfunctions and integrating over the interval $(0,2\pi)$. Verify that we get

\begin{align}\label{eq11-232}
a''_k(r)+\frac{1}{r}a_k'(r)-\frac{k^2}{r^2}a_k(r)&=-A_k(r),\quad k=0,1,2,\dots,\nonumber\\
b''_k(r)+\frac{1}{r}b_k'(r)-\frac{k^2}{r^2}b_k(r)&=-B_k(r),\quad k=0,1,2,\dots,
\end{align}
where $A_k$ and $B_k$ are the Fourier coefficients of $F$,

\begin{align}\label{eq11-233}
& A_0(r)=\frac{1}{\sqrt{2\pi}}\int_0^{2\pi}\mathrm d\theta\,F(r,\theta),\\
&\left.
\begin{aligned}
A_k(r)=\frac{1}{\sqrt{\pi}}\int_0^{2\pi}\mathrm d\theta\,F(r,\theta)\cos k\theta\\
B_k(r)=\frac{1}{\sqrt{\pi}}\int_0^{2\pi}\mathrm d\theta\,F(r,\theta)\cos k\theta
\end{aligned}
\right\}\quad k=1,2,\dots\,\,.
\end{align} 
Note that only solutions of the singular equations (\ref{eq11-232}) that are bounded as $r\to 0$ are acceptable. So our boundary condition at $r=0$ is

\begin{equation}\label{eq11-234}
\lim\limits_{r\to 0}\,\lbrace a_k(r),b_k(r)\rbrace\text{ exists.}
\end{equation}
At $r=R$ the boundary condition is

\begin{equation*}
u(R,\theta)=f(\theta).
\end{equation*}
Taking the Finite Fourier transform of the boundary data $f(\theta)$, we find the the boundary conditions for the functions $a_k(r),b_k(r)$ at $r=R$ aret 

\begin{align}\label{eq11-236}
&a_0(R)= \alpha_0\equiv\frac{1}{\sqrt{2\pi}}\int_0^{2\pi}\mathrm d\theta\,f(\theta),\\
&\left.
\begin{aligned}
a_k(R)=\alpha_k\equiv\frac{1}{\sqrt{\pi}}\int_0^{2\pi}\mathrm d\theta\,f(\theta)\cos k\theta\\
b_k(R)=\beta_k\equiv\frac{1}{\sqrt{\pi}}\int_0^{2\pi}\mathrm d\theta\,f(\theta)\sin k\theta
\end{aligned}
\right\}\quad k=1,2,\dots\,\,.
\end{align} 
Note that both of the equations in (\ref{eq11-232}) are inhomogeneous versions of the homogeneous \emph{Euler equation}

\begin{equation}\label{eq11-237}
C''_k(r)+\frac{1}{r}C_k'(r)-\frac{k^2}{r^2}C_k(r)=0.
\end{equation}
A basis for the solution space of such equations can be found using standard methods. We have

\begin{equation}\label{eq11-238}
C_0(r)=
\left\{
\begin{aligned}
\text{constant}\\
\ln r\\
\end{aligned}
\right.
,\quad
C_k(r)=\left\{
\begin{aligned}
r^k\\
r^{-k}
\end{aligned}
\right.\quad k=1,2,\dots
\end{equation}

The solutions of (\ref{eq11-232}) can now be found using variation of parameters. We have 

\begin{align} \label{eq11-239}
a_0(r)&=\int_0^r\mathrm d t\,\ln\left(\frac{R}{r}\right)A_0(t)t+\int_r^R\mathrm d t\,\ln\left(\frac{R}{t}\right)A_0(t)t+\alpha_k(R),\nonumber\\
a_k(r)&=\frac{1}{2k}\int_0^r\mathrm d t\left[\left(\frac{R}{r}\right)^k-\left(\frac{r}{R}\right)^k\right]\left(\frac{t}{R}\right)^kA_k(t)t\\
&+\frac{1}{2k}\int_r^R\mathrm d t\left[\left(\frac{R}{t}\right)^k-\left(\frac{t}{R}\right)^k\right]\left(\frac{r}{R}\right)^kA_k(t)t+\alpha_k(R)\left(\frac{r}{R}\right)^k,\nonumber\\
b_k(r)&=\frac{1}{2k}\int_0^r\mathrm d t\left[\left(\frac{R}{r}\right)^k-\left(\frac{r}{R}\right)^k\right]\left(\frac{t}{R}\right)^kB_k(t)t\\
&+\frac{1}{2k}\int_r^R\mathrm d t\left[\left(\frac{R}{t}\right)^k-\left(\frac{t}{R}\right)^k\right]\left(\frac{r}{R}\right)^kB_k(t)t+\beta_k(R)\left(\frac{r}{R}\right)^k.\nonumber
\end{align}
This completes the formal solution of (\ref{eq11-209}). 

\noindent Let us now return to the general scheme. As we have indicated, the solutions found using the finite Fourier transform will in general have slow convergence because the eigenfunctions satisfy the homogeneous boundary conditions, whereas in general the solution might satisfies inhomogeneous boundary conditions. 

However,  the problem of slow convergence can in general be solved using an approach which we will now explain. 

Let us first consider the case when the inhomogeneous terms in the equation and boundary conditions are independent of time (we are here thinking about the hyperbolic and the parabolic cases).

Thus we have

\begin{align}\label{eq11-240}
\rho(\mathbf x)K u(\mathbf x,t)+Lu(\mathbf x,t)&=\rho(\mathbf x)F(\mathbf x)\quad\text{ in }G,\text{ for }t>0,\nonumber\\
\alpha(\mathbf x)u(\mathbf x,t)+\beta(\mathbf x)\partial_{\mathbf n}u(\mathbf x,t)&=B(\mathbf x)\quad\text{ on }\partial G,\text{ for }t>0.
\end{align}
Let $v(\mathbf x)$ be a time independent solution to (\ref{eq11-240}). Thus 

\begin{gather}\label{eq11-241}
Lv(\mathbf x)=\rho(\mathbf x)F(\mathbf x)\quad\text{ in }G,\\
\alpha(\mathbf x)v(\mathbf x)+\beta(\mathbf x)\partial_{\mathbf n}v(\mathbf x)=B(\mathbf x)\quad\text{ on }\partial G.\nonumber
\end{gather}
We will assume that we can solve (\ref{eq11-241}). Define $w(\mathbf x,t)=u(\mathbf x,t)-v(\mathbf x)$. Then $w(\mathbf x)$ solves

\begin{align}\label{eq11-242}
&\rho(\mathbf x)Kw(\mathbf x,t)+Lw(\mathbf x,t)\\
&=\rho(\mathbf x)K u(\mathbf x,t)-\rho(\mathbf x)K v(\mathbf x)+Lu(\mathbf x,t)-Lv(\mathbf x)\nonumber\\
&=\rho(\mathbf x)K u(\mathbf x,t)+Lu(\mathbf x,t)-Lv(\mathbf x)\nonumber\\
&=\rho(\mathbf x)(F(\mathbf x)-F(\mathbf x))=0,
\end{align}
and has boundary condition

\begin{align}\label{eq11-243}
&\alpha(\mathbf x)w(\mathbf x,t)+\beta(\mathbf x)\partial_{\mathbf n}w(\mathbf x,t)\\
&=\alpha(\mathbf x)u(\mathbf x,t)-\alpha(\mathbf x)v(\mathbf x)\nonumber\\
&+\beta(\mathbf x)\partial_{\mathbf n}u(\mathbf x,t)-\beta(\mathbf x)\partial_{\mathbf n}v(\mathbf x)\nonumber\\
&=B(\mathbf x)-B(\mathbf x)=0. 
\end{align}
Thus $w(\mathbf x,t)$ satisfies a homogeneous equation and homogeneous boundary conditions. The formal series for $w(\mathbf x,t)$ will therefore have better convergence properties than the formal series for $u(\mathbf x,t)$ . If $v(\mathbf x)$ can be found without using eigenfunction expansions, the solution $u(\mathbf x,t)$ will be represented by a Fourier series with good convergence properties. 

\noindent Finding an expression for $v(\mathbf x)$ might however not be easy. 

There is however a related approach that can be used and that do not rely on finding a stationary solution. 

Let us assume that we can find a function $v(\mathbf x,t)$ such that 

\begin{equation}\label{eq11-244}
\alpha(\mathbf x)v(\mathbf x,t)+\beta(\mathbf x)\partial_{\mathbf n}v(\mathbf x,t)=B(\mathbf x,t).
\end{equation}
The function $v$ must be sufficiently smooth so that it is  possible to substitute into the equation 

\begin{equation}\label{eq11-245a}
\rho(\mathbf x)Ku(\mathbf x,t)+Lu(\mathbf x,t)=\rho(\mathbf x)F(\mathbf x,t).
\end{equation}
Define a function $w(\mathbf x,t)=u(\mathbf x,t)-v(\mathbf x,t)$. Then 

\begin{align} \label{eq11-245b}
&= \alpha(\mathbf x)w(\mathbf x,t)+\beta(\mathbf x)\partial_{\mathbf n}w(\mathbf x,t)\\
&= \alpha(\mathbf x)u(\mathbf x,t)+\alpha(\mathbf x)v(\mathbf x,t)\nonumber\\
&= \beta(\mathbf x)\partial_{\mathbf n}u(\mathbf x,t)- \beta(\mathbf x)\partial_{\mathbf n}v(\mathbf x,t)=B(\mathbf x,t)-B(\mathbf x,t)=0,\nonumber
\end{align}
and 

\begin{align}\label{eq11-246}
&\rho(\mathbf x)K w(\mathbf x,t)+Lw(\mathbf x,t)\\
&=\rho(\mathbf x)K u(\mathbf x,t)-\rho(\mathbf x)K v(\mathbf x,t)\nonumber+ Lu(\mathbf x,t)-Lv(\mathbf x,t)\nonumber\\
&=\rho(\mathbf x)F(\mathbf x,t)-\rho(\mathbf x)K v(\mathbf x,t)-Lv(\mathbf x,t)\nonumber\\
&\equiv \rho(\mathbf x)\hat{F}(\mathbf x,t),\nonumber
\end{align}
where
\begin{equation*}
\hat{F}(\mathbf x,t)=F(\mathbf x,t)-K v(\mathbf x,t)-\frac{Lv(\mathbf x,t)}{\rho(\mathbf x)}.
\end{equation*}
Thus $w(\mathbf x,t)$ solves an inhomogeneous equation with homogeneous boundary conditions. The solution to this problem can be found using Duhamel's principle or the finite Fourier transform. Convergence of the series should now be much faster and thus lead to more efficient numerical calculations. 

Let us illustrate this method for a 1D inhomogeneous heat equation. 

\begin{align}\label{eq11-247}
& u_t(x,t)-c^2u_{xx}(x,t)=g(x,t),\quad 0<x<l,\, t>0,\\
& u(x,0)=f(x),\quad 0<x<l,\nonumber\\
& u(0,t)=g_1(t),\quad t>0,\nonumber\\
& u(l,t)=g_2(t),\quad t>0.\nonumber\\
\end{align}
It is easy to see, using linear interpolation that the function 

\begin{equation}\label{eq11-248}
v(x,t)=\frac{1}{l}(xg_2(t)+(l-x)g_1(t)),
\end{equation}
satisfy the inhomogeneous boundary conditions. Defining $w(x,t)=u(x,t)-v(x,t)$, we find immediately that 

\begin{equation}\label{eq11-249}
w_t(x,t)-c^2w_{xx}(x,t)=g(x,t)-\frac{1}{l}(xg'_2(t)+(l-x)g_1'(t)),
\end{equation}
for $t>0$ and $0<x<l$. The boundary condition are by construction homogeneous 

\begin{align}\label{eq11-250}
w(0,t)&=0,\\
w(l,t)&=0,\nonumber
\end{align}
and the initial value is 

\begin{align}\label{eq11-251}
w(\mathbf x,t)&=u(x,0)-v(x,0)\\
&=f(x)-\frac{1}{l}(xg_2(0)-(l-x)g_1(0)).\nonumber
\end{align}
This problem can now be solved by, for example series expansion, and the convergence of the series should be good. 

\subsection{Nonlinear stability theory}
The main focus of this class is the theory of \emph{linear} partial differential equations. The reason for this is two-fold. Firstly, many systems of scientific and technological interest can be described by such equation. We have seen several modeling examples where this turns out to be the case. Secondly, there is no theory of nonlinear PDEs with a generality anywhere close to the case for linear PDEs. Nonlinear PDEs in general require an individual treatment. In fact some people spend their whole working career studying one or a few such nonlinear equations. This is, however, if one is looking for exact solutions, precise statements about solution spaces etc. If one is satisfied with approximate answers, analytic or numeric, a lot more can be achieved. 
	
The reason for this is that in many situations the nonlinear terms in the PDEs are much smaller than the linear ones. For example, we saw at the beginning of this class that the PDE

\begin{equation}\label{eq11-252}
\partial_tn+a(n)\partial_x n=0,
\end{equation}
is a simple model of traffic flow on a one-lane road. Assuming $n$ is small, thus the density of cars is low, we can Taylor expand $a(n)$ and truncate at first order

\begin{equation}\label{eq11-253}
a(n)\approx a_0+a_1n.
\end{equation}
Inserting (\ref{eq11-253}) into (\ref{eq11-252}) and rearranging the terms we get the nonlinear equation

\begin{equation}\label{eq11-254}
\partial_tn+a_0\partial_xn=-a_1n\partial_xn.
\end{equation}
By assumption (\ref{eq11-253}), the nonlinear term in this equation is smaller that the linear ones. 

The use of truncated expansions like (\ref{eq11-253}) are very common when modeling physical systems by PDEs. We saw this for example for the vibrating string. 

The refractive index, $n$, is a material property that determines the optical response of dielectrics and gases. For inhomogenous materials, like corrugated glass for example, it will depend on position $n=n(\mathbf x)$.

An important equation in the theoretical description of light is 

\begin{equation}\label{eq11-255}
\frac{1}{c^2}\partial_{tt}(n^2E)-\nabla^2E=0.
\end{equation}
This is a equation which model light under normal conditions well. It is clearly a linear equation. For light under more extreme conditions, for example light generated by a laser, one find that the refractive index depends on the  light intensity $I=|E|^2$. We thus in general have have $n=n(\mathbf x,I)$. This functions has been measured for many materials and one find that in many cases a good approximation is 

\begin{equation}\label{eq11-256}
n=n_0+n_2I,
\end{equation}
where both $n_0$ and $n_2$ in general depends on position. Inserting (\ref{eq11-256}) into (\ref{eq11-255}) and keeping only terms that are linear in the intensity $I$, we get the equation

\begin{equation}\label{eq11-257}
\frac{n_0^2}{c^2}\partial_{tt}E-\nabla^2E=-\frac{2n_0n_2}{c^2}\partial_{tt}(|E|^2E).
\end{equation}
This is now a nonlinear equation. The factor $n_2$ is for most materials very small, so ligh of very high intensity $I=|E|^2$ is required  in order for the nonlinearity to have any effect. From your basic physics class you might recall that light tend to bend towards areas of higher refractive index; this property is the basis for Snell's law in optics. Thus the intensity dependent index (\ref{eq11-256}) will tend to focus light into areas of high light intensity. This increases the intensity, creating an even larger index, which focus light even more strongly. This is a runaway effect that quickly creates a local light intensity high enough to destroy the material the light is traveling through.

The study of nonlinear properties of high-intensity light is a large and very active field of research which has given us amazing technologies, like the laser, optical fibers, perfect lenses of microscopic size and even invisibility cloaks! (Well almost anyway)

For some physical systems there are no good linear approximations. One such area is fluid dynamics. The most important equation in this field of science is the \emph{Navier-Stokes} equation. 

\begin{equation}\label{eq11-258}
\rho_0(\partial_t\mathbf u +\mathbf u\cdot \nabla\mathbf u)=-\nabla p+\mu\nabla^2\mathbf u,\quad \nabla\cdot\mathbf u=0,
\end{equation}
where $\rho_0$ is the mass density, $\mathbf u(\mathbf x,t)$ the fluid velocity field, $p(\mathbf x,t)$ the pressure in the fluid and $\mu$ the viscosity (friction). The second term on the left hand side of (\ref{eq11-258}) is clearly nonlinear and in many important applications not small. The essential nonlinearity of the Navier-Stokes equation gives rise to the phenomenon of \emph{turbulence} and other physical phenomena that are extremely hard to describe and predict. In fact, it is not even known, after almost 200 years, if the Cauchy problem for (\ref{eq11-258}) is well posed. There is a possibility that solutions of (\ref{eq11-258}) develop singularities in finite time. If you can solve this problem, the Clay mathematics institute will pay you one million dollars!

For systems where the nonlinear terms are small, one would expect that methods introduced for solving linear PDEs should be useful, somehow. They are, and we will illustrate this using a nonlinear heat equation

\begin{equation}\label{eq11-259}
u_t- u_{xx}=\hat{\lambda}u(1-u^2).
\end{equation}
The parameter $\hat{\lambda}$ measures the strength of the nonlinearity. The equation can model situations where we have a nonlinear heat source. The equation (\ref{eq11-259}) clearly has a solution

\begin{equation}\label{eq11-260}
u_0(x,t)=0.
\end{equation}
This corresponds to a situation where the system has a uniform temperature equal to zero. We know that if the nonlinear heat source is not present, any initial disturbance will decay, and the system will return to a uniform zero temperature state. We have seen this happening using normal modes and also Fourier series solutions. The question we want to ask is if $u_0(x,t)=0$ is stable also in the presence of the nonlinear heat source. We will investigate this problem both for an infinite system, $-\infty<x<\infty$, and a finite one $(0,l)$. We will see that the stability condition is different in the two cases. This is a general fact: Finite and infinite systems display different stability behavior. 

Let us first set the problem up. Since we will be concerned with solutions close to $u_0(x,t)=0$, we introduce a small parameter $0<\epsilon<1$ and write the initial condition as 

\begin{equation}\label{eq11-261}
u(x,0)=\epsilon h(x),
\end{equation}
where $h(x)$ is of normal size, $h(x)=\mathcal{O}(1)$. We are similarly introducing a $\mathcal{O}(1)$ function $w(x,t)$ through 

\begin{equation}\label{eq11-262}
u(x,t)=\epsilon w(x,t).
\end{equation}
The equation (\ref{eq11-259}) now gives

\begin{equation}\label{eq11-263}
w_t-w_{xx}=\hat{\lambda}w-\epsilon^2\hat{\lambda}w^3.
\end{equation}
The first step is to do a \emph{linear stability analysis}. For this we recognize that, if $w=\mathcal{O}(1)$, the last term in (\ref{eq11-263}) is much smaller than the other terms because $\epsilon<<1$. We therefore drop the small term and are left with the equation 

\begin{equation}\label{eq11-264}
w_t-w_{xx}=\hat{\lambda}w.
\end{equation}
(\ref{eq11-264}) is an arbitrary good approximation to (\ref{eq11-263}) if $\epsilon$ is small enough and $w=\mathcal{O}(1)$. Going from (\ref{eq11-263}) to (\ref{eq11-264}) is called \emph{linearization}. Linearization is an idea of outmost importance in applied mathematics. 

We will now investigate the linearized equation (\ref{eq11-264}) for the infinite and finite case separately. 

\subsubsection{The infinite case}
In the infinite case $-\infty<x<\infty$, we can analyze the stability of equation (\ref{eq11-264}) using normal modes as discussed previously in this class. Inserting a normal mode solution 

\begin{equation}\label{eq11-265}
w(x,t)=a(k)e^{ikx+\lambda(k)t},
\end{equation}
into (\ref{eq11-264}), we find the algebraic equation 

\begin{gather}\label{eq11-266}
\lambda(k)+k^2=\hat{\lambda},\\
\Updownarrow\nonumber\\
\lambda(k)=\hat{\lambda}-k^2.
\end{gather} 
From this we see that the stability index for the equation (\ref{eq11-264}) is 

\begin{equation}\label{eq11-267}
\Omega=\sup\limits_k\Re[\lambda(k)]=\hat{\lambda}.
\end{equation}
The conclusion is then that

\begin{enumerate}
	\item $u_0$ is stable if $\hat{\lambda}<0$,
	\item $u_0$ is unstable if $\hat{\lambda}>0$.
\end{enumerate}
In the unstable case the normal mode (\ref{eq11-265}) experience unbounded exponential growth. However, when the growth makes 

\begin{equation}\label{eq11-268}
\hat{\lambda}w\sim\epsilon^2\hat{\lambda}w^3,
\end{equation}
the nonlinear term is of the same size as the linear one and the linearized equation is not valid anymore. From (\ref{eq11-268}) this happend when

\begin{gather}\label{eq11-269}
\epsilon^2w^2\sim 1,\\
\Updownarrow\nonumber\\
w\sim\frac{1}{\epsilon},\nonumber
\end{gather}
which gives us, using (\ref{eq11-265}), a break-down time 

\begin{equation}\label{eq11-270}
t^*=\frac{|\log\epsilon|}{\hat{\lambda}}.
\end{equation}
So what happens after the break-down time? Observe that the nonlinear equation (\ref{eq11-259}) has another simple uniform solution 

\begin{equation}\label{eq11-271}
u_1(x,t)=1.
\end{equation}
We can now ask when this solution is stable. We introduce a small parameter $\epsilon$ and a $\mathcal{O}(1)$ function $v(x,t)$ through

\begin{equation}\label{eq11-272}
u(x,t)=1+\epsilon v(x,t).
\end{equation}
Inserting this into (\ref{eq11-259}) we find 

\begin{align*}
\epsilon(v_t-v_{xx})&=\hat{\lambda}(1+\epsilon v)(1-(1+\epsilon v)^2)\\
&=\hat{\lambda}(1+\epsilon v)(1-1-2\epsilon v-\epsilon^2v^2)\\
&=-\epsilon\hat{\lambda}(1+\epsilon v)(2v+\epsilon v^2)\\
&=-2\epsilon\hat{\lambda}v-2\epsilon^2\hat{\lambda}v^2-\epsilon^2\hat{\lambda}v^2-\epsilon^2\hat{\lambda}v^3,\\
\;\;\;\Downarrow\nonumber\\
v_t-v_{xx}&=-2\hat{\lambda}v-3\epsilon\hat{\lambda}v^2-\epsilon^2\hat{\lambda}v^3.
\end{align*} 
Linearizing, we find the equation
 
\begin{equation}\label{eq11-274}
v_t-v_{xx}=-2\hat{\lambda}v.
\end{equation}
Inserting a normal mode now gives 

\begin{equation}\label{eq11-275}
\lambda=-2\hat{\lambda}-k^2,
\end{equation}
and we can conclude
 
\begin{enumerate}
	\item $u_1$ is stable if $\hat{\lambda}>0$,
	\item $u_1$ is unstable if $\hat{\lambda}<0$.
\end{enumerate}
Thus, the solution $u_1$ is stable when $u_0$ is unstable. A reasonable guess now is that a small perturbation to the solution $u_0$ in the unstable regime $(\hat{\lambda}>0)$ grows and eventually approaches the stable solution $u_1$. Solving equation (\ref{eq11-259}) numerically with a broad selection of small initial data $u(x,0)=\epsilon h(x)$ is needed to support our conjecture. Such situations can not \emph{prove}, in the mathematical sense, that our conjecture is true. However, in real-life applied mathematics, proofs are few and far between,  and we must use the tools at our disposal to do as much as we can to dispel uncertainty as to the validity of our conclusions. 

Another way to support our conjecture is to find solutions for special initial data. Let us consider the special initial data

\begin{equation}\label{eq11-276}
u(x,0)=\epsilon\qquad(\text{independent of }x).
\end{equation}
Let us look for a solution to (\ref{eq11-259}) that is independent of $x$, $u=u(t)$. Such a solution must satisfy

\begin{gather}\label{eq11-277}
u'(t)=\hat{\lambda}u(t)(1-u^2(t),\\
u(0)=\epsilon.\nonumber
\end{gather}
This first order ODE that can be solved exactly 

\begin{equation}\label{eq11-278}
u(x,t)=u(x)=\frac{\epsilon e^{\hat{\lambda}t}}{(1+\epsilon^2(e^{2\hat{\lambda}t}-1))^{\frac{1}{2}}}.
\end{equation}
It is easy to verify that 

\begin{equation}\label{eq11-279}
\lim\limits_{t\to +\infty}u(x,t)=u_t(x,t)=1,
\end{equation}
for the case $\hat{\lambda}>0$. Thus this special solution lends further support to our conjecture. 

\subsubsection{The finite case}
Let us consider the nonlinear heat equation (\ref{eq11-259}) on a finite interval $0<x<\pi$ with zero boundary conditions

\begin{equation}\label{eq11-280}
u(0,t)=u(\pi,t)=0.
\end{equation}
Recall that linearizing (\ref{eq11-256}) around the simple solution $u_0(x,t)=0$ gives us an equation of the form 

\begin{align}\label{eq11-281a} 
& w_t -w_{xx}=\hat{\lambda}w,\quad 0<x<\pi,\\
& w(0,t)=w(\pi,t)=0.
\end{align}
We can solve this problem using the finite Fourier transform. Introduce the operator $L=-\partial_{xx}$. The relevant eigenfunctions satisfy

\begin{align}\label{eq11-281b} 
LM_k(x)=\lambda M_k(x)\\
M_k(0)=M_k(\pi)=0\nonumber
\end{align}
We have analyzed this eigenvalue problem before. The eigenvalues are 

\begin{equation*}
\lambda_k=-k^2,\quad k=1,2,\dots
\end{equation*}
and the normalized eigenfunctions are 

\begin{equation}\label{eq11-282}
M_k(x)=\sqrt{\frac{2}{\pi}}\sin kx.
\end{equation}


\noindent The associated finite Fourier transform of a given function $f(x)$ is

\begin{equation}\label{eq11-282b}  
N_k=(f,M_k),\quad k=1,2,\dots
\end{equation}
and the inverse finite Fourier transform is 

\begin{equation} \label{eq11-283}
f(x)=\sum\limits_kN_kM_k(x).
\end{equation}
We now express an arbitrary solution to (\ref{eq11-281a}) in terms of an inverse finite Fourier transform

\begin{equation}\label{eq11-284}
w(x,t)=\sum\limits_k N_k(t)M_k(x).
\end{equation}
The evolution equations for the $N_k(t)$ are found from (\ref{eq11-281a}) if we multiply by $M_k(x)$ and integrate

\begin{align}\label{eq11-285}
\int^{\pi}_0\mathrm d x\, w_t(x,t)M_k(x)&-\int_0^{\pi}\mathrm d x\,w_{xx}(x,t)M_k(x)\\
&=\int_0^{\pi}\mathrm d x\,\hat{\lambda}w(x,t)M_k(x).\nonumber
\end{align}
For each of the three terms in this equation we have

\begin{equation}\label{eq11-286}
\int^{\pi}_0\mathrm d x\, w_t(x,t)M_k(x)=\partial_t\int^{\pi}_0\mathrm d x\, w(x,t)M_k(x)=N'_k(t),
\end{equation}

\begin{align}\label{eq11-287}
\int^{\pi}_0\mathrm d x\, w_{xx}(x,t)M_k(x)&=w_x(x,t)M_k(x)\big|^{\pi}_0-\int^{\pi}_0\mathrm d x\, w_x(x,t)M_{kx}(x)\\
& -w(x,t)M_{kx}(x)\big|^{\pi}_0+\int^{\pi}_0\mathrm d x\, w(x,t)M_{kxx}(x)\nonumber\\
&= \int^{\pi}_0\mathrm d x\, w(x,t)(-k^2)M_k(x)=-k^2N_k(t),
\end{align}

\begin{equation}\label{eq11-288}
\int^{\pi}_0\mathrm d x\, \hat{\lambda}w(x,t)M_k(x)=\hat{\lambda}N_k(t)
\end{equation}
Inserting (\ref{eq11-280}), (\ref{eq11-287}) and (\ref{eq11-288}) into (\ref{eq11-285}) gives the following uncoupled system of ODEs for the unknown functions $N_k(t)$

\begin{gather}\label{eq11-289}
N'_k(t)+k^2N_k(t)=\hat{\lambda}N_k(t),\\
\Updownarrow\nonumber\\
N'_k(t)=(\hat{\lambda}-k^2)N_k(t),\nonumber\\
\Updownarrow\nonumber\\
N_k(t)=a_ke^{(\hat{\lambda}-k^2)t}.\nonumber
\end{gather}
The general solution to (\ref{eq11-281a}) can thus be written as 

\begin{equation}\label{eq11-290}
w(x,t)=\sum\limits_{k=1}^{\infty}w_k(x,t),
\end{equation} 
where the \emph{normal modes}, $w_k(x,t)$, are 

\begin{equation}\label{eq11-291}
w_k(x,t)=a_ke^{(\hat{\lambda}-k^2)t}\sin kx,\quad k=1,2,\dots 
\end{equation}
We see that 

\begin{enumerate}
	\item $\hat{\lambda}<1\,\Rightarrow$ all normal modes decay $\Rightarrow\,u_0(x,t)=0$ is stable. 
	\item $\hat{\lambda}>1\,\Rightarrow$ at least one normal mode grows exponentially $\Rightarrow\,u_0(x,t)$ is unstable.  
\end{enumerate}
For the infinite domain case the stability condition was $\hat{\lambda}<0$. Thus the stability condition is not the same for the case of infinite and finite domain. This is in general true. 

For the unstable case, $\hat{\lambda}>1$, one or more normal modes experience exponential growth, which after a finite time invalidate the linearized equation (\ref{eq11-281a}). What is the system going to do after this time? This is a questing of \emph{nonlinear stability}. 

Nonlinear stability is a much harder problem than linear stability but powerful techniques exist in the form of asymptotic expansions which in many important cases can provide answer to the question of nonlinear stability. 

Here we will not dive into this more sophisticated approach but will rather introduce an approach based on the finite Fourier transform.

  Observe that the finite Fourier transform can also be applied to the nonlinear equation

\begin{equation}\label{eq11-292}
w_t-w_{xx}=\hat{\lambda}w(1-\epsilon^2w^2)=\hat{\lambda}w-\hat{\lambda}\epsilon^2w^3.
\end{equation}  
Representing the solution using the Finite Fourier transform we have

\begin{equation}\label{eq11-293}
w(x,t)=\sum\limits_{k=1}^{\infty}N_k(t)M_k(x).
\end{equation}
The evolution equations for the $N_k(t)$s are found as for the linearized case. We have 

\begin{gather}\label{eq11-294}
\int_0^{\pi}\mathrm d x\,w_t(x,t)M_k(x)- \int_0^{\pi}\mathrm d xw_{xx}(x,t)M_k(x)\\
=\int_0^{\pi}\mathrm d x\,\hat{\lambda}w(x,t)M_k(x)- \int_0^{\pi}\mathrm d x\,\hat{\lambda}\epsilon^2w^3(x,t)M_k(x).\nonumber
\end{gather}
The only new term as compared with the linearized case is the second term on the right hand side of (\ref{eq11-294}). For this term we have

\begin{align}\label{eq11-295}
&\int_0^{\pi}\mathrm d x\,\hat{\lambda}\epsilon^2w^3(x,t)M_k(x)\\
&= \hat{\lambda}\epsilon^2\sum^{\infty}_{i,j,l=1}\int^{\pi}_0\mathrm d x\,N_i(t)N_j(t)N_l(t)M_i(x)M_j(x)M_l(x)M_k(x)\nonumber\\
&=\hat{\lambda}\epsilon^2\sum^{\infty}_{i,j,l=1}a^{ijl}_{k}N_i(t)N_j(t)N_l(t),\nonumber
\end{align}
where 

\begin{equation}\label{eq11-296}
a^{ijl}_{k}=\int_0^{\pi}\mathrm d x\, M_i(x)M_j(x)M_k(x)M_l(x).
\end{equation}
Using (\ref{eq11-286}), (\ref{eq11-287}), (\ref{eq11-288}) and (\ref{eq11-295}) in (\ref{eq11-294}) gives us 

\begin{equation}\label{eq11-297}
N_k'(t)+(k^2-\hat{\lambda})N_k(t)=-\hat{\lambda}\epsilon^2\sum\limits_{i,j,l=1}^{\infty}a^{ijl}_{k}N_i(t)N_j(t)N_l(t).
\end{equation}
(\ref{eq11-297}) is an infinite system of nonlinear coupled ODEs. This is evidently a very hard problem to solve in any exact sense. (\ref{eq11-297}) shows the complexity of nonlinear PDEs. The very compact notation from calculus makes the nonlinear PDE (\ref{eq11-292}) to \emph{appear} simple, but as (\ref{eq11-297}) shows the problem of solving (\ref{eq11-292}) is actually horrendously difficult.

In order to get anywhere, we must restrict to situations where some of the complexity of (\ref{eq11-297})  is reduced. Finding such situations is what good applied mathematicians are trained to do. 

Let us assume that $\hat{\lambda}>1$ so we are in the unstable regime, but that $\hat{\lambda}<4$ so that only the normal mode $u_1(x,t)$ is unstable. All the other ones decay exponentially. For the $N_k(t)$s this means that 

\begin{equation}\label{eq11-298}
N_1(t)>>N_k(t),\quad\forall k\geq 2.
\end{equation}
Thus, on the right hand side of (\ref{eq11-297}) the term

\begin{equation}\label{eq11-299}
-\hat{\lambda}\epsilon^2a^{111}_{1}N_1^3(t),
\end{equation}
is much larger than the other terms. A good approximation to 
(\ref{eq11-297}) therefore appears to be

\begin{equation}\label{eq11-300} 
N_1'(t)+(1-\hat{\lambda})N_1(t)+\hat{\lambda}\epsilon^2a^{111}_{1}N_1^3(t)=0.
\end{equation}
The coefficient $a^{111}_{1}$ can be calculated using (\ref{eq11-296})

\begin{equation}\label{eq11-301} 
a^{111}_{1}=\int_0^{\pi}\mathrm d x\,(M_1(x))^4=\frac{4}{\pi}\int_0^{\pi}\mathrm d x \sin^4x=\frac{3}{2\pi}.
\end{equation}

\begin{equation}\label{eq11-302} 
\Rightarrow N_1'(t)+(1-\hat{\lambda})N_1(t)+\frac{3\hat{\lambda}}{2\pi}\epsilon^2N_1^3(t)=0
\end{equation}
If we multiply this equation by $N_1(t)$ we obtain

\begin{equation}\label{eq11-303} 
\frac{1}{2}y+(1-\hat{\lambda})y+\frac{3\hat{\lambda}\epsilon^2}{2\pi}y^2=0,
\end{equation}
where $y(t)=N_1^2(t)$. (\ref{eq11-303}) is a so-called Ricatti equation and can be solved exactly.
There is a general way to solve such equations. 

A general Ricatti equation is of the form

\begin{equation}\label{eq11-304}
y'=q_0(x)+q_1(x)y+q_2(x)y^2.
\end{equation} 
First introduce a function $v(x)$ by 

\begin{equation}\label{eq11-305} 
v(x)=q_2(x)y(x).
\end{equation}
Then we get the following equation for $v(x)$

\begin{equation}\label{eq11-306} 
v'=v^2+R(x)v+S(x),
\end{equation}
where 

\begin{align}\label{eq11-307} 
S(x)&=q_2(x)q_0(x),\\
R(x)&=q_1(x)+\frac{q'_2(x)}{q_2(x)}.\nonumber
\end{align}
Next, introduce the function $u(x)$ by 

\begin{equation}\label{eq11-308}
v=-\frac{u'}{u}.
\end{equation}
Then we get the following \emph{linear} equation for $u$

\begin{equation}\label{eq11-309} 
u''(x)-R(x)u'(x)+S(x)u=0.
\end{equation}
Using (\ref{eq11-308}) and (\ref{eq11-305}), any solution of (\ref{eq11-309}) gives a solution of the Ricatti equation (\ref{eq11-304}). Since the general solution to (\ref{eq11-309}) can be found explicitly in some cases, the general solution (\ref{eq11-304}) can also be found for those cases. 

For our equation (\ref{eq11-303}) we have

\begin{align}\label{eq11-310} 
q_0(x)&=0\\
q_1(x)&=2(\hat{\lambda}-1)\nonumber\\
q_2(x)&=-\frac{3\hat{\lambda}\epsilon^2}{\pi}\nonumber
\end{align}
We can find the general solution to this equation which leads to the general solution (\ref{eq11-303}). With initial condition given as in (\ref{eq11-261}) we get

\begin{equation}\label{eq11-311} 
N_1(t)=\frac{(h,M_1)e^{(\hat{\lambda}-1)t}}
{\left[1+\left[3\hat{\lambda}\epsilon^2(h,M_1)^2/(2\pi(\hat{\lambda}-1))\right]\left(e^{2(\hat{\lambda}-1)t}-1\right)\right]^{\frac{1}{2}}}
\end{equation}
This solution tends to a finite stationary value as $t\to\infty$

\begin{equation}\label{eq11-312} 
\lim\limits_{t\to\infty}N_1(t)=\frac{(h,M_1)}{|(h,M_1)|}\sqrt{\frac{2\pi}{3}}\left(\frac{\hat{\lambda}-1}{\hat{\lambda}\epsilon^2}\right)^{\frac{1}{2}}
\end{equation}
rather than growing unbounded as the linearized equation predicts. 

Note that the current method can be extended to study more general nonlinear terms $\hat{\lambda}f(u)$. 

\section{Integral transforms}
Integral transforms are used to solve PDEs on unbounded domains. They have very wide applicability and play a fundamental role in applied mathematics and theoretical physics. 

In these notes I will concentrate on equations of the form 

\begin{equation}\label{eq12-1}
\rho(\mathbf x)K u(\mathbf x,t)+Lu(\mathbf x,t)=\rho(\mathbf x)F(\mathbf x,t),
\end{equation}
where 

\begin{equation}\label{eq12-2}
Lu(\mathbf x,t)=-\nabla\cdot(p\nabla u(\mathbf x,t))+q(\mathbf x)u(\mathbf x,t),
\end{equation}
and where 

\begin{equation}\label{eq12-3}
K=
\left\{
\begin{aligned}
\partial_{tt}\quad& \text{ -- hyperbolic case,}\\
\partial_t\quad&\text{ -- parabolic case,}\\
\partial_{yy}\quad&\text{ -- elliptic case, cylinder symmetry,}\\
0\quad&\text{ -- general elliptic case.}\\
\end{aligned}
\right.
\end{equation}
Recall that the basic idea in the separation of variables is to look for special solution to (\ref{eq12-1}) of the form

\begin{equation}\label{eq12-4}
u(\mathbf x,t)=N(t)M(\mathbf x),
\end{equation}
($t\to y$ in the restricted elliptic case). We found that a separation constant $\lambda$ appeared and that (\ref{eq12-4}) is a solution to (\ref{eq12-1}) only if

\begin{equation}\label{eq12-5}
LM(\mathbf x)=\lambda^2\rho(\mathbf x)M(\mathbf x),
\end{equation}
and 

\begin{equation}\label{eq12-6}
\left\{
\begin{aligned}
N''(t)+\lambda^2N(t)=0\quad& \text{ -- hyperbolic case,}\\
N'(t)+\lambda^2N(t)=0\quad&\text{ -- parabolic case,}\\
N'(y)-\lambda^2N(y)=0\quad&\text{ -- restricted elliptic case.}\\
\end{aligned}
\right.
\end{equation}
Here we have written the separation constant as $\lambda^2$ for convenience. Equation (\ref{eq12-5}) must be supplied with boundary conditions. They must of course be inherited from the boundary conditions on equation (\ref{eq12-1}). The most common ones are that the solution is bounded or that it vanishes as $\|x\|\to\infty$. Thus we will require that $M(\mathbf x)$ is bounded or vanishes at infinity. 

As for the case of bounded domains, (\ref{eq12-5}) with boundary conditions will not have solutions for all $\lambda$. The $\lambda$s that give a solution are part of the \emph{spectrum} of $L$. Things look very similar to the case of bounded domains, so far. There is however a big difference as compared with the bounded case.

	In order to illustrate this difference, let us consider the boundary value problem
	
	\begin{align}\label{eq12-7}
	&-M''(x)=\lambda^2M(\mathbf x)\quad -l<x<l,\\
	& M(l)=M(-l)=0.\nonumber
	\end{align}
	The general solution to the equation is for $\lambda>0$
	
	\begin{equation}\label{eq12-8}
	M(x)=Ae^{i\lambda x}+ Be^{-i\lambda x}.
	\end{equation}
	The boundary conditions imply that 
	
	\begin{align}\label{eq12-9}
	& Ae^{-i\lambda l}+ Be^{i\lambda l}=0,\\
	& Ae^{i\lambda l}+ Be^{-i\lambda l}=0,\nonumber
	\end{align}
	
	\begin{gather}\label{eq12-10}
	\begin{bmatrix}
	e^{-il\lambda} & e^{il\lambda}\\
	e^{il\lambda} & e^{-il\lambda}
	\end{bmatrix} 
	\begin{bmatrix}
	A\\
	B
	\end{bmatrix}
	=0.
	\end{gather}

	\noindent In order for (\ref{eq12-10}) to have a nontrivial solution, and thus for $\lambda$ to be in spectrum of $L=-\partial_{xx}$, the matrix in (\ref{eq12-10}) must have zero determinant
	
	\begin{gather}
	\det
	\begin{bmatrix}
	e^{-il\lambda} & e^{il\lambda}\\
	e^{il\lambda} & e^{-il\lambda}
	\end{bmatrix} =0,\label{eq12-11}\\
	\Updownarrow\nonumber\\
	e^{4il\lambda}=1,\nonumber\\
	\Updownarrow\nonumber\\
	4l\lambda_k=2\pi k,\nonumber\\
	\lambda_k=\frac{\pi k}{2l},\quad k=1,2,\dots\;\;.\label{eq12-12}
	\end{gather}
	There are no eigenvalues for $\lambda\leq 0$. Thus (\ref{eq12-12}) defines the spectrum of $l$.
	
	  Observe that the distance $\Delta \lambda$ between two spectral points is 
	
	\begin{equation}\label{eq12-13}
	\Delta\lambda=\frac{\pi}{2l}.
	\end{equation}
	
	\begin{figure}[h!]
		\centering
		\includegraphics[scale=0.5]{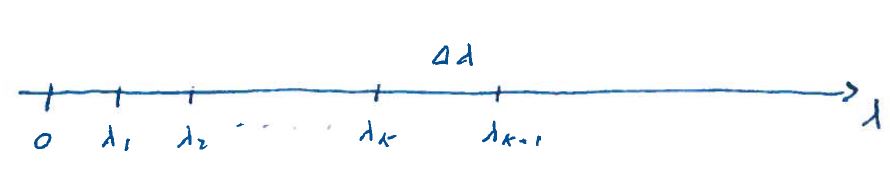}
		\caption{}
	\end{figure}
         \noindent  When the size of the domain increases without bound we observe that the distance between spectral points go to zero 
	
	\begin{equation}\label{eq12-14}
	\lim\limits_{l\to\infty}\Delta\lambda=0.
	\end{equation}
	Thus in this limit the spectrum of $L$ will form a continuum
	
	\begin{equation}\label{eq12-15}
	0<\lambda<\infty.
	\end{equation}

\noindent This is in general the case for unbounded domains; the spectrum of the operators form a continuum. 

For the discrete case, with eigenvalues $\lbrace \lambda_k\rbrace_{k=1}^{\infty}$ and corresponding eigenfunctions $\lbrace M_k\rbrace_{k=1}^{\infty}$, the basic idea was to look for solutions to the PDE in the form of an infinite sum

\begin{equation}\label{eq12-16}
u=\sum\limits_{k=1}^{\infty}N_kM_k.
\end{equation}
For the unbounded case we have eigenfunctions $M_{\lambda}$ and corresponding $N_{\lambda}$, where the $\lambda$s form a continuum that we denote by $D$. By analogy with (\ref{eq12-10}) we seek solutions to the PDE in the form of "sums over the continuum". Such sums are of course integrals. Thus

\begin{equation}\label{eq12-17}
u=\int\limits_D\mathrm d \lambda\,N_{\lambda}M_{\lambda}.
\end{equation}
We recall that in the case of finite domains, the fitting of initial conditions to our solution consisted in finding constants $\lbrace a_k\rbrace_{k=1}^{\infty}$ such that 

\begin{equation}\label{eq12-18}
f(\mathbf x)=\sum\limits_{k=1}^{\infty}a_kM_k(\mathbf x).
\end{equation}
Using orthonormality of the eigenfunctions we found that

\begin{equation}\label{eq12-19}
a_k=(M_k,f)=\int_G \mathrm d V\,\rho(\mathbf x)f(\mathbf x)M_k(\mathbf x).
\end{equation}
Fitting initial conditions for the case of unbounded domains will by analogy involve finding coefficients $a_{\lambda}$, $\lambda\in D$ such that

\begin{equation}\label{eq12-20}
f(\mathbf x)=\int\limits_D\mathrm d \lambda \,a_{\lambda} M_{\lambda}(\mathbf x).
\end{equation} 
In order to determine $a_\lambda$ we need a continuum analog of the orthonormality condition from the discrete case
\begin{equation}
(M_k,M_k')=\delta_{kk'},
\end{equation}
where we recall that $\delta_{kk'}$ is the Kronecker delta.

If we had something like a continuum version of the Kronecker delta one would guess at

\begin{equation}\label{eq12-22}
(M_{\lambda},M_{\lambda'})=\delta_{\lambda\lambda'},
\end{equation} 
This is in fact what \emph{Paul Dirac} did early in the previous century. Dirac was an English genius that was one of the founding fathers of quantum theory. He came to this idea while trying to solve the Schrodinger equation from quantum mechanics. He found that $\delta_{\lambda\lambda''}$ only depended on $\lambda-\lambda'$ and he wrote it as 

\begin{equation}\label{eq12-23}
\delta_{\lambda\lambda'}=\delta(\lambda-\lambda').
\end{equation}
$\delta(\lambda)$ is called the \emph{Dirac delta function} after him. In order to do its intended job he found that it would have to have the following two properties

\begin{equation}\label{eq12-24}
\delta(\lambda)=0,\quad\lambda\neq 0,
\end{equation}

\begin{equation}\label{eq12-25}
\int\limits_D\mathrm d \lambda \,\delta(\lambda)=1.
\end{equation}
It is a trivial exercise to prove that no function exists that have the two properties (\ref{eq12-24}) and \ref{eq12-25}). Mathematicians were quick to point this out but Dirac didn't care, he did not know much mathematics, and he did not have to know it either, because he was in fact a genius and people like him move to a different drum. 

Eventually mathematicians did make sense of the Dirac delta function and in fact discovered a whole new mathematical universe called \emph{the theory of distributions}. We will discuss elements of this theory next semester. In this semester we will use objects like the Dirac delta sparingly and when we do only in a purely formal way like Dirac himself. 

There is an important moral to the story of the Dirac delta. Pure mathematics is important because it gives us abstract frameworks to help organize our mathematical methods. Pure mathematics can also give us more efficient ways of doing our calculations and make precise under which conditions our manipulations will keep us within the boundaries of the known mathematical universe. 

However the framework and bounds imposed by pure mathematics should not stop us from sometimes disregarding the rules and play with the formulas, like Dirac did. In the final analysis, virtually all of mathematics originate from such rule-defying play. 

It is however also appropriate to give a warning; disregarding the rules and playing with the formulas should not be a cover for sloppy thinking. This is unfortunately the case quite frequently in applied science. 

We can't all be a Dirac. Disregarding the rules imposed by pure mathematics will, by most people, most of the time, produce useless nonsense.

In these notes I will sometimes break the rules a la Dirac and I will try to remember telling you when I do. All my rule breaking can be made sense of within the universe of distributions. 

All finite Fourier transforms will in the limit of unbounded domains give rise to integral transforms. There is thus a large set of such transforms available, each tailored to a particular class of problems. In these notes I will focus on the Sine-, Cosine-, Fourier- and Laplace transforms. 

\subsection{The 1D Fourier transform}
Let us consider the eigenvalue problem  

\begin{equation}\label{eq12-26}
M''(x)+\lambda^2M(x)=0,\quad -\infty<x<+\infty.
\end{equation}
The boundary condition is that $M(x)$ is bounded when $|x|\to\infty$. The general solution of the equation (\ref{eq12-26}) is 

\begin{equation}\label{eq12-27}
M_{\lambda}(x)=a(\lambda)e^{i\lambda x}+b(\lambda)e^{-i\lambda x},
\end{equation}
where $\lambda\in\mathbb{C}$. Writing $\lambda=\alpha+i\beta$ we have 

\begin{gather*}
e^{\pm i\lambda x}=e^{\mp \beta x}e^{\pm i\alpha x},\\
\Downarrow\\
 |e^{\pm i\lambda x}|=e^{\mp \beta x},
\end{gather*}
and the functions $e^{\mp \beta x}$ are unbounded on $-\infty<x<\infty$ as long as $\beta\neq 0$. Thus to satisfy the boundary conditions we must restrict $\lambda$ to the real axis $-\infty<x<\infty$. On the real axis impose the conditions 

\begin{equation}\label{eq12-28}
\left.
\begin{aligned}
|a(\lambda)|<M_a,\\
|b(\lambda)|<M_b,
\end{aligned}
\right.\quad \forall\lambda,
\end{equation}
then the boundedness of $M_{\lambda}$ follows from the triangle inequality

\begin{align}\label{eq12-29}
|M_{\lambda}(x)|&=|a(\lambda)e^{i\lambda x}+b(\lambda)e^{-i\lambda x}|\\
&\leq |a(\lambda)||e^{i\lambda x}|+|b(\lambda)||e^{-i\lambda x}|\nonumber\\
&=|a(\lambda)|+|b(\lambda)|\leq M_a+M_b.\nonumber
\end{align}
Therefore the spectrum in this case consists of the whole real axis. We can now use these eigenfunctions to represent functions $f(x)$ in the following way

\begin{align}\label{eq12-30}
f(x)&=\int^{\infty}_{-\infty}\mathrm d \lambda\,\lbrace a(\lambda)e^{i\lambda x}+b(\lambda)e^{-i\lambda x}\rbrace\\
&= \int^{\infty}_{-\infty}\mathrm d \lambda\, a(\lambda)e^{i\lambda x}+\int^{\infty}_{-\infty}\mathrm d \lambda\,b(\lambda)e^{-i\lambda x}\nonumber\\
&=\int^{\infty}_{-\infty}\mathrm d \lambda\, a(-\lambda)e^{-i\lambda x}+\int^{\infty}_{-\infty}\mathrm d \lambda\,b(\lambda)e^{-i\lambda x}\nonumber\\
&=\int^{\infty}_{-\infty}\mathrm d \lambda\,c(\lambda)e^{-i\lambda x},\nonumber
\end{align}
where $c(\lambda)=a(-\lambda)+b(\lambda)$. Note that in (\ref{eq12-30}), $f(x)$ is represented by $e^{i\lambda x}$ in terms of a standard improper integral. In the theory of integral transforms, other, less standard, improper integrals appear. We will discuss these when they arise. 

The formula (\ref{eq12-30}) is by definition the inverse Fourier transform. The key question is how to find the $c(\lambda)$ corresponding to a given function $f(x)$. 

For this purpose, let us introduce a small number $\epsilon>0$ and define

\begin{equation}\label{eq12-31}
\epsilon(x)=
\left\{
\begin{aligned}
-i\epsilon,\quad x>0,\\
i\epsilon,\quad x<0.
\end{aligned}
\right.
\end{equation}
Observe that for all real $\lambda$ we have 

\begin{align}\label{eq12-32}
&\int^{\infty}_{-\infty}\mathrm d x\, e^{-i(\lambda+\epsilon(x))x}=\lim\limits_{L\to\infty}\int_{-L}^{L}\mathrm d x\,e^{-i(\lambda+\epsilon(x))x}\\
&=\lim\limits_{L\to\infty}\left\{\int_{-L}^{0}\mathrm d x\,e^{-i(\lambda+i\epsilon)x}+\int_{0}^{L}\mathrm d x\,e^{-i(\lambda-i\epsilon)x}\right\}\\
& =\lim\limits_{L\to\infty}\left\{\frac{1}{-i(\lambda+i\epsilon)}e^{-i(\lambda+i\epsilon)x}\Big|_{-L}^0+\frac{1}{-i(\lambda-i\epsilon)}e^{-i(\lambda-i\epsilon)x}\Big|_{0}^L\right\} \\
&=\lim\limits_{L\to\infty}\left\{ \frac{1}{-i(\lambda+i\epsilon)}+\frac{1}{i(\lambda-i\epsilon)}\right.\nonumber\\
&\quad+\left.\frac{1}{i(\lambda+i\epsilon)}e^{-\epsilon L}e^{i\lambda L}-\frac{1}{i(\lambda-i\epsilon)}e^{-\epsilon L}e^{-i\lambda L}\right\}\nonumber\\
&= \frac{1}{i(\lambda-i\epsilon)}-\frac{1}{-
	i(\lambda+i\epsilon)}=\frac{2\epsilon}{\epsilon^2+\lambda^2}.
\end{align}
Thus
\begin{equation}\label{eq12-33}
\lim\limits_{\epsilon\to 0}\int^{\infty}_{-\infty}\mathrm d x\, e^{-i(\lambda+\epsilon(x))x}=
\left\{
\begin{aligned}
0,\quad\lambda\neq 0\\
\infty,\quad\lambda=0
\end{aligned}
\right.\;.
\end{equation}

This result makes the following formula plausible (I am playing with formulas now)

\begin{equation}\label{eq12-34}
\int^{\infty}_{-\infty}\mathrm d x\, e^{-\lambda x}\propto\delta(\lambda).
\end{equation}
In the theory of distributions it becomes clear that the constant of proportionality is $2\pi$. So we have the formula 

\begin{equation}\label{eq12-35}
\delta(\lambda)=\frac{1}{2\pi}\int^{\infty}_{-\infty}\mathrm d x\, e^{-\lambda x}
\end{equation}
Since the Dirac delta function is zero for $\lambda\neq 0$ it appears that for all reasonable functions, $f(\lambda)$, we have 

\begin{equation}\label{eq12-36}
\int^{\infty}_{-\infty}\mathrm d \lambda\,f(\lambda)\delta(\lambda-\hat{\lambda})=f(\hat{\lambda}).
\end{equation}
(I am playing now.)

Let us now return to our original question of how to find the $c(\lambda)$  corresponding to a given function $f(x)$ in formula (\ref{eq12-30})

\begin{equation}\label{eq12-37}
f(x)=\int^{\infty}_{-\infty}\mathrm d \lambda\,c(\lambda) e^{-i\lambda x}.
\end{equation}
Multiplying (\ref{eq12-37}) by $e^{i\hat{\lambda}x}$ and integrate with respect to $x$ over $-\infty<x<\infty$ gives

\begin{align}\label{eq12-38a}
&\int^{\infty}_{-\infty}\mathrm d x\,f(x)e^{i\hat{\lambda}x}=\int^{\infty}_{-\infty}\mathrm d x\,e^{i\hat{\lambda}x}\int^{\infty}_{-\infty}\mathrm d \lambda\,c(\lambda)e^{-i\lambda x}\\
&=\int^{\infty}_{-\infty}\mathrm d \lambda\,c(\lambda) \int^{\infty}_{-\infty}\mathrm d x\, e^{-i(\lambda-\hat{\lambda})x}\nonumber\\
&=\int^{\infty}_{-\infty}\mathrm d\lambda\,c(\lambda)2\pi\delta(\lambda-\hat{\lambda})=2\pi c(\hat{\lambda}).
\end{align}
So the formula we are looking for is 

\begin{equation}\label{eq12-38b}
c(\lambda)=\frac{1}{2\pi}\int^{\infty}_{-\infty}\mathrm d x\,f(x)e^{i\lambda x}.
\end{equation}
This formula defines the Fourier transform. Note that there are several conventions about where to place the factor of $2\pi$. A more symmetric choice\cite{textbook}, is

\begin{equation}\label{eq12-39}
F(\lambda)=\frac{1}{\sqrt{2\pi}}\int^{\infty}_{-\infty}\mathrm d x\,f(x)e^{i\lambda x},
\end{equation}

\begin{equation}\label{eq12-40}
f(x)=\frac{1}{\sqrt{2\pi}}\int^{\infty}_{-\infty}\mathrm d x\,F(x)e^{-i\lambda x}.
\end{equation}
Note that not only are there different conventions on where to place the factor of $2\pi$ but also in which exponent to place the minus sign.

In this class (\ref{eq12-39}) defines the Fourier transform and (\ref{eq12-40}) defines the inverse Fourier transform. 

Pure mathematicians have made extensive investigations into for which pair of functions $f(x)$ and $F(\lambda)$ the formulas (\ref{eq12-39}) and (\ref{eq12-40}) are valid. For example if $f(x)$ is piecewise continuously differentiable  on each finite interval and if additionally 

\begin{equation}\label{eq12-41}
\int^{\infty}_{-\infty}\mathrm d x\,|f(x)|<\infty,
\end{equation}
formula (\ref{eq12-39}) and (\ref{eq12-40}) holds pointwise. This result and related ones are important but we will use (\ref{eq12-39}) and (\ref{eq12-40}) in a much wider context than this.

  Let $f(x)=\delta(x)$. Then (\ref{eq12-39}) gives 

\begin{equation}\label{eq12-42}
F(\lambda)=\frac{1}{\sqrt{2\pi}}\int^{\infty}_{-\infty}\mathrm d x\,\delta(x)e^{i\lambda x}=\frac{1}{\sqrt{2\pi}},
\end{equation}
and (\ref{eq12-40}) then gives

\begin{equation}\label{eq12-43}
f(x)=\frac{1}{\sqrt{2\pi}}\int^{\infty}_{-\infty}\mathrm d \lambda\,\frac{1}{\sqrt{2\pi}}e^{-i\lambda x}=\frac{1}{2\pi}\int^{\infty}_{-\infty}\mathrm d \lambda\,e^{-i\lambda x}=\delta(x),
\end{equation} 
where we have used the formula (\ref{eq12-35})  which we argued for previously. More generally, the pair of formulas (\ref{eq12-39}) and (\ref{eq12-40}) holds for any distribution, not just the Delta distribution. More about this next semester. In this semester we will basically use (\ref{eq12-39}) and (\ref{eq12-40}) to take Fourier transform and inverse Fourier transform of anything we want! 

The Fourier transform has many powerful properties. You can find them in handbooks or on the web whenever you need them. For example if $F(\lambda)$ and $G(\lambda)$ are Fourier transforms of $f(x)$ and $g(x)$ then $aF(\lambda)+bG(\lambda)$ is the Fourier transform of $af(x)+bg(x)$. A less trivial, and extremely useful for us, property involves products of Fourier transforms $F(\lambda)G(\lambda)$ of functions $f(x)$ and $g(x)$. Such products often occur when we use Fourier transforms to solve ODEs and PDEs.

Playing with formulas we get the following

\begin{align}\label{eq12-44b}
\frac{1}{\sqrt{2\pi}}\int^{\infty}_{-\infty}\mathrm d\lambda\,F(\lambda)G(\lambda)e^{-i\lambda x}&=\frac{1}{2\pi} \int^{\infty}_{-\infty}\mathrm d\lambda\,G(\lambda)e^{-i\lambda x}\int^{\infty}_{-\infty}\mathrm d t\,f(t)e^{i\lambda t}\nonumber\\
& =\frac{1}{2\pi}\int^{\infty}_{-\infty}\mathrm d t\,f(t)\int^{\infty}_{-\infty}\mathrm d\lambda\,G(\lambda)e^{-i\lambda (x-t)}\nonumber\\
& = \frac{1}{\sqrt{2\pi}}\int^{\infty}_{-\infty}\mathrm d t\,f(t)g(x-t).
\end{align}
The integral on the right hand side of (\ref{eq12-44b}) is known as a \emph{convolution integral}. Canceling the factor $\frac{1}{\sqrt{2\pi}}$ we get the \emph{convolution integral theorem} for Fourier transforms 

\begin{equation}\label{eq12-45}
\int^{\infty}_{-\infty}\mathrm d t\,f(t)g(x-t)=\int^{\infty}_{-\infty}\mathrm d\lambda \,F(\lambda)G(\lambda)e^{-i\lambda x}.
\end{equation}
Evaluating at $x=0$ we get 

\begin{equation}\label{eq12-46}
\int^{\infty}_{-\infty}\mathrm d t\,f(t)g(-t)=\int^{\infty}_{-\infty}\mathrm d\lambda\, F(\lambda)G(\lambda).
\end{equation}
Now observe that

\begin{gather}\label{eq12-47}
F(\lambda)=\frac{1}{\sqrt{2\pi}}\int^{\infty}_{-\infty}\mathrm d t\,f(t)e^{i\lambda t},\nonumber\\
\Downarrow\nonumber\\
\overline{F}(\lambda)=\frac{1}{\sqrt{2\pi}}\int^{\infty}_{-\infty}\mathrm d t\,\overline{f}(t)e^{-i\lambda t}=\frac{1}{\sqrt{2\pi}}\int^{\infty}_{-\infty}\mathrm d t\,\overline{f}(-t)e^{i\lambda t},
\end{gather}
where $\overline{F}(\lambda)$ is the complex conjugate of $F(\lambda)$.

Thus the fourier transform of $\overline{f}(-t)$ is $\overline{F}(\lambda)$. Putting $g(t)=\overline{f}(-t)$ into (\ref{eq12-46}) we get 

\begin{equation}\label{eq12-48}
\int^{\infty}_{-\infty}\mathrm d t\,|f(t)|^2=\int^{\infty}_{-\infty}\mathrm d t\,f(t)\overline{f}(t)=\int^{\infty}_{-\infty}\mathrm d t\,f(t)g(-t)=\int^{\infty}_{-\infty}\mathrm d \lambda\,|F(\lambda)|^2.
\end{equation}
The identity 

\begin{equation}\label{eq12-49}
\int^{\infty}_{-\infty}\mathrm d t\, |f(t)|^2=\int^{\infty}_{-\infty}\mathrm d \lambda\,|F(\lambda)|^2.
\end{equation}
is called the Parseval equation and is clearly a generalization of the Parseval equation for generalized Fourier series to continuous  spectrum. 

Another important property of Fourier transforms involve derivatives. Observe that 

\begin{align}\label{eq12-50}
\frac{1}{\sqrt{2\pi}}\int^{\infty}_{-\infty}\mathrm d xf'(x)e^{i\lambda x}=\frac{1}{\sqrt{2\pi}}f(x)e^{i\lambda x}\bigg|^{\infty}_{-\infty}-\frac{1}{\sqrt{2\pi}}\int^{\infty}_{-\infty}\mathrm d xf(x)(i\lambda)e^{i\lambda x}\nonumber\\
= -i\lambda \frac{1}{\sqrt{2\pi}} \int^{\infty}_{-\infty}\mathrm d x f(x)e^{i\lambda x}=-i\lambda F(\lambda).
\end{align}
Thus, derivatives are converted into multiplication by the Fourier transform. (\ref{eq12-50}) can be easily generalized to

\begin{equation}\label{eq12-52}
\frac{1}{\sqrt{2\pi}}\int^{\infty}_{-\infty}\mathrm d x\, f^{(n)}(x)e^{i\lambda x}=(-i\lambda)^nF(\lambda).
\end{equation}
Formulas (\ref{eq12-50}) and (\ref{eq12-52}) are the main sources for the great utility of Fourier transforms in the theory of PDEs and ODEs.

\subsubsection{A boundary value problem for an ODE}
Let us consider the following boundary value problem 

\begin{equation}\label{eq12-53}
y''(x)-k^2y(x)=-f(x),\quad-\infty<x<\infty,
\end{equation}
where $k$ is a fixed  constant and $f$ a given function. The boundary conditions attached to (\ref{eq12-53}) are 

\begin{equation}\label{eq12-54}
y(x),\,y'(x)\to 0,\text{ when }|x|\to\infty.
\end{equation}
Let $F(\lambda)$ and $Y(\lambda)$ be the Fourier transforms of $f(x)$ and $y(x)$. Then the transform of the equation (\ref{eq12-53}) is, upon using (\ref{eq12-52}), given by

\begin{equation}\label{eq12-55}
-(\lambda^2+k^2)Y(\lambda)=-F(\lambda)
\end{equation}
This is a simple algebraic equation whose solution is

\begin{equation}\label{eq12-56}
Y(\lambda)=\frac{F(\lambda)}{\lambda^2+k^2}.
\end{equation}
Using the inverse Fourier transform and the convolution theorem we get the following representation of $y(x)$:

\begin{align}\label{eq12-57}
y(x)&=\frac{1}{\sqrt{2\pi}} \int^{\infty}_{-\infty}\mathrm d \lambda\,\frac{1}{\lambda^2+k^2}F(\lambda)e^{-i\lambda x}\nonumber\\
& =\frac{1}{\sqrt{2\pi}} \int^{\infty}_{-\infty}\mathrm d t\,g(x-t)f(t),
\end{align}
where 

\begin{equation}\label{eq12-58}
g(\xi)=\frac{1}{\sqrt{2\pi}} \int^{\infty}_{-\infty}\mathrm d \lambda\frac{e^{-i\lambda\xi}}{\lambda^2+k^2}.
\end{equation}
Using a table of transforms or residue calculus from complex function theory we have 

\begin{equation}\label{eq12-59}
g(\xi)=\frac{\sqrt{2\pi}}{2k}e^{-k|\xi|},
\end{equation}
so that the solution $y(x)$ is 

\begin{equation}\label{eq12-60}
y(x)=\frac{1}{2k}\int^{\infty}_{-\infty}\mathrm d t \,e^{-k|x-t|}f(t).
\end{equation}
Because of the exponentially decaying factor, $y(x)$ and $y'(x)$ will both go to zero as $|x|\to\infty$ if $f(t)$ decays fast enough at infinity. 

Actually, computing $y(x)$ for a given $f(t)$ using formula (\ref{eq12-60}) is no simple matter; tables of integrals, complex function theory, numerical methods and asymptotic methods are all of help here. 

For some simple choices for $f(t)$ the integral in (\ref{eq12-60}) can be evaluated easily. Let us first choose $f(t)=1$. Formula (\ref{eq12-60}) implies 

\begin{align}\label{eq12-61}
y(x)&=\frac{1}{2k} \int^{\infty}_{-\infty}\mathrm d t \,e^{-k|x-t|}\\
&=\frac{1}{2k}\left[\int^{x}_{-\infty}\mathrm d t \,e^{-k(x-t)}+ \int_{x}^{\infty}\mathrm d t \,e^{k(x-t)}\right]\nonumber\\
&=\frac{1}{2k}\left[e^{-kx}\left(\frac{1}{k}e^{kt}\right)\Big|_{-\infty}^x+e^{kx} \left(-\frac{1}{k}e^{-kt}\right)\Big|_{x}^{\infty}\right]\nonumber\\
&=\frac{1}{2k}\left[\frac{1}{k}+\frac{1}{k}\right]=\frac{1}{k^2}.\nonumber
\end{align} 
Observe that 

\begin{equation}\label{eq12-62}
y''(x)-k^2y(x)=-k^2\left(\frac{1}{k^2}\right)=-1=-f(x),
\end{equation}
so $y(x)$ defined in (\ref{eq12-61}) is a solution to the equation (\ref{eq12-53}), but it does clearly not satisfy the boundary conditions (\ref{eq12-54}). The reason why this occurs is that in deriving the formula (\ref{eq12-60}) we assumed  that $f(x)$ had a Fourier transform $F(\lambda)$. However $f(x)=1$ leads to the divergent integral 

\begin{equation}\label{eq12-63}
F(\lambda)=\frac{1}{\sqrt{2\pi}} \int^{\infty}_{-\infty}\mathrm d x \,e^{i\lambda x}
\end{equation}
for $F(\lambda)$. Using (\ref{eq12-43}) we conclude that $F(\lambda)$ is not a function but a distribution

\begin{equation}\label{eq12-64}
F(\lambda)=\sqrt{2\pi}\delta(\lambda).
\end{equation}
Calculations like the one leading up to the "solution" (\ref{eq12-61}) are very common whenever transforms are used to solve PDEs. They are \emph{ideal} calculations that is meant to represent the \emph{limit} of a sequence of problems, each producing an honest to God solution to the boundary value problem (\ref{eq12-53}),(\ref{eq12-54}). Such ideal calculations are often possible to do in closed form and does often lead to some insight into the general case, but in what sense the result of the calculation is a limit of a sequence of actual solution is for the most part left unsaid. This will be clarified when we discuss distributions next semester. However, let us preview next semester by saying exactly what we mean by a limit here.

Define a sequence of functions $\lbrace f_N(x)\rbrace$

\begin{equation} \label{eq12-65}
f_N(x)=\left\{
\begin{aligned}
1,\quad |x|\leq N,\\
0,\quad |x|> N.
\end{aligned}
\right.
\end{equation}

\begin{figure}[h!]
	\centering
	\includegraphics[scale=0.5]{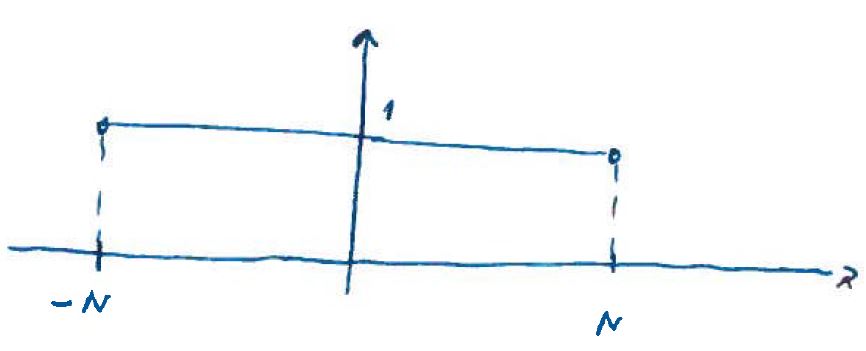}
	\caption{}
\end{figure}
This function has a Fourier transform for each $N$ and our formula (\ref{eq12-60}) should give a solution to the boundary value problem (\ref{eq12-53}),(\ref{eq12-54})

\begin{align} \label{eq12-66}
y_N(x)&=\frac{1}{2k}\int^{\infty}_{-\infty}\mathrm d t \,e^{-k|x-t|}f_N(t)\\
&=\frac{1}{2k}\int^N_{-N}\mathrm d t\,e^{-k|x-t|}.\nonumber
\end{align}
Let us first assume that $|x|\leq N$

\begin{align}\label{eq12-67}
\Rightarrow y_N(x)&=\frac{1}{2k}\left[ \int^{x}_{-N}\mathrm d t \,e^{-k(x-t)}+ \int_{x}^{N}\mathrm d t \,e^{k(x-t)}\right]\\
&=\frac{1}{2k}\left[e^{-kx}\left(\frac{1}{k}e^{kt}\right)\Big|_{-N}^x+e^{kx} \left(-\frac{1}{k}e^{-kt}\right)\Big|_x^N\right]\nonumber\\
&=\frac{1}{2k^2}\left[e^{-kx}\left(e^{kx}-e^{-kN}\right)-e^{kx}\left(e^{-kN}-e^{-kx}\right)\right]\nonumber\\
&=\frac{1}{2k^2}\left[1-e^{-kN}e^{-kx}-e^{-kN}e^{kx}+1\right]\nonumber\\
&=\frac{1}{k^2}\left[1-e^{-kN}\cosh kx\right]. 
\end{align}
The calculations for $|x|>N$ is completed in a similar way. This gives us 

\begin{equation} \label{eq12-68}
y_N(x)=\left\{
\begin{aligned}
\frac{1}{k^2}\left(1-e^{-kN}\cosh kx\right),\quad |x|\leq N,\\
\frac{1}{k^2}\left(e^{-k|x|}\sinh kN\right),\quad |x|\geq N.
\end{aligned}
\right.
\end{equation}
$y_N(x)$ is a solution to (\ref{eq12-53}) and it is easy to verify that $y_N(x)$, $y'_N(x)\to 0$ when $|x|\to\infty$, so we have an infinite sequence of solutions $\lbrace y_N(x)\rbrace$ to the boundary value problem defined by (\ref{eq12-53}),(\ref{eq12-54}). 

For any $x$ there evidently exists a $N$ so large that $N\geq|x|$. Therefore for such an $x$ we have $f_N(x)=f(x)$ for all $N>|x|$ and consequently the sequence of functions $\{f_N\}_N^\infty$ converge pointwise to the functions $f(x)=1$. This convergence is however \emph{not} uniform.

Observe that for any $x$ there obviously exists inifinitely many $N$ such that $N\geq|x| $. For such $N$ we have from (\ref{eq12-68}) that 

\begin{equation}\label{eq12-69}
y_N(x)=\frac{1}{k^2}\left(1-e^{-kN}\cosh kx\right)\to\frac{1}{k^2}=y(x),
\end{equation}
when $N\to\infty$

Thus, this convergence, like the one for the source $f(x)$,  is also pointwise but not uniform. This non-uniformity of the limit is why the limiting function, $y(x)$, does not satisfy the boundary condition (\ref{eq12-54}) even though all $y_N(x)$ does. 

As a second example of an ideal calculation let us choose $f(t)=\delta(t)$. Then formula (\ref{eq12-60}) gives 

\begin{equation}\label{eq12-71}
y(x)=\frac{1}{2k}\int^{\infty}_{-\infty}\mathrm d t\,e^{-k|x-t|}\delta(t)=\frac{1}{2k}e^{-k|x|}.
\end{equation}
This function satisfies the boundary condition and the equation for all $x\neq 0$. It is however \emph{not} differentiable at $x=0$.

 Introduce a sequence of functions $\lbrace g_N(x)\rbrace$ by

\begin{equation}\label{eq12-72}
g_N(x)=\left\{
\begin{aligned}
\frac{N}{2},\quad |x|\leq \frac{1}{N},\\
0,\quad |x|> \frac{1}{N}.
\end{aligned}
\right.
\end{equation}

Using formula (\ref{eq12-60}) with $f(x)=g_N(x)$, it is straight forward to show that the corresponding solution is 
\begin{equation}\label{eq12-72b}
y_N(x)=\left\{
\begin{aligned}
\frac{N}{2k^2}\left(1-e^{-\frac{k}{N}}\cosh(kx)\right),\quad |x|\leq \frac{1}{N},\\
\frac{N}{2k^2}e^{-|x|}\sinh\left(\frac{k}{N}\right),\quad |x|> \frac{1}{N}.
\end{aligned}
\right\}
\end{equation}
For any $x$ there exist infinitely many $N$ such that $|x|\geq\frac{1}{N}$, and for such $N$ we have from (\ref{eq12-72b}) that
\begin{equation*}
y_N(x)=\frac{N}{2k^2}e^{-|x|}\sinh\left(\frac{k}{N}\right)\to \frac{1}{2k} e^{-|x|},
\end{equation*}
when $N$ approach infinity. Furthermore, for $x=0$ we have from (\ref{eq12-72b}) that
\begin{equation*}
y_N(0)=\frac{N}{2k^2}\left(1-e^{-\frac{k}{N}}\right)\to\frac{1}{2k},
\end{equation*}
when $N$ approach infinity. 

We have thus established that the infinite sequence of solutions (\ref{eq12-72b}) of equation (\ref{eq12-53}) corresponding to the sequence $\{g_N\}$ converge pointwise to the solution $y(x)=\frac{1}{2k}e^{-k|x|}$ corresponding to the Dirac delta function $f(x)=\delta(x)$. The sequence (\ref{eq12-72b}) does not converge uniformly to the function $y(x)$ and this explains the fact that all functions from the sequence  (\ref{eq12-72b}) are differentiable at $x=0$ whereas the limiting function $y(x)=\frac{1}{2k}e^{-|x|}$ is not.

In the previous example the sequence $\{f_N(x)\}$ converged pointwise to the limiting function $f(x)=1$ that generated the solution $y(x)=\frac{1}{k^2}$. Is it similarly true that in the current example the sequence $\{g_N(x)\}$ converge to the Dirac delta function that generated the solution $f(x)=\frac{1}{2k}e^{-k|x|}$?

In order to answer this question we start by observing that for all $N$
\begin{equation}\label{eq12-73}
\int^{\infty}_{-\infty}\mathrm d x\,g_N(x)= \int^{\frac{1}{N}}_{-\frac{1}{N}}\mathrm d x\,\frac{N}{2}=\frac{N}{2}(\frac{1}{N}+\frac{1}{N})=1,
\end{equation}
so $g_N(x)$ satisfies the same integral identity (\ref{eq12-25}) as the Dirac delta function (Here $D=(-\infty,\infty)$). Furthermore for any $x\neq 0$ there exists a $N_x$ such that

\begin{equation}\label{eq12-74}
|x|>\frac{1}{N_x}.
\end{equation}
 Therefore for all $N\geq N_x$, we have $g_N(x)=0$. Therefore 

\begin{equation}\label{eq12-75}
\lim\limits_{N\to \infty}g_N(x)=0,
\end{equation}
so the limit of the sequence $g_N$ when $N\to 0$ also satisfies the identity (\ref{eq12-24}). So in fact 

\begin{equation}\label{eq12-76}
\lim\limits_{N\to \infty} g_N(x)=\delta(x).
\end{equation}
I am playing with formulas here; the limit in (\ref{eq12-76}) does not mean the usual thing. After all, $\delta(x)$ is not a function but something different, namely a distribution. It would perhaps be more fair to say that the limit of the sequence $\lbrace g_N(x) \rbrace$, when $N\to \infty$ does \emph{not} exist, since there can be \emph{no} function satisfying (\ref{eq12-24}) and (\ref{eq12-25}).

What you see here is an example of a very subtle and profound mathematical idea called \emph{completion}. 
You have actually used this idea constantly since grade school, but you may not have been told. 

When we learn children arithmetic we start with whole numbers and eventually move on to fractions or \emph{rational numbers} as they are called. The universe of numbers are \emph{rational numbers} at grade school level. At some point the children are confronted with the equation

\begin{equation}\label{eq12-77}
x^2=2.
\end{equation}  
If the teacher is any good he will by trial and error convince the children that (\ref{eq12-77}) has no solutions. So what is there to do? He will say that (\ref{eq12-77}) actually has a solution that we denote by $\sqrt{2}$, and he will show, perhaps using a calculator, how to find a sequence of rational numbers that more and more closely approximate $\sqrt{2}$. In a similar way the teacher introduces $\sqrt{3}$, $\sqrt{5}$ etc. and he develops rules for calculations with such things. 

But the children  has been tricked. When the teacher argues that calculated sequences of rational numbers approximate $\sqrt{2}$ better and better and in the limit equals $\sqrt{2}$, he has been pulling wool over their eyes. What he has done is to \emph{complete} the set of rational numbers, in the process \emph{creating} the \emph{real numbers}. There is \emph{no} sequence of rational numbers converging to $\sqrt{2}$, in reality $\sqrt{2}$ is defined to \emph{be} a certain set of sequences of rational numbers. The existence of $\sqrt{2}$ is \emph{postulated}, pure and simple, through the completion process. After the completion process has lead to the set of real numbers, pure and applied mathematics can give them additional intuitive meaning by mapping them onto other physical and mathematical systems. In this way real numbers have come to be interpreted as geometrical points on a line etc. 

The relation between functions and distributions is exactly the same as the one between rational and real numbers. For the case of distributions Dirac is the teacher and the rest of us the children. Next semester we will map distributions onto a system derived from linear algebra and in this way gain more insight into the nature of distributions. 

The particular "solution" (\ref{eq12-71}) found by choosing $f(x)=\delta(x)$ is an example of what is called a \emph{Green's function}. Writing 

\begin{equation}\label{eq12-78}
g(x-y)=\frac{1}{2k}e^{-k|x-y|}. 
\end{equation} 
We have seen that $g(x-y)$ is a solution to the equation 

\begin{equation}\label{eq12-79}
g''-k^2g=-\delta(x-y),
\end{equation}
where we are taking the derivative with respect to the variable $x$, the variable $y$ is a parameter here.
Green's functions are of \emph{outmost} importance in applied mathematics and theoretical physics and will be discussed extensively next semester. As you can see from (\ref{eq12-79}), Green's functions are solutions to differential equations formulated in the universe of distributions. For this reason Green's functions are not actually functions but distributions. 

Let us now proceed to consider a slightly different ODE

\begin{equation}\label{eq12-80}
y''(x)+k^2y(x)=-f(x),\quad\ -\infty<x<\infty.
\end{equation}
For this equation the boundary conditions can not be that $y(x)$ and $y'(x)\to 0$ as $|x|\to\infty$. The reason is that (\ref{eq12-80}) does not have solutions decaying at infinity. In fact for $f(x)=0$ the solution consists of $\sin kx$ and $\cos kx$, which are bounded but does not vanish. We will therefore only require that $y$ and $y'$ are bounded at infinity. This implies, by an argument identical to the one on page 11, that the spectrum is real. We apply the Fourier transform like for the previous example and find

\begin{equation}\label{eq12-81}
Y(\lambda)=\frac{F(\lambda)}{\lambda^2-k^2}.
\end{equation} 
We observe that, compared to the previous case, a problem has appeared: $Y(\lambda)$ is not defined at $\lambda=\pm k$. The proper way to handle this is to use complex function theory and for more complicated situations complex function theory becomes a necessity. Here we will proceed in a more heuristic way by observing that (\ref{eq12-80}) becomes (\ref{eq12-53})  if we let $k\to ik $ or $k\to -ik $. Making the same replacement in the solution formula (\ref{eq12-60}) we get two solutions 

\begin{equation}\label{eq12-82}
y_{\pm}(x)=\pm\frac{i}{2k}\int^{\infty}_{-\infty}\mathrm d t\,e^{\pm ik|x-t|}f(t).
\end{equation}
Both of $y_{\pm}(x)$ satisfy (\ref{eq12-80}) and are bounded if $f(t)$ is absolutely integrable $\left(\int^{\infty}_{-\infty}\mathrm d t|f(t)|<\infty \right)$. Thus the boundary value problem (\ref{eq12-80}) does not have a unique solution. Why is that?

Let us start with the wave equation 

\begin{equation}\label{eq12-83}
-\frac{1}{c^2}\varphi_{tt}+\varphi_{xx}=-f(x),\quad-\infty<x<\infty,
\end{equation}
and let the source, $f(x)$, be supported in some bounded domain $G$. Outside this domain $G$, (\ref{eq12-83}) is the homogeneous wave equation whose general solution consists of left- and right going waves.

\begin{figure}[h!]
	\centering
	\includegraphics[scale=0.5]{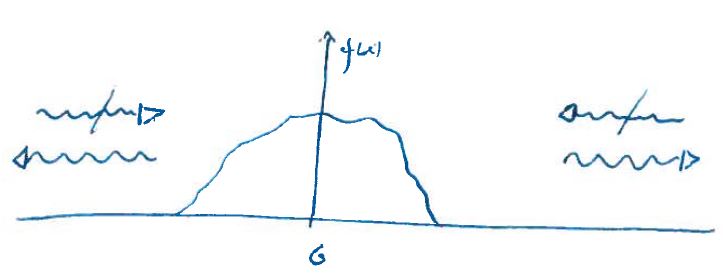}
	\caption{}
\end{figure}
We assume that the source is the originator of all waves. Thus, no source, no waves. Then it is evident that there can be no left moving waves to the right of $G$ and no right-moving waves to the left of $G$. This is a physical condition satisfied by all known wave phenomena; waves move \emph{away} from the source that creates them. This is in fact part of the definition of what a physical source is! 

Let us separate variables in (\ref{eq12-83}) 

\begin{gather}
\varphi(x,t)=e^{-i\omega t}y(x)\label{eq12-84}\\
\Rightarrow y''(x)+\left(\frac{\omega}{c}\right)^2y(x)=-f(x). \label{eq12-85}
\end{gather}
Defining $k=\frac{\omega}{c}$ we find our equation (\ref{eq12-80}). The homogeneous  equation (\ref{eq12-85})  has a solution proportional to $e^{\pm ikx}$. Taking these solutions together with (\ref{eq12-84}) we observe that

\begin{align}\label{eq12-86a}
& e^{ikx}\to\ e^{i(kx-\omega t)},\quad\text{ right travelling wave,}\\
& e^{-ikx}\to e^{-i(kx+\omega t)},\quad\text{ left travelling wave.}\nonumber
\end{align}  
The physical boundary condition thus consists of killing terms containing $e^{-ikx}$ to the right of the support of $f$ and terms containing $e^{ikx}$ to the left of the support of $f$. This can be expressed by saying that

\begin{align}\label{eq12-86b}
&\partial_x y(x)=iky(x)\quad  x>\text{Supp(f)}\\
&\partial_x y(x)=ikx\quad  x<\text{Supp(f)}.\nonumber
\end{align}
where $\text{Supp(f)}=\{x;f(x)\neq0\}$. Note that Supp(f) by definition includes its boundary points and is thus a closed set.

We have here assumed that the support of $f$ is bounded. However it is enough that $f(x)\to 0$, $|x|\to\infty$. For this more general case we pose the conditions

\begin{align}\label{eq12-87}
&\lim\limits_{x\to\infty}(\partial_x y(x)-iky(x))=0,\\
&\lim\limits_{x\to\infty}(\partial_x y(x)+iky(x))=0.\nonumber
\end{align} 
(\ref{eq12-86b}) or more generally (\ref{eq12-87}) is called \emph{the radiation condition at infinity}, or \emph{outgoing waves at infinity} and is the correct boundary condition for equations describing wavelike phenomena. (\ref{eq12-87}) picks out the unique solution $y_{+}(x)$ to our boundary value problem (\ref{eq12-80}),(\ref{eq12-87}). Choosing $f(x)=\delta(x)$ gives us the following Green's function for this case

\begin{equation}\label{eq12-88}
G(x-\xi)=\frac{i}{2k}e^{ik|x-\xi|}.
\end{equation}

\subsubsection{The Cauchy problem for the heat equation}
Let us next apply the Fourier transform to the Cauchy problem for the heat equation

\begin{align} \label{eq12-89}
u_t(x,t)-c^2u_{xx}(x,t)=0,\quad &-\infty<x<\infty,\,t>0,\\
u(x,0)=f(x),\quad  &-\infty<x<\infty.\nonumber
\end{align}
In order to apply the Fourier transform we must assume that $f$ has a Fourier transform and that $u$, $u_x$ goes to zero at $x=\pm\infty$.

Taking the Fourier transform with respect to $x$ of the heat equation (\ref{eq12-89}) we get

\begin{gather}\label{eq12-94}
\partial_t U(\lambda,t)+c^2\lambda^2U(\lambda,t)=0,\\
U(\lambda,0)=F(\lambda),\nonumber
\end{gather}
where 

\begin{eqnarray}\label{eq12-90}
U(\lambda,t)=\frac{1}{\sqrt{2\pi}}\int_{-\infty}^{\infty}\mathrm d x\, e^{i\lambda x}u(x,t),\\
F(\lambda)=\frac{1}{\sqrt{2\pi}}\int_{-\infty}^{\infty}\mathrm d x\, e^{i\lambda x}f(x).
\end{eqnarray}

\noindent The equations (\ref{eq12-94}) are uncoupled for different values of $\lambda$. It is the linearity of the heat equation that produce an \emph{uncoupled} system of ODEs for the set of functions $U(\lambda,t)$ indexed by the spectral parameter $\lambda$. Note how similar this is to what we got when we applied the finite Fourier transform to linear PDEs on a finite domain. 

The initial value problems (\ref{eq12-94}) are easy to solve. We get

\begin{equation}\label{eq12-96}
U(\lambda,t)=F(\lambda)e^{-c^2\lambda^2t}.
\end{equation}
Using the inverse Fourier transform we get the solution

\begin{align}\label{eq12-97}
u(x,t)&=\frac{1}{\sqrt{2\pi}}\int^{\infty}_{-\infty}\mathrm d \lambda\, e^{-i\lambda x}U(\lambda,t)\\
&=\frac{1}{\sqrt{2\pi}}\int^{\infty}_{-\infty}\mathrm d \lambda\, e^{-i\lambda x}F(\lambda)e^{-c^2\lambda^2 t}\nonumber\\
&=\frac{1}{\sqrt{2\pi}}\int^{\infty}_{-\infty}\mathrm d \lambda\, \frac{1}{\sqrt{2\pi}}\int^{\infty}_{-\infty}\mathrm d x'\,e^{-i\lambda x}e^{i\lambda x'}f(x')e^{-c^2\lambda^2t}\nonumber\\
&=\frac{1}{2\pi}\int_{-\infty}^{\infty}\mathrm d x'\,\left[\int_{-\infty}^{\infty}\mathrm d\lambda\,e^{-i\lambda(x-x')-c^2\lambda^2t}\right]f(x'),\nonumber
\end{align}
where we have changed order of integration. Defining 
\begin{equation}\label{eq12-98}
g(\alpha)=\frac{1}{2\pi}\int_{-\infty}^{\infty}\mathrm d\lambda\,e^{-i\lambda\alpha-c^2\lambda^2t},
\end{equation}
\noindent the solution can be written in the form 
\begin{equation}\label{eq12-98.1}
u(x,t)=\int_{-\infty}^\infty\;dx'\;g(x-x')f(x').
\end{equation}
\noindent Note that we could have gotten this expression more directly by using the convolution theorem.
Observe that 

\begin{align}\label{eq12-98}
&g(\alpha)=\frac{1}{2\pi}\int_{-\infty}^{\infty}\mathrm d\lambda\,e^{-i\lambda\alpha-c^2\lambda^2t}\\
&=\frac{1}{2\pi} \int_{-\infty}^{0}\mathrm d\lambda\,e^{-i\lambda\alpha-c^2\lambda^2t}+\frac{1}{2\pi}\int_{0}^{\infty}\mathrm d\lambda\,e^{-i\lambda\alpha-c^2\lambda^2t}\nonumber\\
&=\frac{1}{2\pi} \int_{0}^{\infty}\mathrm d\lambda\,e^{i\lambda\alpha}e^{-c^2\lambda^2t}+\frac{1}{2\pi}\int_{0}^{\infty}\mathrm d\lambda\,e^{-i\lambda\alpha}e^{-c^2\lambda^2t}\nonumber\\
&=\frac{1}{\pi}\int^{\infty}_0\mathrm d\lambda\cos\lambda\alpha e^{-c^2\lambda^2t}.\nonumber
\end{align}
This integral, because of the exponential factor present for $t>0$, converges very fast, so we can differentiate under the integral sign 

\begin{align}\label{eq12-100}
\frac{\mathrm d g}{\mathrm d\alpha}(\alpha)&=-\frac{1}{\pi}\int_0^{\infty}\mathrm d\lambda\,\lambda\sin\lambda\alpha e^{-c^2\lambda^2t}\\
&=\frac{1}{2\pi c^2t}\sin\lambda\alpha e^{-c^2\lambda^2t}\Big|_{\lambda=0}^{\lambda=\infty}-\frac{\alpha}{2\pi c^2t}\int_0^{\infty}\mathrm d\lambda\cos\lambda\alpha e^{-\lambda^2c^2},\nonumber
\end{align}
so 

\begin{equation}\label{eq12-101}
\frac{\mathrm d g}{\mathrm d \alpha}=-\frac{\alpha}{2 c^2t}g.
\end{equation}
This ODE is simple to solve. The solution is

\begin{equation}\label{eq12-102}
g(\alpha)=g(0)e^{-\frac{\alpha^2}{4 c^2t}},
\end{equation}
and 

\begin{equation}\label{eq12-103}
g(0)=\frac{1}{\pi}\int_0^{\infty}\mathrm d\lambda\,e^{-c^2\lambda^2t}=\frac{1}{\sqrt{4\pi c^2t}}. 
\end{equation}
(This is a Gaussian integral.) Thus 

\begin{equation}\label{eq12-104}
g(\alpha)=\frac{1}{\sqrt{4\pi c^2t}}e^{-\frac{\alpha^2}{4c^2t}}.
\end{equation}
The formula for $u(x,t)$ thus becomes 

\begin{equation}\label{eq12-106}
u(x,t)=\frac{1}{\sqrt{4\pi c^2t}}\int_{-\infty}^{\infty}\mathrm d x'\,e^{-\frac{(x-x')^2}{4c^2t}}f(x').
\end{equation} 
If we formally use the initial condition

\begin{equation}\label{eq12-107}
f(x)=\delta(x-x'),\quad -\infty<x<\infty,
\end{equation}
we get the solution 

\begin{equation}\label{eq12-108}
u(x,t;x')\equiv G(x-x',t)=\frac{1}{\sqrt{4\pi c^2t}}e^{-\frac{(x-x')^2}{4c^2t}}.
\end{equation}
$G$ is yet another example of a Green's function. In this context it is called the \emph{fundamental solution} to the heat equation. The fundamental solution is the temperature profile generated by a point-like heat source located at $x=x'$ for $t=0$. 

  Note that it is an actual classical solution (verify!) for all $t>0$ even if the initial condition is extremely singular. This is yet another example illustrating the extreme smoothening properties of the heat equation.

From the solution formula (\ref{eq12-106}) it is evident that even if the initial heat profile is zero beyond a finite interval $[-L,L]$, the solution is nonzero for \emph{any} point $x\in\mathbb R$ for any $t>0$, however small $t$ is.

  The heat equation clearly does not respect causality.

\begin{figure}[h!]
	\centering
	\includegraphics[scale=0.4]{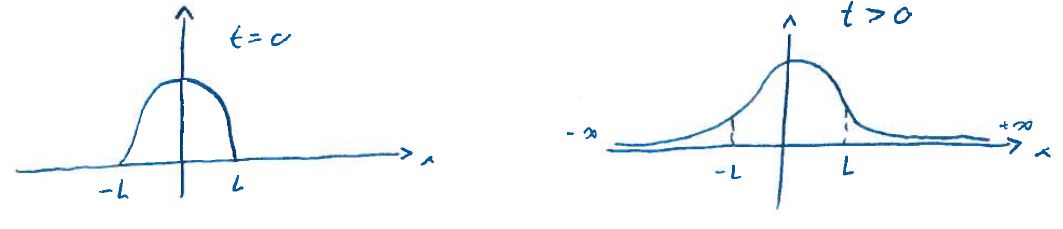}
	\caption{}
\end{figure}
\noindent{Also note that for the fundamental solution we have}

\begin{equation}\label{eq12-109}
\int^{\infty}_{-\infty}\mathrm d x'\, G(x-x',t)g(x')\to g(x),\quad t\to 0^{+}.
\end{equation}
Thus 

\begin{equation}\label{eq12-110}
G(x-x',t)\to\delta(x-x'),\quad t\to 0^{+}.
\end{equation}
The fundamental solution is thus a sequence of smooth functions converging to the Dirac delta distribution. Which is exactly what it needs to do in order to satisfy the initial condition!

We can argue for the correctness of (\ref{eq12-109}) as follows. Observe that as $t\to0^{+}$ the function $G(x-\xi,t)$ becomes more and more localized around $x=\xi$. Thus for $t$ close enough to zero and for all reasonable functions $g(x)$ we have

\begin{equation}\label{eq12-111}
\int^{\infty}_{-\infty}\mathrm d x'\, G(x-x',t)g(x')\approx g(x)\int^{\infty}_{-\infty}\mathrm d x'\, G(x-x',t)=g(x),
\end{equation}
because

\begin{equation}\label{eq12-112}
\int^{\infty}_{-\infty}\mathrm d x'\, G(x-x',t)=1.
\end{equation}
Verify this!


\subsubsection{The Cauchy problem for the wave equation}
Let us consider the problem 

\begin{gather}\label{eq12-113}
u_{tt}(x,t)-c^2u_{xx}(x,t)=0,\quad -\infty<x<\infty,\,t>0,\\
u(x,0)=f(x),\nonumber\\
u_t(x,0)=g(x).\nonumber
\end{gather}
Taking the Fourier transform of this problem, proceeding like for the heat equation, we get

\begin{gather}\label{eq12-114}
\partial_{tt}u(\lambda,t)-c^2\lambda^2u(\lambda,t)=0,\quad t>0,\\
u(\lambda,0)=F(\lambda),\nonumber\\
\partial_tu(\lambda,0)=G(\lambda).\nonumber
\end{gather}
The initial value problem (\ref{eq12-114}) is entirely standard and the solution is

\begin{align}\label{eq12-115}
u(\lambda,t)&=\left(\frac{1}{2}F(\lambda)+\frac{1}{2i\lambda c}G(\lambda)\right)e^{i\lambda ct}\\
&+\left(\frac{1}{2}F(\lambda)-\frac{1}{2i\lambda c}G(\lambda)\right)e^{-i\lambda ct}.\nonumber
\end{align}
Taking the inverse transform we get 

\begin{align}\label{eq12-116}
u(x,t)&=\frac{1}{2\sqrt{2\pi}}\int_{-\infty}^{\infty}\mathrm d\lambda\, e^{-i\lambda(x-ct)}F(\lambda)+\frac{1}{2\sqrt{2\pi}}\int_{-\infty}^{\infty}\mathrm d\lambda\, e^{-i\lambda(x+ct)}F(\lambda)\\
&+\frac{1}{2c\sqrt{2\pi}}\int_{-\infty}^{\infty}\mathrm d\lambda\, e^{-i\lambda(x-ct)}\frac{G(\lambda)}{i\lambda}-\frac{1}{2c\sqrt{2\pi}}\int_{-\infty}^{\infty}\mathrm d\lambda\, e^{-i\lambda(x+ct)}\frac{G(\lambda)}{i\lambda}.\nonumber
\end{align}
By the definition of the inverse Fourier transform we have

\begin{equation}\label{eq12-117}
g(x)=\frac{1}{\sqrt{2\pi}}\int_{-\infty}^{\infty}\mathrm d\lambda\, e^{-i\lambda x}G(\lambda).
\end{equation}
But then 

\begin{equation}\label{eq12-118}
\int_0^x\mathrm d s\, g(s)=-\frac{1}{\sqrt{2\pi}}\int_{-\infty}^{\infty}\mathrm d\lambda\, e^{-i\lambda x}\frac{G(\lambda)}{i\lambda}+D,
\end{equation}
where D is a constant. Given (\ref{eq12-118}), (\ref{eq12-116}) can be written as 

\begin{equation}\label{eq12-119}
u(x,t)=\frac{1}{2}\left(f(x-ct)+f(x+tc)\right)+\frac{1}{2c}\int^{x+ct}_{x-ct}\mathrm d s\,g(s).
\end{equation}
This is the d'Alembert formula, now derived using the Fourier transform. 

There is however something rather fishy with what we did. Let us say that we did not realize that the identity (\ref{eq12-118}) holds. Then, for given $F(\lambda)$, $G(\lambda)$, we would have to find $u(x,t)$ by solving the integrals in (\ref{eq12-116}) analytically or numerically. We should however worry about the last two integrals. They are of the general form

\begin{equation}\label{eq12-120}
I(\alpha)=\int_{-\infty}^{\infty}\mathrm d\lambda\, e^{-i\lambda \alpha}\frac{H(\lambda)}{i\lambda},\quad\alpha\in\mathbb R,
\end{equation}
where in general $H(0)\neq 0$. The integral is singular at $\lambda=0$ and the question how the integral should be understood must be raised.

Away from $\lambda=0$ there are no problems, so the question we need to ask is if we can make sense of integrals of the type

\begin{equation}\label{eq12-121}
\int_{-L}^{L}\mathrm d\lambda\,\frac{1}{\lambda},
\end{equation} 
where $L$ is arbitrary but fixed. In calculus we have seen that at least some integrals with a singular integrand in the domain of integration can be made sense of by considering it to be an improper integral. Let us try to use this approach on (\ref{eq12-121}):

\begin{align}\label{eq12-122}
\int_{-L}^{L}\mathrm d\lambda\,\frac{1}{\lambda}&=\int_{-L}^{0}\mathrm d\lambda\,\frac{1}{\lambda}+\int_{0}^{L}\mathrm d\lambda\,\frac{1}{\lambda}\\
&=\lim\limits_{\epsilon\to 0} \int_{-L}^{-\epsilon}\mathrm d\lambda\,\frac{1}{\lambda}+\lim\limits_{\delta\to 0} \int_{\delta}^{L}\mathrm d\lambda\,\frac{1}{\lambda}\nonumber\\
&=\lim\limits_{\epsilon\to 0} \int_{L}^{\epsilon}\mathrm d\mu\,\frac{1}{\mu}+\lim\limits_{\delta\to 0} \int_{\delta}^{L}\mathrm d\lambda\,\frac{1}{\lambda}\nonumber\\
&=-\lim\limits_{\epsilon\to 0} \int_{\epsilon}^{L}\mathrm d\mu\,\frac{1}{\mu}+\lim\limits_{\delta\to 0} \int_{\delta}^{L}\mathrm d\lambda\,\frac{1}{\lambda}\nonumber\\
&=-\lim\limits_{\epsilon\to 0}(\ln L-\ln\epsilon)+\lim\limits_{\delta\to 0}(\ln L-\ln\delta)\nonumber\\
&=\lim\limits_{\epsilon\to 0}\ln\epsilon-\lim\limits_{\delta\to 0}\ln\delta.\nonumber
\end{align}
This limit does not exist. Thus (\ref{eq12-121}) can \emph{not} be given meaning as an improper integral. Observe however that if we in (\ref{eq12-122}) let $\epsilon=\delta$ \emph{and} take the difference \emph{before} we take the limit we get zero, which is a finite answer. Thus what we are proposing is to understand (\ref{eq12-121}) in the following way

\begin{align}\label{eq12-123}
\int_{-L}^{L}\mathrm d\lambda\,\frac{1}{\lambda}&=\lim\limits_{\epsilon\to 0}\left\{\int_{-L}^{-\epsilon}\mathrm d\lambda\,\frac{1}{\lambda}+\int_{\epsilon}^{L}\mathrm d\lambda\,\frac{1}{\lambda}\right\}\\
&=\lim\limits_{\epsilon\to 0} \left\{\int_{L}^{\epsilon}\mathrm d\mu\,\frac{1}{\mu}+\int_{\epsilon}^{L}\mathrm d\lambda\,\frac{1}{\lambda}\right\}\nonumber\\
&=\lim\limits_{\epsilon\to 0}\left\{\ln\epsilon -\ln L+\ln L-\ln\epsilon \right\}=0.\nonumber
\end{align}
When we use the definition (\ref{eq12-123}) we understand the integral as a \emph{Cauchy principal value} integral. This is now a new way for you to calculate integrals with a singular integrand. In general for a function $f(x)$, with a singularity at $x=x_0$, the definition is 

\begin{equation}
PV\int_{-L}^{L}\mathrm d x\,f(x)=\lim\limits_{\epsilon\to 0}\left\{\int_{-L}^{x_0-\epsilon}\mathrm d x\,f(x)+\int^{L}_{x_0+\epsilon}\mathrm d x\,f(x)\right\},
\end{equation}
where $PV$ signifies that we are calculating the principal value. 

The principal value is finite only if the singularity in $f$ is of the right type. For example for $f(x)=\frac{1}{x^2}$ we get 

\begin{align}\label{eq12-125}
&PV\int_{-1}^{1}\mathrm d x\,\frac{1}{x^2}=\lim\limits_{\epsilon\to 0}\left\{\int_{-1}^{-\epsilon}\mathrm d x\, \frac{1}{x^2}+\int_{\epsilon}^{1}\mathrm d x\,\frac{1}{x^2}\right\}\\
&=\lim\limits_{\epsilon\to 0}\left\{\left(-\frac{1}{x}\right)\Big|^{-\epsilon}_{-1}+\left(-\frac{1}{x}\right)\Big|_{\epsilon}^1\right\}\nonumber\\
& \lim\limits_{\epsilon\to 0}\left\{ \frac{1}{\epsilon}+1-1+\frac{1}{\epsilon}\right\}=\lim\limits_{\epsilon\to 0}\frac{1}{\epsilon}=\infty!\nonumber
\end{align}
If we interpret the two integrals containing $\dfrac{1}{\lambda}$ in (\ref{eq12-110}) as Cauchy principal value integrals, we can use the formula to calculate the solution to the wave equation.

However, there is still a weak fishy smell here; why is interpreting the singular integrals in (\ref{eq12-110}) as Cauchy principal value integrals the right thing to do? How do we know that this will give the same solution as, say, the d'Alembert formula? Could one not imagine some other funny way of getting a finite number out of integrals like (\ref{eq12-121})? And would this not perhaps produce a different solution; which should not be possible because the Cauchy problem for the wave equation has a unique solution? What is going on?

Answering this question, which is an important one, would take us too far from the main focus of the course, so I will refrain from giving an answer. Here I will only ensure you that this \emph{is} the right way to understand the integral (\ref{eq12-121}) and it \emph{does} give the same solution as the d'Alembert formula. It is also worth noting that questions like the one we have been considering frequently turn up when we use transforms to solve PDEs, and solving them can often be a deep and vexing problem. 

\subsubsection{The Laplace equation in a half-plane}
Let us consider the Laplace equation in the upper halfplane:

\begin{align}\label{eq12-126}
& u_{xx}(x,y)+u_{yy}(x,y)=0,\quad-\infty<x<\infty,\,y>0,\\
& u(x,0)=f(x),\nonumber\\
& u(x,y) \text{ is bounded as }y\to\infty.\nonumber
\end{align}
We introduce the Fourier transform in $x$

\begin{equation}\label{eq12-127}
U(\lambda,y)=\frac{1}{\sqrt{2\pi}}\int_{-\infty}^{\infty}\mathrm d x\,e^{i\lambda x}u(x,y).
\end{equation}
Using this Fourier transform on (\ref{eq12-126}) gives us the following infinite system of uncoupled ODEs:

\begin{align}\label{eq12-128}
& \partial_{yy}U(\lambda,y)-\lambda^2U(\lambda,y)=0,\quad 0<y<\infty,\\
& u(\lambda,0)=F(\lambda)\nonumber,\\
& u(\lambda,y) \text{ is bounded as }y\to\infty.\nonumber
\end{align}
This problem is standard and the solution is

\begin{equation}\label{eq12-129}
U(\lambda,y)=F(\lambda)e^{-|\lambda|y}.
\end{equation}
The inverse transform gives us the solution 

\begin{align}\label{eq12-130a}
u(x,y)&=\frac{1}{\sqrt{2\pi}}\int_{-\infty}^{\infty}\mathrm d \lambda\,e^{-i\lambda x}F(\lambda)e^{-|\lambda|y}\\
&=\frac{1}{2\pi}\int_{-\infty}^{\infty}\mathrm d \lambda\,e^{-i\lambda x}e^{-|\lambda|y}\int_{-\infty}^{\infty}\mathrm dx'\,e^{i\lambda x'}f(x')\nonumber\\
&=\frac{1}{2\pi}\int_{-\infty}^{\infty}\mathrm dx'\,f(x')\left[\int_{-\infty}^{\infty}\mathrm d\lambda\,e^{-i\lambda(x-x')-|\lambda|y}\right].\nonumber
\end{align}
Furthermore, observe that 

\begin{align}\label{eq12-131a}
&\int_{-\infty}^{\infty}\mathrm d\lambda\,e^{-i\lambda(x-x')-|\lambda|y}\nonumber\\
&=2\int_{0}^{\infty}\mathrm d\lambda\,\cos\lambda(x-x')e^{-\lambda y}\nonumber\\
&=\frac{2y}{(x-x')^2+y^2}\;,
\end{align}
and therefore we have the following expression for the solution of the Laplace equation in a half-space
\begin{equation}\label{eq12-130b}
u(x,y)=\frac{y}{\pi}\int_{-\infty}^{\infty}\mathrm dx'\,\frac{f(x')}{(x-x')^2+y^2}.
\end{equation}
By direct substitution it can be verified that (\ref{eq12-130b}) is a solution to the Laplace equation for $y>0$. 

When $y\to 0$ the factor $\dfrac{y}{\pi}$ goes to zero, but the integrand approaches

\begin{equation}\label{eq12-131b}
\frac{f(x')}{(x-x')^2},
\end{equation}
which is singular at $x'=x$. A careful study of the limit shows that the limit of (\ref{eq12-130b}) when $y\to 0$ is actually $f(x)$. We will not do this here.

For simple choices of $f(x)$ we can solve the integral (\ref{eq12-130b}). For example if $f(x)=H(x)$ where $H(x)$ is the Heaviside step function 

\begin{equation}\label{eq12-132}
H(x)=\left\{
\begin{aligned}
1,\quad x>0,\\
0,\quad x<0,
\end{aligned}
\right.
\end{equation}
we find 

\begin{align}\label{eq12-133}
u(x,y)&=\frac{y}{\pi}\int_{-\infty}^{\infty}\mathrm dx'\,\frac{H(x')}{(x-x')^2+y^2}\\
&=\frac{y}{\pi} \int_{0}^{\infty}\mathrm dx'\,\frac{1}{(x-x')^2+y^2}\nonumber\\
&=\frac{1}{\pi y} \int_{0}^{\infty}\mathrm dx'\,\frac{1}{1+\left(\frac{x-x'}{y}\right)^2}.\nonumber
\end{align}
Substitute  $u=\dfrac{x'-x}{y}$ and $\mathrm du=\dfrac{1}{y}\mathrm dx'$ in this integral. This gives us the following explicit expression for $u(x,y)$

\begin{align}\label{eq12-134}
u(x,y)&=\frac{1}{\pi}\int_{-\frac{x}{y}}^{\infty}\mathrm du\frac{1}{1+u^2}=\frac{1}{\pi}tg^{-1}(u)\Big|^{\infty}_{-\frac{x}{y}}\\
&=\frac{1}{2}+\frac{1}{\pi}tg^{-1}\left(\frac{x}{y}\right).\nonumber	
\end{align}
This function satisfies the Laplace equation and

\begin{equation}\label{eq12-135a}
\lim\limits_{y\to 0}u(x,y)=\frac{1}{2}+\frac{1}{\pi}\lim\limits_{y\to 0}tg^{-1}\left(\frac{x}{y}\right)=\left\{
\begin{aligned}
&\frac{1}{2}+\frac{1}{\pi}\cdot\frac{\pi}{2}=1,\,x>0,\\
&\frac{1}{2}+\frac{1}{\pi}\cdot\left(-\frac{\pi}{2}\right)=0,\,x<0,
\end{aligned}
\right.
\end{equation} 
so $u(x,y)$ satisfies the boundary condition at $y=0$. Also

\begin{equation}\label{eq12-135b}
\lim\limits_{y\to\infty}u(x,y)=\frac{1}{2},
\end{equation}
so $u(x,y)$ is bounded at $y=+\infty$. 

However, $u(x,y)$ does not go to zero as $|x|\to\infty$ and so does not actually have a Fourier transform, which we assumed. The reason for this is that our boundary function $f(x)=H(x)$ does not have a Fourier transform in the conventional sense. Observe that $u(x,y)$ is infinitely smooth for $y>0$ even if it is discontinuous in $x$ for $y=0$. We are rediscovering the fact that the Laplace equation can not support solutions with singularities. As we have seen, this is a general feature of elliptic equations. 

Let us say that we were rather to solve the Laplace equation with Neumann conditions

\begin{align}\label{eq12-136}
u_{xx}(x,y)+u_{yy}(x,y)=0,\quad&-\infty<x<\infty,\,y>0,\\
u_y(x,0)=g(x).\nonumber
\end{align}
We can reduce this to a problem of the previous type using a trick, invented by an English mathematician, called \emph{Stoke's rule}.

\noindent Introduce $v(x,y)=u_y(x,y)$. Differentiating (\ref{eq12-130b}) with respect to $y$ gives

\begin{align}\label{eq12-137}
v_{xx}+v_{yy}=0,\quad&-\infty<x<\infty,\,y>0,\\
v(x,0)=g(x),\quad&-\infty<x<\infty.\nonumber
\end{align}
The solution to this problem is given by (\ref{eq12-130b}) and since

\begin{equation}\label{eq12-138}
u(x,y)=\int_0^y\mathrm d s\, v(x,s),
\end{equation}
we get 

\begin{align}\label{eq12-139}
u(x,y)&=\frac{1}{\pi}\int_0^y\mathrm ds\, s\int^{\infty}_{-\infty}\mathrm dx'\,\frac{f(x')}{(x-x')^2+s^2}\nonumber\\
&=\frac{1}{\pi} \int^{\infty}_{-\infty}\mathrm dx'\,f(x') \int_0^y\mathrm ds\,\frac{s}{(x-x')^2+s^2}\nonumber\\
&=\frac{1}{2\pi} \int^{\infty}_{-\infty}\mathrm dx'\,\ln\left[(x-x')^2+y^2\right]f(x').
\end{align}
Observe that $u(x,y)$ defined by  (\ref{eq12-139}) is not the only solution to (\ref{eq12-136}).

  For any solution $v(x,y)$ to (\ref{eq12-137}) we can add an arbitrary function $C(x)$ to the right hand side of (\ref{eq12-138}) and in this way define a whole infinite familty of functions $u(x,y)$
\begin{equation}\label{eq12-138b}
u(x,y)=\int_0^y\mathrm d s\, v(x,s)+C(x).
\end{equation}
All such functions satisfy the boundary condition $u_y(x,0)=g(x)$, but in order to satisfy the Laplace equation (\ref{eq12-136}) the added function $C(x)$ must be linear. Furthermore, in order for the solution to be bounded in $x$ this linear function must be a constant. 

\subsection{The Fourier sine- and cosine transform}
Let us consider the following eigenvalue problem

\begin{equation}\label{eq12-141}
M''(x)+\lambda^2M(x)=0,\quad0<x<\infty. 
\end{equation}
on a semi-infinite 1D domain.
We will consider two sets of boundary conditions

\begin{align}
M(0)=0,\quad M(x)\text{ bounded as }x\to+\infty,\label{eq12-142}\\
M'(0)=0,\quad M(x)\text{ bounded as }x\to+\infty.\label{eq12-143}
\end{align}
Let us first consider (\ref{eq12-141}) with boundary conditions (\ref{eq12-142}) . 

We will only consider (\ref{eq12-141}) for real $\lambda$. For $\lambda=0$ the general solution of (\ref{eq12-141}) is 

\begin{equation}\label{eq12-144}
M(x)=Ax+B,
\end{equation}
and

\begin{align*}
& M(0)=0\,\Rightarrow\, B=0\Rightarrow\,M(x)=Ax,\\
& M(x) \text{ bounded }\Rightarrow\,A=0\,\Rightarrow M(x)=0.
\end{align*}
Thus $\lambda=0$ is not an eigenvalue. For $\lambda\neq 0$ the general solution to (\ref{eq12-141}) is 
\begin{equation}\label{eq12-145}
M(x)=A\cos\lambda x+B\sin\lambda x,
\end{equation}
and 
\begin{align*}
& M(0)=0\,\Rightarrow\,A=0\,\Rightarrow\,M(x)=B\sin\lambda x.
\end{align*}
This function is clearly bounded for all real $\lambda$ and therefore all such $\lambda$ are in the spectrum. Observe however that 

\begin{equation}\label{eq12-146}
M_{-\lambda}(x)=-M_{\lambda}(x).
\end{equation}
Therefore, the eigenfunctions corresponding to $\lambda$ and $-\lambda$ are linearly dependent. We therefore loose nothing by restricting to the range
\begin{equation}\label{eq12-147}
0<\lambda<\infty,
\end{equation}
and eigenfunctions 
\begin{equation}\label{eq12-148}
M_{\lambda}(x)=\sin\lambda x.
\end{equation}
For (\ref{eq12-141}) with the boundary conditions (\ref{eq12-143}) we find in a similar way the spectrum
\begin{equation}\label{eq12-149}
0\leq\lambda<\infty,
\end{equation} 
and eigenfunctions

\begin{equation}\label{eq12-150}
M_{\lambda}(x)=\cos\lambda x.
\end{equation}
In order to find the transform corresponding to (\ref{eq12-147})  and (\ref{eq12-148}) it is useful to consider the finite domain

\begin{equation}\label{eq12-151}
0<x<l,
\end{equation} 
and then let $l$ approach infinity. On the bounded domain (\ref{eq12-151}), the boundary condition at $x=+\infty$ is replaced by the condition $M(l)=0$. This gives us the following normalized eigenfunctions

\begin{equation}\label{eq12-152}
M_k(x)=\sqrt{\frac{2}{l}}\sin\lambda_kx,\quad k=1,2,\dots,
\end{equation}
where

\begin{equation}\label{eq12-153}
\lambda_k=\frac{k\pi}{l}.
\end{equation}
A function $\tilde{f}(x)$ can now be expanded in a Fourier series with respect to (\ref{eq12-152}) according to 

\begin{equation}\label{eq12-154}
\tilde{f}(x)=\sqrt{\frac{2}{l}}\sum\limits_{k=1}^{\infty}F_k\sin\lambda_k x,
\end{equation}
\begin{equation}\label{eq12-155}
F_k=\sqrt{\frac{2}{l}}\int_0^l\mathrm d x\,\sin\lambda_k x\tilde{f}(x).
\end{equation}
Observe that the distance between the points $\lambda_k$ on the $\lambda$-axis is 

\begin{equation}\label{eq12-156}
\Delta \lambda=\lambda_{k+1}-\lambda_k=\frac{\pi}{l}.
\end{equation}
We use this and rewrite (\ref{eq12-154}), (\ref{eq12-155}) in the following way

\begin{equation}\label{eq12-157}
\tilde{f}(x)=\frac{\sqrt{2}}{\pi}\sum\limits_{k=1}^{\infty}\Delta\lambda\sin\lambda_k x\tilde{F}_k,
\end{equation}

\begin{equation}\label{eq12-158}
\tilde{F_k}(x)=\sqrt{2}\int_0^l\mathrm d x\,\sin\lambda_k x\tilde{f}(x),
\end{equation}
where $\tilde{F}_k=\sqrt{l}F_k$. If we formally let $l\to\infty$ in (\ref{eq12-157}) and (\ref{eq12-158}) we get 

\begin{equation}\label{eq12-159}
\tilde{f}(x)=\frac{\sqrt{2}}{\pi}\int_0^{\infty}\mathrm d\lambda\,\sin\lambda x\tilde{F}(\lambda),
\end{equation}

\begin{equation}\label{eq12-160}
\tilde{F}(\lambda)=\sqrt{2}\int_0^{\infty}\mathrm dx\,\sin\lambda x\tilde{f}(x).
\end{equation}
In order to get a more symmetric and pleasing form of the transform (\ref{eq12-159}), (\ref{eq12-160}) we introduce scalings 

\begin{equation}\label{eq12-161}
\tilde{F}=\alpha F,\quad \tilde{f}=\beta f,
\end{equation}

\begin{equation}\label{eq12-162}
\Rightarrow f(x)=\frac{\sqrt{2}\alpha}{\pi\beta}\int_0^{\infty}\mathrm d\lambda\,\sin\lambda xF(\lambda),
\end{equation}

\begin{equation}\label{eq12-163}
F(\lambda)=\frac{\sqrt{2}\beta}{\alpha}\int_0^{\infty}\mathrm dx\,\sin\lambda xf(x).
\end{equation}
Choose $\alpha=\sqrt{\pi}$, $\beta=1$. Then we get the \emph{Fourier sine transform}

\begin{align}
f(x)&=\sqrt{\frac{2}{\pi}}\int_0^{\infty}\mathrm d\lambda\,\sin\lambda xF_s(\lambda),\label{eq12-164}\\
F_s(\lambda)&=\sqrt{\frac{2}{\pi}}\int_0^{\infty}\mathrm dx\,\sin\lambda xf(x).\label{eq12-165}
\end{align}
Our friendly analysis can tell us for which pair of nice functions the improper integrals converge and for which (\ref{eq12-164}) and (\ref{eq12-165}) holds. 

However we will apply it to much more general things and will proceed formally and hope that in the end everything will be ok in the universe of distributions. Here it is only worth noting that even though $\sin\lambda x$ is zero for $x=0$, $f(x)$ in (\ref{eq12-164}) might not approach zero as $x$ approach zero. This is because in general the convergence of the improper integral will not allow us to interchange integral and limit. 

A similar derivation for the boundary condition (\ref{eq12-143}) gives us the \emph{Fourier cosine transform} 

\begin{align}
f(x)&=\sqrt{\frac{2}{\pi}}\int_0^{\infty}\mathrm d\lambda\,\cos\lambda x\,F_c(\lambda),\label{eq12-166}\\
F_c(\lambda)&=\sqrt{\frac{2}{\pi}}\int_0^{\infty}\mathrm dx\,\cos\lambda x\,f(x).\label{eq12-167}
\end{align}
The Fourier sine- and cosine transforms have several interesting and useful properties\cite{textbook}. Here we will discuss how they behave with respect to derivatives. Observe that

\begin{equation}\label{eq12-168}
\int_0^{\infty}\mathrm d x\,\sin\lambda xf'(x)=\sin\lambda xf(x)\Big|^{\infty}_0-\lambda\int_0^{\infty}\mathrm dx\,\cos\lambda xf(x).
\end{equation}
If we assume that $f(x)$ vanishes at $x=+\infty$ we get the identity

\begin{equation}\label{eq12-169}
\sqrt{\frac{2}{\pi}}\int_0^{\infty}\mathrm d x\,\sin\lambda xf'(x)=-\lambda F_c(\lambda). 
\end{equation}
In a similar way we find, assuming that also $f'(x)$ vanishes at infinity,

\begin{equation}\label{eq12-170}
\sqrt{\frac{2}{\pi}}\int_0^{\infty}\mathrm d x\,\sin\lambda xf''(x)=\lambda\sqrt{\frac{2}{\pi}}f(0)-\lambda^2F_s(\lambda).
\end{equation}
Thus the Fourier sine transform of $f'$ is expressed in terms of the Fourier cosine transform. The Fourier sine transform of $f''$ is expressed in terms of the Fourier sine transform of $f$ and a term depending on $f(0)$. 

For the Fourier cosine transform we get in a similar way 

\begin{equation}\label{eq12-171}
\sqrt{\frac{2}{\pi}}\int_0^{\infty}\mathrm d x\,\cos\lambda xf'(x)=-\sqrt{\frac{2}{\pi}}f(0)+\lambda F_s(\lambda),
\end{equation}
\begin{equation}\label{eq12-172}
\sqrt{\frac{2}{\pi}}\int_0^{\infty}\mathrm d x\,\cos\lambda xf''(x)=-\lambda\sqrt{\frac{2}{\pi}}f'(0)-\lambda^2F_c(\lambda).
\end{equation}
We observe that in order to apply the Fourier sine- and cosine transforms to solve differential equations on a semi infinite domain the boundary conditions must be of a particular type. For second order equations the Fourier sine transform requires Dirichlet condition and the cosine transform requires Neumann condition. Since one typically solve time dependent problems on a semi-infinite domain $t>0$, one might think that the Fourier Sine- and Cosine transforms should be useful for this, but they are in fact not. 

If the equation contains a first derivative with respect to time we will find that the Fourier Sine transform of $u_t$ is expressed in terms of the Fourier Cosine transform for $u$ and vice versa. If the equation contains two time derivatives, $u_{tt}$, there are two data given at $t=0$, $u$ and $u_t$. But both the Fourier Sine and the Fourier Cosine transform can only include one condition at $t=0$, $u$ for the Sine transform and $u_t$ for the Cosine transform. 

We will see that for initial value problems,  the \emph{Laplace transform} is the appropriate one. 

Another useful property of the Fourier Sine- and Cosine transform is 

\begin{align}
&\sqrt{\frac{2}{\pi}}\int_0^{\infty}\mathrm d x\,x\sin\lambda xf(x)=-\frac{\partial F_c}{\partial \lambda},\label{eq12-173}\\
& \sqrt{\frac{2}{\pi}}\int_0^{\infty}\mathrm d x\,x\cos\lambda xf(x)=\frac{\partial F_s}{\partial \lambda}.\label{eq12-174}
\end{align}

\subsubsection{The heat equation on a semi-infinite interval}
We consider the equation 

\begin{equation}\label{eq12-175}
u_t(x,t)-c^2u_{xx}(x,t)=0,\quad 0<x<\infty,\,t>0,
\end{equation}
with initial condition

\begin{equation}\label{eq12-176}
u(x,0)=f(x),\quad 0<x<\infty,
\end{equation}
and boundary condition

\begin{equation}\label{eq12-177}
u(0,t)=g(t),\quad t>0. 
\end{equation}
Because of the Diricelet boundary condition at $x=0$  it is appropriate to use the Fourier sine transform on the $x$-variable. 

Taking the Fourier Sine transform of (\ref{eq12-175}) and make use  properties  (\ref{eq12-170}) and (\ref{eq12-177}), we get

\begin{equation}\label{eq12-178}
\partial_tU_s(\lambda,t)+(\lambda c)^2U_s(\lambda,t)=\lambda c^2\sqrt{\frac{2}{\pi}}g(t),
\end{equation}
with initial condition

\begin{equation}\label{eq12-179}
U_s(\lambda,0)=F_s(\lambda),
\end{equation}
where 

\begin{align}
U_s(\lambda,t)=\sqrt{\frac{2}{\pi}}\int_0^{\infty}\mathrm d x\,\sin\lambda x\,u(x,t)\label{eq12-180},\\
F_s(\lambda)=\sqrt{\frac{2}{\pi}}\int_0^{\infty}\mathrm d x\,\sin\lambda x\,f(x). \label{eq12-181}
\end{align}
The initial value problem (\ref{eq12-178}), (\ref{eq12-179}) is standard with solution

\begin{align}\label{eq12-182}
u_s(\lambda,t)&=F_s(\lambda)e^{-\lambda^2c^2t},\\
&+\lambda c^2 \sqrt{\frac{2}{\pi}} \int_0^{t}\mathrm d \tau\,e^{-\lambda^2c^2(t-\tau)}g(\tau).\nonumber
\end{align}
Taking the inverse Fourier Sine transform (\ref{eq12-104}) we get the solution in the form 

\begin{align}\label{eq12-183}
u(x,t)&= \sqrt{\frac{2}{\pi}}\int_0^{\infty}\mathrm d \lambda\,\sin\lambda x \,F_s(\lambda)e^{-\lambda^2c^2 t}\\
&+\frac{2c^2}{\pi} \int_0^{\infty}\mathrm d \lambda\,\sin\lambda x\,\lambda\int^t_0\mathrm d\tau\,e^{-\lambda^2c^2(t-\tau)}g(\tau)\nonumber\\
&=\frac{2}{\pi}\int_0^{\infty}\mathrm d \lambda\,\sin\lambda x \int_0^{\infty}\mathrm d s\,\sin\lambda s\,f(s) e^{-\lambda^2c^2 t}\nonumber\\
&+\frac{2c^2}{\pi}\int_0^{\infty}\mathrm d\lambda\,\sin\lambda x\,\lambda\int^t_0\mathrm d\tau\,e^{-\lambda^2c^2(t-\tau)}g(\tau)\nonumber\\
&=\frac{2}{\pi} \int_0^{\infty}\mathrm d s\,\left[\int_0^{\infty}\mathrm d\lambda\,\sin\lambda x\sin\lambda s\,e^{-\lambda^2c^2t}\right]f(s)\nonumber\\
&+\frac{2c^2}{\pi}\int^t_0\mathrm d\tau\,\left[\int_0^{\infty}\mathrm d\lambda\,\lambda\sin\lambda x\,e^{-\lambda^2c^2(t-\tau)}\right]g(\tau).\nonumber
\end{align}

Using the addition formulas for cosine we have

\begin{equation}\label{eq12-185}
2\sin\lambda x\sin\lambda s=\cos\lambda(x-s)-\cos\lambda(x+s).
\end{equation}
But then 

\begin{align}\label{eq12-186}
&\frac{2}{\pi}\int_0^{\infty}\mathrm d\lambda\,\sin\lambda x\sin\lambda se^{-\lambda^2c^2t}\\
&=\frac{1}{\pi}\int_0^{\infty}\mathrm d\lambda\,\lbrace\cos\lambda(x-s)-\cos\lambda(x+s)\rbrace e^{-\lambda^2c^2t}\nonumber\\
& =\frac{1}{\pi}\int_0^{\infty}\mathrm d\lambda\,\cos\lambda(x-s) e^{-\lambda^2c^2t}-\frac{1}{\pi}\int_0^{\infty}\mathrm d\lambda\,\cos\lambda(x+s) e^{-\lambda^2c^2t}\nonumber\\
&=\frac{1}{\sqrt{4\pi c^2t}}e^{-\frac{(x-s)^2}{4c^2t}}-\frac{1}{\sqrt{4\pi c^2t}}e^{-\frac{(x+s)^2}{4c^2t}}\nonumber\\
&=G(x-s,t)-G(x+s,t),\nonumber
\end{align}
where we have used quantities defined in (\ref{eq12-98}) and (\ref{eq12-104}). We also have

\begin{align}\label{eq12-187}
& \frac{1}{\pi}\int_0^{\infty}\mathrm d\lambda\,\lambda\sin\lambda x\,e^{-\lambda^2c^2(t-\tau)}\\
&=-\frac{1}{\pi}\partial_x\int_0^{\infty}\mathrm d\lambda\,\cos\lambda x \,e^{-\lambda^2c^2(t-\tau)}\nonumber\\
&=-\partial_xG(x,t-\tau), \nonumber
\end{align}
so the solution can compactly be written as 

\begin{align}\label{eq12-188}
u(x,t)&=\int_0^{\infty}\mathrm d s\,\lbrace G(x-s,t)-G(x+s,t)\rbrace f(s)\\
&-2c^2\int^t_0\mathrm d\tau\,\partial_xG(x,t-\tau)g(\tau).
\end{align}
Let us look at a couple of simple cases. 

For the first case we assume $g(t)=0$, so the boundary condition is homogeneous, and that the initial temperature profile is uniform $u(x,0)=f(x)=u_0$. Then the second term in (\ref{eq12-188}) vanishes and

\begin{align}\label{eq12-189}
u(x,t)&=\int_0^{\infty}\mathrm ds\,\left\{G(x-s,t)-G(x+s,t)\right\}u_0\\
&=\frac{u_0}{\sqrt{4\pi c^2t}}\left\{\int_0^{\infty}\mathrm ds\,e^{-\frac{(x-s)^2}{4c^2t}}-\int_0^{\infty}\mathrm ds\,e^{-\frac{(x+s)^2}{4c^2t}}\right\}\nonumber\\
&\qquad\quad r=\frac{s-x}{2c\sqrt{t}}\qquad\qquad r=\frac{s+x}{2c\sqrt{t}}\nonumber\\
&=\frac{u_0}{\sqrt{\pi}}\left\{\int^{\infty}_{-\frac{x}{2c\sqrt{t}}}\mathrm dr\,e^{-r^2} - \int^{\infty}_{\frac{x}{2c\sqrt{t}}}\mathrm dr\,e^{-r^2} \right\}\nonumber\\
&=\frac{2u_0}{\sqrt{\pi}} \int^{\frac{x}{2c\sqrt{t}}}_{0}\mathrm dr\,e^{-r^2}=u_0 \erf\left(\frac{x}{2c\sqrt{t}}\right),\nonumber
\end{align}
where $\erf(\xi)$ is the \emph{error function}

\begin{equation}\label{eq12-190}
\erf(\xi)=\frac{2}{\sqrt{\pi}}\int_0^{\xi}\mathrm dr\,e^{-r^2}.
\end{equation}
We observe that 

\begin{equation}\label{eq12-191}
\erf(0)=0,\quad \erf(\infty)=1.
\end{equation}
Using the properties of the error function we find that

\begin{equation}\label{eq12-192}
u(x,0)=u_0\erf(\infty)=u_0,
\end{equation}

\begin{equation}\label{eq12-193}
u(0,t)=u_0\erf(0)=0,
\end{equation}

\begin{equation}\label{eq12-194}
u(x,t)\to0\text{ as }t\to\infty.
\end{equation}
Thus $u(x,t)$ satisfies \emph{all} requirements of our problem even if $u(x,0)=u_0$ does not have a Fourier sine transform. One can show that the solution (\ref{eq12-189}) is the limit of a sequence of problems that all have initial data for which the transform exist. 

For the second case we don't assume that $f(x)$ is constant, but merely that it  is bounded

\begin{equation}\label{eq12-195}
|f(x)|\leq M<\infty,\quad\text{for all }x>0.
\end{equation}
We let the boundary data at $x=0$ be given by an arbitrary time dependent function $g(t)$. Observe that for this case, the first term in the solution formula (\ref{eq12-188}) approach zero for large $t$

\begin{gather*}
\Big|\int_0^{\infty}\mathrm ds\,\lbrace G(x-s,t)-G(x+s,t)\rbrace f(s)\Big|\\
\leq M \int_0^{\infty}\mathrm ds\,\lbrace G(x-s,t)-G(x+s,t)\rbrace=M\erf\left(\frac{x}{2c\sqrt{t}}\right)\to 0.
\end{gather*}
Therefore as $t\to\infty$ the solution $u(x,t)$ approaches

\begin{equation}\label{eq12-196}
u(x,t)\approx-2c^2\int^t_0\mathrm d\tau\,\partial_x G(x,t-\tau)g(\tau).
\end{equation}
This solution is called a \emph{steady state}. It is not time invariant but represent the solution for times so large that the initial data has been dissipated. For such large times the solution only depends on the boundary condition, or \emph{drive at the boundary}. Recall that an analogy behavior is displayed by a driven, damped harmonic oscillator. 

\subsection{The Laplace transform}
The last transform we will discuss in this class is the \emph{Laplace transform}. It is, as we will see, well suited for solving initial value problems. 
 We will introduce the Laplace transform using the Fourier transform. 
 
 The Laplace transform is concerned with functions $f(t)$ which are zero for $t<0$.
The function $f(t)$ can be quite wild for $t>0$, but we will set a bound on how fast it can grow when $t\to\infty$. 

We will therefore require $f(t)$ to be of \emph{exponential order} as $t\to\infty$, This means that there are constants $c>0$, $k$ and $M$ such that 

\begin{equation}\label{eq12-197}
|f(t)|\leq Me^{kt},\quad\text{for }t>c.
\end{equation}
Let now $\alpha$ be a constant such that $\alpha>k$. Then the Fourier transform of $e^{-\alpha t}f(t)$ will exist and is given by

\begin{equation}\label{eq12-198}
\hat{F}(\hat{\lambda})=\frac{1}{\sqrt{2\pi}}\int^{\infty}_0\mathrm dt\,e^{i\hat{\lambda}t}e^{-\alpha t}f(t),
\end{equation}
where the fact that $f(t)=0$ for $t\leq 0$ has been taken into account. For example if $f(t)=1$ we get

\begin{equation}\label{eq12-199}
\hat{F}(\hat{\lambda})=\frac{1}{\sqrt{2\pi}}\int^{\infty}_0\mathrm dt\,e^{(i\hat{\lambda}-\alpha)t}=\frac{-1}{\sqrt{2\pi}(i\hat{\lambda}-\alpha)}.
\end{equation}
In the Fourier transform we know that $\hat{\lambda}$ is a real number related to the spectrum of a certain differential operator. Observe, however, that the right hand side of (\ref{eq12-199}) can be evaluated for complex $\hat{\lambda}$ also. For example

\begin{equation}\label{eq12-200}
\hat{F}(i)=\frac{-1}{\sqrt{2\pi}(i^2-\alpha)}=\frac{1}{\sqrt{2\pi}(1+\alpha)}.
\end{equation}
Also observe that (\ref{eq12-199}) can not be evaluated for \emph{all} complex numbers, in fact for $\hat{\lambda}=-i\alpha$ the denominator in (\ref{eq12-199}) is zero. However, this is the only problematic point and $\hat{F}(\hat{\lambda})$ can be extended to the whole complex plane otherwise

\begin{equation}\label{eq12-201}
\hat{F}(\hat{\lambda})=\frac{-1}{\sqrt{2\pi}(i\hat{\lambda}-\alpha)},\quad\hat{\lambda}\in\mathbb C-\lbrace -i\alpha\rbrace.
\end{equation}
In a similar way we will now prove that $\hat{F}(\lambda)$ can be extended to complex values for any $f(t)$ of exponential order.

Evaluating (\ref{eq12-198})  at $\hat{\lambda}=\hat{\lambda}_r+i\hat{\lambda}_i\in\mathbb C$ gives us

\begin{equation}\label{eq12-202}
\hat{F}(\hat{\lambda})=\frac{1}{\sqrt{2\pi}}\int^{\infty}_0\mathrm dt\,e^{(i\hat{\lambda}-\alpha)t}f(t).
\end{equation}
For this to make sense the integral has to converge. But

\begin{align}\label{eq12-203}
|\hat{F}(\hat{\lambda})|&\leq\frac{1}{\sqrt{2\pi}}\int^{\infty}_0\mathrm dt\,e^{-(\hat{\lambda_i}+\alpha)t}|f(t)|\\
& \leq\frac{1}{\sqrt{2\pi}}\int^{c}_0\mathrm dt\,e^{-(\hat{\lambda_i}+\alpha)t}|f(t)|\nonumber\\
&+ \frac{1}{\sqrt{2\pi}}\int^{\infty}_c\mathrm dt\,Me^{(k-\hat{\lambda_i}-\alpha)t}\nonumber
\end{align}
and the last integral is finite only if

\begin{equation}\label{eq12-204}
\hat{\lambda_i}>-(\alpha-k)\quad\text{(Recall that }\alpha-k>0\text{).}
\end{equation}
So $\hat{F}(\lambda)$ can not be defined in the whole complex plane but only in a domain defined by (\ref{eq12-204}). 

\begin{figure}[h!]
	\centering
	\includegraphics[scale=0.5]{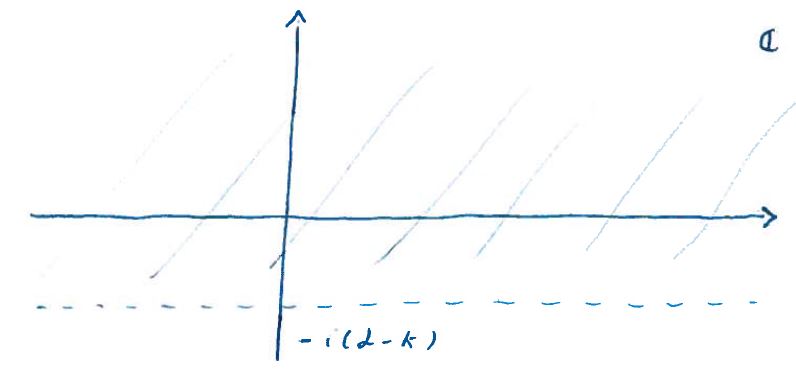}
	\caption{}
\end{figure}
Below the horizontal line $\hat{\lambda}=-i(\alpha-k)$ in the complex plane, we expect the function $\hat{F}(\hat{\lambda})$ to have singularities. For the special choice 

\begin{equation}\label{eq12-205}
f(t)=\left\{
\begin{aligned}
1,\quad t\geq 0,\\
0,\quad t<0,
\end{aligned}
\right.
\end{equation}
this was exactly what occured.

  Change variables in the complex plane to $\lambda$ through 

\begin{equation}\label{eq12-206}
\hat{\lambda}=i(\lambda-\alpha)
\end{equation} 
and define a function $F(\lambda)$ by 

\begin{equation}\label{eq12-207}
F(\lambda)=\sqrt{2\pi}\hat{F}(i(\lambda-\alpha)).
\end{equation}
From (\ref{eq12-202}) we then get 

\begin{equation}\label{eq12-208}
F(\lambda)=\int_0^{\infty}\mathrm dt\,e^{-\lambda t}f(t).
\end{equation}
For a function of exponential order $k$, $F(\lambda)$ is defined for complex arguments $\lambda=\lambda_r+i\lambda_i$, with 

\begin{equation}\label{eq12-209}
\lambda_r>k.
\end{equation}
\begin{figure}[h!]
	\centering
	\includegraphics[scale=0.6]{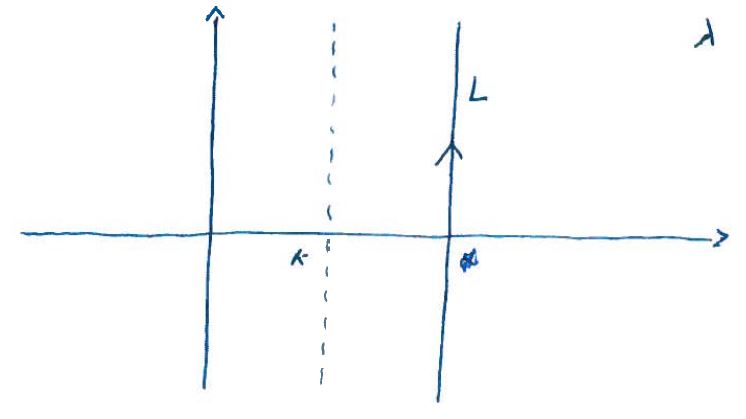}
	\caption{}
	\label{Fig12-6}
\end{figure}

$F(\lambda)$ is the \emph{Laplace transform} of $f(t)$. Using the Fourier inversion formula we have 

\begin{equation}\label{eq12-211}
e^{-\alpha t}f(t)=\frac{1}{\sqrt{2\pi}}\int_{-\infty}^{\infty}\mathrm d\hat{\lambda}\,e^{-\hat{\lambda}t}\hat{F}(\hat{\lambda}).
\end{equation}
Introducing $\hat{\lambda}=i(\lambda-\alpha)$ gives us $\mathrm d\hat{\lambda}=i\,\mathrm d\lambda$ so that (\ref{eq12-211}) implies

\begin{equation}\label{eq12-212}
f(t)=\frac{1}{2\pi i}\int\limits_L\mathrm d\lambda\,e^{\lambda t}F(\lambda),
\end{equation}
where $L$ is a parameterized curve 

\begin{equation}\label{eq12-213}
L:\quad\lambda=\alpha+iu,\quad -\infty<u<\infty.
\end{equation}
The curve $L$ is indicated in Figure \ref{Fig12-6}. 

The formula (\ref{eq12-212}) is the \emph{inversion formula} for the Laplace transform. The inversion formula is a very effective and flexible way to calculate the inverse Laplace transform, but using it requires \emph{residue calculus} from complex analysis, so we will not make further use of formula (\ref{eq12-212}) in this class. Inverse Laplace transforms can in many common situations be found using a combination of formal properties of the inverse Laplace transform and tables of inverse transforms. Lists of such formal properties can be found in mathematical tables or on the web. Here we will only display the most important property we need when we solve initial value problems for differential equations using the Laplace transform

\begin{equation}\label{eq12-214}
\int_0^{\infty}\mathrm dt\,e^{-\lambda t}f^{(n)}(t)=\lambda^n \int_0^{\infty}\mathrm dt\,e^{-\lambda t}f(t)-\lambda^{n-1}f(0)-\dots f^{(n-1)}(0).
\end{equation}
This formula is easy to prove using integration by parts. From (\ref{eq12-214}) it is clear that the $n$ initial conditions for a $n^{\text{th}}$ order differential equation can be included in a natural way. The Laplace transform is therefore very useful for solving initial value problem for both ODEs and PDEs. Here we will apply it to initial boundary value problems for PDEs. 

We will now apply the Laplace transform to two initial value problems for PDEs.

\subsubsection{The wave equation on a half interval}
\begin{align}\label{eq13-1}
& u_{tt}-c^2u_{xx}=0,\quad & 0<x<\infty,\,t>0,\\
& u(x,0)=u_t(x,0)=0,&\nonumber\\
& u(0,t)=f(t),&\nonumber\\
& u(x,t)\to 0,\quad & x\to\infty.\nonumber
\end{align}
Taking the Laplace transform of the equation we get

\begin{equation}\label{eq13-2}
\lambda^2U(x,\lambda)-\lambda u(x,0)-\lambda^2u_t(x,0)-c^2U_{xx}(x,\lambda)=0
\end{equation}
Using the initial conditions we get 

\begin{gather}\label{eq13-3}
\lambda^2U(x,\lambda)-c^2U_{xx}(x,\lambda)=0,\\
U(0,\lambda)=F(\lambda).\nonumber
\end{gather}
The general solution to the family of ordinary differential equations (\ref{eq13-3}) is

\begin{equation}\label{eq13-4}
U(x,\lambda)=a(\lambda)e^{-\frac{\lambda}{c}x}+b(\lambda)e^{\frac{\lambda}{c}x}.
\end{equation}
Recall that we evaluate (\ref{eq13-4}) in the right halfspace only, so $\lambda_r>0$. In order to satisfy the boundary condition 

\begin{equation}\label{eq13-5}
u(x,t)\to 0,\,x\to\infty\Leftrightarrow U(x,\lambda)\to 0,\,x\to\infty.
\end{equation}
We must clearly assume that 

\begin{equation}\label{eq13-6}
b(\lambda)=0.
\end{equation}
Finally applying the boundary condition at $x=0$ we get 

\begin{equation}\label{eq13-7}
U(x,\lambda)=F(\lambda)e^{-\frac{\lambda}{c}x}.
\end{equation}
We now must take the inverse Laplace transform of (\ref{eq13-7}). 

Observe that if $H(\xi)$ is the Heaviside function

\begin{equation}\label{eq13-8}
H(\xi)=\left\{
\begin{aligned}
1,\quad\xi\geq 0,\\
0,\quad\xi<0.
\end{aligned}
\right.
\end{equation}
Then 

\begin{align}\label{eq13-9}
&\int_0^{\infty}\mathrm dt\,H(t-b)f(t-b)e^{-\lambda t}\\
&= \int_b^{\infty}\mathrm dt\,f(t-b)e^{-\lambda t},\quad\xi=t-b,\,\mathrm d\xi\,=\mathrm dt,\nonumber\\
&=\int_0^{\infty}\mathrm d\xi\,f(\xi)e^{-\lambda(\xi+b)}=e^{-\lambda b} F(\lambda).\nonumber
\end{align}
This is equal to (\ref{eq13-7}) with $b=\frac{x}{c}$. Thus the inverse Laplace transform of (\ref{eq13-7}) must be

\begin{equation}\label{eq13-10}
u(x,t)=H\left(t-\frac{x}{c}\right)f\left(t-\frac{x}{c}\right)
\end{equation}
and we have solved the initial value problem.

\subsubsection{The heat equation on a finite interval}
Consider the problem 

\begin{align}\label{eq12-215}
& u_t(x,t)-c^2u_{xx}(x,t)=0,\quad 0<x<l,\, t>0,\\
& u(0,t)=u(l,t)=0,\quad t>0,\nonumber\\
& u(x,0)=f(x),\quad 0<x<l.
\end{align} 
Taking the Laplace transform of (\ref{eq12-215}) we get

\begin{align}\label{eq12-216}
&\lambda U(x,\lambda)-f(x)-c^2U_{xx}(x,\lambda)=0,\quad 0<x<l,\\
& U(0,\lambda)=U(l,\lambda)=0.\nonumber
\end{align}
This is now a boundary value problem. Equation (\ref{eq12-216}) can be solved using variation of parameters from the theory of ODEs. Fitting the boundary conditions to the variation of parameter solution gives

\begin{align}\label{eq12-217}
U(x,\lambda)&=\frac{1}{c\sqrt{\lambda}\sinh\left(\frac{l\sqrt{\lambda}}{c}\right)}\int_0^x\mathrm ds\,\sinh\left({\frac{\sqrt{\lambda}}{c}s}\right) \sinh\left(\frac{\sqrt{\lambda}}{c}(l-x)\right)f(s)\\
& +\frac{1}{c\sqrt{\lambda}\sinh\left(\frac{l\sqrt{\lambda}}{c}\right)}\int_x^l\mathrm ds\,\sinh\left(\frac{\sqrt{\lambda}}{c}x\right) \sinh\left(\frac{\sqrt{\lambda}}{c}(l-s)\right)f(s).\nonumber
\end{align}
Finding the inverse transform of (\ref{eq12-217}) is not easy, but we have a secret weapon. In the inversion formula for the Laplace transform (\ref{eq12-212}) we can place the integration curve $L$ anywhere as long as $\lambda_r>k$. Thus we may assume in (\ref{eq12-217}) that $|\lambda|$ is as large as we want. 

For such $\lambda$ we have

\begin{equation}\label{eq12-218}
\frac{1}{\sinh\left(\frac{\sqrt{\lambda}l}{c}\right)}=\frac{2}{e^{\frac{l\sqrt{\lambda}}{c}}-e^{-\frac{l\sqrt{\lambda}}{c}}}=\frac{2e^{-\frac{l\sqrt{\lambda}}{c}}}{1-e^{-2\frac{l\sqrt{\lambda}}{c}}}=2e^{-\frac{l\sqrt{\lambda}}{c}}\sum\limits_{j=0}^{\infty}e^{-\frac{2jl\sqrt{\lambda}}{c}}.
\end{equation}
The last equality sign in (\ref{eq12-218})  is valid because for $|\lambda|$ large enough $\big|e^{-2l\sqrt{\lambda}/c}\big|<1$ and therefore 

\begin{equation}\label{eq12-219}
\frac{1}{1-x}=\sum\limits_{j=0}^{\infty}x^j.
\end{equation} 
The hyperbolic sine under the integral signs can also be expressed in terms of exponentials. Thus we end up with an infinite sum where each term is an exponential.

From a table of Laplace transforms we have

\begin{equation}\label{eq12-220}
\int_0^{\infty}\mathrm dt\,e^{-\lambda t}\frac{1}{\sqrt{\pi t}}e^{-\frac{\alpha^2}{4t}}=\frac{1}{\sqrt{\lambda}}e^{-\alpha\sqrt{\lambda}},\quad\alpha>0.
\end{equation}
Inserting (\ref{eq12-220}) into the infinite sum, formally 
interchanging summation and integration gives us finally

\begin{equation}\label{eq12-221}
u(x,t)=\int_0^l\mathrm ds\,\sum\limits_{j=-\infty}^{\infty}\lbrace G(x-s-2jl,t)-G(x+s-2jl,t)\rbrace f(s),
\end{equation}
where $G(y,t)$ is the fundamental solution to the heat equation. We have previously found the solution to (\ref{eq12-215}) using a Fourier sine series

\begin{equation}\label{eq12-222}
u(x,t)=\sqrt{\frac{2}{l}}\sum\limits_{j=1}^{\infty}a_je^{-\left(\frac{\pi jc }{l}\right)^2t}\sin\left(\frac{\pi j}{l}x\right).
\end{equation}
Problem (\ref{eq12-215}) has a unique solution so even if (\ref{eq12-221}) and (\ref{eq12-222}) \emph{appear} to be very different they must in fact be identical. 

The formulas can be proven to be identical using the \emph{Poisson summation formula}\cite{textbook}.

The formulas are mathematically equivalent, but are useful for evaluating $u(x,t)$ for different ranges of $t$. Because of the exponential term in (\ref{eq12-222}), evaluating $u(x,t)$ for large $t$ using this formula requires only a few terms. 

However for small $t$, $j$ must be large before the exponential term kicks in. Thus we require many terms in (\ref{eq12-222})  in order to evaluate $u(x,t)$ for small $t$. 

On the other hand, for small $t$, the fundamental solution $G(y,t)$ is sharply peaked around $y=0$. Thus for small $t$, only a few terms in (\ref{eq12-221}) contribute to $u(x,t)$. Thus for small $t$ (\ref{eq12-221}) is fast to evaluate.

In general it is not easy to invert integral transforms. There are however powerful analytic approximation methods that can be brought to bear on this problem.

Observe that the formula for the Laplace transform

\begin{equation}\label{eq12-223}
F(\lambda)=\int_{0}^{\infty}\mathrm dt\, e^{-\lambda t}f(t)
\end{equation}
clearly shows that $F(\lambda)$ for large $|\lambda|$ depends mainly on $f(t)$ for small $t$. This can be made more precise. From the body of results called \emph{Abelian asymptotic theory} we have the following result:\\ 
If

\begin{equation}\label{eq12-224}
f(t)\approx\sum\limits_{n=0}^{\infty}\alpha_nt^n,\quad t\to 0,
\end{equation}
then 

\begin{equation}\label{eq12-225}
F(\lambda)\approx\sum\limits_{n=0}^{\infty}\alpha_n\frac{n!}{\lambda^{n+1}},\quad|\lambda|\to\infty.
\end{equation}
The converse in also true, (\ref{eq12-225}) implies (\ref{eq12-224}). Thus if we are only interested in $f(t)$ for small $t$ we can find the expression for $f(t)$ by expanding $F(\lambda)$ in inverse powers of $\lambda$. This expansion gives us the coefficients $\alpha_n$ which we then insert into (\ref{eq12-224}).

There are also results relating the behavior of $f(t)$ for large $t$ to the behavior of $F(\lambda)$ for small $\lambda$. This body of results is called \emph{Tauherian asymptotic theory}. Describing the simplest results from this theory requires however complex analysis so we will not do it here.

\subsection{Group velocity 1} 
Let us consider the Klein-Gordon equation

\begin{equation}\label{eq12-226}
u_{tt}(x,t)-\gamma^2u_{xx}(x,t)+c^2u(x,t)=0,\quad t>0,\,-\infty<x<\infty.
\end{equation}
Initial data is of the form

\begin{align}\label{eq12-227}
u(x,0)&=f(x),\\
u_t(x,0)&=g(x).\nonumber
\end{align}
Recall that the Klein-Gordon equation is of dispersive type. We have found this previously using normal mode analysis. Taking the Fourier transform of the Klein-Gordon equation we get 

\begin{equation}\label{eq12-228}
\partial_{tt}U(\lambda,t)+(\gamma^2\lambda^2+c^2)U(\lambda,t)=0.
\end{equation}
The general solution to this equation is 

\begin{equation}\label{eq12-229}
U(\lambda,t)=F_+(\lambda)e^{i\omega t}+F_{-}(\lambda)e^{-i\omega t},
\end{equation}
where $\omega=\omega(\lambda)$ is the dispersion relation for the Klein-Gordon equation 

\begin{equation}\label{eq12-230}
\omega=(\gamma^2\lambda^2+c^2)^{\frac{1}{2}}.
\end{equation}
Using the inverse Fourier transform we get the following general solution to (\ref{eq12-226}) 

\begin{align}\label{eq12-231}
u(x,t)&=\frac{1}{\sqrt{2\pi}}\int^{\infty}_{-\infty}\mathrm d\lambda\,F_+(\lambda)e^{i(\omega(\lambda)t-\lambda x)}\\
&+ \frac{1}{\sqrt{2\pi}}\int^{\infty}_{-\infty}\mathrm d\lambda\,F_-(\lambda)e^{-i(\omega(\lambda)t+\lambda x)}.\nonumber
\end{align}
The functions $F_+(\lambda)$ and $F_-(\lambda)$ are determined from initial data.

From (\ref{eq12-227}) the Fourier transform of the initial data is

\begin{align}\label{eq12-232}
U(\lambda,0)&=F(\lambda)\equiv\frac{1}{\sqrt{2\pi}}\int_{-\infty}^{\infty}\mathrm dx\,f(x)e^{i\lambda x},\\
U_t(\lambda,0)&=G(\lambda)\equiv \frac{1}{\sqrt{2\pi}}\int_{-\infty}^{\infty}\mathrm dx\,g(x)e^{i\lambda x}.\nonumber
\end{align}
Using (\ref{eq12-229}) we get 

\begin{align}\label{eq12-233}
F_+(\lambda)+F_-(\lambda)&=F(\lambda),\\
i\omega F_+(\lambda)-i\omega F_-(\omega)&=G(\lambda).\nonumber
\end{align}
Solving (\ref{eq12-233}) gives us 

\begin{align}\label{eq12-234}
F_+(\lambda)&=\frac{i\omega F(\lambda)+G(\lambda)}{2i\omega},\\
F_-(\lambda)&=\frac{i\omega F(\lambda)-G(\lambda)}{2i\omega}.\nonumber
\end{align}
(\ref{eq12-234}) together with (\ref{eq12-231}) solves the initial value problem (\ref{eq12-227}) for the Klein-Gordon equation (\ref{eq12-226}).

The solution (\ref{eq12-231}) is clearly a superposition of normal modes

\begin{equation}\label{eq12-235}
F_{\pm}(\lambda)e^{i(\pm\omega(\lambda)t-\lambda x)}.
\end{equation}
The phase-speed for the normal modes are

\begin{equation}\label{eq12-236}
v_f(\lambda)=\frac{\mathrm d x}{\mathrm d t}=\pm\frac{\omega(\lambda)}{\lambda}=\pm\sqrt{\gamma^2+\left(\frac{c}{\lambda}\right)^2}.
\end{equation}
Recall that the phase speed of a normal mode is the speed we have to move at in order for the phase

\begin{equation}\label{eq12-237}
\theta_{\pm}(x,t)=\pm\omega(\lambda)t-\lambda x,
\end{equation}
to stay fixed. Thus we have 

\begin{gather}\label{eq12-238}
\theta_{\pm}(x(t),t)=\theta_{\pm}^0,\\
\Downarrow\nonumber\\
\partial_x\theta_{\pm}\frac{\mathrm dx}{\mathrm dt}+\partial_t\theta_{\pm}=0,\nonumber\\
\Downarrow\\
-\lambda\frac{\mathrm dx}{\mathrm dt}\pm\omega(\lambda)=0,\nonumber\\
\Downarrow\nonumber\\
\frac{\mathrm dx}{\mathrm dt}=\pm\frac{\omega(\lambda)}{\lambda}.\nonumber
\end{gather}
The Klein-Gordon equation has the same principle part as the wave equation. It thus have the same characteristic speed $\pm\gamma$.

As we have seen, discontinuities move at characteristic speed $\pm\gamma$ and $\pm\gamma$ represents the upper speed limit for the Klein-Gordon equation. 

But from (\ref{eq12-236}) we see that the phase speed is bigger than $\gamma$ for all modes. The physical relevance of the phase speed is therefore in question. 
Furthermore, noting that the solution to the Klein-Gordon equation is a superposition of normal modes which all move at different phase velocities, it is not at all clear that the phase speed is a relevant quantity for this equation. 

In order to investigate this problem we will consider special initial conditions of the form

\begin{align}\label{eq12-239}
u(x,0)&=p(\epsilon x)e^{-ik_0x},\\
u_t(x,0)&=q(\epsilon x)e^{-ik_0x},\nonumber
\end{align}
where $k_0$ is fixed. The function $e^{ik_0x}$ is clearly periodic of period

\begin{equation}\label{eq12-240}
L_0=\frac{2\pi}{k_0}.
\end{equation}
In (\ref{eq12-239}) we assume that $0<\epsilon<<1$ is a number so small that 

\begin{equation}\label{eq12-241}
\epsilon L_0<<1.
\end{equation}
Then the amplitudes $p(\epsilon x)$ and $q(\epsilon x)$ vary very little over a period $L_0$. The initial conditions (\ref{eq12-239}) are what we call \emph{wave packets}. These are functions that locally look like plane waves but that has a (complex) amplitude that varies over scales much larger than $L_0$.

\begin{figure}[h!]
	\centering
	\includegraphics[scale=0.5]{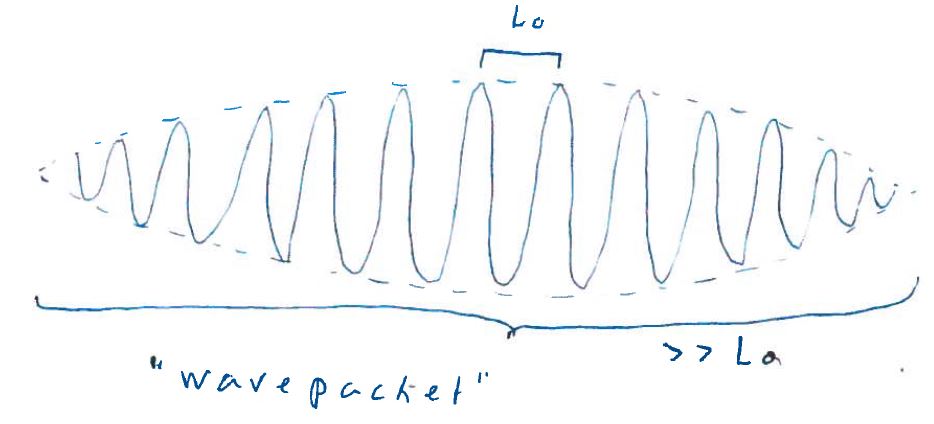}
	\caption{}	
\end{figure} 
Using simple properties of the Fourier transform we find easily from (\ref{eq12-231}) that 

\begin{align}\label{eq12-242}
F(\lambda)&=\frac{1}{\epsilon}P\left(\frac{\lambda-k_0}{\epsilon}\right),\\
G(\lambda)&=\frac{1}{\epsilon}Q\left(\frac{\lambda-k_0}{\epsilon}\right),\nonumber
\end{align}
where $P$ and $Q$ are the Fourier transform of $p$ and $q$.
\begin{figure}[h!]
	\centering
	\includegraphics[scale=0.5]{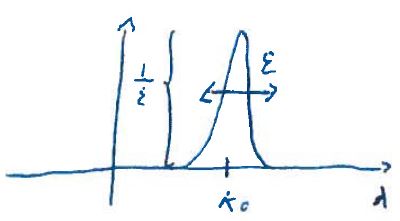}
	\caption{}	
	\label{fig12-8}
\end{figure} 
 From Figure \ref{fig12-8} it is clear that in Fourier space, $F(\lambda)$ and $G(\lambda)$ are narrow "spikes" of width $\epsilon$, height $\frac{1}{\epsilon}$ and centered at $\lambda=k_0$. 

The functions $F_+(\lambda)$, $F_-(\lambda)$ will also be of this form

\begin{equation}\label{eq12-243}
F_{\pm}(\lambda)=\frac{1}{\epsilon}A_{\pm}\left(\frac{\lambda-k_0}{\epsilon}\right),
\end{equation}
for some functions $A_{\pm}$.

  The solution of the Klein-Gordon equation consists of two terms of a very similar structure. The first term is 

\begin{align}\label{eq12-244}
u_+(x,t)&=\frac{1}{\sqrt{2\pi}}\int^{\infty}_{-\infty}\mathrm d\lambda\,F_+(\lambda)e^{i(\omega(\lambda)t-\lambda x)}\\
&=\frac{1}{\sqrt{2\pi}}\int^{\infty}_{-\infty}\mathrm d\lambda\,\frac{1}{\epsilon}A_+\left(\frac{\lambda-k_0}{\epsilon}\right)e^{i\omega(\lambda)t}e^{-i\lambda x}.\nonumber
\end{align}
Consider the following picture

\begin{figure}[h!]
	\centering
	\includegraphics[scale=0.5]{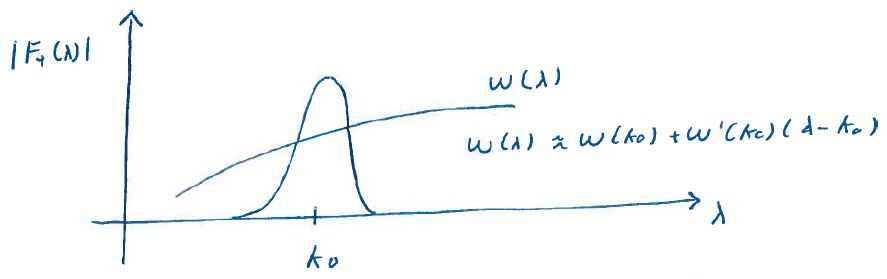}
	\caption{}
	\label{fig12-81}
\end{figure}
\noindent{From figure \ref{fig12-81} it is clear that if the wavepacket is narrow enough, the dispersion relation $\omega(\lambda)$ is well approximated by its first order Taylor polynomial within the support of $F_+(\lambda)$.} Thus we have 

\begin{align}\label{eq12-245}
u_+(x,t)&\approx\frac{1}{\sqrt{2\pi}}\int_{-\infty}^{\infty}\mathrm d\lambda\, \frac{1}{\epsilon}A_+\left(\frac{\lambda-k_0}{\epsilon}\right)e^{i(\omega(k_0)+\omega'(k_0)(\lambda-k_0))t}e^{-i\lambda x}\\
&=e^{i(\omega(k_0)t-k_0x)}\frac{1}{\sqrt{2\pi}}\int^{\infty}_{-\infty}\mathrm d\lambda\,\frac{1}{\epsilon}A_+\left(\frac{\lambda-k_0}{\epsilon}\right)e^{-i((\lambda-k_0)x-\omega'(k_0)(\lambda-k_0)t}.\nonumber
\end{align}
Changing integration variable in the integral 

\begin{equation}\label{eq12-246}
\eta=\frac{\lambda-k_0}{\epsilon}\,\Rightarrow\,\mathrm d\eta=\frac{1}{\epsilon}\mathrm d\lambda,
\end{equation}
gives us

\begin{equation}\label{eq12-247}
u_+(x,t)\approx e^{i(\omega(k_0)t-k_0x)}\frac{1}{\sqrt{2\pi}}\int_{-\infty}^{\infty}\mathrm d\eta\,f_+(\eta)e^{-i\eta\epsilon(x-\omega'(k_0)t)}.
\end{equation}
Defining 

\begin{equation}\label{eq12-248}
a_+(\xi)=\frac{1}{\sqrt{2\pi}} \int_{-\infty}^{\infty}\mathrm d\eta\,A_+(\eta)e^{-i\eta\xi},
\end{equation}
we get the final expression

\begin{equation}\label{eq12-249}
u(x,t)=a_+(x-\omega'(k_0)t)e^{i(\omega(k_0)t-k_0x)}.
\end{equation}
We observe that the \emph{envelope} of the wave packet, $a_+$, translates without changing shape at a speed

\begin{equation}\label{eq12-250}
v_g\equiv\frac{\mathrm d\omega}{\mathrm d\lambda}(k_0).
\end{equation}
$v_g$ is the \emph{group velocity} for the wave packet. The group velocity, and not the phase velocity, is the relevant speed for a wave packet. 
This conclusion holds true for all equations of dispersive type. 

For the Klein-Gordon equation

\begin{gather}\label{eq12-251}
\omega=(\gamma^2\lambda^2+c^2)^{\frac{1}{2}},\\
\Downarrow\nonumber\\
v_g=\frac{\gamma^2\lambda}{(\gamma^2\lambda^2+c^2)^\frac{1}{2}}=\frac{\gamma}{\sqrt{1+\left(\frac{c}{\gamma\lambda}\right)^2}}\,<\gamma.\nonumber
\end{gather}
The group velocity is always smaller than the characteristic velocity $\gamma$ for the Klein-Gordon equation

\begin{figure}[h!]
	\centering
	\includegraphics[scale=0.4]{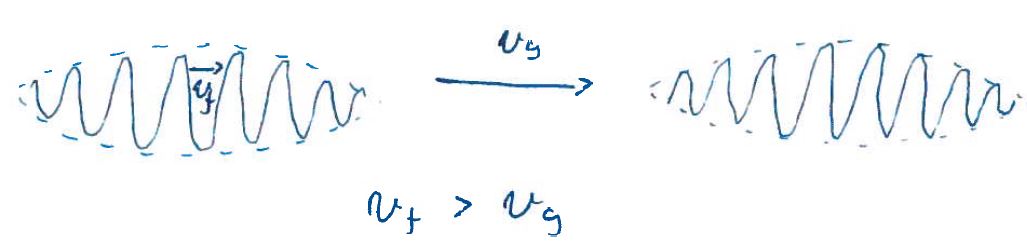}
	\caption{}
\end{figure}
\noindent Note that the statement that the wave packet move at the group velocity is only approximately true. A closer investigation shows that eventually the wave packet spreads out and disperse. 

We will shortly see that the group velocity plays a pivotal role in the description of general solutions to dispersive equations. However, before  we do this we must describe a tool of outmost importance. This is the \emph{method of stationary phase}. 

\subsection{Method of stationary phase}
The method of stationary phase is part of a very large and important domain of applied mathematics called \emph{asymptotic approximations}, or \emph{perturbation methods}. 

The method of stationary phase seeks to approximate integrals of complex integrands, where the \emph{phase} of the integrand has a component that varies very quickly.

Let us consider the integral 

\begin{equation}\label{eq12-252}
I(k)=\int_{-\infty}^{\infty}\mathrm dt\,f(t)e^{ik\varphi(t)},
\end{equation}
where $\varphi(t)$ is a real function and $f(t)$ can be a complex function. Both functions are assumed to be smooth. 

We want to investigate the integral for large $k$. The formulas we derive will hold in the limit $k\to\infty$. In any actual application of the formulas, $k$ is a fixed number, that is large, but not infinite. The errors we make in applying the asymptotic formulas using a given finite $k$, are not so easy to estimate. That is why we call it \emph{method} of stationary phase rather than \emph{theory} of stationary phase: Proofs are few and far between in this domain of applied mathematics. The method of stationary phase and asymptotic methods in general tend however to be surprisingly accurate even for values of $k$ that are nowhere "close" to infinity. 

Observe that when $k$ is very large the phase term 

\begin{equation}\label{eq12-253}
e^{ik\varphi(t)},
\end{equation}  
performs very fast oscillations. Since $f(t)$ is assumed to be smooth, the phase (\ref{eq12-253}) will oscillate on a scale much faster than the scale at which the amplitude, $f(t)$, varies. 

\begin{figure}[h!]
	\centering
	\includegraphics[scale=0.5]{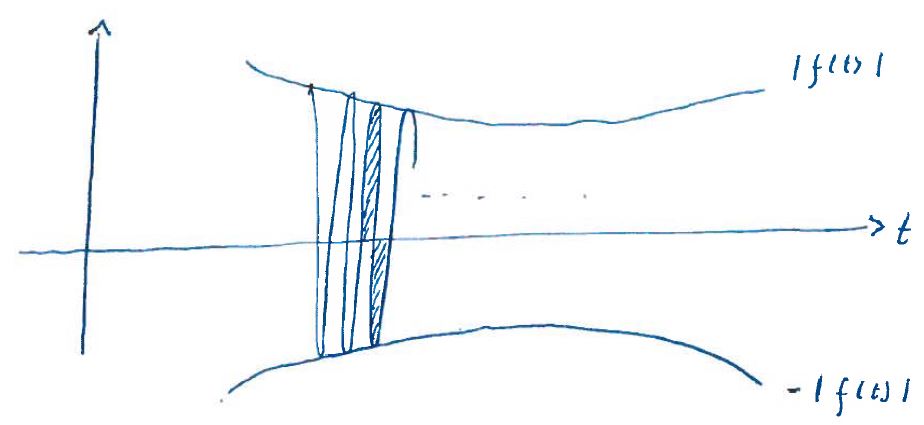}
	\caption{}
\end{figure}
\noindent{In this situation, as the figure illustrates, the positive and negative contributions in the integral tends to cancel each other so that we get}
\begin{equation}\label{eq12-254}
I(k)\to 0,\quad k\to\infty.
\end{equation}
This argument can be made into precise mathematics using the \emph{Riemann-Lebesques} lemma from complex analysis, but we will not pursue this any further here.

Let $t=t_0$. Close to $t_0$ we have 

\begin{equation}\label{eq12-255}
\varphi(t)=\varphi(t_0)+\varphi'(t_0)(t-t_0)+\frac{1}{2}\varphi''(t_0)(t-t_0)^2+\dots
\end{equation}
The largest contribution to $I(k)$ comes close to points where 

\begin{equation}\label{eq12-256}
\varphi'(t_0)=0,
\end{equation}
because around such points the phase varies more slowly than around points where $\varphi'(t_0)\neq 0$ for $k$ large. 

A point $t_0$ such that (\ref{eq12-256}) hold is called a \emph{stationary point}. If $t=t_0$ is a stationary point, the contribution to $I(k)$, which we call $I(k;t_0)$,  is approximately given by 

\begin{align}\label{eq12-257}
I(k;t_0)&\approx\int_{t_0-\epsilon}^{t_0+\epsilon}\mathrm dt\,f(t)e^{i(\varphi(t_0)+\frac{1}{2}\varphi''(t_0)(t-t_0)^2)k}\\
&\approx f(t_0)e^{i\varphi(t_0)k}\int_{t_0-\epsilon}^{t_0+\epsilon}\mathrm dt\,e^{i\frac{k}{2}\varphi''(t_0)(t-t_0)^2},\nonumber
\end{align}
where we have used the fact that for a small enough interval $(t_0-\epsilon,t_0+\epsilon)$, $f(t)$ varies very little around the value $f(t_0)$. This is where we use the smoothness of $f(t)$. If $f(t)$ has a singularity at $t=t_0$ the approximation (\ref{eq12-257}) is not valid. 

\noindent If $k$ is large enough the integrand will oscillate very fast beyond $(t_0-\epsilon,t_0+\epsilon)$, so by the same argument as the one leading up to (\ref{eq12-254}), we can conclude that the contribution from the integral 

\begin{equation}\label{eq12-258}
\int\limits_{|t-t_0|>\epsilon}\mathrm dt\,e^{i\frac{k}{2}\varphi''(t_0)(t-t_0)^2}
\end{equation}
is negligible with respect to the contribution from the integral

\begin{equation}\label{eq12-259}
\int\limits_{|t-t_0|<\epsilon}\mathrm dt\,e^{i\frac{k}{2}\varphi''(t_0)(t-t_0)^2}.
\end{equation} 
Thus we can without loss of generality write (\ref{eq12-257}) in the form

\begin{equation}\label{eq12-260}
I(k;t_0)\approx f(t_0)e^{i\varphi(t_0)k}\int_{-\infty}^{\infty}\mathrm dt\,e^{i\frac{k}{2}\varphi''(t_0)(t-t_0)^2}.
\end{equation}
The integral in (\ref{eq12-260}) can be evaluated exactly using complex analysis

\begin{equation}\label{eq12-261}
I(k;t_0)\approx f(t_0)e^{i\varphi(t_0)k}\sqrt{\frac{2\pi}{k|\varphi''(t_0)|}}e^{i\frac{\pi}{4}\frac{\varphi''(t_0)}{|\varphi''(t_0)|}}.
\end{equation}
If you have several stationary points, each point gives a contribution of type (\ref{eq12-261}) so that finally 

\begin{equation}\label{eq12-262}
I(k)\approx\sum\limits_{j=1}^n I(k;t_j),
\end{equation}
where $\varphi'(t_j)=0,\,j=1,2,\dots,n.$

The method of stationary phase can be extended to integrals over finite domains and to multidimensional integrals.

The importance of the method of stationary phase can hardly be overstated. For example, it is by using this method that we understand how our macroscopic everyday world, ruled by classical physics, arise from a substratum ruled by the laws of quantum physics. 

\subsection{Group velocity 2} 
Let us now return to the Klein-Gordon equation (\ref{eq12-226}), whose solution satisfying the general initial data (\ref{eq12-227}) is 

\begin{align}\label{eq12-263}
u(x,t)&=\frac{1}{\sqrt{2\pi}}\int^{\infty}_{-\infty}\mathrm d\lambda\,F_+(\lambda)e^{i(\omega(\lambda)t-\lambda x)}\\
& +\frac{1}{\sqrt{2\pi}}\int^{\infty}_{-\infty}\mathrm d\lambda\,F_-(\lambda)e^{-i(\omega(\lambda)t+\lambda x)},\nonumber
\end{align}
where $F_+$ and $F_-$ are defined in terms of the initial data in (\ref{eq12-234}). Both terms in (\ref{eq12-263}) have a similar structure, so let us focus our attention on the first of them 

\begin{equation}\label{eq12-264}
I(x,t)\equiv\frac{1}{\sqrt{2\pi}}\int^{\infty}_{\infty} \mathrm d\lambda\, F_+(\lambda)e^{+i(\omega(\lambda)t-\lambda x)}.
\end{equation}
We want to investigate this expression for both $x$ and $t$ large. For such large values of $x$ and $t$ the phase varies very quickly as a function of $\lambda$ and we will, according to the method of stationary phase, get the largest contributions to $I(x,t)$ close to values of $\lambda$  where the phase is stationary. Thus we seek $\lambda$ such that 

\begin{gather}\label{eq12-265}
\partial_\lambda(\omega(\lambda)t-\lambda x)=0,\\
\Updownarrow\nonumber\\
\omega'(\lambda)t-x=0.\nonumber
\end{gather}
Thus for given $x$ and $t$ large we find the stationary points by solving the equation 

\begin{equation}\label{eq12-266}
\omega'(\lambda)=\frac{x}{t}.
\end{equation}
For the Klein-Gordon equation we have from (\ref{eq12-230}) 

\begin{equation}\label{eq12-267}
\omega(\lambda)=(\gamma^2\lambda^2+c^2)^{\frac{1}{2}}. 
\end{equation}
Thus the equation we must solve is 

\begin{gather}
\frac{\gamma^2\lambda}{(\gamma^2\lambda^2+c)^{\frac{1}{2}}}=\frac{x}{t},\label{eq12-268}\\
\Downarrow\nonumber\\
\gamma^4\lambda^2=\left(\frac{x}{t}\right)^2(\gamma^2\lambda^2+c^2),\nonumber\\
\Updownarrow\nonumber\\
\left(\gamma^4-\gamma^2\left(\frac{x}{t}\right)^2\right)\lambda^2=c^2\left(\frac{x}{t}\right)^2,\\
\Updownarrow\nonumber\\
\gamma^2\left(\gamma^2-\left(\frac{x}{t}\right)^2\right)\lambda^2=c^2\left(\frac{x}{t}\right)^2. \label{eq12-269}
\end{gather}
Clearly this equation only has solutions if

\begin{equation}\label{eq12-270}
\big|\frac{x}{t}\big|<\gamma.
\end{equation}
But the light cone for the Klein-Gordon equation is $|x|=\gamma|t|$,

\begin{figure}[h!]
	\centering
	\includegraphics[scale=0.5]{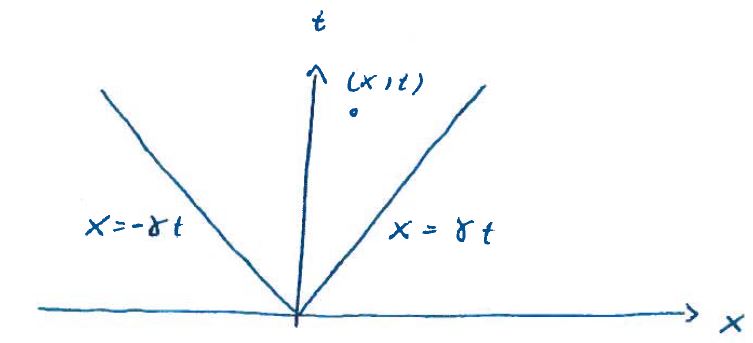}
	\caption{}
\end{figure}
\noindent{so stationary points only exists for $(x,t)$ inside the light cone. Outside the light cone there are no stationary points and according to the Riemann Lebesques lemma $I(x,t)\approx 0$.}

Assuming $\big|\frac{x}{t}\big|<\gamma$ we solve (\ref{eq12-269}) and find 

\begin{equation}\label{eq12-271}
\lambda=\frac{c}{\gamma}\frac{x}{t}\left(\gamma^2-\left(\frac{x}{t}\right)^2\right)^{-\frac{1}{2}}.
\end{equation}
We do not get $\pm$ in this formula because according to (\ref{eq12-268}) $\lambda$ and $\frac{x}{t}$ have the same sign. 

Observe that 

\begin{equation}\label{eq12-272}
\omega''(\lambda)=\gamma^2c^2\left(\gamma^2\lambda^2+c^2\right)^{-\frac{3}{2}}>0,
\end{equation}
so that 

\begin{equation}\label{eq12-273}
\varphi''(\lambda)=\omega''(\lambda)>0\quad\text{(see \ref{eq12-266})},
\end{equation}
and therefore 

\begin{equation}\label{eq12-274}
\frac{\varphi''(\lambda)}{|\varphi''(\lambda)|}=1.
\end{equation}
Thus the stationary phase formula (\ref{eq12-261}) gives us

\begin{equation}\label{eq12-275}
I(x,t)\approx\frac{1}{\sqrt{t}}H(\lambda)e^{i(\omega(\lambda)t-\lambda x+\frac{\pi}{4})},
\end{equation}
where $H(\lambda)=F_+(\lambda)(\omega''(\lambda))^{-\frac{1}{2}}$ and where the phase in (\ref{eq12-264}) has been written as 

\begin{equation}\label{eq12-276}
\omega(\lambda)t-\lambda x=t(\omega(\lambda)-\frac{x}{t})=t\varphi,
\end{equation}
so $k$ in the stationary phase formula is identified with $t\gg1$.

  So what have we found? Observe that 

\begin{equation}\label{eq12-277}
\frac{x}{t}=\omega'(\lambda)=\frac{\gamma}{\left(1+\left(\frac{c}{\gamma\lambda}\right)\right)^{\frac{1}{2}}}<\gamma.
\end{equation}
We have the picture

\begin{figure}[h!]
	\centering
	\includegraphics[scale=0.5]{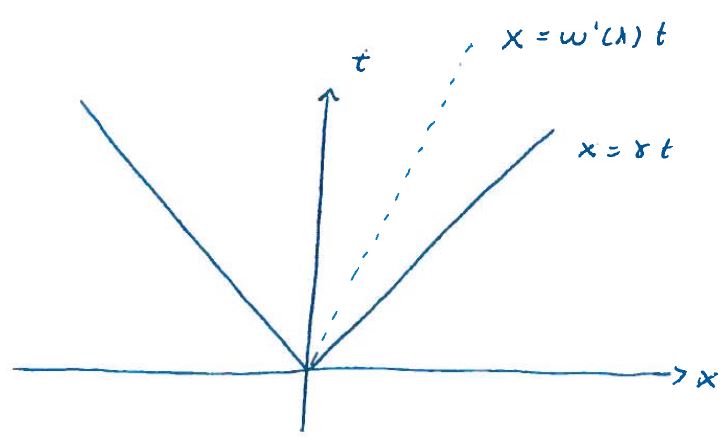}
	\caption{}
\end{figure}
According to (\ref{eq12-275}) and (\ref{eq12-260}) for large $x$ and $t$ along the line $x=\omega'(\lambda)t$, the wave field take the form of a plane wave with frequency $\omega(\lambda)$, wavelength $\frac{2\pi}{\lambda}$ and amplitude

\begin{equation}\label{eq12-278}
|I(x,t)|\sim\frac{1}{\sqrt{t}}|H(\lambda)|.
\end{equation}
The amplitude of the plane wave will decay as $\frac{1}{\sqrt{t}}$. Or to put it another way; if we shift to a moving frame moving at speed $\omega'(\lambda)$ with respect to the lab frame we would see a plane wave whose amplitude decays as $\frac{1}{\sqrt{t}}$ and whose frequency $\omega(\lambda)$ and wavelength $\frac{2\pi}{\lambda}$ will depend on our speed. Observe that $\omega'(\lambda)$ is just the group velocity that we introduced in (\ref{eq12-250}).

One can prove that for the Klein-Gordon equation, and also for many other dispersive equations, the energy in the wave field is transported at the group velocity. 

The group velocity is the relevant velocity for dispersive wave equations. It is also a concept whose importance for the description of wave phenomena in general hardly can be overstated.  

\section{Projects}
\subsection{Project1}

In this project the aim is to develop and implement on a computer a finite
difference algoritm for the solution of the 1D wave equation. Deriving the
finite difference algoritm and implementing it in a computer code is not
really the challenging and time consuming step here. The real challenge is
to make a convincing case for the correctness of the computer
implementation. As I have stressed in my lectures we should never, ever,
trust anything that comes out of a computer without extensive testing of the
implemented code.

So how do we proceed to build the needed trust in our implementation? The
best way to do this is to run the code in a situation where an analytic
solution is known. Most of the time an analytic solution also needs to be
implemented numerically, it typically include infinite sums of functions and
integrals that can not be solved in terms of elementary functions. So what
we are really doing is to compare the output of two different computer
implementations for the same solution. If this is to make sense and build
confidence, it is of paramount importance that the two codes implementing
the finite difference algoritm and the analytical solution should be as
independent as humanly possible. Under such circumstances closely similar
output from the two codes will build confidence that our finite difference
solution is correct.

So the challenge then is to find analytical solutions to the problem at
hand, in our particular case this is the 1D wave equation subject to some
boundary conditions that we will specify shortly. It is typically the case
that we will not be able to find analytic solutions to the actual problem we
are trying to solve but to a related problem or a special case of the
problem where the domain, boundary conditions and/or equation is changed
somewhat or fixed by some specific choises. After all, if we could find an
analytic solution to the actual problem it would not be much point in
spending a lot of effort in trying to solve it using finite differences. It
is however important that changing from the actual problem to the related
problem does not involve major modifications to the core of the finite
difference code. The reason for this is that such major modifications are
their own source of error and as a consequence our test will not be very
convincing as a test of the actual problem.

So what is the actual problem then? Here it is:%
\begin{eqnarray*}
u_{tt}-u_{xx} &=&0\text{\ \ \ \ }-l<x<l \\
u(-l,t) &=&h(t) \\
u(l,t) &=&k(t) \\
u(x,0) &=&\varphi (x) \\
u_{t}(x,0) &=&\psi (x)
\end{eqnarray*}

a) \ \ Implement a finite difference method for this problem. The code
should be flexible enough to handle arbitrary choises for the functions $h$,$%
k\,$,$\varphi $ and $\psi $.

b) \ \ Show that 
\[
u(x,t)=\frac{1}{2}(\varphi (x+t)+\varphi (x-t)) 
\]

is a solution to the problem 
\begin{eqnarray*}
u_{tt}-u_{xx} &=&0\text{\ \ \ \ }-\infty <x<\infty \\
u(x,0) &=&\varphi (x) \\
u_{t}(x,0) &=&0
\end{eqnarray*}

In your finite difference code choose $l$ and $\varphi $ in such a way that
the interval is much larger than the support of $\varphi $.( By the support
of $\varphi $ I mean the set of points where the function is nonzero). A
Gaussian function of the form 
\[
\varphi (x)=ae^{-bx^{2}} 
\]%
is a good choise if the parameters $a$, $b$ and $x_{0}$ are choosen in the
right way. Also choose $h=k=0$. As long as the support of the solution $u$
is well away from the boundary of the interval $[0,l]$ our exact solution
and the finite difference solution should match. This is our first test.

c) \ \ \ Our second test involve a rather general idea called an "artificial
source". We consider the modified problem 
\begin{eqnarray*}
u_{tt}-u_{xx} &=&\rho (x,t)\text{\ \ \ \ }-l<x<l \\
u(-l,t) &=&h(t) \\
u(l,t) &=&k(t) \\
u(x,0) &=&\varphi (x) \\
u_{t}(x,0) &=&\psi (x)
\end{eqnarray*}

This equation would appear as the evolution equation for a string clamped at
both ends and where there is a precribed force acting along the string. This
kind of prescribed influence is called a "source". In the artificial source
test we turn the problem around and rather assume a solution of a particular
type and then compute what the corresponding source must be.

We thus pick a function $u(x,t)$ and calculate the source $\rho (x,t)$,
boundary functions $h(t),k(t)$ and initial data $\varphi (x)$ and $\psi (x)$
so that $u(x,t)$ is in fact a solution to our modified problem. Thus 
\begin{eqnarray*}
\rho (x,t) &=&u_{tt}(x,t)-u_{xx}(x,t) \\
\varphi (x) &=&u(x,0) \\
\psi (x) &=&u_{t}(x,0) \\
h(t) &=&u(-l,t) \\
k(t) &=&u(l,t)
\end{eqnarray*}

We now have a analytic solution to our modified problem and can compare this
solution to the solution produced by our finite difference code. You are
free to choose the function $u(x,t)$ for yourself but here is a possible
family of choises%
\[
u(x,t)=ae^{-b(x-x_{0}\cos (\omega t))^{2}} 
\]

This is a family of Gaussian functions whose center move back and forth in a
periodic manner determined by the parameters $x_{0}$ and $\omega $. Compare
the numerical solution to the exact solution for several choises of
parameters. Make sure that at least one of the test involve the boundary
conditions in a significant way. By this I mean that $h(t)$ and $k(t)$
should be nonzero for some values of $t$.

d) \ \ Now use you code to solve the problem 
\begin{eqnarray*}
u_{tt}-u_{xx} &=&0\text{\ \ \ \ }-l<x<l \\
u(-l,t) &=&\sin (\omega t) \\
u(l,t) &=&0 \\
u(x,0) &=&0 \\
u_{t}(x,0) &=&0
\end{eqnarray*}

This model what happend when you periodically move the left end of the
string,whose right end is clamped, up and down. What do you expect to
happend? Do your code verify your intuition?

\subsection{Project2}
In this project the aim is to develop a stable numerical scheme for the
initial value problem%
\begin{eqnarray}
au_{t}+bu_{x} &=&\rho ,\text{ \ \ \ \ }-L<x<L,\text{ \ }0<t<T
\label{equation} \\
u(-L,t) &=&0  \nonumber \\
u(L,t) &=&0  \nonumber \\
u(x,0) &=&\varphi (x)  \nonumber
\end{eqnarray}%
where 
\begin{eqnarray*}
a &=&a(x,t) \\
b &=&b(x,t,u) \\
\rho  &=&\rho (x,t)
\end{eqnarray*}%
Remark: Since the equation is first order in time, one can actually not pose
boundary conditions at both endpoints. Thus we should only use the condition
at one end point, for example $u(-L,t)=0$. However since we are using center
difference we really also need a boundary condition at $x=L$. In this
problem I ask you to use the condition $u(-L,t)=0$. Using these two
conditions you can construct a finite difference scheme using center
difference in space. However, and this is of outmost importance, you must
make sure that the solution does not become nonzero in the region close to
the boundaries, if it does the scheme will go unstable no matter what we do.
But this restriction is really ok since we really are trying to simulate the
dynamics on an infinite domain, and thus the computations boundaries at $%
x=-L,L$ should be "invisible" at all time.

\begin{description}
\item[a)] Derive a numerical scheme for equation (\ref{equation}) using
center differences for the time and space derivatives. Find the Von-Neuman
condition for numerical stability of your scheme for the case when $%
a=a_{0},b=b_{0}$ are constants and $\rho =0$.

\item[b)] Verify the numerical stability condition by running the your
numerical scheme for some initial condition of your choosing. Run with grid
parameters both in the Von Neumann stable and unstable regime Since our
numerical scheme is meant to be a good approximation to solutions of
equation (\ref{equation}) on an infinite interval, you must make sure that
the initial condition is very small or zero close to the end points of the
interval $[-L,L]$. A well localized Gaussian function will do the job.

\item[c)] Solve the initial value problem%
\begin{eqnarray*}
(3x^{2}+1)u_{t}+2tu_{x} &=&0,\text{ }-\infty <x<\infty \\
u(x,0) &=&e^{-\gamma x^{2}}\text{\ \ \ \ \ \ }
\end{eqnarray*}%
using the method of characteristics. Plot some of the base characteristic
curves in the $x-t$ plane and also plot the solution $u(x,t)$ as a function
of $x$ for several $t>0$.

\item[d)] Solve the initial value problem from c) using your finite
difference scheme. Compare with the solution exact solution from c) by
plotting them on top of each other for several times.

\item[e)] Let the following quasi linear initial value problem be given%
\begin{eqnarray*}
u_{t}+uu_{x} &=&0,\text{ \ \ \ \ }-\infty <x<\infty \\
u(x,0) &=&\varphi (x)
\end{eqnarray*}%
Find exact formulas for the shock time for solutions to the initial value
problem from e) corresponding to the initial data

\begin{enumerate}
\item 
\[
\varphi (x)=\left\{ 
\begin{array}{cc}
0, & x>1 \\ 
1-x, & 0<x<1 \\ 
x+1, & -1<x<0 \\ 
0, & x<0%
\end{array}%
\right. 
\]

\item 
\[
\varphi (x)=\frac{1}{1+x^{2}} 
\]

\item 
\[
\varphi (x)=sech(x) 
\]
\end{enumerate}

 Test the validity of your finite difference code, and exact formulas
for the shock time, by running your code with the given initial conditions
and verifying that the graph of the numerical solutions become approximately
vertical for some $x$ as the shock time is approached.
\end{description}

\subsection{Project3}

In this project we are going to solve a heat equation numerically using the
finite difference method and the finite Fourier transform.

\begin{eqnarray}
u_{t}-u_{xx} &=&\rho (x,t)\text{ },0<x<l\text{ },\text{ }t>0  \label{problem}
\\
u(0,t) &=&f(t)  \nonumber \\
u(l,t) &=&g(t)  \nonumber \\
u(x,0) &=&\varphi (x)  \nonumber
\end{eqnarray}

\begin{description}
\item[a)] \ \ Implement a finite difference method for this problem. The
code should be flexible enough to handle arbitrary choises for the functions 
$f$, $g\,$, $\rho $ and $\varphi $.

\item[b)] The function

\item 
\[
u_{e}(x,t)=\frac{1}{l}(xh(x,t)(\frac{g(t)}{h(l,t)})+(l-x)h(x,t)(\frac{f(t)}{%
h(0,t)}) 
\]

\item satisfy the boundary conditions. Calculate the corresponding
artificial source $\rho $ and initial condition $\varphi $. Now run your
finite difference code with the calculated $\rho $ and $\varphi $ and
compare the numerical solution to the exact solution $u_{e}(x,t)$. Do this
with a couple of choises for the functions $f(t),g(t)$ and $h(x,t)$ and for
several values of the parameters in the problem. For each choise plot the
numerical and exact solution in the same plot and thereby verify that they
coincide.

\item[c)] Solve (\ref{problem}) using the finite Fourier transform. Before
you apply the finite Fourier transform you should convert (\ref{problem})
into a problem with homogenous boundary conditions like in equation 4.6.45
in our textbook. Use the same choises for $f(t)$ and $g(t)$ that you used in
problem b) and choose your own source. Compare the finite Fourier transform
solution with the finite difference solution for some choises of initial
conditions. Make sure that the initial conditions is consistent with the
boundary conditions%
\begin{eqnarray*}
\varphi (0) &=&f(0) \\
\varphi (l) &=&g(0)
\end{eqnarray*}

\item[d)] Let the source in problem (\ref{problem}) be the nonlinear source
from equation 4.7.1 in our textbook.%
\[
\rho (x,t)=\widehat{\lambda }u(x,t)(1-u^{2}(x,t)) 
\]

\item Use the finite difference code to verify the linear and nonlinear
stability results found in section 4.7 of our textbook. In particular you
should shown that
\end{description}

\begin{enumerate}
\item \ \ $u(x,t)=0$ is stable for $\widehat{\lambda }<1$

\item $u(x,t)$ is unstable for $1<$ $\widehat{\lambda }<4$ and approach the
solution 4.7.24 as time increase to infinity.
\end{enumerate}

Run the finite difference code for values of $\widehat{\lambda }$ higher
than 4 and plot the resulting solutions. Any conclusions?


\begin{thebibliography}{99}
\bibitem{textbook} ``Partial Differential Equations of Applied Mathematics'', Erich Sauderer, 2006, Wiley.
\bibitem{Weinberger} ``A first Course in Partial Differential Equations'',H. F. Weinberger, 1965, John Wiley and Sons.
\bibitem{Folland} ``Introduction to Partial Differential Equations'', Gerald B. Folland,1976, Princeton University Press.
\bibitem{Strauss} ``Partial Differential Equations, An Introduction'',Walter A. Strauss,2007, John Wiley and Sons.
\bibitem{Courant} ``Methods of mathematical Physics, Volume II'', Courant and Hilbert,1962, Interscience Publishers.
\bibitem{John} ``Partial Differential Equations'', Fritz John,1982,Springer.
\bibitem{Keener} ``Methods of Applied Mathematics'', James P. Keener, 1988,Springer.
\bibitem{Stakgold} ``Green's functions and Boundary value problems'', Ivar Stakgold and Michael J. Holst,2011, Wiley.
\end{thebibliography}

\end{document}